\documentclass[a4paper,12pt]{article}

\usepackage{amssymb,amsthm, graphicx,latexsym,amsmath}
\newtheorem{theorem}{Theorem}[section]
\newtheorem{lemme}{Lemma}[section]
\newtheorem{remark}{Remark}[section]
\newtheorem{corollary}{Corollary}[section]
\newtheorem{proposition}{Proposition}[section]

\usepackage{color}  

\makeatletter

\@addtoreset{equation}{section}

\makeatother

\begin{document}

\centerline {\bf ROBIN FUNCTIONS FOR COMPLEX MANIFOLDS}
\centerline {\bf  AND APPLICATIONS  \footnote{\  This paper is dedicated
to Professor John Wermer on the occasion of his 80th birthday.}}

\vskip6pt

\centerline {\bf Kang-Tae Kim, Norman Levenberg and Hiroshi Yamaguchi}

\vskip10pt

\section{Introduction} \setcounter{section}{1} In \cite{Y} and \cite{LY}
we analyzed the second variation of
the Robin function associated to a smooth variation of domains in ${\mathbb{C}}^n$
for $n\geq 2$; i.e., ${\cal D}=\cup_{t\in B}(t,D(t))\subset
B\times {\mathbb{C}}^n$ is a variation of domains $D(t)$ in ${\mathbb{C}}^n$ each
containing a fixed point $z_0$ and with $\partial D(t)$ of class
$C^{\infty}$ for $t\in B:=\{t\in {\mathbb{C}}:|t|<\rho\}$. For such $t$ and
for $z\in \overline {D(t)}$ we let $g(t,z)$ be the ${\mathbb{R}}^{2n}$-Green
function for the domain $D(t)$ with pole at $z_0$; i.e., $g(t,z)$ is
harmonic in $D(t)\setminus \{z_0\}$, $g(t,z)=0$ for $z\in \partial
D(t)$, and $g(t,z)-{1\over ||z-z_0||^{2n-2}}$ is harmonic near $z_0$. We
call $$\lambda (t):=\lim_{z\to z_0} [g(t,z)-{1\over ||z-z_0||^{2n-2}}]$$
the {\it Robin constant} for $(D(t),z_0)$. Then
 \begin{align}
 \label{eqn:ff}{\partial ^2\lambda \over \partial t \partial \overline
t}(t)=-c_n\int_{\partial
D(t)}k_2(t,z)||\nabla_zg||^2d S_z
-4c_n\int_{D(t)}\sum_{a=1}^n|{\partial ^2g\over \partial t
\partial \overline z_a}|^2 dV_z.
\end{align}

\noindent Here, $c_n={1\over (n-1)\Omega_n}$ is a positive
dimensional constant where $\Omega_n$ is the
area of the unit sphere in ${\mathbb{C}}^n$, $dS_z$ and $dV_z$ are the Euclidean area
element on $\partial D(t)$ and volume element on $D(t)$,
$\nabla _zg=(\frac{\partial
 g}{\partial z_1}, \cdots ,  \frac{\partial g}{\partial z_n}  ) $ and
$$k_2(t,z):=||\nabla_z \psi||^{-3}\bigl[{\partial ^2 \psi\over \partial t
\partial \overline
t}||\nabla_z \psi||^2-2\Re\{{\partial  \psi\over \partial t}\sum_{a=1}^n
{\partial  \psi\over
\partial \overline z_a}{\partial ^2\psi \over \partial \overline t
\partial  z_a}\}+ |{\partial
\psi\over \partial t}|^2\Delta_z\psi \bigr] $$
is the so-called {\it Levi-curvature} of
$\partial {\cal D}$ at $(t,z)$. The function
$\psi(t,z)$ is a defining function for ${\cal D}$ and the numerator is
the sum of the Levi-form of $\psi$ applied to the
$n$ complex tangent vectors $(-{\partial \psi\over \partial
z_j},0,...,{\partial \psi\over \partial
t},0,...,0)$. In particular, if
${\cal D}$ is pseudoconvex (strictly pseudoconvex) at a point $(t,z)$
with $z\in \partial D(t)$, it
follows that
$k_2(t,z)\geq 0 \ (k_2(t,z)> 0)$ so that
$-\lambda (t)$ is subharmonic (strictly subharmonic) in $B$. Given a bounded domain
$D$ in ${\mathbb{C}}^n$, we let $\Lambda (z)$ be the Robin constant for $(D,z)$.
If we fix a point
$\zeta_0 \in D$, for $\rho>0$ sufficiently small and $a\in {\mathbb{C}}^n$, the
disk
$\zeta_0+aB:=\{\zeta=\zeta_0+at, \ |t|< \rho\}$ is contained in $D$.
Using the biholomorphic
mapping $T(t,z)=(t,z-at)$ of $B\times {\mathbb{C}}^n$, we get the variation of domains
${\cal D}=T(B\times
D)$ where each domain $D(t):=T(t,D)=D-at$ contains $\zeta_0$. Letting
$\lambda (t)=\Lambda
(\zeta_0+at)$ denote the Robin constant for $(D(t),\zeta_0)$ and using
(\ref{eqn:ff}) yields part of the
following surprising result (cf., ~\cite{Y} and ~\cite{LY}).

\begin{theorem}\label{thm:01}\   Let $D$ be a bounded pseudoconvex
domain in ${\mathbb{C}}^n$ with $C^{\infty}$ boundary.
Then $\log {(-\Lambda(z))}$ and $-\Lambda(z)$ are real-analytic, strictly
plurisubharmonic
exhaustion functions for $D$. \par
 \end{theorem}

We now study a generalization of the second variation formula
(\ref{eqn:ff}) to complex manifolds $M$ equipped with a Hermitian metric $ds^2$ and a smooth, nonnegative function $c$. Our purpose is that, with this added flexibility, {\it we are able to give a criterion for a bounded, smoothly bounded, pseudoconvex domain $D$ in a complex homogeneous space to be Stein}. In particular, we are able to do the following:
\smallskip

\begin{enumerate}
\item Describe concretely all the non-Stein pseudoconvex domains $D$ in the complex torus of Grauert (section 5).
\item Give a description of all the non-Stein pseudoconvex domains $D$ in the special Hopf manifolds $\mathbb{H}_n$ (section 6).
\item Give a description of all the non-Stein pseudoconvex domains $D$ in the complex flag spaces ${\cal F}_n$ (section 7).
\item Give another explanation as to why all pseudoconvex subdomains of complex projective space, or, more generally, of complex Grassmannian manifolds, are Stein (Appendix A).
\end{enumerate}

\smallskip
The metric $ds^2$ and the function $c$ give rise to a $c$-Green function
and
$c$-Robin constant associated to an open set $D\subset M$ and a point $p_0\in D$. We then take a variation ${\cal D}=\cup_{t\in B}(t,D(t))\subset B\times M$ of
domains $D(t)$ in $M$ each containing a fixed point $p_0$ and define a $c$-Robin function $\lambda(t)$. The precise definitions of these notions and the new variation formula (\ref{eqn:f2}) will be given
in the next section. In section 3 we impose a natural condition (see (\ref{eqn:psm})) on the metric $ds^2$ which will be useful for applications. K\"ahler metrics, in particular, satisfy (\ref{eqn:psm}).
After discussing conditions which insure that the function
$-\lambda $ is subharmonic, we will use (\ref{eqn:f2}) to develop a
``rigidity lemma'' (Lemma \ref{lem:rigid}) which will imply, if $-\lambda $ is not strictly subharmonic, the existence of a nonvanishing, holomorphic vector field on $M$ with certain properties (Corollary \ref{lem:rigid-2}). This will be a key tool in constructing strictly
plurisubharmonic exhaustion functions for pseudoconvex
subdomains $D$ with smooth boundary in certain complex Lie groups and in
certain complex homogeneous spaces; i.e., we use these $c$-Robin
functions to verify that $D$ is Stein.

Specifically, in section 5 we study pseudoconvex domains $D$ in a
complex Lie group $M$. In a sense which will be made precise in
Corollary \ref{thm:4.3} the functions $-\Lambda$ we construct in this setting are  the ``best
possible'' plurisubharmonic exhaustion functions: {\it
if a $c$-Robin function $\Lambda$ for $D$ is such that $-\Lambda$ is {\bf not} strictly plurisubharmonic, then $D$ does not admit a strictly plurisubharmonic exhaustion function}. We characterize the smoothly bounded, relatively compact pseudoconvex domains $D$ in a complex Lie group $M$ which are Stein in Theorem \ref{thm:fund-lie-case}. Then we apply this result to describe all of the non-Stein pseudoconvex domains $D$ in the complex torus example of Grauert.

In section 6 we let $M$ be an $n$-dimensional complex homogeneous space with
an associated connected complex Lie group
$G\subset$Aut$M$ of complex dimension $m \ge n$.  We set up our $c$-Robin function machinery to discuss sufficient conditions on the pair $(M,G)$ such that for every
smoothly bounded, relatively compact pseudoconvex domain $D$ in $M$, the function $-\lambda$ is strictly plurisubharmonic on $D$ (Theorem \ref{3ptspan}). In particular, the Grassmann manifolds $M=G(k,n)$ with $G=GL(n,{\mathbb{C}})$ satisfy one of these conditions.

In Theorems \ref{thm:mainth} and \ref{thm:main-nonconnected} we give an
analogue of Theorem \ref{thm:fund-lie-case} to characterize the smoothly
bounded, relatively compact pseudoconvex domains $D$ in a complex
homogeneous space $M$ which are Stein. We immediately apply this result
to special Hopf manifolds $\mathbb{H}_n$. Then in section 7 we apply the result to
describe all of the non-Stein pseudoconvex domains $D$ in complex flag
spaces ${\cal F}_n$ (Theorem \ref{lasttheorem}).

Some of the
results in this paper were
announced without proof in \cite{LY2} and \cite{KLY}; in this paper, we provide complete
proofs and illustrate the
significance of this generalization of the second variation formula with
applications and concrete examples. The material presented is completely self-contained; in particular, all concepts pertaining to Lie theory and homogeneous spaces are explained.

We thank Professor T. Ueda for his helpful advice in our study of Levi problems for flag spaces. We also
thank Professor T. Morimoto  for his useful comments regarding Lie algebras.

\vskip8pt

\section{The variation formula.}
\setcounter{section}{2}
\  Our general set-up is this: let $M$ be
an
$n$-dimensional complex manifold (compact or not) equipped with a
Hermitian metric
$$ds^2=\sum_{a,b=1}^n g_{a\overline b}dz_a\otimes d{\overline z}_b$$
and let $\omega:=i \sum_{a,b=1}^n g_{a\overline b}dz_a\wedge d{\overline
z}_b$ be the associated real  $(1,1)$ form. As in the introduction, we
take $n\geq 2$. We write $g^{\overline a
b}:=(g_{a\overline b})^{-1}$ for the elements of the inverse matrix to
$(g_{a\overline b})$ and
$G:=\hbox{det}(g_{a\overline b})$. Note that $\omega^n=2^nn!Gdx_1\wedge
\cdots \wedge dx_{2n}$ locally where
$z_k=x_{2k-1}+ix_{2k}$. For a domain $W\subset M$, we let ${\cal
L}^{p,q}(W)$ denote the $(p,q)$ forms
on $W$ with complex-valued, $C^{\infty}(W)$-functions as
coefficients. We have the standard linear operators

$$*:{\cal L}^{p,q}(W)\to {\cal L}^{n-q,n-p}(W),$$

$$\partial:{\cal L}^{p,q}(W)\to {\cal L}^{p+1,q}(W),$$

$$\overline \partial: {\cal L}^{p,q}(W)\to {\cal L}^{p,q+1}(W),$$

$$\delta:=-*\partial *:{\cal L}^{p,q}(W)\to {\cal L}^{p,q-1}(W),$$

$$\overline{  \delta}:=-* \overline{  \partial} *:{\cal L}^{p,q}(W)\to
{\cal L}^{p-1,q}(W),$$
and $d=\partial + \overline \partial$.  We get the box Laplacian operator
$$\delta \overline \partial + \overline \partial \delta:{\cal
L}^{p,q}(W)\to {\cal L}^{p,q}(W)$$
and its conjugate
$$\overline \delta \partial + \partial  \overline \delta:{\cal
L}^{p,q}(W)\to {\cal L}^{p,q}(W).$$
Adding these, we obtain the Laplacian operator
$$\Delta= \delta \overline \partial + \overline \partial \delta +
\overline \delta \partial +
 \partial  \overline \delta$$
which is a real operator; in local coordinates acting on functions this
has the form
\begin{align}
 \Delta u &=-2\bigl[\sum_{a,b=1}^ng^{\overline b a}{\partial ^2 u \over
\partial \overline z_b \partial
z_a  }+ {1\over 2}\sum_{a,b=1}^n({1\over G}{\partial (Gg^{\overline b
a})\over \partial
z_a}{\partial u \over \partial \overline z_b}+ {1\over G}{\partial
(Gg^{\overline a b})\over
\partial \overline z_a}{\partial u \over \partial z_b})\bigr] \nonumber \\
&=:-2[Pu +Ru]. \label{eqno(1.1)}
\end{align}

We call $Pu$ the {\it principal part} of $\Delta u$. Also, we remark that
if $ds^2$ is K\"ahler, i.e., if $d\omega =0$,
then $\Delta u=-2Pu=-2\sum_{a,b=1}^ng^{\overline
b a}{\partial ^2 u \over \partial \overline z_b
\partial  z_a}$.
 As usual, we set, for any $\alpha \in {\cal L}^{p,q}
(W)$, $ \|\alpha\|^2_W = \int_{ W} \alpha \wedge \overline{  *\alpha} \
\ge 0.$

\smallskip
Given a nonnegative $C^{\infty}$ function $c=c(z)$ on $M$, we call a
$C^{\infty}$ function $u$ on an open set $D\subset M$ {\it $c$-harmonic}
on $D$ if $\Delta u +cu=0$ on $D$. Choosing local coordinates near a
fixed point $p_0\in M$ and a coordinate neighborhood $U$ of $p_0$ such
that $[g_{a\overline
b}(p_0)]_{a,b=1,...,n}=[\delta_{ab}]_{a,b=1,...,n}$, the Laplacian
$\Delta$ corresponds to a second-order elliptic operator $\tilde \Delta$
in ${\mathbb{C}}^n$. In particular, we can find a $c$-harmonic function $Q_0$ in
$U\setminus \{p_0\}$ satisfying $$\lim_{p\to p_0} {Q_0(p)
d(p,p_0)^{2n-2}}=1$$ where $d(p,p_0)$ is the geodesic distance between
$p$ and $p_0$ with respect to the metric $ds^2$. We call $Q_0$ a {\it
fundamental solution} for $\Delta$ and $c$ at $p_0$. Fixing $p_0$ in a
smoothly bounded domain $D\Subset M$ and fixing a fundamental
solution $Q_0$, the {\it $c$-Green function} $g$ for $(D,p_0)$ is the
$c$-harmonic function in $D\setminus \{p_0\}$ satisfying $g=0$ on
$\partial D$ ($g$ is continuous up to $\partial D$) with $g(p)-Q_0(p)$
regular at $p_0$. The $c$-Green function always and uniquely exists (cf., \cite{NS}) and  is nonnegative on $D$.  Then $$\lambda:=\lim_{p\to p_0} [g(p)-Q_0(p)]$$
is called the {\it $c$-Robin constant} for $(D,p_0)$. Thus we have
$$g(p)=Q_0(p)+\lambda +h(p)$$ for $p$ near $p_0$, where $h(p_0)=0$. In
case $M$ is compact, if $c\not \equiv 0$ on $M$, then the $c$-Green
function $g$ for $(M, p_0)$ exists and is positive on $M$, hence the $c$-Robin
constant is finite. But if $c\equiv 0$ on $M$, a $c$-harmonic function is harmonic and
cannot attain its minimum; thus, in this case, $g(z)\equiv +\infty$ on
$M$ (cf., \cite{NS}). In this case we set $\lambda =+\infty$.

Now let ${\cal D}=\cup_{t\in B}(t,D(t))\subset B\times M$ be a
$C^{\infty}$ variation of
domains $D(t)$ in $M$ each containing a fixed point $p_0$ and with
$\partial D(t)$ of class
$C^{\infty}$ for $t\in B$. This means that there exists $\psi(t,z)$ which
is $C^{\infty}$ in
a neighborhood $N\subset B\times M$ of $\{(t,z): t\in B, \ z\in \partial D(t)\}$, negative
in $N\cap \{(t,z): t\in B, \ z\in D(t)\}$, and for each
$t\in B, \ z\in \partial D(t)$,
we require that $\psi(t,z)=0$ and ${\partial \psi\over \partial
z_i}(t,z)\not = 0$ for some
$i=1,...,n$. We call $\psi(t,z)$ a defining function for
${\cal D}$. Assume that $B \times \{p_0\}\subset {\cal D}$. Let
$g(t,z)$ be the
$c$-Green function for
$(D(t),p_0)$ and $\lambda (t)$ the corresponding $c$-Robin constant. The
hypothesis that ${\cal D}$ be a $C^{\infty}$ variation implies that for
each $t\in B$, the $c$-Green function $g(t,z)$ extends
of class $C^\infty$  beyond $\partial D(t)$; this follows from the general
theory of partial differential equations. Most of the
calculations and the subsequent results in this paper remain valid under
weaker ($C^2$ or $C^3$) regularity assumptions on ${\cal D}$.

\smallskip
Our formulas are the following:

\begin{align}
\label{eqn:f1}
&\frac{\partial \lambda }{\partial t} (t) = -c_n \int_{ \partial D(t)}
k_1(t,z)\sum_{a,b=1}^n(g^{\overline a b}{\partial g\over \partial
\overline z_a}{\partial g\over
\partial  z_b})d \sigma_z ,\\
 &{\partial ^2\lambda \over \partial t \partial \overline
t}(t)=-c_n\int_{\partial
D(t)}k_2(t,z)\sum_{a,b=1}^n(g^{\overline a b}{\partial g\over \partial
\overline z_a}{\partial g\over
\partial  z_b})d \sigma_z \nonumber  \\
& \quad - {c_n \over 2^{n-2}} \,  \bigl\{\,||\overline \partial {\partial
g\over \partial
t}||_{D(t)}^2+{1\over 2}||\sqrt c {\partial g\over \partial
t}||_{D(t)}^2+
 \frac{ 1}2 \Re
\int_{D(0)}
{\partial g \over
\partial  t}\bigl[{1\over i}\partial *\omega\wedge \overline \partial {
\partial g \over
\partial \overline t} +{1\over i}\overline \partial *\omega\wedge
\partial {\partial g \over \partial \overline t}\bigr]
 \bigr\} \nonumber\\
&\ \  :=-c_nI-{c_n\over 2^{n-2}}J \label{eqn:f2}
\end{align} where  $d \sigma_z$ is the area element on
$\partial D(t)$ with respect to the Hermitian metric, $c_n={1\over
(n-1)\Omega_n}$ and
\begin{align}\nonumber
&k_1(t,z):= [\sum_{a,b=1}^ng^{\overline a b}{\partial \psi\over \partial
\overline z_a}{\partial
\psi\over \partial  z_b}]^{-1/2} {\partial \psi \over \partial t}, \\
\nonumber
& k_2(t,z):=
[\sum_{a,b=1}^ng^{\overline a b}{\partial \psi\over \partial \overline
z_a}{\partial
\psi\over \partial  z_b}]^{-3/2} \times\\
& \bigl[{\partial ^2 \psi\over \partial t \partial \overline
t}(\sum_{a,b=1}^ng^{\overline a b}{\partial \psi\over \partial \overline
z_a}{\partial
\psi\over \partial  z_b})-2\Re\{{\partial  \psi\over \partial
t}(\sum_{a,b=1}^ng^{\overline a b}
{\partial \psi\over \partial \overline z_a}{\partial ^2 \psi\over
\partial  z_b \partial
\overline t})\} -\frac{ 1}2   |{\partial
\psi\over \partial t}|^2  \Delta_z \psi \bigr] ,
\label{eqn:k2}
\end{align}
$\psi(t,z)$ being a defining function for ${\cal D}$. Here, $k_i(t,z)\
(i=1,2)$
is a real-valued function for $(t,z)\in \partial {\cal D}$ which is
independent of both the
choice of defining function for ${\cal D}$ and of the choice of local
parameter $z$ in the manifold $M$. We call $k_2(t,z)$ the Levi scalar
curvature with respect to the metric $ds^2$.

\smallskip
Formula (\ref{eqn:f1}) is a generalization of the classical Hadamard
variation formula. For the study of several complex variables
the variation formula (\ref{eqn:f2}) of the second order is
fundamental and we now give the proof. First, for each $t \in
B$, the equation
$$g(t,p)=Q_0(p)+\lambda(t) +h(t,p),$$
where $h(t,p_0)=0$ for all $t\in B$, holds; hence we have, for $p\ne p_0$,
$${\partial ^2g\over \partial t \partial \overline{  t}} (t,p)
= {\partial^2 \lambda\over \partial t \partial \overline{  t}}(t)+
{\partial ^2 h\over \partial t \partial \overline{  t}}(t,p).
$$
Since ${\partial ^2h \over \partial t \partial \overline{  t}}(t, p_0)
=0$ for all $t \in B$, if we set ${\partial^2 g\over
\partial t \partial \overline{  t}}(t, p_0)= {\partial^2 \lambda\over
\partial t \partial \overline{  t}}(t)$, it follows that
${\partial^2 g\over \partial t \partial  \overline{  t}}(t,p)$ is a
$c$-harmonic function on all of $D(t)$ even though $g (t,p)$ has a
singularity at $p_0$. Fix $t_0\in B$. Using Green's formula for
$c$-harmonic functions, we thus obtain the
formula

\begin{align}
 \label{eqn:lt0}
{\partial ^2 \lambda\over \partial t \partial \overline t}(t_0)
 ={-c_n\over 2^ {n}}\int_{\partial D(t_0)}{\partial ^2g\over \partial t
\partial
\overline t}(t_0,z)*dg(t_0,z).
\end{align}
We note that, in  ${\mathbb{C}}^n$, from the  definition of
the $*$-operator, we have $$*dg(t_0,z) =-2^{n}\|\nabla_z g(t_0, z)\|dS_z
=2^{n-1}
 {\partial g \over \partial n_z}(t_0, z)dS_z$$ for $z \in \partial
 D(t_0)$, where
$dS_z$ and $n_z $ are  the Euclidean area element and
 the unit outer normal vector for $\partial D(t_0)$ at $z$.

We may assume $t_0=0$. Now using the fact that $-g(t,z)$ is a defining
function for ${\cal D}$, we can write, from (\ref{eqn:k2}),
\begin{align}
 & \nonumber k_2(0,z)\cdot(\sum_{a,b=1}^ng^{\overline ab}{\partial g\over
\partial \overline z_a}{\partial g\over
\partial  z_b})^{3/2}=-{\partial ^2g\over \partial t \partial
\overline t}(\sum_{a,b=1}^ng^{\overline ab}{\partial g\over \partial
\overline z_a}{\partial g\over
\partial  z_b})\\
 & \ \ \ +2\Re[(\sum_{a,b=1}^ng^{\overline ab}{\partial^2 g\over \partial
z_b\partial \overline t}{\partial g\over
\partial \overline z_a}){\partial g\over
\partial t}]+ \frac{ 1}2 |{\partial g\over
\partial t}|^2 \Delta g.
\label{eqno(1.5)}
\end{align}
Since $\Delta g+cg=0$ on $\partial D(0)$
 (here we use the fact that $g$
is of class $C^{\infty}$ on $\overline{  D(0)}\setminus \{p_0\})$ and $g=0$ on $\partial D(0)$,
we have from (\ref{eqn:lt0}) and  (\ref{eqno(1.5)})
 \begin{align*}
 {\partial ^2 \lambda\over \partial t \partial \overline t}(0)
&={-c_n\over 2^{n}}\int_{\partial D(0)}
\bigl\{\,-k_2(0,z)\sum_{a,b=1}^n(g^{\overline a b}{\partial g\over
\partial \overline z_a}{\partial g\over
\partial  z_b})^{1/2}\\
& \qquad \qquad \qquad  \quad   +{
2\Re\  \sum_{a,b=1}^n
  ( g^{\overline ab}{\partial^2 g\over \partial z_b\partial \overline t}
{\partial g\over
\partial \overline z_a}){\partial g\over
\partial t}
\over
\sum_{a,b=1}^ng^{\overline ab}{\partial g\over \partial \overline
z_a}{\partial g\over
\partial  z_b}}\, \bigr\}*dg(0,z)
\\
&\equiv (I) + (II).
\end{align*}

In general, for any $C^{\infty}$ defining function $v$ on $\overline{D(0)}$
($v=0$ on $\partial D(0)$ and $v<0$ on $D(0)$),
$$
*dv=2^{n}(\sum_{a,b=1}^ng^{\overline ab}{\partial v\over \partial
\overline z_a}{\partial v\over
\partial  z_b})^{1/2}d \sigma_z\ge 0
$$
for $z\in \partial D(0)$, where $d\sigma_z$ is the area element on
$\partial D(0)$ at $z$ with respect to the metric $ds^2$. We apply this
to $v=-g(0,z)$, plug this into ({\it I}\,), and obtain
the formula
\begin{align*}
 (I)&=  -c_n\int_{\partial
D(0)}k_2(0,z)\sum_{a,b=1}^n(g^{\overline a b}{\partial g\over \partial
\overline z_a}{\partial g\over
\partial  z_b})d \sigma_z.
\end{align*}

Now we work with ({\it II}\,). We need to calculate
$${\partial g\over \partial \overline z_a}*dg(0,z)$$
on $\partial D(0)$. To this end, we first note that for a function $u$,
we have the following formulas for $*\partial u$
and $*\overline \partial u$:
$$*\partial u = -i^n\sum_{a,b=1}^nGg^{\overline ab}{\partial u \over
\partial z_b}dz_a\wedge dz_1\wedge d\overline z_1 \wedge
\cdots\widehat{dz_a \wedge d \overline z_a}\cdots\wedge
dz_n\wedge d\overline z_n;$$
$$*\overline \partial u = i^n\sum_{a,b=1}^nGg^{\overline ba}{\partial u
\over \partial \overline z_b}d\overline z_a\wedge  dz_1\wedge d\overline
z_1 \wedge
\cdots\widehat{dz_a \wedge d \overline z_a}\cdots\wedge dz_n\wedge
d\overline z_n.$$ Now if $u=0$ on $\partial D(0)$, $du=\partial u +
\overline \partial u=0$ along $\partial D(0)$ and we obtain
$$*du=-2i^n\sum_{a,b=1}^nGg^{\overline ab}{\partial u \over \partial
z_b}dz_a\wedge  dz_1\wedge d\overline z_1 \wedge \cdots
\widehat{dz_a \wedge d \overline z_a}\cdots\wedge dz_n\wedge d\overline
z_n.$$ Applying this to $u=g(0,z)$ on $\partial D(0)$, we get
\begin{align*}
{\partial g \over \partial \overline z_a} * dg(0,z)&= -2i^n
\sum_{i,j=1}^nGg^{\overline ij}{\partial g \over \partial
\overline{z}_a}{\partial g \over \partial z_j}dz_i\wedge  dz_1\wedge
d\overline z_1
\wedge \cdots\widehat{dz_i \wedge d
\overline z_i}\cdots\wedge dz_n\wedge
d\overline z_n.
\end{align*}
 Again we use the fact that $g(0,z)=0$ on $\partial D(0)$ so that
$dg(0,z)=0$ along $\partial D(0)$. Then, if $i \ne a$, ${\partial g\over
\partial \overline{  z} _a}d \overline{  z}_a$ in the above formula can
be replaced by $ - {\partial g \over
\partial \overline{z} _i }d \overline{  z}_i$. It follows that
 \begin{align}\nonumber
 {\partial g \over \partial \overline z_a} * dg(0,z)=
-2i^n\bigl(G\, \sum_{i,j=1}^n g^{\overline ij}{\partial g \over \partial
z_j}{\partial g \over \partial
\overline z_i}\bigr)\, dz_a\wedge  dz_1\wedge d\overline z_1 \wedge
\cdots\widehat{dz_a \wedge d
\overline z_a}\cdots\wedge dz_n\wedge
d\overline z_n,
\end{align}
where the term in parenthesis is independent of $a=1,...,n$. We plug this
into ({\it II}\,) where
$$(II)={-c_n\over 2^{n}}\int_{\partial D(0)}F$$
and
\begin{align*}
 F&:=\bigl\{\,
{2\Re\sum_{a,b=1}^n (g^{\overline ab}{\partial^2 g\over \partial
z_b\partial \overline t}{\partial g\over
\partial \overline z_a}){\partial g\over
\partial t}
\over \sum_{a,b=1}^ng^{\overline ab}{\partial g\over \partial \overline
z_a}{\partial
g\over
\partial  z_b}}\, \bigr\}*dg(0,z)\\
&=-4\Re \ \bigl\{i^n
\bigl(\sum_{a,b=1}^n Gg^{\overline ab}{\partial^2 g\over \partial
z_b\partial \overline t}{\partial g\over
\partial t} \bigr)\\
 &\qquad \qquad \qquad  \ \ \ \
dz_a\wedge  dz_1\wedge d\overline z_1 \wedge \cdots\widehat{dz_a \wedge d
\overline z_a}\cdots\wedge dz_n\wedge d\overline z_n \bigr\}.
\end{align*}
Note that the denominator cancels.

Next we use the relation
$$*\partial ({\partial g \over \partial \overline t})=-i^n\sum_{a,b=1}^n
Gg^{\overline ab}{\partial ^2g \over \partial z_b \partial \overline
t}\,dz_a\wedge  dz_1\wedge d\overline z_1 \wedge \cdots\widehat{dz_a \wedge
d \overline z_a}\cdots\wedge dz_n\wedge d\overline z_n$$
to obtain
$$F=4\Re \bigl\{{\partial g \over \partial t}(*\partial ({\partial g
\over \partial \overline t}))\bigr\}$$
on $\partial D(0)$.

Thus we obtain
\begin{align*}
 {\partial ^2 \lambda\over \partial t \partial \overline t}(0)= &-c_n
\, \int_{\partial D(0)}k_2(0,z)
[\sum_{a,b=1}^ng^{\overline ab}{\partial g\over \partial \overline
z_a}{\partial
g\over
\partial  z_b}]d \sigma_z \\
&- {c_n \over 2^{n-2}}\Re\int_{\partial D(0)}\, {\partial g \over
\partial t}(*\partial ({\partial g \over \partial \overline t}))
 =:-I-J
\end{align*}
where
$$J:= { c_n \over 2^{n-2} }\ \Re\int_{\partial D(0)}f$$
 and
$$
  f:={\partial g \over \partial t}(*\partial ({\partial g \over \partial
\overline t}))
$$
is an $(n,n-1)$ form. We want to convert $J$ into a volume integral; to
do so we must compute $\overline \partial f$ since, $f$ being $(n,n-1)$
form,
$$J={c_n\over 2^{n-2}}\, \Re\, \int_{D(0)}df={c_n\over 2^{n-2}}\,
\Re\int_{D(0)}\overline \partial f.$$
We get two terms whose sum comprise $\overline \partial f$:
$$\hbox{(i)}=(\overline \partial {\partial g \over \partial t})\wedge
(*\partial ({\partial g \over \partial \overline t}));$$
$$
\hbox{(ii)}={\partial g \over \partial t}\overline \partial (*\partial
({\partial g \over \partial \overline t})).$$

Now
$$\int_{D(0)}\hbox{(i)}=\int_{D(0)}(\overline \partial {\partial g \over
\partial t})\wedge (*\partial ({\partial g \over \partial \overline
t}))=||\overline \partial
{\partial g \over \partial t}||^2_{D(0)}.$$
We note that in ${\mathbb{C}}^n $, we have
$||\overline \partial
{\partial g \over \partial t}||^2_{D(0)}=2^n \int_{ D(0)}
(\sum_{ i=1}^n |{\partial ^2 g \over \partial t \partial \overline{
z}_i}|^2)dV_z$, where $dV_z$ is the Euclidean volume element of ${\mathbb{C}}^n$,
so that the first term ${c_n \over 2^{n-2}}\,{\Re\,
\int_{D(0)}\hbox{(i)}}$ of
$J$ coincides with the last term of formula (\ref{eqn:ff}).

Next, for (ii),
$$\overline \partial (*\partial ({\partial g \over \partial \overline
t}))=i^n\sum_{a,b=1}^n\bigl[ {\partial (Gg^{\overline ab})\over \partial
\overline z_a}
{\partial ^2g \over \partial z_b \partial \overline t}+Gg^{\overline
ab}{\partial ^3g \over \partial \overline z_a \partial z_b \partial
\overline
t}\bigr]dz_1\wedge d\overline z_1 \wedge \cdots \wedge dz_n\wedge
d\overline
z_n.$$
Using the relation
\begin{align}
 \partial *\omega= -i^{n+1}\sum_{a,b=1}^n {\partial (Gg^{\overline
ab})\over \partial z_b } dz_1\wedge d\overline z_1 \wedge
\cdots dz_a \wedge \widehat{ d
\overline z_a}\cdots\wedge dz_n\wedge d\overline z_n,
\label{eqn:som}
\end{align}
with an analogous formula for $\overline \partial *\omega$, if
$u$ is a complex-valued  function, we obtain
$$(\overline \partial *\omega)\wedge \partial u = i^{n+1}\sum_{a,b=1}^n
{\partial (Gg^{\overline ab})\over \partial \overline z_a }{\partial u
\over
\partial z_b}\,dz_1\wedge d\overline z_1
\wedge \cdots\wedge dz_n\wedge d\overline z_n.$$
Also
 \begin{align*}
  \partial * \overline \partial u &= i^{n}\sum_{a,b=1}^n [ Gg^{\overline
ba}{\partial ^2 u \over \partial z_a \partial \overline z_b }+{\partial
(Gg^{\overline
ba})\over \partial  z_a }{\partial u \over
\partial \overline z_b}]\,dz_1\wedge d\overline z_1
\wedge \cdots\wedge dz_n\wedge d\overline z_n \\
&=:D_2u+D_1u.
\end{align*}
 Similarly,
$$\overline \partial *  \partial u = D_2u+\overline {D_1\overline u}.$$
Finally, it is
straightforward to verify that
$$(\overline \partial *\omega)\wedge \partial u =i\, \overline
{D_1\overline u}.$$
Using these relations, we obtain
$$(\overline \partial *\omega)\wedge \partial {\partial g \over \partial
\overline t}=i\, \overline {D_1 ({\partial g \over \partial  t})}$$
and
$$( \partial *\omega)\wedge \overline \partial { \partial g \over
\partial \overline t}=-i {D_1 ({\partial g \over \partial \overline
t})}.$$

Now since $g$ is a function, $\delta g = \overline \delta g =0$; thus we
rewrite $\Delta g+cg=0$ as
$$(\overline \delta \partial + \delta \overline \partial)g +cg= -(*
\overline \partial * \partial +* \partial * \overline \partial)g +cg=0.$$
Applying $*$ to this relation and using $**K=K$, we
have
$$-(\overline \partial * \partial + \partial * \overline
\partial)g+*cg=0.$$
This last equation can be written using $D_1$ and $D_2$ as
$$-(D_2g+\overline {D_1g} +D_2g +D_1g)+*cg=0.$$
Hence
\begin{align}
 \label{eqn:d2g}
D_2g={-1\over 2}[D_1g + \overline {D_1g}-*cg].
\end{align}
Equation (\ref{eqn:d2g}) is valid in $D(0)\setminus \{p_0\}$ -- indeed,
in each $D(t)\setminus \{p_0\}$ -- as an equality of $(n,n)$ forms.
Differentiating (\ref{eqn:d2g}) with
respect to $\overline t$, we thus obtain
\begin{align}
 \label{eqn:dd2}
D_2({\partial g \over \partial \overline t})={-1\over 2}[D_1({\partial g
\over \partial \overline t}) +
\overline {D_1({\partial g \over \partial t})}-*c{\partial g \over
 \partial \overline t}]
\end{align}
in $D(t)$.

We use (\ref{eqn:dd2}) in (ii):
$$\overline \partial (*\partial ({\partial g \over \partial \overline
t}))=D_2({\partial g \over \partial \overline t})+\overline
{D_1({\partial g \over \partial
 t})}={-1\over 2}[D_1({\partial g \over \partial \overline t}) -
\overline {D_1({\partial g \over \partial  t})}-*c{\partial g \over
\partial \overline t}]$$
$$={-1\over 2}\bigl[{-1\over i}\partial *\omega\wedge \overline \partial
{ \partial g \over \partial \overline t}-{1\over i}\overline \partial
*\omega\wedge
\partial {\partial g \over \partial \overline t}-*c{\partial g \over
\partial \overline t}\bigr].$$

Inserting parts (i) and (ii) of $\overline \partial f$ back in to
$J = {c_n \over 2^{n-2}} \ \Re\int_{D(0)}\overline \partial f$,
we obtain
$$J={c_n \over 2^{n-2}}\Re\int_{D(0)}\bigl\{\overline \partial
{\partial g \over \partial t}\wedge * \partial
{\partial g \over \partial \overline t}+{1\over 2}{\partial g \over
\partial  t}\bigl[{1\over i}\partial *\omega\wedge \overline \partial {
\partial g \over
\partial \overline t}+{1\over i}\overline \partial *\omega\wedge
\partial {\partial g \over \partial \overline t}+*c{\partial g \over
\partial \overline t}\bigr] \bigr\}.$$
 Hence,
$$J={c_n \over 2^{n-2}}\ \int_{D(0)}\bigl\{\overline \partial
{\partial g \over \partial t}\wedge * \partial
{\partial g \over \partial \overline t}+\Re\, (\,{1\over 2}{\partial g \over
\partial  t}\bigl[{1\over i}\partial *\omega\wedge \overline \partial {
\partial g \over
\partial \overline t} +{1\over i}\overline \partial *\omega\wedge
\partial {\partial g \over \partial \overline t}\bigr]\,) +\, {1\over
2}c|{\partial g \over \partial  t}|^2{\omega^n\over n!}\bigr\}.$$
Thus we obtain
$${\partial ^2\lambda \over \partial t \partial \overline
t}(0)=-c_n\int_{\partial
D(0)}k_2(0,z)\sum_{a,b=1}^n(g^{\overline a b}{\partial g\over \partial
\overline z_a}{\partial g\over
\partial  z_b})d \sigma_z $$
$$\    -{c_n\over 2^{n-2}}  \bigl\{\,||\overline \partial {\partial
g\over \partial
t}||_{D(0)}^2+{1\over 2}||\sqrt c {\partial g\over \partial
t}||_{D(0)}^2+ \frac{ 1}2 \Re
\int_{D(0)}
{\partial g \over
\partial  t}\bigl[\,{1\over i}\partial *\omega\wedge \overline \partial {
\partial g \over
\partial \overline t} +{1\over i}\overline \partial *\omega\wedge
\partial {\partial g \over \partial \overline t}\,\bigr]\,
 \bigr\}$$
which is (\ref{eqn:f2}) for $t=0$. The second order variation formula
(\ref{eqn:f2}) is thus proved. \hfill $\Box$
\begin{remark} \  {\rm Note that
(\ref{eqn:f2}) reduces to (\ref{eqn:ff}) if $M={\mathbb{C}}^n, \ ds^2=|dz|^2$ is
the Euclidean metric, and $c\equiv 0$.}
\end{remark}
\section{Subharmonicity of
 $-\lambda$.} \setcounter{section}{3}\
We impose  the
following condition on the Hemitian metric $ds^2$ on $M$:
\begin{align}
 \partial *\omega =0  \quad  \mbox{  on } \ M.  \label{eqn:psm}
\end{align}

\begin{remark} \label{re:3-1} \  {\rm A  K\"ahler  metric $ds^2$  on
$M$ (i.e., $d\omega=0$) satisfies (\ref{eqn:psm}). Indeed, since  $\displaystyle{
*\omega= \frac{ \omega^{n-1}}{n-1}}$ on $M$ and $\omega$ is real, it
follows that
$$
\partial *\omega = \partial \omega \wedge \omega^{n-2} = \frac{ 1}2 \,
d\omega \wedge \omega^{n-2} =0  \quad  \mbox{  on }  M.
$$ }
\end{remark}

\begin{theorem}\label{thm:kansya}\
 Assume that $ds^2$ satisfies condition (\ref{eqn:psm}).
Then:
\begin{enumerate}
 \item  The Levi scalar curvature $k_2(t,z)$ of $\partial {\cal D}$
 reduces to:
\begin{align}
 &K_2(t,z):= \nonumber
[\sum_{a,b=1}^ng^{\overline a b}{\partial \psi\over \partial \overline
z_a}{\partial
\psi\over \partial  z_b}]^{-3/2} \times\\
& \bigl[{\partial ^2 \psi\over \partial t \partial \overline
t}(\sum_{a,b=1}^ng^{\overline a b}{\partial \psi\over \partial \overline
z_a}{\partial
\psi\over \partial  z_b})-2\Re\{{\partial  \psi\over \partial
t}(\sum_{a,b=1}^ng^{\overline a b}
{\partial \psi\over \partial \overline z_a}{\partial ^2 \psi\over
\partial  z_b \partial
\overline t})\}
+ |{\partial
\psi\over \partial t}|^2(\sum_{a,b=1}^ng^{\overline a b}{\partial ^2
\psi\over \partial \overline
z_a \partial z_b})\bigr] .
\nonumber
\end{align}
 \item
 The second variation formula (\ref{eqn:f2}) of $\lambda(t)$ reduces to:
\begin{align}
 {\partial ^2\lambda \over \partial t \partial \overline
t}(t)&=-c_n\int_{\partial
D(t)}K_2(t,z)\sum_{a,b=1}^n(g^{\overline a b}{\partial g\over \partial
\overline z_a}{\partial g\over
\partial  z_b})d \sigma_z  \nonumber  \\
& \quad \qquad \qquad  - {c_n \over 2^{n-2}} \,  \bigl\{\,||\overline \partial {\partial
g\over \partial
t}||_{D(t)}^2+{1\over 2}||\sqrt c {\partial g\over \partial
t}||_{D(t)}^2 \, \bigr\}. \label{omega0}
\end{align}
\end{enumerate}
\end{theorem}

\noindent {\bf  Proof}.
 In local coordinates $z=(z_1,...,z_n)$ we see from formula (\ref{eqn:som}) that the condition $\partial *\omega=0$ on
 $M$ is equivalent to
\begin{eqnarray}
 I_a:=\sum_{a=1}^n {\partial (Gg^{\overline ab})\over \partial z_b}=0, \
a=1,...,n,  \quad  \mbox{  on } \ M,  \label{eqn:is}
\end{eqnarray}
so that
 the Laplacian $\Delta$ on functions $u$ has the form
\begin{eqnarray}\label{eqn:is2}
 \Delta u= -2
\sum_{a,b=1}^ng^{\overline a b}{\partial ^2
u\over \partial \overline
z_a \partial z_b}.
\end{eqnarray}
Thus 1. and 2. follow from formulas (\ref{eqno(1.1)}) and
 (\ref{eqn:f2}). \hfill $\Box$

\begin{remark}\label{re:3-2}  \ {\rm \begin{enumerate}
 \item  Under the condition in Theorem 3.1
in case $\dim M=n=2$, if  $\frac{\partial^2 \lambda }{\partial t \partial
\overline { t}}(t_0)=0 $, then $\partial D(t_0)$ is Levi flat. If
$n>2$, this conclusion is not necessarily true.
 \item   As mentioned in section 2, $k_2(t,z)$ is  a well-defined function on
 $\partial {\cal D}$ for any Hermitian metric $ds^2$ on $M$; i.e., it is
 independent of the local coordinates $z$. This is not the case for
 $K_2(t,z)$. From (\ref{eqno(1.1)}) and (\ref{eqn:is}) it follows that
 $K_2(t,z)$ is a well-defined function on $\partial {\cal D}$ if and
 only if $ds^2$ satisfies condition (\ref{eqn:psm}).
 \end{enumerate}}
\end{remark}

We now turn to the
question as to when ${\partial ^2 \lambda\over \partial t \partial
\overline t}\leq 0$; i.e.,
the  subharmonicity of $-\lambda(t)$. We begin with the following:

\begin{theorem}\label{thm:psm}\
Assume that $ds^2$ satisfies condition (\ref{eqn:psm}).
 If ${\cal D}$ is pseudoconvex in $B\times M$, then $- \lambda (t)$ is
 subharmonic on $B$.
\end{theorem}

\noindent {\bf Proof}.  From the second variation formula for $\lambda(t)$
in 2. of Theorem \ref{thm:kansya} it suffices
to prove that $K_2(t,z)\ge 0$ on $\partial
{\cal D}.$

Fix $(t_0,z_0)\in
\partial {\cal D}$; we may assume $(t_0, z_0)=(0,0)$.
We choose coordinates $z=(z_1,...,z_n)$ in a
neighborhood $V$ of $0$ so
that $g_{\overline a b}(0)=\delta_{ab}$. Let $\psi(t,z)$ be a defining
function for ${\cal D}$. We take a sufficiently small disk
$B_0\subset B$ centered at $t=0$ so that $\psi(t,z)$ is defined in $B_0
\times V$. Then from Theorem \ref{thm:kansya}  we have
\begin{equation}
\label{k200}
K_2(0,0)= ||\nabla_z \psi||^{-3}{L\psi(0,0)} \\
\end{equation}
$$:=||\nabla_z
\psi||^{-3}\bigl[{\partial ^2 \psi\over \partial t \partial \overline
t}||\nabla_z \psi||^2-2\Re\{{\partial  \psi\over \partial t}\sum_{a=1}^n
{\partial  \psi\over
\partial \overline z_a}{\partial ^2\psi \over \partial \overline t
\partial  z_a}\}+ |{\partial
\psi\over \partial t}|^2\Delta_z\psi \bigr].$$
Here  $\Delta_z \psi  =\sum_{ \alpha=1  }^{ n  }\frac{\partial ^2\psi}{\partial
z_\alpha \partial \overline { z}_\alpha }   $, \  $\nabla_z \psi=({\partial \psi \over\partial z_1}, \cdots , {\partial
\psi\over \partial z_n})$ and the terms on the right hand side are
evaluated at
$(t,z)=(0,0)$.
\  Let
${\cal D}_0:={\cal D}\cap (B_0\times V)$. Then $\psi(t,z)$ is  a defining
function for ${\cal D}_0$ near $\partial {\cal D}\cap {\cal D}_0$;
i.e., $\psi$ is negative on ${\cal D}_0$ and $\nabla_z \psi(t,z)\not = 0$
on $\partial {\cal D}\cap {\cal D}_0$.
 For $j\in \{1,...,n\}$, let
$${\cal D}_{0,j}:={\cal D}_0\cap \{z_1=...=\widehat {z_j} = ...=z_n=0\}.$$
Since ${\cal D}_0$ is pseudoconvex in ${\mathbb{C}}^{n+1}$, it follows that
${\cal D}_{0,j}$ is pseudoconvex in ${\mathbb{C}}^2$ in the variables
$(t,z_j)$;
hence, for $j=1,...,n$,
$${\partial ^2 \psi\over \partial t \partial \overline t}|{\partial
\psi\over \partial z_j}|^2-
2\Re\{{\partial  \psi\over \partial t} {\partial  \psi\over
\partial \overline z_j}{\partial ^2\psi \over \partial \overline t
\partial  z_j}\}+|{\partial
\psi\over \partial t}|^2{\partial ^2\psi \over \partial \overline z_j
\partial  z_j}\geq 0$$
at the point $(0,0)$. Summing up from $j=1,...,n$, we get
$$L\psi(0,0)\geq 0.$$
Using (\ref{k200}) yields the result. \hfill $\Box$

\begin{remark}
 \  {\rm We consider the following condition on the metric $ds^2$ on $M$:
 $$
{\cal W}= {\cal W}_{ds^2}:={1\over i}\overline \partial \partial * \omega
- ||\partial * \omega||^2{\omega^n\over n!} \ge 0  \eqno(W)
$$
as an $(n,n)$ form on $M$. This is a
weaker condition than (\ref{eqn:psm}). We can prove from
 (\ref{eqn:f2}) that, if $ds^2$ satisfies condition
 $(W)$ and if ${\cal D}$ is pseudoconvex in $B\times M$, then
 $-\lambda(t)$ is subharmonic on $B$. We put ${\cal W}=W(z) \frac{
 \omega^n}{n!}$. In local coordinates, $W(z)$ has the form
$$
W(z)=\frac{ 1}G
\sum_{ \alpha,\beta=1  }^{ n  }
\bigl[
\frac{\partial ^2(Gg ^{\overline  \alpha \beta})}
{\partial  \overline {z}_\alpha \partial z_\beta }
 -
\frac{ g_{\alpha \overline { \beta}}} G\
\overline {
\bigl(
\sum_{l=1   }^{n   } \frac{\partial (G g^{\overline  \alpha l})}
{\partial   z_l }
\bigr)  }
\bigl(
\sum_{k=1   }^{n   } \frac{
\partial (G g^{\overline  \beta
 k})}
{\partial  z_k }
\bigr)
\bigr].
$$
Using standard notation from classical differential geometry,
we define complex Christoffel symbols
$$\Gamma_{\lambda \beta}^{\alpha}:=\sum_{\gamma=1}^n g^{\bar \gamma \alpha}{\partial g_{\beta \bar \gamma}\over \partial z_{\lambda}}$$
and complex torsion
$$T_{\lambda \beta}^{\gamma}:=\Gamma_{\lambda \beta}^{\gamma}- \Gamma_{\beta \lambda}^{\gamma}.$$
Summing, we define $T_{\alpha}:= \sum_{\lambda =1}^n T_{\lambda \alpha}^{\lambda}$. Then
$$W(z)=\sum_{\alpha, \beta =1}^n g^{\bar \beta \alpha}\,{\partial T_{\alpha} \over \partial \bar z_{\beta}},$$
which is a type of scalar curvature function on $M$. Thus the condition
 $(W)$ means that this scalar curvature $W(z)$ is nonnegative. As an example,  on the unit ball
$M=\{z\in {\mathbb{C}}^n: \|z\| <1\}$ in ${\mathbb{C}}^n, \
n\ge 2$, take
$ds^2:={|dz|^2\over (1-\|z\|^2)^2}, $
where $|dz|^2$ is the Euclidean metric in ${\mathbb{C}}^n$. Then
 $W(z)=2(n-1)(n-(n-1)\|z\|^2)\ge 2(n-1)>0$ on  $M$.
On the other hand, the special Hopf manifold $\mathbb{H}_n$, $n\ge 2$ (see definition (\ref {eqn:special-hoph}) in section 6) admits the Hermitian metric
$ds^2:= {|dz|^2\over \|z\|^2}$. A calculation shows that
$W(z)= -(n-1)^2<0$ on $\mathbb{H}_n$.
We plan to investigate condition $(W)$ in a future
 paper. }
\end{remark}

J-C. Joo \cite{joo} has a generalization of our second variation formula (\ref{eqn:f2}) to certain
almost complex manifolds.

\section{Rigidity.}\setcounter{section}{4}\  We continue under the same hypotheses: $M$ is
an
$n$-dimensional complex manifold equipped with a Hermitian metric
$ds^2$; ${\cal D}=\cup_{t\in B}(t,D(t))\subset B\times M$ is a
$C^{\infty}$ variation of
domains $D(t)$ in $M$ each containing a fixed point $p_0$ and with
$\partial D(t)$ of class
$C^{\infty}$ for $t\in B$; and $c$ is a nonnegative $C^{\infty}$ function
on $M$.

Throughout this section we will assume that {\it
\begin{enumerate}
 \item [{\rm({\bf  1})}] \ $ds^2$ satisfies condition (\ref{eqn:psm}) on $M$;
 \item [{\rm({\bf  2})}] \ ${\cal D}$ is pseudoconvex in $B\times M$.
\end{enumerate}}

\begin{lemme}[Rigidity]\label{lem:rigid} \
If there exists $t_0\in B$ at which ${\partial ^2 \lambda\over \partial t
\partial \overline t}(t_0)=0$,
then ${\partial g\over \partial t}(t_0,z)\equiv 0$ for $z\in D(t_0)$
provided
at least one of the following conditions hold:
\begin{enumerate}
\item [{i)}] \ $c(z)\not \equiv 0$ on $D(t_0)$;
\item [{ii)}] \ $\partial D(t_0)$ is not Levi flat.
\end{enumerate}
\end{lemme}

\noindent {\bf Proof}.  Since ${\cal D}$ is pseudoconvex in $B\times M$,
Theorems \ref{thm:kansya} and \ref{thm:psm} imply that $k_2(t,z)\geq 0$ on $\partial {\cal D}$ and we
obtain the following
estimate from (\ref{omega0}):
\begin{align}\label{abc}
 {\partial ^2 \lambda\over \partial t \partial \overline t}(t_0)\leq
{-c_n\over 2^{n-1}}\ \{\ 2||\overline{  \partial} {\partial g\over \partial
t}||^2_{D(t_0)}+  ||\sqrt
c{\partial g\over \partial
t}||^2_{D(t_0)}\}.
\end{align}
Thus if the left hand side ${\partial ^2 \lambda\over \partial t \partial
\overline t}(t_0)=0$, then (see (\ref{eqn:f2}))
$
k_2(t_0,z) \equiv 0 $ on $\partial D(t_0),$
and each term on the right in (\ref{abc}) must vanish. Hence
\begin{enumerate}
 \item [{1.}] $ \overline{  \partial}{ \partial g\over \partial
t}(t_0,z)\equiv 0$ on
$D(t_0)$;
\item [{2.}] $\sqrt c{\partial g\over \partial t}(t_0,z)\equiv 0$ on
$D(t_0)$.
\end{enumerate}
In particular, 1. says that ${\partial g\over \partial t}(t_0,z)$ is
holomorphic on $D(t_0)$. First assume i) holds. Since $c$ is
 of class $C^\infty$,
$c(z)\not \equiv 0$ on $D(t_0)$ together with 2. implies that the
holomorphic function ${\partial g\over \partial t}(t_0,z)$ vanishes on an
open set in $D(t_0)$
and hence vanishes identically on $D(t_0)$.

Now assume ii) holds but i) does not; i.e., $c(z)\equiv 0$ on $D(t_0)$.
 We want to show that $ {\partial g\over \partial t}(t_0,\cdot)\equiv 0$
on $
  \partial D(t_0)$. To see this, fix a point $q_0\in \partial D(t_0)$.
Using
local coordinates $z=(z_1, \ldots , z_n)$ in a neighborhood $V$ of $q_0$
where $z=0$ corresponds to $q_0$,  we can assume that
$g_{a \overline {  b}}(0)=\delta_{ab}$. We show that ${\partial g\over
\partial t}(t_0,0)=0$.

Since ${\partial g\over \partial
t}(t_0,z)$ is holomorphic in $D(t_0)$ and of class $C^{\infty}$ on
$\overline{
  D(t_0)} $, we have $ \overline{  \partial}
  {\frac{\partial g }{\partial t} }=0$ on $\partial D(t_0)$. Then
  formula (\ref{eqno(1.5)})
$$k_2(t_0,z)\cdot(\sum_{a,b=1}^ng^{\overline ab}{\partial g\over \partial
\overline z_a}{\partial g\over
\partial  z_b})^{3/2}=-{\partial ^2g\over \partial t \partial
\overline t}(\sum_{a,b=1}^ng^{\overline ab}{\partial g\over \partial
\overline z_a}{\partial g\over
\partial  z_b})$$
$$+2\Re[(\sum_{a,b=1}^ng^{\overline ab}{\partial^2 g\over \partial
z_b\partial \overline t}{\partial g\over
\partial \overline z_a}){\partial g\over
\partial t}]-|{\partial g\over
\partial t}|^2(\sum_{a,b=1}^ng^{\overline ab}{\partial^2 g\over \partial
\overline z_a\partial  z_b})$$
simplifies at the point $(t_0,0)$ to
$$k_2(t_0,0)\cdot ||\nabla g||^3+{\partial ^2g\over \partial t \partial
\overline t}||\nabla g||^2 +|{\partial g\over
\partial t}|^2\sum_{i=1}^n {\partial^2 g\over \partial z_i\partial
\overline z_i}=0.$$
Since $k_2(t_0,z)=0$ for $z\in \partial D(t_0)$,  we
have
\begin{align}
 \label{eqn:(I)}
{\partial ^2g\over \partial t \partial \overline t}||\nabla g||^2
+|{\partial g\over
\partial t}|^2 \sum_{i=1}^n{\partial^2 g\over \partial z_i\partial
\overline z_i}=
0 \quad  \mbox{ at} \ (t_0, 0).
\end{align}
Note that $-g$ is a defining function  for ${\cal D}$. Using
pseudoconvexity of ${\cal D}$ at $(t_0,z), \ z\in \partial D(t_0)$,
and using $\overline{  \partial} \, {\partial g\over \partial t}=0$ on
$\partial D(t_0)$,
it follows  that
$${\partial ^2g\over \partial t \partial \overline t}|{\partial g\over
\partial z_i}|^2+|{\partial g\over
\partial t}|^2 {\partial^2 g\over \partial z_i\partial \overline z_i}
\le 0  \quad \mbox{ on } \ \partial D(t_0)\cap V$$
for each $i=1,...,n$.
Now since $\partial *  \omega \equiv 0$ on $M$,
it follows  from (\ref{eqn:is2}) and $g_{a\overline {
b}}(0)=\delta_{ab}$ that  the Laplacian $\Delta$ on functions $u$ at
$z=0$ has the form
$$
\Delta u(0)= -2 \sum_{i=1}^n{\partial ^2u\over \partial z_i \partial
\overline z_i}(0).
$$
We are assuming that $c\equiv 0$ on $D(t_0)$
and hence on $\partial D(t_0)$; thus
$\Delta g(t_0,z)=0$ for $z\in D(t_0)$ except  at the pole $p_0$. We thus
have
$$\sum_{i=1}^n {\partial^2 g\over \partial z_i\partial \overline
z_i}(t_0,0)=0. $$
Using (\ref{eqn:(I)}), we conclude that ${\partial ^2g\over \partial t
\partial
\overline t}||\nabla g||^2=0$ at $(t_0,0)$. Hence
$$|{\partial g\over
\partial t}|^2 {\partial^2 g\over \partial z_i\partial \overline z_i}\le
0 \ \  \mbox{ at}\ \,  (t_0,0), \ \ \  i=1,...,n.
$$
From the assumption ii), ${\partial^2 g\over \partial z_i\partial
\overline z_i}(t_0,0) >0$ for some $i$, and we conclude that
${\partial g\over \partial t}=0$ at $(t_0,0)$. Thus the holomorphic
function
${\partial g\over \partial t}$ vanishes identically  on $\partial
D(t_0)$, and hence on $D(t_0)$.
\hfill $\Box$

\begin{corollary} [Trivial
variation]\label{cor:trivial-v}
  \ If $\lambda(t)$ is harmonic in $B$ and i) or ii) in Lemma \ref{lem:rigid} holds, then ${\cal D}=B\times D(0)$.
\end{corollary}
\noindent {\bf Proof}.  By  Lemma \ref{lem:rigid},  ${\partial g\over
\partial t}=0$ on ${\cal D}$, i.e., $g(t,z)$ does not depend on $t\in
B$. By analytic continuation it suffices to prove that there exists a
small disk $B_0\subset B$ centered at $t=0$ such that ${\cal D}|_{B_0} =
B_0 \times D(0)$, where ${\cal D}|_{B_0}= \cup _{t\in B_0}(t, D(t))$. If
not, there exists a sequence $\{t_n\}$
such that
$t_n \ne 0$, $t_n \to 0$ as $n\to \infty$, and $D(t_n)\ne D(0)$. Thus we
can find a point $z_n\in
\bigl (\partial D(0) \cap D(t_n)\bigr ) \cup \bigl (D(0) \cap \partial
D(t_n)\bigr )$ for
each $n$. If $z_n\in \partial D(0) \cap D(t_n)$, then, $g(0,z_n)=0$
and $g(t_n, z_n)>0$, which is a contradiction; if  $z_0\in D(0) \cap
\partial D(t_n)$,  $g(0,z_n)>0$
and $g(t_n, z_n)=0$ which is again a contradiction.
\hfill $\Box$

\smallskip
We next consider the following set-up. Let $F:B\times M \to M$ be a
holomorphically varying, one-parameter family of automorphisms of $M$;
i.e., $F(t,z)$ is
holomorphic in $(t,z)$ with ${\partial F\over \partial t}\not = 0$ for
$(t,z)\in B \times M$, and, for each $t\in B$, the mapping $F_t:M\to M$ via $F_t(z):=
F(t,z)$ is an
automorphism of $M$. Then the mapping $T:B\times M \to B\times M$ defined
as
$$T(t,z)=(t,w):=(t,F(t,z))$$
provides a fiber-wise automorphism of $M$; i.e., for each $t\in B$, the
map $w=F(t,z)$ is an
automorphism of $M$.

Let $ds^2$ be a Hermitian metric on $M$ satisfying condition (\ref{eqn:psm})
and let $c(z)\ge 0$ be a $C^\infty$ function on $M$.
 Fix a pseudoconvex domain $D\Subset M$ and let ${\cal D}:=T(B\times D)$.
This is
a variation of pseudoconvex domains $D(t)=F(t,D)$ in the image space
$B\times M$ of
$T$. Assume
there exists a common point $w_0$ in each domain $D(t), \ t\in B$.
 Let $g(t,w)$ and $\lambda(t)$ denote
the $c$-Green
function and $c$-Robin constant of $(D(t),w_0)$ for the Laplacian $\Delta$
associated to $ds^2$. We obtain the following fundamental result, utilizing the rigidity lemma, Lemma \ref{lem:rigid}.
\begin{corollary} \label{lem:rigid-2}
If
$$
{\partial ^2 \lambda\over \partial t \partial \overline{t}} (t_0)=0
$$
at some $t_0\in B$, and either $\{z\in D(t_0) : c(z)\ne 0\}\ne
 \emptyset$ or
$\partial D(t_0)$ (or equivalently $\partial D$) is not Levi flat,
then
there exists a nonvanishing holomorphic vector field $\Theta$ on $M$
 which is tangential on $\partial D$;\,i.e., the entire integral curve $I(z_0)
$ associated to $\Theta$
  for any initial point $z_0\in \partial D $ lies on $\partial D$.
\end{corollary}

\noindent {\bf Proof}. We may take $t_0=0$. For fixed $t$, we let
$\phi(t,w):=F_t^{-1}(w)$ so that $z=\phi(t,w)$ if $w=F(t,z)$. Let
$z_0(t):=\phi(t,w_0)$ for $t\in B$; thus, in the $(t,z)$-coordinates in
$B\times M$
we have a constant variation of domains $t\to D(t)\equiv D$ with a
varying family of poles $z_0(t)$; letting $c(t,z):= c(F(t,z))$ and
$ds_t^2(z):=F(t,\cdot)^*(ds^2(w))$ (the pull-back of $ds^2$ by $F(t,\cdot)$)-- note this is a varying family of
Hermitian metrics on the fixed domain $D$ -- we have the $c(t,z)$-Green function ${\cal G}(t,z)$
for $(D,z_0(t))$ for the Laplacian
 $\tilde \Delta_t$ associated to $ds_t^2(z)$.
 Indeed,
$${\cal G}(t,z):=g(t,F(t,z))$$
satisfies $\Delta g(t,w)= \tilde \Delta_t g(t,F(t,z))$
  by the invariance of the Hodge-* operator under holomorphic mappings.

The assumption that ${\partial^2 \lambda\over \partial t \partial
\overline{ t}}(0)=0$ implies ${\partial g\over \partial t}(0,w)\equiv
0$ on $D(0)$ by the rigidity lemma.
For $p_0\in M$, we let $q_0=F(0,p_0)$ and we
choose local coordinates $z$ on a coordinate neighborhood $U$ of $p_0$
and
$w$ on a coordinate neighborhood $V$ of $q_0$ and write
$$z=\phi(t,w)=(z_1,...,z_n)=(\phi_1(t,w),...,\phi_n(t,w)).$$
Using these coordinates and the relation ${\cal G}(t,\phi(t,w))=g(t,w)$,
and noting that each $\phi_j, \ j=1,...,n$ is holomorphic, the condition
${\partial g\over
\partial t}(0,w)\equiv 0$ on $D(0)$ becomes
$${\partial {\cal G}\over \partial t}(0,z)+\sum_{k=1}^n{\partial {\cal
G}\over \partial z_k}(0,z){\partial \phi_k \over \partial t}(0,F(0,z))=0$$
for $z\in D$.
Since $g(t,w)$ is of class $C^{\infty}$ on $\overline{  {\cal D}}$,
${\cal  G}(t, \phi(t,w))$ is of class $C^{\infty}$ on $B \times \overline
D$. Thus, the above
equality is valid for $z\in \overline{  D}$.
Since ${\cal G}(t,z)=0$ for $z\in \partial D$ for all $t\in B$,
 we have  ${\partial
{\cal G}\over
\partial t}(0,z)=0$ for $z\in \partial D$, and hence
\begin{align}
 \sum_{k=1}^n{\partial {\cal G}\over \partial z_k}(0,z){\partial \phi_k
\over \partial t}(0,F(0,z))=0,  \quad  \mbox{  for $z\in \partial D$.}
\label{eqn:tyokou}
\end{align}

The vector field defined by
\begin{align}
 \Theta(z):=\sum_{k=1}^n {\partial \phi_k \over \partial
t}(0,F(0,z)){\partial \over \partial z_k}
\label{thetaz}
\end{align}
in each coordinate
neighborhood $U$ with coordinates $z$ (and coordinate neighborhood $V$ for $w=\phi(t,z)$, $z\in U$)
is easily seen
to be a holomorphic vector field on all of $M$. Furthermore, since
${\partial F\over \partial t}\not = 0$ on $B\times M$ and $\phi^{-1}=F$,
we have $\Theta(z)\not = 0$ on $M$.

We shall verify the property that
the entire integral curve $I(z_0)$ associated to $\Theta$
for any initial point $z_0\in \partial D$ lies on $\partial D$, i.e.,
writing $\Theta(z)=\sum_{k=1}^n \zeta_k(z){\partial \over \partial
z_k}$,  if we consider the ordinary
differential equation
\begin{align}\label{eqn:moku}
 {dZ_j\over dt}=\zeta_j(Z_1, \ldots , Z_n),  \quad j=1, \ldots , n,
\end{align}
on $M$, then the solution curve $Z(t)=(Z_1(t), \ldots , Z_n(t)),\  t\in
{\mathbb{C}}_t$ with
initial value $Z(0)=z_0$ always lies on $\partial D$.

To see this, we fix
such a $z_0$ which we take to be $0$. Since the $c(0,z)$-Green function
${\cal G}(0,z)$ vanishes on
$\partial D$, the vector
$$({\partial {\cal G}\over \partial z_1}(0,z),...,{\partial {\cal G}\over
\partial z_n}(0,z))$$
is a nonvanishing complex normal vector to $\partial D$ at $z$. By
(\ref{eqn:tyokou}) we have
$\sum_{k=1}^n \zeta_k(z){\partial {\cal G}(0,z)\over \partial z_k}=0
$ on $\partial D$. Thus, at each point $z\in \partial D$, the vector
$(\zeta_1(z), \ldots ,  \zeta_n(z))$ is a holomorphic tangent vector.
This observation gives the intuitive reason for the property
of the solution $Z(t)$ of the system of equations (\ref{eqn:moku}) in the
statement of the corollary.

To verify this property of $Z(t)$, we set up local coordinates $z=(z_1, \ldots , z_n)$ near the point $z_0$ where $z_k=x_k+iy_k$ so that
$$z=(x,y)=(x_1,y_1, \ldots , x_n,y_n) \in {\mathbb{R}}^{2n}.$$
We let $z_0$ correspond to $z=0$ and define $\Psi(z)=\Psi(x,y):= {\cal
G}(0,z)$. Then $\Psi$ is a $C^\infty$ real-valued function near
$\partial D$ with $
\Psi=0 $  and $ \frac{\partial
\Psi }{\partial n_{z} }<0$  on $ \partial D,$
where $n_{z}$ is the unit outer normal vector to $\partial D$ at
       $(x,y)$; i.e.,
$$
{\rm  Grad}\,\Psi_{(x,y)}=(\frac{\partial \Psi}{\partial x_1},
\frac{\partial \Psi}{\partial y_1}, \ldots ,
\frac{\partial \Psi}{\partial x_n}, \frac{\partial \Psi}{\partial y_n})
$$
is a normal vector to $\partial D$ at $(x,y)$.

The holomorphic vector field
$$
\Theta: = \sum_{k=1   }^{ n   } \zeta_k (z)\frac{\partial }{\partial z_k}
 \quad  \mbox{  on} \ M
$$
satisfies
$$
\sum_{ k=1  }^{ n  }\zeta_k(z) \frac{\partial \Psi}{\partial z_k}(z)=0
     \ \quad     \mbox{  on } \ \partial D.
$$
We set
$$
\zeta_k(z):=a_k(x,y)+ib_k(x,y), \ \ \ k=1,\ldots , n,
$$
where $a_k,\ b_k$ are real-valued on $M$. Since $\Psi$ is real-valued, it follows that
$$
\sum_{ k=1  }^{n   } (a_k \frac{\partial \Psi }{\partial x_k}+ b_k
\frac{\partial \Psi}{\partial y_k }) +i\,
\sum_{ k=1  }^{n   } (b_k
\frac{\partial \Psi }{\partial x_k }- a_k
\frac{\partial \Psi}{\partial y_k }) =0  \quad    \mbox{  on } \ \partial D.
$$
Consequently, if we set
\begin{align*}
\Theta_1&:=\Re \, \Theta = \sum_{ k=1  }^{ n  }(a_k \frac{\partial }{\partial x_k}+
b_k \frac{\partial }{\partial y_k})\\
\Theta_2& :=\Im \, \Theta= \sum_{ k=1  }^{ n  }(b_k \frac{\partial }{\partial x_k}
- a_k \frac{\partial }{\partial y_k}),
\end{align*}
then $\Theta_1,\ \Theta_2$ are
real vector fields on $M$ which are tangential on $\partial D$.

\smallskip
We consider the integral curve $Z(t)=(Z_1(t), \ldots, Z_n(t)), \  t\in {\mathbb{C}}$, with initial value $Z(0)=z_0$ as mentioned above, so that
$$
\frac{ dZ_k(t)}{dt}= \zeta_k(Z(t)),  \quad  t\in {\mathbb{C}}, \ \ \ k=1,\ldots ,n.
$$
We write $t=t_1+it_2$ where $(t_1,t_2)\in \mathbb{R}^2$;
 $Z_k(t)=x_k(t_1,t_2)+i y_k(t_1,t_2)$ and
$$
x(t_1,t_2)=(x_k(t_1,t_2))_{k=1,\ldots ,n}, \quad
y(t_1,t_2)=(y_k(t_1,t_2))_{k=1,\ldots ,n}
$$
for $(t_1,t_2)\in \mathbb R^2$, so that $Z(t)=x(t_1,t_2)+iy(t_1,t_2)$.
Since $Z_k(t)$ is holomorphic, we
have $\frac{
dZ_k(t)}{dt}= \frac{\partial Z_k(t)}{\partial t_1} $ and hence
$$
(1) \quad   \left\{
\begin{array}{llll}
\displaystyle{  \frac{\partial x_k(t_1,t_2)}{\partial t_1}}
 &=a_k(x(t_1,t_2), y(t_1,t_2)) \\
\displaystyle{  \frac{\partial y_k(t_1,t_2)}{\partial t_1} } &=
b_k(x(t_1,t_2),y(t_1,t_2))
\end{array}
\right.
$$
for $(t_1,t_2)\in {\mathbb{R}}^2$. Analogously,
$\frac{
dZ_k(t)}{dt}= \frac{ 1}{i} \frac{\partial Z_k(t)}{\partial t_2} $
implies
$$
(2)  \quad  \left\{
\begin{array}{llll}
\displaystyle{  \frac{\partial x_k(t_1,t_2)}{\partial t_2}}
 &=-b_k(x(t_1,t_2),y(t_1,t_2)) \\
 \displaystyle{ \frac{\partial y_k(t_1,t_2)}{\partial t_2}}
 &=a_k(x(t_1,t_2),y(t_1,t_2))
\end{array}
\right.
$$
for $(t_1,t_2)\in {\mathbb{R}}^2$. If we set $t_2=0$ in (1), then $
(x(t_1,0), y(t_1,0))$
satisfies
$$
 \left\{
\begin{array}{llll}
\displaystyle{  \frac{d x_k(t_1,0)}{d t_1}} &=a_k(x(t_1,0),y(t_1,0)) \\
 \displaystyle{ \frac{d y_k(t_1,0)}{d t_1}} &=b_k(x(t_1,0),y(t_1,0))
\end{array}
\right.
$$
for $t_1\in (-\infty, \infty)$ and $(x(0,0),y(0,0))=z_0$. Thus
$(x(t_1,0), y(t_1,0))$ is an integral curve for $\Theta_1$ with
initial vaue $z_0$. Since $\Theta_1$ is
tangential on $\partial D$, it follows that
$$
(x(t_1,0), y(t_1,0))\in \partial D, \qquad t_1\in (-\infty, \infty).
$$
Next, considering equation (2) with fixed $t_1\in (-\infty, +\infty)$
and variable $t_2\in (-\infty, \infty)$, we see that
$(x(t_1,t_2), y(t_1,t_2)) $
 is an integral curve for $\Theta_2$ with
initial vaue $(x(t_1,0),y(t_1,0))\in \partial D$.
 Since
$\Theta_2$ is tangential on $\partial D$, we have
$$
(x(t_1,t_2), y(t_1,t_2)) \in \partial D, \qquad t_2\in (-\infty, \infty).
$$
Hence
$$
Z(t)=x(t_1,t_2)+iy(t_1,t_2)\in \partial D,  \quad  t\in {\mathbb{C}},
$$
as required.  \hfill $\Box$

\smallskip
We remark that the hypothesis of the corollary  is satisfied if, for example,
$c>0$ in $M$.

\begin{remark}\label{re:4-1} \ {\rm \begin{enumerate}
\item By an elementary topological argument, the corollary immediately
implies that
the integral curve $I(z_0)$ of the vector field $\Theta$
with initial value $z_0\in D$ ($\overline{  D}^c$) always lies in $D$
($\overline{  D}^c$). However, this does not necessarily imply that
$I(z_0)\Subset D$
( $\overline{  D}^c)$.
\item If $\partial D$ contains at least one strictly
pseudoconvex point, then ${\partial ^2 \lambda\over \partial t \partial
\overline{  t}}(t)>0 $ on $B$. \
\item In the case where dim$M=n=2$, ${\partial ^2 \lambda\over \partial t
\partial
\overline{  t}}(t_0)=0$ for some  $t_0\in B$ implies that $\partial D$ is
Levi flat.
\end{enumerate}}
\end{remark}

\noindent Corollary \ref{lem:rigid-2} is fundamental for the
applications in this paper. We now prove a partial converse.
First we remark that for any fixed $t\in B$,
\begin{align}
 \label{eqn:bigtheta}
\Theta (t,z):= \sum_{ k=1  }^{ n  } \frac{\partial \phi _k }{\partial t
} (t, F(t,z)) \frac{\partial }{\partial z_k }
\end{align}
is a holomorphic vector field on all of $M$. Thus
$$
\mathbf{\Theta} : t\in B \to \Theta (t,z)
$$
is a variation of holomorphic vector fields on $M$. We have
the following result.

\begin{corollary}\label{cor:converse} \  Assume that either $\{z\in M:c
(z)\ne
 0\}$ is a dense, open set in $M$ or that $\partial D\subset M$ is not Levi flat. Then
$\frac{\partial^2 \lambda}{\partial t \partial \overline { t}
 }(t_0)=0$ for some $t_0\in B$ if and only if $\Theta (t_0,z)$ is a tangential vector field on  $\partial D$.
 \end{corollary}

\noindent {\bf Proof}.  We verified the sufficiency in the proof of  the
previous corollary. We verify the necessity. From the relation ${\cal G}(t,\phi(t,w))=g(t,w)$, for any fixed $t\in B$,
$$
\frac{\partial  g}{\partial t}(t,w)=
 {\partial {\cal G}\over \partial t}(t,z)+\sum_{k=1}^n{\partial {\cal
G}\over \partial z_k}(t,z){\partial \phi_k \over \partial t}(t,F(t,z))
$$
for $z\in \overline { D}$ and $w=F(t,z)$. Since ${\cal G}(t,z)=0$ on $\partial D$ for all $t\in B$,
 it follows that
$$
\frac{\partial  g}{\partial t}(t,w)=
\sum_{k=1}^n{\partial {\cal
G}\over \partial z_k}(t,z){\partial \phi_k \over \partial t}(t,F(t,z)),
  \ \   \  \mbox{ for $z\in \partial D$.}
$$
Suppose that for some $t_0\in B$,
$ \Theta (t_0, z)$ is tangential on $\partial D$. Since $-{\cal
 G}(t_0,z)$ is a defining function for $\partial D$, this  is equivalent to
$$
\sum_{k=1}^n{\partial {\cal
G}\over \partial z_k}(t_0,z){\partial \phi_k \over \partial
t}(t_0,F(t_0,z))=0,
  \ \   \  \mbox{ for $z\in \partial D$,}
$$
and hence
$$
\frac{\partial g}{\partial t}(t_0, w)=0,   \ \   \   \mbox{ for $w\in
 \partial D(t_0)$.}
$$
Since $\frac{\partial g }{\partial t}
(t_0, w)$ is a continuous function on $\overline { D(t_0)} $ and is
 $c$-harmonic on all of $D(t_0)$, it follows from the uniqueness theorem that $\frac{\partial g }{\partial t} (t_0, w)\equiv 0$ on $D(t_0)$.
From our second variation formula  (\ref{omega0}) we have
\begin{align*}
 &{\partial ^2\lambda \over \partial t \partial \overline
t}(t_0)=-c_n\int_{\partial
D(t_0)}K_2(t_0,w)\sum_{a,b=1}^n(g^{\overline a b}{\partial g\over \partial
\overline w_a}{\partial g\over
\partial  w_b})d \sigma_w \nonumber  \\
& \quad - {c_n \over 2^{n-2}}  \bigl\{||\overline \partial {\partial
g\over \partial
t}||_{D(t_0)}^2+{1\over 2}||\sqrt c {\partial g\over \partial
t}||_{D(t_0)}^2  \bigr\}
\\
& \quad  \qquad \ \ \ =-c_n\int_{\partial
D(t_0)}K_2(t_0,w)\sum_{a,b=1}^n(g^{\overline a b}{\partial g\over \partial
\overline w_a}{\partial g\over
\partial  w_b})d \sigma_w.
 \end{align*}
Thus to finish the proof it suffices to show that
$K_2(t_0, w)=0$ on $\partial D(t_0)$.

\smallskip
Recall that ${\cal D}=T(B\times D)$ where $T(t,z)=(t,w)=(t,
F(t,z))$ and we write $z=\phi(t, w)$. Let $\psi(z)$ be a
defining function for $\partial D$ and set
\begin{eqnarray*}
 \begin{array}{lll}
  & \varphi (t,w):=\psi(\phi (t,w))   \quad  &\mbox{on} \ {\cal D}; \
 \hbox{i.e.}, \vspace{2mm} \\
& \psi(z)=
\varphi (t,F(t,z))  \quad &\mbox{on } \ B \times D.
\end{array}\end{eqnarray*}
Then $\varphi (t,w)$ is a defining function for $\partial {\cal D}$.
From the definition  of $K_2(t,w)$, we can choose local coordinates
$w$ so that
\begin{align*}
 K_2(t_0,w)||\nabla_w \varphi ||^{3}=
{\partial ^2 \varphi \over \partial t
\partial \overline
t}||\nabla_w \varphi ||^2-2\Re\{{\partial  \varphi \over \partial t}\sum_{a=1}^n
{\partial  \varphi \over
\partial \overline w_a}{\partial ^2\varphi  \over \partial \overline t
\partial  w_a}\}+ |{\partial
\varphi \over \partial t}|^2\Delta_w\varphi .
\end{align*}
Writing
$$
\phi(t,w)=(\phi_1(t,w), \ldots , \phi_n(t,w)),  \quad  \mbox{  for $(t,w)\in {\cal D}$,}
$$
we have
\begin{align*}
\frac{\partial \varphi }{\partial t}(t,w)=
\sum_{ i=1  }^{n   } \frac{\partial \psi}{\partial z_i}(t,\phi (t,w))
\frac{\partial  \phi_i}{\partial t}(t,w)  \quad  \mbox{  near } \ \partial
{\cal D}.
\end{align*}
Since $\phi(t,w)$ is holomorphic in $(t,w)$, it follows that
\begin{align*}
 \frac{\partial^2 \varphi }{\partial \overline { t} \partial t}=
\sum_{ i,j=1  }^{n   } \frac{\partial^2 \psi}{\partial z_i  \partial
\overline { z_j}}(t,\phi (t,w)) \overline { (\frac{\partial
\phi_j}{\partial t } )}\frac{\partial  \phi_i}{\partial t}  \quad  \mbox{  near } \ \partial
{\cal D}.
\end{align*}
By assumption, $\Theta (t_0,z )$ is tangential on $\partial D$ so that
\begin{align}
 \label{eqnn2}
 \sum_{ i=1  }^{n   } \frac{\partial \phi_i}{\partial t}(t_0, F(t_0, z))
\frac{\partial \psi}{\partial z_i}(z)=0  \quad  \mbox{  on } \ \partial
D.
\end{align}
This implies
$$
\frac{\partial \varphi }{\partial t}(t_0, F(t_0, z))=0  \quad  \mbox{
on } \ \partial D;
$$
i.e., $\frac{\partial \varphi }{\partial t}(t_0,w)=0 $ on
$\partial D(t_0)$.  Hence
\begin{align}
 \label{eqnn3}
K_2(t_0,w)\|\nabla _w \varphi\|= \frac{\partial^2 \varphi}{\partial t \partial
\overline { t}}(t_0, w)   \quad  \mbox{  on } \ \partial D(t_0).
\end{align}
Again using (\ref{eqnn2}) and the fact that $\psi=0$ on $\partial D$, we have
$$
\sum_{ i=1  }^{n   } \frac{\partial \phi_i }{\partial t}(t_0, F(t_0, z)) \frac{\partial \psi
}{\partial z_i } (z)=f(z)\psi(z)
\quad  \mbox{near } \ \partial D,
$$
where $f(z)$ is a $C^\infty$ function near $\partial D$.
 Since $F(t_0,z)$ is holomorphic in $z$ and $\psi =0$ on $\partial D$, if we differentiate both
sides of the previous equation with respect to $\overline z_j$ we obtain
$$
\sum_{ i=1  }^{ n  } \frac{\partial \phi_i }{\partial t}(t_0, F(t_0, z)) \frac{\partial^2 \psi
}{\partial z_i \partial \overline { z}_j} (z)=f(z)
\frac{\partial
\psi }{\partial \overline { z}_j} (z) \quad  \mbox{  on } \ \partial D.
$$
Multiplying by $\overline { (\frac{\partial \phi_j }{\partial  {
t}})(t_0, F(t_0,z))}$
 and adding for $j=1,\ldots , n$ gives
$$
\sum_{i, j=1  }^{ n  } \frac{\partial \phi_i }{\partial t}
\overline { (\frac{\partial \phi_j }{\partial { t}})}
 \frac{\partial^2 \psi
}{\partial z_i \partial \overline { z}_j} =f \sum_{ j=1  }^{n   }
\frac{\partial
\psi }{\partial \overline { z}_j} (z)
\overline { (\frac{\partial \phi_j }{\partial { t}})}   \quad  \mbox{
on } \ \partial D.
$$
Since $\psi$ is real-valued, formula (\ref{eqnn2}) implies that the right-hand side of this equation vanishes on $\partial D$; thus (\ref{eqnn3})
implies that $K_2(t_0, w)=0$ on $\partial D(t_0)$. \hfill $\Box$

\section{Complex Lie Groups.}\setcounter{section}{5}
We apply Corollary
\ref{lem:rigid-2} of the
previous section to the study of complex Lie
groups. Let
$M$ be a complex Lie group of complex dimension
$n$  with identity
$e$. We always take $M$ to be connected.
By a theorem of Kazama, Kim and Oh ~\cite{K2}, there always exists a
K\"ahler metric on $M$; thus using Remark \ref{re:3-1}, we conclude that $M$
is equipped with a Hermitian metric $ds^2$ satisfying condition
(\ref{eqn:psm}). We fix such a Hermitian metric and a strictly positive
$C^{\infty}$ function $c=c(z)$ on $M$ throughout this section.

We recall some general properties of finite-dimensional complex Lie groups. Let $M$ be a complex Lie group of complex
dimension $n$. We let $\mathfrak X $ denote the set of all left-invariant
holomorphic vector fields $X$ on $M$ and we write $\exp tX$ for the integral
curve of $X$ with initial value  $e$. We note that for any $\zeta\in M$
 the integral curve of $X $ with
 initial value $\zeta$ is given by $\zeta \exp tX$. Moreover, for $s,t\in {\mathbb{C}}$, we have the relations
$$(\exp tX)^{-1}=\exp (-tX),   \quad (\exp sX)
(\exp tX)=\exp (s+t)X.
$$
For $z\in M$, we let $T_zM$ denote the complex tangent space of
 $M$ at $z$ and we define
$$
X(z):= [\frac{d (z\exp tX)}{dt } ] |_{t=0}\in T_ zM.
$$
We identify $\mathfrak X $ with the Lie algebra $\mathfrak g$ for $M$ with the usual Lie bracket $[X,Y]$ for vector fields $X,Y\in \mathfrak X $. Note that $\dim \mathfrak g =n$. Finally, by the connectedness of $M$, given $g\in M$, there exists $t_i\in {\mathbb{C}}$ and $X_i\in
\mathfrak g, \ i=1,\ldots ,N=N(g)$ such that $g= {\textstyle \prod_{i=1}^N \exp t_iX_i}$; i.e.,

\begin{align}
\label{eqn:liesub} {\textstyle  M=\{\prod_{i=1}^\nu \exp t_iX_i \ \in M \ : \ \nu\in {\boldsymbol{Z}}^+, \ t_i\in {\mathbb{C}}, \ X_i\in
 \mathfrak g \}}.\end{align}

Given a Lie subalgebra $\mathfrak g _0$ of $\mathfrak g $ with $\dim
 \mathfrak g _0= d< n$, there is a corresponding connected complex Lie subgroup $H$,
 called the {\it integral manifold} of $\mathfrak g _0$, with the following properties:
 \begin{enumerate}
\item There is a holomorphic isomorphism  from a neighborhood  $V$ of
      $e$ in $H$ onto a neighborhood  $U$ of $0$ in $\mathfrak g_0 $
      whose inverse is given by the exponential map: $X \in U \mapsto
      g=\exp X \in V$.
\item  $M $ is {\it   foliated}  by cosets $\{zH : z\in M\}$ of $H$; i.e., given $z_0\in M$, there exist local
coordinates $\varphi : V \mapsto \Delta={\textstyle  \prod_{i=1}^n\{|w_i|<1\}\subset
{\mathbb{C}}_w^n}$ defined on a neighborhood  $V$ of $z_0$ in $M$ with $\varphi
(z_0)=0$ such that for each $|b_j|<1, \ j=d+1,\ldots ,n$, the $d$-dimensional
polydisk  $\{w\in  \Delta: w_j=b_j, \  j=d+1, \ldots ,n\}$ corresponds to
the connected component of $ (zH)\cap V$ where $z=\varphi^{-1}
(0,\ldots ,0, b_{d+1}, \ldots , b_n)$. We call $(\varphi; V,\Delta)$
      {\it $\mathfrak g _0$-local coordinates} at $z_0$ in $M$.
\item If $H $
is closed in $M$;\,i.e., if $H$ is a
 $d$-dimensional, non-singular, complex analytic set in $M$, the
      quotient space
$M/ H$ is an $(n-d)$-dimensional complex manifold with canonical projection
\begin{align*}
\pi:\,z\in M \mapsto z H \in  M/H.
\end{align*}
\end{enumerate}

We remark that $H$ is obtained as in (\ref{eqn:liesub}) with $\mathfrak g _0$ replacing $\mathfrak g$. Note that in property 2.,
regardless of the size of the neighborhood  $V$ of $z_0$,
we may have $(zH \cap V)\cap (z'H\cap V)\ne \emptyset$  for
      some $z'=\varphi ^{-1}(0,\ldots ,0, b'_{d+1}, \ldots , b'_n)\ne z$.
If $H$ is closed in $M$, however, we can always find a
 neighborhood  $V$ of $z_0$ with $(zH \cap V)\cap (z'H\cap V)=
 \emptyset$ if $z\ne z'$.

These three items concerning the subalgebra $\mathfrak g_0$ of
$\mathfrak g $ form part of a holomorphic {\it Frobenius theorem} for
complex Lie subgroups and subalgebras. Unless otherwise noted, all of our Lie subgroups will be complex Lie subgroups.

\smallskip
Let $D \Subset M$
be a pseudoconvex domain in $M$ with $C^{\infty }$ boundary. Fix $X \in
\mathfrak X \setminus \{0\}$;  $\zeta
\in D$ and
$$B=\{t\in {\mathbb{C}}: |t| < \rho\}$$
with $\rho >0$ sufficiently small so that $\zeta \exp tX \in D$ for $t\in
B$. Given $w\in M$, let $R_w:M\to M$ via right multiplication by $w$:
$R_w(z):=zw$.
We consider the holomorphic map $T:B\times M\to B\times M$ defined as
$$T(t,z)=(t,w):=(t,F(t,z)) \ \ \ \hbox{where} \ F(t,z):=z(\zeta \exp
tX)^{-1}.$$
As in Corollary \ref{lem:rigid-2}, we write $z=\phi(t,w)$ if $w=F(t,z)=z(\zeta \exp
tX)^{-1}$ so that $\phi(t,w)=w(\zeta \exp tX)$. Let ${\cal D}=T(B\times D)$.
Fixing a value of $t\in B$, we will write $D(t):= F(t,D)=D\cdot (\zeta \exp
tX)^{-1}$.
 Note for $t\in B$
fixed,
\begin{align}
 F(t,\cdot)=R_{(\zeta \exp tX)^{-1}}\in \hbox{Aut}(M). \label{eqn:righttras}
\end{align}
Furthermore, for each $t\in B$, we have $e\in D(t)$ since, by hypothesis,
$\zeta
\exp tX\in D$. Thus we can construct the $c$-Robin constant $\lambda (t)$
for $(D(t),e)$; and we have the following.
\begin{lemme}
 \label{thm:lie} \ Suppose
$${\partial ^2 \lambda\over \partial t \partial \overline t}(t_0)=0$$
for some $t_0\in B$. Then the integral curve
$z\,\exp tX, \ t\in{\mathbb{C}}$,
with initial value $z$, of the holomorphic vector field
$X$ satisfies, for $A=D, \ \partial D$, or $\overline D^c$,
\begin{enumerate}
 \item [1.]\ $z\in A$ implies \ $z\,\exp tX\in A$ for all $t\in {\mathbb{C}}$;

\item [2.]  \ $A\cdot z^{-1}=A \cdot (z\,\exp  tX)^{-1}$ for all $t\in {\mathbb{C}}$
       and all $z\in M$.
\end{enumerate}
\end{lemme}
\noindent {\bf Proof}. From the proof of Corollary \ref{lem:rigid-2}
 the nonvanishing holomorphic vector field $\Theta=\Theta(t_0,\cdot)$ on $M$ defined in
 (\ref{thetaz}) has the property that the
entire integral
curve $I(z_0)$ associated to $\Theta$
for any initial point $z_0\in \partial D$ lies on $\partial D$. We claim
that $X=\Theta$. To verify the claim, write,
in local coordinates,
$$X=\sum_{k=1}^n\eta_k {\partial  \over \partial  z_k}$$
and
$$\Theta = \sum_{k=1}^n\zeta_k {\partial \over \partial z_k}.$$
By left-invariance of $X$, for $z\in M$,
\begin{align}\label{eqn:lefti}
 [{d (z\exp tX)_k\over d t}]|_{t=t_0}=\eta_k(z \exp t_0X),  \quad \
k=1,...,n.
\end{align}
On the other hand,  by (\ref{thetaz}),
$$\zeta_k(z) = { \partial  \phi_k \over \partial  t}(t_0,F(t_0,z))
=[{
\partial
\phi_k (t,w)\over \partial  t}]|_{t=t_0, w= z(\zeta \exp t_0 X) ^{-1}  }.$$
But $z=\phi(t,w)= w(\zeta \exp tX)$; in particular,
$$z_k=\phi_k(t,w)= (w\zeta \exp tX)_k, \quad  \ k=1,\ldots , n. $$
Hence, by (\ref{eqn:lefti})
\begin{align*}
 \zeta_k(z) &=
\bigl[{
\partial
\phi_k (t,w)\over \partial  t}\bigr]|_{t=t_0, w= z(\zeta t_0 X) ^{-1}  }\\
&=\bigl[\eta_k(w\zeta \exp t_0X)\bigr]|_{w= z(\zeta \exp t_0 X) ^{-1}}\\
&= \eta _k(z),
\end{align*}
which proves $X=\Theta$.

By Corollary \ref{lem:rigid-2}, this proves 1. for the case $A=\partial
D$.
Using this result, as mentioned in 1. of
Remark \ref{re:4-1}, \, 1. holds for $A=D$ and $A=\overline D^c$.

To prove 2., we first observe that
$$z_1z_2^{-1}=[z_1(\exp tX)][z_2(\exp tX)]^{-1},  \quad  t\in \mathbb C.$$
Take $A=\partial D$ and let $z\in M$ be fixed. For $z_1\in \partial D$, by 1., $z_1(\exp tX)\in
\partial D$ for all $t\in {\mathbb{C}}$. Thus
$$z_1z^{-1}=[z_1(\exp tX)][z(\exp tX)]^{-1}\in (\partial D)
\cdot[z (\exp tX)]^{-1}$$
so that $ (\partial D)\cdot{z^{-1}} \subset (\partial D)\cdot [z(\exp tX)]^{-1}$.
For the reverse inclusion, take
$z_1\in \partial D$; then $z_1(\exp (-tX))\in \partial D$ for all $t\in
{\mathbb{C}}$ by 1.; hence
$$z_1[z(\exp tX)]^{-1}=z_1(\exp (-tX)z^{-1})\in (\partial D)\cdot z^{-1}$$
so that $ (\partial D)[z(\exp tX)]^{-1}\subset (\partial D)\cdot z^{-1}
$.
Similar arguments show that 2. holds for $A=D$ and $A=\overline
D^c$.
 \hfill $\Box$

\smallskip
We make another observation based on Lemma \ref{thm:lie}.
From the proof, we see by using (\ref{eqn:righttras}) that the holomorphic vector field $\Theta(t, z)$
on $M$ defined in (\ref{eqn:bigtheta}) is independent of
$t\in B$; indeed, it coincides with $X$. Suppose that
$\frac{\partial ^2\lambda}{\partial t \partial
\overline { t  }} (t_0)=0 $ for some $t_0\in B$.
From the necessity condition in Corollary
\ref{cor:converse}, $X=\Theta(t_0,z)$ is thus a tangential vector field
 on $\partial D$. Fix $t\in B$. Since $\Theta(t,z)=X$, the
 suffciency condition of the corollary implies that $\frac{\partial ^2\lambda}{\partial t \partial
\overline { t  }} (t)=0 $. From Theorem
\ref{cor:trivial-v} we conclude that $D(t)=D(0)$ for $t\in B$. Thus,
$D( \zeta \exp
tX) ^{-1}  = D \zeta   ^{-1}  $ for $t\in B$, and, aposteriori, $D( \zeta \exp
tX) ^{-1}  = D \zeta   ^{-1}  $ for all $t\in
 \mathbb C$.

\smallskip
We continue in the setting of a complex Lie group $M$ as before. Consider
the following automorphism $T$ of $M\times M$:
$$T(z,w)=(z,W):=(z,wz^{-1})=(z,R_{z^{-1}}(w)).$$
Let $D\Subset M$ be a domain with $C^{\infty}$ boundary and let
$${\cal D}:=T(D\times D)={\textstyle  \bigcup_{z\in D}}\, (z,D(z))$$
where
$$D(z):=D\cdot z^{-1}=R_{z^{-1}}(D)=\{wz^{-1}:w\in D\}.$$
This is a variation of domains $D(z)$ in $M$ with parameter space
$D\subset M$. Note that $e\in D(z)$ for all $z\in D$ since, for $z\in D$,
$zz^{-1}\in D(z)$. Let
$G(z,W)$ be the $c$-Green function for $(D(z),e)$ and $\Lambda(z)$ the
$c$-Robin constant. Then $\Lambda(z)$ is a $C^\infty$ function on $D$, called the {\bf $c$-Robin function} for $D$. Finally, let $\{a_k(z)\}_{k=1,...,n}$ be the
eigenvalues (repeated with
multiplicity) of the complex Hessian matrix $[{\partial^2(-\Lambda)\over
\partial \overline z_j \partial z_k}]_{j,k}$ at $z\in D$.

The following elementary result, which we state without proof, will be used repeatedly in this paper.
\begin{proposition}
 \label{prop:often} \
Let $V$ be a domain in ${\mathbb{C}}^n$ and let $s(z)$ be a $C^2$  function on $V$. Assume that $\bigl[\frac{ \partial^2 s(z_0+
 {\bf  a}t)}{\partial t \partial \overline { t}}\bigr]|_{t=0}=0$ for some
 $z_0\in V$ and some ${\bf  a}\in {\mathbb{C}}^n\setminus \{0\}$. Then for any
holomorphic curve $z=\zeta(t),
 |t| \ll 1$ in $V$ with $\zeta(0)=z_0$
 and $\bigl[\frac{d\zeta(t) }{\partial t}\bigr]|_{t=0}= {\bf  a}$,  we have
$\bigl[\frac{ \partial^2 s(\zeta(t))
 }{\partial t \partial \overline { t}}\bigr]|_{t=0}=0$.
\end{proposition}

We use Lemma \ref{thm:lie} to prove the following result, which will be crucial in all that follows.
\begin{lemme}
 \label{prop:crucial} \ Suppose $D$ is pseudoconvex. Then $-\Lambda
(z)$ is a plurisubharmonic exhaustion function for $D$. Furthermore, if
$a_i(\zeta)=0$ for
$i=1,...,k\leq n$ for some $\zeta\in D$, then there exist $k$
 linearly independent  holomorphic
vector fields $X_1,...,X_k$ in $\mathfrak X $
which satisfy the following conditions:
\  For each $z\in D$ and each $X=X_i\, (i=1,\ldots , k)$,
\begin{enumerate}
 \item [o.]\ $[\frac{\partial^2 \Lambda(z \exp tX)}{\partial t \partial
\overline { t}}] |_{t=0}=0$;
 \item [i.]\ $z\exp tX\in D$, $t\in \mathbb C$ \   and \    $D\cdot(z \exp tX )^{-1}=D\cdot
 z^{-1}, \   t\in {\mathbb{C}};$
 \item [ii.]\ $ \Lambda( z \exp tX)= \Lambda(z), \ t\in {\mathbb{C}};$
 \item [iii.]\ $\{z \exp tX \in M: t \in {\mathbb{C}}\} \Subset D$.
\end{enumerate}
\end{lemme}

\noindent {\bf Proof}. We know that $-\Lambda$ is plurisubharmonic by Theorem \ref{thm:psm} since ${\cal D}$ is pseudoconvex in $D\times M$. To show that
$-\Lambda$ is an exhaustion function for $D$, since $D(z):=D\cdot z^{-1}$,
if $z \to \eta \in \partial D$, then $\eta z^{-1} \in \partial D(z) $ and
$\eta z^{-1} \to e$; hence we have
$\Lambda (z) \to -\infty$ from standard potential theory arguments.

Now suppose $a_i(\zeta)=0$ for
$i=1,...,k\leq n$ for some $\zeta\in D$; i.e.,
$$[{\partial ^2\Lambda(\zeta +b_it)\over \partial t \partial \overline
t}]|_{t=0}=0, \quad  i=1,...,k$$
for linearly independent vectors $b_1,...,b_k\in {\mathbb{C}}^n$ (in local
coordinates).
We take any $X_i \in \mathfrak X $  such that
$$
X_i(\zeta)=  [\frac{d (\zeta \exp tX_i)}{dt } ]|_{t=0} =b_i.
$$
Then $X_1, \ldots ,X_k$ are linearly
independent in $\mathfrak X $ and satisfy, by Proposition
\ref{prop:often},
$$
[\frac{\partial^2 \Lambda(\zeta \exp tX_i)}{\partial t \partial
\overline { t}}] |_{t=0}=0.
$$

For sufficiently small $|t|$, say for $t$ in
$B=\{t:|t| < \rho\}$, we can assume that
$\zeta\exp tX_i \subset D, \ i=1,...,k$.
Consider the transformation
$$T_i(t,z)=(t,w)=(t,F_i(t,z))$$
where $w=F_i(t,z)=z(\zeta\exp tX_i)^{-1}$. Then
${\cal D}_i:=T_i(B\times D)$ defines a variation of domains
$D_i(t):=F_i(t,D)=\{z(\zeta\exp tX_i)^{-1}\in
M:z\in D\}=D\cdot(\zeta\exp tX_i)^{-1}$;
moreover, $e\in D_i(t), \ t\in B$. Let $g_i(t,w)$ be the $c$-Green
function for $(D_i(t),e)$ and let $\lambda_i(t)$ be the
$c$-Robin constant for $(D_i(t),e)$.  By definition, we have $\lambda
_i(t)=\Lambda (\zeta\exp tX_i)$, so that
$${\partial ^2 \lambda_i\over \partial t \partial \overline
t}(0)=[{\partial ^2\Lambda(\zeta \exp tX_i)\over \partial t \partial \overline
t}]|_{t=0}=0, \quad  i=1,...,k.$$
We now make direct use of the previous result for $A=D$ to conclude that
 $z\in D$ if and only if $z\exp tX_i\in D$ for all $t\in {{\mathbb{C}}}$; also
$D\cdot z^{-1}=D\cdot
(z\exp tX_i)^{-1}$ for all $ t \in {\mathbb{C}}$, so that i. is proved. Assertion
i. and the definition of $\Lambda(z)$ imply ii. Furthermore, since $-\Lambda$ is an exhaustion function for $D$, for a fixed $z\in D$ the
level set $S_z:=\{z'\in D:\Lambda(z')=\Lambda(z)\}$ is  relatively
compact in $D$.
Together with ii. this implies iii. Finally,
 since  $-\Lambda$ is
a plurisubharmonic exhaustion function for $D$, $-\Lambda$ is
subharmonic and bounded on the complex curve
$\{z\exp tX, \ t\in {\mathbb{C}} \}$. This curve is conformally  equivalent
to ${\mathbb{C}}$, ${\mathbb{C}}\setminus \{0\}$ or a one-dimensional compex  torus and hence $-\Lambda$
is constant on this curve, which implies $o.$
\hfill $\Box$

\begin{remark}\label{remark:one-inv} \
{\rm Condition $i.$ for $X \in \mathfrak X $ is clearly equivalent to
\begin{align}
 \label{eqn:eqwot}
 D\cdot \exp tX =D, \ t\in {\mathbb{C}};
\end{align}
i.e., $D $ is invariant under the one parameter
transformation group $\exp tX$ of $M$. Indeed, the proof of the lemma shows that $o.-iii.$ for $X \in \mathfrak X $ are all equivalent; and hence equivalent to (\ref{eqn:eqwot}).}
\end{remark}
\begin{corollary}
\label{oldremark5.1} \ Suppose $D \Subset M$ is pseudoconvex with smooth
 boundary. Then $-\Lambda
(z)$ is not strictly plurisubharmonic in $D$ if and only if  there exists $X\in
 \mathfrak X \setminus \{0\}$ for which
 one of the four  equivalent conditions $o.\sim iii.$ in Lemma \ref
{prop:crucial} hold.
 Moreover, for such $X$, $\{ z \exp tX \in
            M : t\in C\} \Subset A$ for $z\in A$ where $A =D, \ \partial D \mbox{  or }
            \ \overline { D}^c$.
\end{corollary}

\noindent {\bf  Proof}. We have already verified the first part; and iii. verifies the last statement for $A=D$.  To prove the last part for
 $A=\partial D$ or $\overline { D}^c$, from Lemma
\ref{thm:lie} we have  $z\in A$ if and only if $ z \exp tX \in A$ for all $t \in {\mathbb{C}}$; and
$A\cdot z^{-1}= A \cdot (z\exp tX)^{-1}$. Since $\partial D$ is compact,
 this verifies the corollary
for  $A=\partial D$. To verify the corollary for  $A=\overline D^c$, fix $z\in \overline D^c$ and $z_0\in
\partial D$. Then
$$z\exp tX=zz_0^{-1}(z_0\exp tX)$$
and $\{z_0 \exp tX : t\in {\mathbb{C}}\}\Subset \partial D$ since $z_0\in \partial D$ so
that $\{z \exp tX: t \in {\mathbb{C}}\}\Subset M$. Thus it
suffices to show that
\begin{align}
 \label{eqn:4.1}
\overline {z\exp tX}\cap \partial D=\emptyset.
\end{align}
If (\ref{eqn:4.1}) does not hold, there exists a point $ \eta \in
\partial D $ and a  sequence
$\{t_\nu\}_\nu$ in ${\mathbb{C}}$ such that
 $z(\exp t_\nu X) \to \eta $ as  $\nu \to \infty$.
 We write $D_1(z):=\overline D^c \cdot z^{-1}$ for $z\in
\overline D^c$. Then $e\in D_1(z)$ for such $z$,
so that we can consider  the $c$-Green function $G_1(z, \cdot)$ and the
$c$-Robin constant $\Lambda_1(z)$ for $(D_1(z),e)$. (Note that
although $D_1(z)$ is not necessarily relatively compact in $M$, we can
define the $c$-Green function $G_1(z,\cdot)$ and the $c$-Robin constant
$\Lambda_1 (z)$ for $(D_1(z),
e)$ by a standard exhaustion method.)
 If $ Z \in \overline D^c \to \eta \in \partial D$, then $\eta Z^{-1} \in
 \partial D_1(Z) $ and $ \eta Z^{-1} \to e$; thus
 $\Lambda_1(Z) \to -\infty$ as $Z \to \eta$. In particular,
$\Lambda_1(z(\exp t_\nu X))
 \to - \infty$ as $\nu \to \infty$. On the other hand, from the previous
results, $D_1(z)=D_1(z\exp t_\nu X), \nu=1,2,\ldots.$
Thus, $\Lambda_1(z)=\Lambda_1(z\exp t_\nu X) $ for
all $\nu$, and since $\Lambda_1(z)> -\infty$, we obtain a contradiction,
verifying (\ref{eqn:4.1}).
 \hfill $\Box$
\begin{corollary}\label{thm:4.3} \ Suppose $D \Subset M$ is pseudoconvex
with $C^{\infty}$ boundary. Let $c>0$ be a $C^{\infty}$ function on $M$
and let $\Lambda (z)$
be the $c$-Robin function for $D$. If the complex
Hessian matrix $[{\partial
^2(-\Lambda)\over \partial z_j \partial \overline z_k}(\zeta)]$ has a zero
eigenvalue with multiplicity
$k\geq 1$ at {\bf some} point $\zeta \in D$, then the complex Hessian
matrix of {\bf any}
plurisubharmonic exhaustion function $s(z)$ for $D$ has a zero eigenvalue
with
multiplicity at least $k$ at {\bf each} point $z \in D$.
\end{corollary}
\par
\noindent {\bf  Proof}. Take  $X_i\in \mathfrak X \setminus \{0\}, \ i=1,...,k$
as in the lemma. Given $z\in D$, $\{z\exp tX_i:t\in {\mathbb{C}}\}\Subset D$. Thus if $s$ is any plurisubharmonic exhaustion function for $D$, $s$ is
subharmonic and bounded on this complex curve and hence, as in the proof of Corollary \ref{oldremark5.1}, $s$ is constant on the curve. Thus $[{\partial ^2s\over
\partial z_a
\partial \overline z_b}(z)]$ has a zero eigenvalue of multiplicity at
least $k$ (cf., ~\cite{ST}\  Proposition 2.1). This is valid for each $z\in D$.
\hfill $\Box$

\smallskip
We remark that the conclusion of the corollary implies, in
particular, that $D$ is not Stein. We now address the following question: {\it for a complex Lie group $M$, when is a pseudoconvex domain
$D\Subset M$ with $C^{\infty}$ boundary Stein}?
An answer is provided in the following result.
\begin{theorem}\label{thm:fund-lie-case}\ Let $D \Subset M$ be a pseudoconvex domain
 with smooth boundary which is not Stein. {Then
\begin{enumerate}
 \item [I.]
there exists a unique connected complex Lie
subgroup $H$ of $M$   such that
\begin{enumerate}
 \item [1.]\  $\dim \, H \ge 1$;
 \item [2.]\ $D$ is  foliated by cosets of $H$; $D=\bigcup_{z\in D}zH$
 with $zH \Subset D$;
 \item  [3.] \ any holomorphic curve  $ \ell  :=\{z=z(t) \in D: t\in \mathbb{C} \} $  with $\ell \Subset D$ is necessarily contained in some coset $zH$ in $D$.
\end{enumerate}
\item [II.]\ Furthermore, if $H$ is closed in $M$ and
 $\pi: M
 \mapsto M/H$ is the canonical projection, then $H$ is a  complex
  torus and there exists a
 Stein domain $D_0 \Subset M/H$ with smooth boundary such that
$D= \pi^{-1}(D_0)$.
\end{enumerate}}
\end{theorem}

\noindent {\bf Proof}. Let $\Lambda(z)$ denote the $c$-Robin function on
$D$; this is a $C^\infty$ plurisubharmonic exhaustion function on
$D$. We utilize $\Lambda(z)$ to define a
Lie subalgebra $\mathfrak g _0$ of $\mathfrak g$ and a candidate Lie subroup $H$ using (\ref{eqn:liesub}) with $\mathfrak g _0$ replacing $\mathfrak g$. We begin with the construction of $\mathfrak g _0$, then of $H$; finally we verify the properties of $H$. We divide this procedure into several steps.

\smallskip
{\it $1^{st}$ step.} \ {\it   Let $\zeta\in D$ and define
$$
\mathfrak X _\zeta: =\{X\in \mathfrak X : [ \frac{\partial^2 \Lambda(\zeta
\exp tX) }{\partial t \partial \overline { t }} ]|_{t=0}=0\}.
$$
Then $\mathfrak X _\zeta $ is independent of $\zeta\in D$ and $\dim \,
\mathfrak X_\zeta=d \ge 1$. }

\smallskip Indeed, since $D$ is not Stein, $-\Lambda(z)$ is not strictly
plurisubharmonic on $D$; i.e., in local coordinates there exists a point $\zeta\in D$ and
$b \in {\mathbb{C}}^n \setminus \{0\}$
 such that $[\frac{\partial^ 2\Lambda(\zeta+bt)}{\partial t \partial \bar t}]|_{t=0}=0 $. The $1^{st}$ step now follows from Lemma \ref{prop:crucial} and Corollary \ref{oldremark5.1}. \hfill
$\Box$

\smallskip
Thus we write  $\mathfrak g _0:=\mathfrak
X _\zeta$. We will show in the $3^{rd}$ step that  $\mathfrak g _0$ is, indeed, a
Lie subalgebra of  $\mathfrak g$.

\smallskip
{\it $2^{nd}$ step.} \ {\it   Fix $\zeta\in D$ and an integer
$m\ge 1$. Then
$$
{\textstyle  \Lambda( \zeta \, \prod _{k=1}^m \exp t_kX_k)=\Lambda(\zeta),}  \quad  \hbox{for all} \ t_k\in {\mathbb{C}}, \ X_k \in \mathfrak g _0 .
$$}

\smallskip To see this, from Corollary \ref{oldremark5.1}  and  (\ref{eqn:eqwot}) we have
$$
D\cdot (\exp t_mX_m)^{-1}=D,  \quad \hbox{for all} \ t_m\in {\mathbb{C}}, \
X_m\in \mathfrak g _0.
$$
Thus for $t_{m-1}\in {\mathbb{C}}$,
\begin{align*}
D\cdot ( \exp t_{m-1}X_{m-1} \exp t_mX_m)^{-1} &=
D\cdot (\exp t_mX_m)^{-1} (\exp t_{m-1}X_{m-1})^{-1}\\
&= D\cdot (\exp t_{m-1}X_{m-1})^{-1}\\
&=D.
\end{align*}
Inductively, we have
$$
{\textstyle  D\cdot [\prod _{k=1}^m \exp t_kX_k]^{-1}=D,}  \quad \hbox{for all} \ t_k\in {\mathbb{C}},
\   X_k\in \mathfrak X_0.
$$
This is equivalent to
$$
D\cdot[\zeta \,{\textstyle   \prod _{k=1}^m \exp t_kX_k]^{-1}=D\cdot \zeta^{-1},}\quad \hbox{for all} \ t_k\in {\mathbb{C}},  \ X_k\in \mathfrak g _0.
$$
The $2^{nd}$ step now follows immediately  from the definition of $\Lambda(z)$. \hfill $\Box$

\smallskip
{\it   $3^{rd}$ step.} \ {\it  The space   $\mathfrak g _0$ is a Lie subalgebra  of
$\mathfrak g$; i.e., $\mathfrak g_0$ is a linear
subspace of $\mathfrak g $ and $X,\, Y\in \mathfrak g _0$ implies
$[X,Y]\in \mathfrak g _0$.}

\smallskip
To prove that $ \mathfrak g_0$ is a linear subspace of $\mathfrak
g $, fix $\zeta\in D$; let $X,\ Y\in \mathfrak g_0$; and let
$\alpha , \, \beta \in {\mathbb{C}}$. Since it is clear that $\alpha \,X, \ \beta\, Y \in
\mathfrak g_0 $, from the $2^{nd}$ step we have
$$
\Lambda(\zeta\, \exp (\alpha t)X \, \exp (\beta t)Y)= \Lambda(\zeta),
\quad \hbox{for all} \ t\in {\mathbb{C}}.
$$
Thus
$$
[\frac{\partial^2
\Lambda(\zeta\, \exp (\alpha t)X \, \exp (\beta t)Y)}{\partial t \partial
\overline { t}}]|_{t=0}= 0.
$$
It is a standard result that for any $X,Y\in  \mathfrak g$
\begin{align}
 \label{eqn:well-known-1}
\exp (\alpha t)X\, \exp (\beta t)Y=
\exp [t\,(\alpha X+ \beta Y)+
O(t^2)\,]
\end{align}
near $t=0$ and hence it follows that
{ \begin{align*}
[\frac{d(\zeta \exp t(\alpha X+ \beta
 Y)) }{dt } ]|_{t=0}
&= [\frac{d (\zeta \exp (\alpha t)X\,\exp (\beta t)Y)}{dt } ]|_{t=0}.
\end{align*}}
Using Proposition  \ref{prop:often} we see that $\alpha X + \beta Y
\in \mathfrak
g_0$; i.e.,
\begin{align}
 \label{eqn:tri-but}
[\frac{\partial^2
\Lambda(\zeta \exp t(\alpha X+ \beta
 Y))}{\partial t \partial
\overline { t}}]|_{t=0}=[\frac{\partial^2
\Lambda(\zeta\, (\exp (\alpha t)X \, \exp (\beta t)Y))}{\partial t \partial
\overline { t}}]|_{t=0}= 0.
\end{align}
To prove that $X,\, Y\in \mathfrak g _0$ implies $ [X,Y]\in
\mathfrak g_0$, fix $\zeta \in D$. From the $2^{nd}$ step  we have
$$
\Lambda \bigl(\zeta \exp ( - \sqrt t X)\, \exp (-\sqrt{t}Y) \, \exp
       \sqrt tX \, \exp \, \sqrt tY \bigr)
= \Lambda (\zeta), \  \quad \hbox{for all} \ t\in {\mathbb{C}}.
$$
On the other hand, it is known that
\begin{align}
 \label{eqn:well-known-2}
 \exp ( - \sqrt t X)\, \exp (-\sqrt{t}Y) \, \exp
       \sqrt tX \, \exp \, \sqrt tY &=
 \exp \bigl(t[X,Y])+O(t^{3/2}) \bigr)
\end{align}
near $t=0$ so that
$$
[\frac{ d  (\zeta\,\exp ( - \sqrt t X)\, \exp (-\sqrt{t}Y) \, \exp
       \sqrt tX \, \exp \, \sqrt tY)}{dt}
]|_{t=0}
=[ \frac{ d( \zeta \, \exp t[X,Y])}{dt}]|_{t=0}.
$$
It follows from Proposition \ref{prop:often}
 that
\begin{eqnarray}
&& [\frac{\partial^2 \Lambda (\zeta\exp \bigl(t[X,Y])} { \partial t
\partial \overline { t}}]| _{t=0} \nonumber \\
&=&[\frac{\partial^2 \Lambda
(\zeta\, \exp ( - \sqrt t X)\, \exp (-\sqrt{t}Y) \, \exp
       \sqrt tX \, \exp \, \sqrt tY)}{\partial t \partial \overline
       { t}}]|_{t=0} \nonumber \\
&=&
[\frac{\partial^2 \Lambda (\zeta) }{\partial t
\partial \overline { t} }]|_{t=0}=0.  \label{eqn:tri-but2}
\end{eqnarray}
Since  $[X,Y]\in \mathfrak X $, we see from Lemma \ref{prop:crucial}
and Corollary \ref{oldremark5.1} that $[X,Y]\in \mathfrak g _0$.
 \hfill $\Box$

\smallskip
 From the Frobenius theorem described at the beginning of this section and the fact that $\mathfrak g _0$ is a Lie subalgebra  of
$\mathfrak g$, we have  the
integral manifold $H$ in $M$ for $\mathfrak g _0$ passing through the
identity $e$,\, i.e., from (\ref{eqn:liesub}),
\begin{align}
 \label{eqn:defofH}
{\textstyle  H=\{\prod_{i=1}^\nu \exp t_iX_i  \in M \,:\,\nu\in {\boldsymbol{Z}}^+, \ t_i\in {\mathbb{C}}, \ X_i\in
 \mathfrak g _0\}}.
\end{align}

\smallskip
{\it $4^{th}$ step.} \ {\it  $H$ satisfies conditions 1.,  2. and 3.
 in Theorem
\ref{thm:fund-lie-case}.}

\smallskip
Since $\dim \, H=\dim \, \mathfrak g _0\ge 1 $, $H$ satisfies 1.
 Fix $\zeta\in D$. We consider the level set
$$
S_\zeta:=\{z\in D : \Lambda(z)= \Lambda(\zeta)\}.
$$
Since $-\Lambda$ is an exhaustion function for $D$, we have $S_\zeta \Subset D$.
Fix $X\in \mathfrak g_0 $. We have
$$
\Lambda(\zeta \exp tX)=\Lambda(\zeta),  \quad \hbox{for  all} \ t\in {\mathbb{C}}
$$
so that $\zeta \exp tX \subset S_\zeta$. Since $t\in {\mathbb{C}}$ and $X\in
\mathfrak g _0$ were arbitrary, it follows from
 (\ref{eqn:defofH}) that we have $\Lambda( \zeta g)= \Lambda(\zeta)$ for all
       $g\in H$. Hence $\zeta H \subset S_\zeta \Subset D$.  The conclusion that $D$ has a
       foliation $D=\cup_{z\in D}zH$ follows from this observation together with the fact that $M$ itself is foliated by cosets $M=\cup_{z\in M}zH$ (the second item of the Frobenius theorem). Thus $H$ satisfies 2.

 We shall prove that $H$ satisfies 3. by contradiction. Thus we assume
 that there exists a holomorphic curve $ \ell  :=\{z=z(t) \in D: t\in \mathbb{C} \} $  with $\ell \Subset D$ and with the property that there exists a coset $z_0H$ in $D$  with $\ell \cap z_0H\ne \emptyset
$ and $\ell \not \subset z_0H$. Since $D$ is foliated by cosets $zH$, we
 may assume that $\ell$ and $z_0H$ intersect transversally at some point
 $\zeta$ in $D$.
Let $a\in {\mathbb{C}}\setminus \{0\}$
be a tangent vector of $\ell $ at $\zeta$ and choose $X\in \mathfrak X
$ such that $X(\zeta)=a$.  By assumption $\Lambda=const.\
=\Lambda(\zeta)$ on $\ell$. It follows from Lemma
\ref{thm:fund-lie-case} that $X\in \mathfrak g _0$, so that $C:=\{\zeta \exp
tX: \ t\in {\mathbb{C}}\} \subset \zeta H$. However, the curve $C$ intersects $z_0H$ transversally; in particular, $C \not \subset z_0H$. From property 2., since $z_0H\cap \zeta H\not = \emptyset$, we have $z_0H=\zeta H$, which yields a contradiction.
 \hfill $\Box$

\smallskip
{\it $5^{th}$ step.} \ {\it  $H$ satisfying  1., 2. and 3. is uniquely
determined. }
\smallskip

Let $\widehat{ H}$ be another Lie subgroup of $M$ which satisfies
1., 2. and 3. in the theorem. We write $\widehat{ \mathfrak g }_0$ for
the connected Lie subalgebra of $\mathfrak g
$ which corresponds to $\widehat{ H}$. Thus from (\ref{eqn:liesub}),
$$
{\textstyle  \widehat{ H}=\{\,\prod_{i=1}^\nu \exp t_iX_i  \in M  :  \nu\in {\boldsymbol{Z}}^+, \  t_i\in {\mathbb{C}}, \ \widehat {X}_i\in
 \widehat {\mathfrak g} _0\,\}}.$$
Fix $\widehat{ X}\in \widehat{ \mathfrak g
 }_0$ and $\zeta\in D$. Since $\zeta\widehat{ H} \Subset D$, it follows from
  condition 3. for $H$ that  the analytic curve $\{\zeta\exp\, t \widehat X,\
 t\in {\mathbb{C}}\}$ passing through $\zeta$
is contained in $\zeta H$. By the above description of  $\widehat{ H}$ we
inductively have $\zeta \widehat{ H} \subset \zeta H$, so that
 $\widehat{ H}\subset H$. The converse inclusion is proved in a similar fashion.
 \hfill $\Box$

\smallskip
 {\it $6^{th}$ step.} \ {\it Assertion $II$ is true. }

\smallskip
Suppose $H$ is closed in $M $. By 2. $H$ is compact in $M$, so that
$H$ is a compact complex Lie group, i.e., a complex torus. Let $\pi:M\mapsto M/H$ be the
canonical projection. We set $D_0:=\pi(D) \subset M/H$. Since $D$ is foliated by cosets of $H$ (property 2.),
$\pi^{-1}(D_0)=D$; moreover, $\pi^{-1}(\partial D_0)=\partial D$, and since $\partial D$ is smooth in $M$ we have $\partial D_0$ is smooth in $M/H$. Given
$\xi\in D_0$, we choose a point $\zeta\in \pi  ^{-1}(\xi)\subset D$ and define
$$
\Lambda_0(\xi):=\Lambda(\zeta).
$$
Again from property 2., this does not depend on the choice of $\zeta\in D$ with
$\pi(\zeta)=\xi$; i.e., $\Lambda_0$ is a well-defined function on
$D_0$. We show that $\Lambda_0$ is a $C^\infty$
 strictly plurisubharmonic exhaustion function for $D_0$.

To this end, let $\xi_n\in D_0 \to \xi\in \partial D_0\subset M/H$ as $n \to
\infty$. We can choose $\zeta_n\in D$ and $\zeta\in \partial D$ with $\pi
(\zeta_n)= \xi_n$ and $\pi (\zeta)= \xi$ in such a way that $\zeta_n \to \zeta$ in $M$ as $n \to \infty$. Since
$-\Lambda$ is an exhaustion function for $D$, it follows
that $-\Lambda_0(\xi_n)= -\Lambda(\zeta_n) \to \infty$, so that
$-\Lambda_0$ is an exhaustion function for $D_0$. Next, fix $\xi_0\in D_0$ and
take a point $\zeta_0\in D$ with $\pi ^{-1} (\zeta_0)=\xi_0$. From 2. of the
Frobenius theorem we have $\mathfrak g_0$-local coordinates
$\varphi :V \Subset D \mapsto \Delta=\prod _{i=1}^n\{|w_i|<1\}$ with
$\varphi (\zeta_0)=0$.
We write
$$
w=(w_1, \ldots , w_d;w_{d+1}, \ldots , w_n)=(w';w'')\in \Delta' \times
\Delta''=\Delta$$
where $\Delta'=\prod_{i=1}^d \{|w_i|<1\}$ and
$\Delta''=\prod_{i=d+1}^n \{|w_i|<1\}$; we also write $0=(0';0'')\in
\Delta' \times \Delta''$. We define, for $w''\in \Delta''$,
$s(w''):= \varphi ^{-1} (0';w'') \subset V$. Then, since $H$ is closed
in $M$, $V= \cup_{w''\in
\Delta''}(s(w'')H)\cap V$  and this is a disjoint union (if necessary,
take a smaller  neighborhood  $V$ of $\zeta_0$). Setting $\pi (V)=\delta
\subset D_0$, for any $\xi\in \delta$ there exists a unique point
$w''\in \Delta''$ such that $\pi(s(w''))=\xi$. Consequently,
$$
\varphi _0: \xi\in \delta\mapsto w''\in \Delta''
$$
gives local coordinates for the neighborhood  $\delta$ of $\xi_0$ in
$D_0$. Identifying $s(w''  )$ with $(0';w'')$ via the mapping $\varphi $ and identifying $\xi$ with $w''$ via the mapping $\varphi _0$, from the definition of
$\Lambda_0$ for $w''=(w_{d+1}, \ldots , w_n)\in \Delta''$ we have
\begin{align}
 \label{eqn:ddashpro}
\Lambda_0(w_{d+1}, \ldots , w_{n})=\Lambda(0,\ldots ,0;w_{d+1},\ldots ,
w_n).
\end{align}
Since $-\Lambda$ is plurisubharmonic on $\Delta$, it follows that $-\Lambda_0$ is plurisubharmonic on $\Delta''$, and hence $-\Lambda_0$ is plurisubharmonic
on $D_0$. It remains to verify that $-\Lambda_0$ is strictly
plurisubharmonic on $D_0$. If not, there exists a point $\xi_0\in D_0$ at
which $-\Lambda_0$ is not strictly plurisubharmonic. From formula
(\ref{eqn:ddashpro}) $-\Lambda_0(w_{d+1}, \ldots , w_n)$ is not strictly plurisubharmonic
at the origin $0''$. Thus $-\Lambda(0,\ldots ,0;w_{d+1},
\ldots , w_n)$ is not strictly plurisubharmonic at $w''=0''$. This means that  there exists $(b_{d+1}, \ldots , b_n)\in {\mathbb{C}}^{n-d}\setminus \{0''\}$  such
that
$$
[\frac{\partial ^2\Lambda(0,\ldots ,0;b_{d+1}t, \ldots ,
b_{n}t)}{\partial t \partial \overline { t}} ]|_{t=0}=0.
$$
On the other hand, from the previous steps we know that $\Lambda(w';0'')\equiv  \Lambda(0';0'')$ on
$\Delta'$; thus it follows that there exist at least $d+1$ linearly independent eigenvectors with eigenvalue $0$ for the matrix $[\frac{\partial^2 (-\Lambda)}{\partial z_i
\partial \overline { z}_j} ]_{i,j=1,\ldots ,n}$ at $\zeta_0=s(0'')\in D$.
This contradicts the definition of the dimension $d$ of $\mathfrak g_0$. \hfill $\Box$

\begin{remark} \label{re:5-2} \
{\rm In Theorem \ref{thm:fund-lie-case} the complex Lie
 group $M$
 is a real Lie group of dimension $2n$ with real Lie algebra
$\mathfrak g$. If $H $ is not closed in $M$
the closure $\overline { H}$ of $H$ in $M$ is a closed real Lie
subgroup of $M$ whose real dimension $m$
is less than $2n$.  In this case, we have the
real Lie subalgebra $\mathfrak g_0 \subset \mathfrak g $ corresponding
to $\overline { H}$, and the projection $\widetilde { \pi} : M\to M/
\overline { H}$, where $M/ \overline { H}$ is a real
manifold of dimension $2n-m$. From properties 1. and
2., $\overline { H} $ is a compact real submanifold in $M$  with
 $\overline { H} \subset D$ and
$D$ is foliated by cosets of $\overline { H}$, so that $\Lambda= \
 const.$ on each $z \overline { H}, \ z\in D$. Furthermore, $D_0:=
\widetilde { \pi}(D)$ is a relatively compact domain
with smooth boundary in $M/\overline { H}$ and $D=\widetilde
 {\pi}^{-1}( D_0)$. }
\end{remark}
\begin{corollary}\  The subgroup $H$ of $M$ in Theorem
 \ref{thm:fund-lie-case}
is  normal.
\end{corollary}
\noindent {\bf  Proof}. We may assume that $D$ contains the identity,
$e$. For if $z_0\in D$ and we define $D_1:=z_0^{-1}D$, then $D$ is a
pseudoconvex domain with smooth boundary; $D_1$ is not Stein; and $e\in
D_1$. Let $\Lambda_1$ be the $c$-Robin function for $D_1$.
 Assume that $X\in \mathfrak X $ satisfies
$$[{\partial^2  { \Lambda}( z_0\exp t  { X}   ) \over \partial t
\partial \overline { t}}]|_{t=0}=0.\eqno(*)$$
Then by Corollary
\ref{thm:4.3} we have $C:=\{z_0 \exp tX\in M: t\in {\mathbb{C}}\} \Subset D$.
 We thus have
$C_1:= \{\exp tX\in M: t\in {\mathbb{C}}\} \Subset z_0^{-1}D=D_1$, so that
$\Lambda_1=\ const. $ on  $C_1$. Thus,
$$[{\partial^2  { \Lambda_1}( \exp t  { X}   ) \over \partial t
\partial \overline { t}}]|_{t=0}=0.\eqno(**)$$
A similar argument shows that (**) implies (*). It follows that the Lie subroup $H_1$ in Theorem
 \ref{thm:fund-lie-case} corresponding to $D_1$ coincides with $H$. Thus we may indeed assume that $e\in D$.

We consider the isomorphism $\widetilde {
T}(\zeta,z)$ of $M$ for $z\in M$ given by left multiplication by $z^{-1}$:
$$\widetilde { T}(\zeta,z):\  \zeta\in M \mapsto w=  z^{-1}\zeta \in M.$$
This induces the variation of domains $\widetilde { D}(z)$ in $M$:
$$
\widetilde { {\cal D}}: \ z\in D \mapsto \widetilde { D}(z)= z^{-1}D.
$$
Since each $\widetilde { D}(z), \ z\in D$, contains the identity element $e$, we can
form
 the $c$-Robin function $\widetilde \lambda(z)$ on $D$.  We let
 $\widetilde { \mathfrak g} $ denote the Lie algebra consisting of all
 right-invariant holomorphic vector fields on $M$. For a fixed $\zeta\in D$ we define
$$
\widetilde { \mathfrak g }_0:=\{ \widetilde {  X}\in \widetilde { \mathfrak g }:
[\frac{\partial^2 \widetilde { \lambda}( (\exp t \widetilde { X} ) \zeta ) }{\partial t
\partial \overline { t}}]|_{t=0}=0\,\}.
$$
As with $\mathfrak g _0$, the set  $\widetilde { \mathfrak g
}_0$ is a Lie subalgebra of $\widetilde { \mathfrak g}$.
 We construct the connected  Lie subgroup
$\widetilde { H}$ corresponding to $\widetilde { \mathfrak g }_0$:
$${\textstyle  \widetilde { H}:=
\{\,\prod_{i=1}^\nu \exp t_iX_i \, \in M  :  \nu\in {\boldsymbol{Z}}^+, \ t_i\in {\mathbb{C}},  \ X_i\in
 \widetilde {\mathfrak g} _0\,\}}$$
with
$\dim \widetilde H=\widetilde { d}\ge 1$ and satisying the
 corresponding properties $1., 2.$ and $3.$ in Theorem
 \ref{thm:fund-lie-case}.
We first prove that
\begin{align} \label{eqn:1}
z H = \widetilde { H}z, \quad \mbox{  for any }  z\in D.
\end{align}
To prove this, fix $z\in D$ and $p\in
 \widetilde { H}z$ with $p=(\prod_{i=1}^\nu \exp t_i \widetilde { X}_i)\, z $ for some  $t_i\in
 {\mathbb{C}}$ and $\widetilde { X}_i\in \widetilde { H}$. Condition 2. for
 $\widetilde { H}$ implies that the analytic curve $C_\nu:=\{(\exp t \widetilde {
 X}_\nu) z : t\in {\mathbb{C}}\}$ is relatively compact in $D$. Condition 3. for
 $H$ thus implies that $C_\nu \Subset zH$. In paricular,   $z_\nu:=(\exp
 t_\nu \widetilde { X}_\nu)z\in zH$  and hence  $z_\nu H\subset zH$. If we consider
 the analytic curve $C_{\nu-1}:=\{ (\exp t \widetilde { X}_{\nu-1}) z_\nu\in
 M: t\in {\mathbb{C}}\}$, then a similar argument shows that $C_{\nu-1}
 \Subset z_1 H$,  so that $z_{\nu-1}:=(\exp t_{\nu-1} \widetilde {
 X}_{\nu-1}) z_\nu  \in z_\nu H \subset zH$, and hence $z_{\nu-1}H\subset zH$. Repeating this argument we obtain $(\exp t_1 \widetilde { X}_1) z_{2}
\in z_{2}H\subset zH$, i.e., $p\in zH$. We thus have
 $\widetilde {
 H}z \subset zH$. The reverse inclusion is proved similarly and
 equality (\ref{eqn:1}) holds for $z \in D$.

Since $e\in D$ it follows that $\widetilde { H}=H$. Thus $zH=Hz$ for $z\in D$. To extend this equality to all $z\in M$, we consider
$${\mathcal A}:= \{z\in M: z^{-1}Hz=H\}\subset M.$$
Then $D\subset {\mathcal A}$ and we want to show that ${\mathcal
A}=M$. Clearly $H\subset {\mathcal A}$ and ${\mathcal A}$ is a subgroup of
$M$. Since $e\in D\subset {\mathcal A}$, we can choose a neighborhood
$V\subset {\mathcal A}$ of $e$. We
utilize this neighborhood $V$ to show that ${\mathcal A}$ is open and
closed in $M$. To show that ${\mathcal A}$ is open, fix $z_0\in
{\mathcal A}$. Then $Vz_0$ is a neighborhood of $z_0$ in
$M$. Since $V \subset {\mathcal A}$, it follows that $Vz_0 \subset
{\mathcal A}$ (since ${\mathcal A}$ is a group), and hence ${\mathcal
A}$ is open. To show that ${\mathcal A}$ is closed, fix a point $\zeta
\in \partial {\mathcal A}$. Choose $z\in {\mathcal A}$ sufficiently
close to $\zeta$ so that $a:=z^{-1}\zeta\in V\subset {\mathcal
A}$. Hence $\zeta= za \in \mathcal A$, and ${\mathcal A}$ is closed.
\hfill $\Box$

\smallskip
We shall give two concrete examples of complex Lie groups $M$ in Theorem
\ref{thm:fund-lie-case} or Remark \ref{re:5-2}.

\smallskip
{\bf  1.}\ Consider
$$M=\hbox{Aut}\,{\mathbb{C}}=\{az+b:a\in {\mathbb{C}}\setminus \{0\}, \ b\in {\mathbb{C}} \}.$$
This is a complex Lie group of complex dimension $2$ with the group
operation of composition. This composition gives a group multiplication on
${\mathbb{C}}\setminus
\{0\}\times  {\mathbb{C}}$ which is noncommutative: $(a,b)\cdot (c,d)=(ac,ad
+b)$. That is, if $f(z)=az+b$ and $g(z)=cz+d$, then
$$(f\circ g)(z)=(ac)z + (ad+b).$$
 Note that any left-invariant holomorphic vector field  on $M$ is of the
form
$$
X= \alpha X_1 + \beta X_2, \quad \alpha , \beta \in {\mathbb{C}},
$$
where
$$
{\textstyle  X_1= a  \frac{\partial }{\partial a }, \quad X_2= a {\partial \over
       \partial b}.}
$$
We have
$$
\exp tX_1=(e^t,0) \qquad \hbox{and}  \qquad \exp tX_2 =(1,t).
$$
Hence,  $\exp tX= (e^{\alpha t}, \, \beta \, \frac{ e^{\alpha t}-1}{\alpha})$,
so that the integral curve $I_{(a_0,b_0)}$ for $X$ with given initial value
$(a_0,b_0)\in {\mathbb{C}} \setminus \{0\} \times {\mathbb{C}}$ is
$$
{\textstyle  (a_0,b_0)\cdot \exp tX =(a_0 \exp {\alpha t},\  a_0 \beta\,  \frac{ e^{\alpha t}-1}{\alpha}+b_0),}
$$
so that $I_{(a_0,b_0)}$ is a complex line in ${\mathbb{C}} \setminus \{0\} \times {\mathbb{C}}$
defined by $\beta a-\alpha b=\beta a_0 -\alpha b_0$.
 Thus a non-constant, entire integral curve for $X$ is never
relatively compact in   $M$ (this fact also follows from the Liouville theorem). It follows that any smoothly bounded pseudoconvex domain $D\Subset M$ is Stein.

\smallskip
{\bf  2.}\
Grauert gave an example of a pseudoconvex domain $D$ with smooth
boundary which is not Stein. Moreover, $D$ admits
no nonconstant holomorphic functions. This domain lies in a complex torus ${\bf
T}$ of complex dimension $2$ (cf., \cite{grauert} and Example 3 in p. 324 in \cite
{nishino}).
 Our goal is to describe {\bf all} pseudoconvex subdomains $D$ of ${\bf
 T}$ with smooth boundary which are not Stein. The key tools we will use
 are Theorem \ref{thm:fund-lie-case} and Remark \ref{re:5-2}.

We begin with real $4$-dimensional Euclidean
space ${\mathbb{R}}^4$ with
coordinates $x=(x_1,x_2,x_3,x_4)$. Let
$$
e_1=(1,0,0,0), \ e_2=(0,1,0,0), \ e_3=(0,0,1,0), \ e_4=(0,0,\xi,1)
\quad  \mbox{
in }  {\mathbb{R}}^4,
$$
where $\xi$ is an irrational number. Initially we consider the real $4$-dimensional torus:
$$
T:= {\mathbb{R}}^4/[e_1,e_2,e_3,e_4] = T_1 \times T_2,
$$
where $T_1={\mathbb{R}}_{x_1,x_2}/[e_1', e_2']$ and $T_2={\mathbb{R}}_{x_3, x_4}/
[e_3', e_4']$, with
\begin{align*}
 &e'_1=(1,0), \ e'_2=(0,1)  \quad   \mbox{ in } {\mathbb{R}}_{x_1}\times {\mathbb{R}}_{x_2};\\
&e'_3=(1,0), \ e'_4=(\xi,1)  \quad  \mbox{  in }{\mathbb{R}}_{x_3}\times {\mathbb{R}}_{x_4}.
\end{align*}
This torus $T$ is a real Lie group with real  Lie
algebra
$$
{\textstyle  \mathfrak h =\{   \sum  {}_{ j=1  }^{4  } c_j\frac{\partial
}{\partial x_j}:  c_j\in {\mathbb{R}} \}.}
$$
Following Grauert, we impose the complex structure
$$
z=x_1+ix_3,  \quad  w=x_2+ix_4
$$
on $T$. Then $T$, equipped with this complex structure, becomes a
complex  torus ${\bf  T}$ of complex dimension $2$. Note that
$e_1,e_2,
e_3,e_4$ correspond to $(1,0),(0,1),(i,0),(i\xi,i)$ in ${\mathbb{C}}^2$.
 Grauert showed
that
\begin{align}
 D=D(c_1,c_2):=\{c_1<\Re \ z < c_2\} \subset {\bf  T}, \nonumber
\end{align}
where $0\le c_1<c_2<1$, is a pseudoconvex domain which admits
no non-constant holomorphic functions.

\smallskip
The complex Lie algebra of the complex Lie group ${\bf  T}$ is
$$
{\textstyle
\mathfrak g = \{\alpha \frac{\partial }{\partial z} + \beta  \frac{\partial
}{\partial w} : \alpha, \beta \in {\mathbb{C}} \}. }
$$

\begin{remark}
 \label{rem:6-3} ${}$

  {\rm 1. \  A point $(x_1,x_2,x_3,x_4)$ in ${\mathbb{R}}^4$ (and, more generally, a subset $K \subset
 {\mathbb{R}}^4$) will often be identified  with the point $(x_1,x_2,x_3,x_4)/[e_1,e_2,e_3,e_4]$ in
 ${\bf  T} $ (and $K/[e_1,e_2,e_3, e_4]$ in ${\bf  T} $). A similar
 remark applies to points and subsets in ${\mathbb{R}}_{x_1}\times {\mathbb{R}}_{x_2}  $
 (or ${\mathbb{R}}_{x_3}\times {\mathbb{R}}_{x_4}  $) being identified with points and
 subsets in $T_1$ (or $T_2$). In particular, let $L_1$ be a line in ${\mathbb{R}}_{x_1}\times {\mathbb{R}}_{x_2}  $ containing the
 origin. Then $L_1$ defines a simple closed curve $l_1=L_1/[e'_1,e'_2]$ in $T_1 $ if and only if $L_1$ contains a point in ${\mathbb{R}}_{x_1}\times {\mathbb{R}}_{x_2}  $ of the form $(m,n)$ where $m,n \in {\boldsymbol{Z}}$. In this case,
there exist unique such $m, n\in {\boldsymbol{Z}}$ up to a sign with $(m,n)=\pm 1$, and $l_1$ may be identified with the segment joining $(0,0)$ and
$(m,n)$ in ${\mathbb{R}}_{x_1}\times {\mathbb{R}}_{x_2}$. If we need to consider a directed line $L_1$ and an oriented curve $l_1$, we fix the direction of $L_1$ as that given by the directed line segment from $(0,0)$ to $(m,n)$, and this induces what we will call the positive orientation on $l_1$. Moreover, if we fix
 $m_1,n_1\in {\boldsymbol{Z}}$  such that $m_1n- m n_1 =1$, then the parallelogram
 $\Delta_1 $ with vertices $(0,0), (m_1,n_1), (m,n),(m_1+m,n_1+n)$ is a fundamental domain of $T_1$, which we call
a  {\sl standard fundamental domain associated to $L_1$}. {(For future
 convention we choose $m_1n-n_1m=1$.)}
The segment joining $(0,0)$ and
 $(m_1,n_1)$ in ${\mathbb{R}}_{x_1}\times {\mathbb{R}}_{x_2}$ defines a simple closed
 curve in $T_1$, which we denote by $l_1^*$. Note that
 $(m_1,n_1)\in L_1+(1/n,0)$.

\smallskip
2. \ Let $L_2$ be a line in ${\mathbb{R}}_{x_3}\times {\mathbb{R}}_{x_4}  $ containing the origin. Then $L_2$ defines a simple closed curve $l_2=L_2/[e'_3,e'_4]$ in  $T_2 $ if and only if $L_2$ contains a point
 $(M',n')$ where $M'=m'+n' \xi$ with $m',n' \in {\boldsymbol{Z}}$. Then there exist unique such $m'\in {\boldsymbol{Z}},\  n'\in {\boldsymbol{Z}}^+$ with
 $(m',n')=\pm 1$, and $l_2$ may be identified with the segment joining
 $(0,0)$ and $(M',n')$ in ${\mathbb{R}}_{x_3}\times {\mathbb{R}}_{x_4} $. We fix the direction of $L_2$ as that given by the directed line segment from $(0,0)$ to $(m'+n' \xi,n')$, and this induces what we will call the positive orientation on $l_2$. If we fix
 $m'_1,n'_1\in {\boldsymbol{Z}}$  such that $m'_1n'- m' n'_1 =1$, then the
 parallelogram $\Delta_2$ with vertices $(0,0), (M'_1,n'_1),
 (M_1,n_1),(M_1'+M_1,n_1'+n_1)$ (where $M_1'=m_1'+n_1' \xi$) is a fundamental domain of $T_
2$, which we call
a  {\sl standard fundamental domain associated to $L_2$}.
The segment joining $(0,0)$ and
 $(M'_1,n'_1)$ in ${\mathbb{R}}_{x_3}\times {\mathbb{R}}_{x_4}$ defines a simple closed
 curve in $T_2$, which we denote by $l_2^*$. Note that
 $(M_1',n_1')\in L_2+(1/n',0)$.

\smallskip
3. \ Finally, we mention that if $L_1,L_1'$ ($L_2,L_2'$) are distinct
 parallel lines in ${\mathbb{R}}_{x_1}\times {\mathbb{R}}_{x_2}  $ (${\mathbb{R}}_{x_3}\times
 {\mathbb{R}}_{x_4}  $), then for the corresponding curves $l_1,l_1'$
 ($l_2,l_2'$) in $T_1$ ($T_2$) either $l_1\cap l_1'=\emptyset$ ($l_2\cap
 l_2'=\emptyset$) or $l_1 =l_1'$ ($l_2 = l_2'$). }
\end{remark}

 Let $D \Subset {\bf  T}$ be a pseudoconvex domain with smooth boundary which is not Stein. We consider the $c$-Robin function
$\Lambda (z,w)$ on $D$, where $c\equiv 1$ on ${\bf  T}$. By
Theorem \ref{thm:fund-lie-case} there exists $X= \alpha \frac{\partial }{\partial z} + \beta  \frac{\partial
}{\partial w }\in \mathfrak g $ with $(\alpha, \beta)\ne (0,0)$
such that $D(\exp tX) ^{-1} = D, \ t\in {\mathbb{C}}$. Since the integral
curve $\exp tX, \ t\in {\mathbb{C}}$ for $X$ passing through ${\mathbf{0}}$ in ${\bf  T}$ is
$$
(z, w)= (\alpha t,\beta t) \in {\mathbb{C}}_z \times {\mathbb{C}}_w,  \quad  t\in {\mathbb{C}},
$$
the above formula means simply that $D+(\alpha t,\beta t)=D, \ t\in
{\mathbb{C}}$. Since $\dim\ {\bf  T}=2$, the Lie subalgebra $\mathfrak g_0 $
associated to $D$ from Theorem \ref{thm:fund-lie-case} and the corresponding Lie subgroup $H$
are of the form
\begin{align}
 \label{eqn:start}
\mathfrak g _0=\{cX  \in \mathfrak g : c\in {\mathbb{C}}\},  \quad  H=\{(\alpha t,
\beta t)\in {\bf  T}: t\in {\mathbb{C}}\}.
\end{align}

We consider three cases in (\ref{eqn:start}):
$$
(1)\ \alpha=0,  \quad  \quad    (2)\ \beta=0,  \quad  \quad (3)\ \alpha, \ \beta\ne 0.
$$

In case $(1)$, the subalgebra $ \mathfrak g _0$ is $\{c \frac{\partial }{\partial w} :
 c\in {\mathbb{C}}\}$
and the Lie subgroup
$$
H=\{(0, \beta t) :t\in {\mathbb{C}}\}= \{0\}\times \mbox{  the $w$-plane }=\{0 \}\times {\mathbb{C}}_w;
$$
i.e.,
$$
H= \{(0, w)/
[(1,0),(0,1),(i,0),(i\xi,i)]
$$
We shall show that $H$ is not closed in ${\bf  T}$ and the closure $ \bar {
 H}$ in ${\bf  T}$  is diffeomorphic
 to $S^1 \times T_2$.
In fact, as a set in ${\bf  T}=T_1 \times T_2$, $H$ is equal to
$$
A \times B: =(\  (0,{\mathbb{R}}_{x_2})/[(1,0), (0,1)]\  ) \times ( \
(0, {\mathbb{R}}_{x_4})/[(1,0), (\xi,1)] \ ) \subset
T_1 \times T_2.
$$
Clearly $A$ is closed, and $A=\{0\}\times {\mathbb{R}}_{x_2}/[1]= \{0\}\times
S^1$. Since $\xi$ is irrational, it is standard that $B$ is not closed in $T_2$ and $\bar {
B}= T_2$.
 Thus the closure $\bar {  H}$ of $H$ in ${\bf  T} $ can be identified with
$$\bar H =(\{0\}\times S^1)\times T_2=\{ x_1=0 \}/
[(1,0),(0,1),(i,0),(i\xi,i)],
$$
as claimed.  We note that $\bar { H}$ is a  closed real
Lie subgroup of $T$ with
real Lie subalgebra $\{ \sum_{j=2}^4 c_j\frac{\partial  }{\partial
x_j} \}$. The quotient space $
 {\bf    T} / \bar { H}= S^1:=\{ x_1\in {\mathbb{R}}_{x_2} / [1]\}$
and the projection
$\widetilde { \pi}: {\bf  T} \mapsto {\bf  T} / \bar { H}$  is
$\widetilde {  \pi}(z,w)= x_1$. From Remark \ref{re:5-2}, $D$ is foliated by cosets of $\bar H$; thus there exists
a subdomain $D_0=(c_1,c_2)\subset S^1$  such that $D=\widetilde {
\pi}^{-1}(D_0)$; i.e., $D=\{c_1<x_1<c_2\}\subset {\bf T}$. This is the Grauert example.

In case $(2)$, the subalgebra $ \mathfrak g _0$ is $\{c \frac{\partial }{\partial z} :
 c\in {\mathbb{C}}\}$
and the Lie subgroup
$$
H=\{(\alpha t,0) :t\in {\mathbb{C}}\}= \mbox{  the $z$-plane } \times \{0\}={\mathbb{C}}_z \times \{0\};
$$
i.e.,
$$
H= \{(z, 0)/ [(1,0),(0,1),(i,0),(i\xi,i)]
 \in {\bf  T}: z\in {\mathbb{C}}_z\}.
$$
We shall show that $H$ is closed in ${\bf  T}$ and is holomorphically
 equivalent to the standard one-dimensional torus ${\bf  T}_1={\mathbb{C}}/[1,i]$.
In fact, since
$z_1\sim z_2 $ modulo $[1,i]$  in ${\mathbb{C}}_z$ if and only if  $(z_1,0) \sim
(z_2,0)$ modulo $
[(1,0),(0,1),(i,0),(i\xi,i)] $ in ${\mathbb{C}}^2$,
it follows that
$$
H= {\mathbb{C}}_z /[1,i] \times \{0\}
$$
so that $H= {\bf  T}_1 \times \{0\}$, as claimed.
Moreover, we have
$
(z_1,w_1)+H= (z_2,w_2)+H \mbox{\ \  in ${\bf  T}$}$ if and ony if  $w_1 \sim w_2 \mbox{ \ modulo $[1,i]$}.
$
The ``only if'' direction is clear. To
 prove the reverse direction, fix $(z_1,w_1), (z_2,w_2)\in {\mathbb{C}}^2$ with
$w_1 \sim w_2$ modulo $[1,i]$, i.e., $w_2-w_1=p+iq$ for some $p,q\in
 {\boldsymbol{Z}}$. Since $(z,p+qi)-(z-q\xi i,0)=p(0,1)+ q(i\xi,i)=(0,0)$ in ${\bf
 T}$  for any $z\in {\mathbb{C}}$,
  it follows that
\begin{align*}
 (z_2,w_2)&= (z_1,w_1)+ (z_2-z_1, p+qi)\\
&= (z_1,w_1)+ (z_2-z_1-q\xi i,0) \in (z_1,w_1)+ H  \quad  \mbox{  in
 ${\bf  T}$},
\end{align*}
as claimed. In particular, $(z,w)+H= (0,w)+H$ in ${\bf  T}$ for any $z\in {\mathbb{C}} $,  and
 $(0,w)+H = (0,w')+H$ in ${\bf  T}$ if and only  if $w \sim w'$ modulo $[1,i]$.
Hence ${\bf  T}/H= \{0\} \times {\mathbb{C}}_w/[1,i] = \{0\} \times {\bf  T}_1$, so that
 ${\bf  T}= {\bf  T}_1 \times {\bf  T}_1$. By II in Theorem
 \ref{thm:fund-lie-case}, using the projection
$\pi: {\bf  T} \mapsto {\bf T}/H= {\bf  T}_1$, we see that
 there exists a domain $D_0 \subset {\bf  T}_1$  with
smooth boundary such that $D=\pi^{-1}(D_0)$.

\smallskip  We consider case $(3)$. Write $\beta/\alpha=a+ib$.
 Then the integral curve $\{\exp tX, \ t\in {\mathbb{C}}\}$ starting at $(0,0)$ can be written as
$$
   S : \quad  w= (a+ib)z    \quad  \mbox{  in } {\mathbb{C}}_z \times {\mathbb{C}}_w,
$$
so that $S/[(1,0),(0,1),(i,0),(i\xi,\xi)]$ is the Lie
subgroup $H$ of ${\bf  T}$. For future use we note that
\begin{align}
 \label{eqn:ell}
S: \ x_2+ix_4= (a+ib)(x_1+ix_3)  \quad  \mbox{  in } {\mathbb{R}}^4
\end{align}
defines a real two-dimensional plane passing through the origin in ${\mathbb{R}}^4$; moreover we may identify $H$ with $S/[e_1,e_2,e_3,e_4]$.

We consider two subcases:
\begin{equation} \label{bcases}
{\rm  (i)}\  b =0 ;  \qquad \quad   {\rm  (ii)}  \  b \ne 0.
\end{equation}
\noindent Note that in case (i) $S=\{x_2=ax_1,\ x_4=ax_3\}$ is the
product of two (not necessarily closed) curves $H_1=\{x_2=ax_1\} \subset
T_1$ and $H_2=\{x_4=ax_3\} \subset T_2$ in ${\bf  T}$.

We also make the distinction between the two options:
\begin{center}
 (i')\ $a$ is rational; $a=q/p$, $p,q$ relatively prime;   \quad  (ii')\ $a$ is
 irrational.
\end{center}

In case (i'), $H_1$ defines  a closed curve $\gamma_1$ in $T_1$ determined by
the segment $[(0,0), (p,q)]$ in ${\mathbb{R}}^2$, while $H_2$ is dense in
$T_2$. Hence $\bar { H}=\gamma_1 \times T_2 $ is a closed real
Lie subgroup in ${\bf  T}$ with corresponding real Lie subalgebra given by
$$
{\textstyle  \mathfrak h  _0= \{ c_1(\frac{\partial }{\partial x_1}+a \frac{\partial
}{\partial x_2})+ c_2 \frac{\partial }{\partial x_3}+ c_3 \frac{\partial
}{\partial x_4 } : c_1,c_2,c_3\in {\mathbb{R}}\}.}$$
Using the projection
$\widetilde { \pi}: {\bf  T} \to {\bf  T}/
\bar { H} $ yields that ${\bf  T}/
\bar { H}$ is diffeomorphic to $ S^1={\mathbb{R}}/[1/q]$. From Remark \ref{re:5-2} there exist
$0\le a_1<a_2< 1/q$  such that
$D=\pi^{-1}(a_1,a_2)=
\cup _{a_1<s<a_2}((s,0;0,0)+ \gamma_1 \times T_2)
$. Thus, in case (i'),
$$D \approx (a_1,a_2)\times \gamma_1 \times T_2$$
as a two-dimensional complex manifold.

Case (ii') will be divided further into two subcases depending on whether $\frac{ 1}a- \xi$ is rational (i.e., $\frac{ 1}a- \xi=p/q$) or irrational. The former case is similar to case (i'):
since $q=a(q\xi+p)$, we have
 $(q \xi+p, q ) \in H_2$. This point is equivalent to $(0,0) \in T_2$, so that
 $H_2$ defines a simple, closed curve $\gamma_2$ in $T_2$ (see Remark
 \ref{rem:6-3}). On the other hand, since $a$ is
 irrational, $H_1$ is dense in $T_1$. Hence $\bar { H}= T_1 \times
 \gamma_2$, which is a real Lie subgroup of ${\bf  T}$. The
 quotient space ${\bf  T}/\bar { H}= {\mathbb{R}}/[1/q]=:S^1$. Again from Remark
 \ref{re:5-2} there exist $0\le a_1<a_2< 1/q$  such that
$$D={\textstyle  \bigcup_{a_1<s<a_2}}\ (T_1 \times \gamma_2+(0,0;s,0)).$$
In case $\frac{ 1}a- \xi$ is irrational, $H_1$ and $H_2$ are dense in $T_1$ and $T_2$. But $H_1
\times H_2 \Subset D$, so that $T_1 \times T_2 \Subset D$, which
yields a contradiction since $D \ne {\bf  T}$. In other words, the case $\frac{ 1}a- \xi$ irrational cannot occur.

\smallskip
We turn to case (ii) in our previous dichotomy (\ref{bcases}): $b\ne 0$. Here $a$ may be either rational or irrational. From
(\ref{eqn:ell}) $ S \subset {\mathbb{R}}^4$ may be written as
\begin{eqnarray*}
\left\{   \begin{array}{ll}
 x_2&=ax_1-bx_3, \medskip \\
x_4& =bx_1+ax_3,
\end{array}\right.   \qquad  (x_1,x_3)\in {\mathbb{R}}_{x_1}\times {\mathbb{R}}_{x_3}
\end{eqnarray*}
or, equivalently,
\begin{eqnarray}
   \label{eqn:fundam}
 \left\{   \begin{array}{ll}
 x_3&=A x_1 +B x_2,   \medskip  \\
x_4&= C x_1 -A x_2,
\end{array}\right.   \qquad  (x_3,x_4)\in {\mathbb{R}}_{x_3}\times {\mathbb{R}}_{x_4},
\end{eqnarray}
where $$
 A= \frac ab ,  \quad  B= - \frac 1b,  \quad  C=
\frac{ a^2+b^2}b \ne 0.
$$
Note that $A,B,C$ are related via the Jacobian:
$$
{   \frac{\partial (x_3,x_4) }{\partial (x_1,x_2)}= -A^2-BC=1.}
$$

We will show the following:

\smallskip
\noindent  {\bf  1.\,(i)}  \ {\it There exists a unique set of six integers $m,n,m'\in {\boldsymbol{Z}};n',p,q\in {\boldsymbol{Z}}^+$ where $(m,n)$,
$(m',n')=\pm 1$, $(p,q)=1$, such that $a,\, b$ can be written in the following
form:
\begin{eqnarray}
 \label{eqn:a-b}
   a = { \frac{ p_1p_2+q_1q_2 }{p_1^2+q_1^2}} ,  \quad
  b =  { \frac{ p_2q_1-p_1q_2}{p_1^2+q_1^2}}\  ,
\end{eqnarray}
where
\begin{align} \label{eqn:defMprime}
 M' &:= m'+n' \xi; \ \ \ p_1:=M'p, \ \  \ p_2:=n'p, \ \  \ q_1:=mq, \ \
 \ q_2: =nq.
\end{align}   }
\noindent {\bf 1.\,(ii)}  \ {\it The integral curve $S$ in (\ref{eqn:ell}) contains the following two  points:
$$
(\, q(m,n),p(M',n')\, ),  \quad  (\, p(1/n,0),q(1/n',0)+ \eta (M',n')\,),
$$
where $\eta=\frac{ p}{nn'} \frac{ p_2^2+q^2_2}{p_2q_1-p_1q_2}$. Note
that $\eta$ is irrational.}

 \bigskip
\noindent {\bf  2.\,(i)}  \  {\it The closure $\bar { H}$ of $H$ in ${\bf  T}$ is
a closed real Lie subgroup of ${\bf  T}$ whose corresponding real Lie subalgebra $\mathfrak h _0$ of $\mathfrak
h$ is generated by
\begin{align}
 \label{eqn:hzero}
\bigl\{& q_1\frac{\partial }{\partial x_1 }+ q_2 \frac{\partial
}{\partial x_2 }, \  p_1 \frac{\partial }{\partial x_3}+ p_2 \frac{\partial
}{\partial x_4},  \nonumber \\
 & \ \ \ ({p_2q_1-p_1q_2})\,  \frac{\partial }{\partial x_1 }+
 ({p_1p_2+q_1q_2}) \frac{ \partial }{\partial x_3} +
 ({p_2^2+q_2^2})\frac{ \partial }{\partial  x_4}     \bigr\}.
\end{align}}

We proceed to give a more precise description of $\bar H$. Assuming {\bf 1.}, let
\begin{eqnarray*}
\begin{array}{llll}
 L_1: &\{(x_1,x_2):mx_2=nx_1\}=\{t(m,n): t\in {\mathbb{R}}\} \  \quad  &\mbox{  in } {\mathbb{R}}_{x_1}\times {\mathbb{R}}_{x_2} ;\\
L_2: &\{(x_3,x_4):M'x_4=n'x_3\}=\{t(M',n'):t\in {\mathbb{R}}\} \ &\mbox{ in }
 {\mathbb{R}}_{x_3}\times {\mathbb{R}}_{x_4}.
\end{array}
\end{eqnarray*}
Since $m,n \in {\boldsymbol{Z}}$, $L_1$ defines a simple closed curve $l_1$
with positive orientation in the real torus $T_1$ (see 1. in  Remark
\ref{rem:6-3}), and from (\ref{eqn:defMprime}), specifically, the
relation $M' = m'+n'\xi$, $L_2$ defines a  simple closed curve $l_2$
with positive orientation in the real torus $T_2$ (see 2. in  Remark
\ref{rem:6-3}). Given $0\le s\le 1$, define
\begin{eqnarray*}
\begin{array}{llll}
 {\cal L}_1(s)&:= L_1+ ps({ 1}/n,0) =\{t(m,n)+ps({ 1}/n,0): t\in {\mathbb{R}}\}
   &\mbox{  in } {\mathbb{R}}_{x_1}\times
 {\mathbb{R}}_{x_2};\\
{\cal L}_2(s)&:=  L_2+ qs({ 1}/{n'},0) =\{t(M',n')+qs({ 1}/{n'},0):t\in {\mathbb{R}}\}   &\mbox{  in }
 {\mathbb{R}}_{x_3}\times {\mathbb{R}}_{x_4}.
\end{array}
\end{eqnarray*}
Then ${\cal L}_1(s)$ and ${\cal L}_2(s)$ also define simple closed curves
$l_1(s)$ and $l_2(s)$ in $T_1$ and $T_2$;
$l_1(s)$ is a  curve in $T_1$, which, in ${\mathbb{R}}_{x_1}\times {\mathbb{R}}_{x_2}$, is parallel to $l_1$ translated by the vector $ps(\frac{1}n,0
)$. Similarly, $l_2(s)$ is parallel to $l_2$ in ${\mathbb{R}}_{x_3}\times {\mathbb{R}}_{x_4}$ translated by the vector $qs(\frac{ 1}{n'},0)$. We have $l_i=l_i(0)=l_i(1)$ for $i=1,2$ and
\begin{equation}
 \label{diff}
(l_1(s')\times l_2(s'))\cap
(l_1(s'')\times l_2(s''))= \emptyset \ \mbox{if} \  s'\ne s''.
\end{equation}
Note that we may have, e.g., $l_1(s')=l_1(s'')$ if $s'\ne s''$ but then
$l_2(s')\ne l_2(s'')$ from $(p,q)=1$ so that (\ref{diff}) holds (cf.,
3. in Remark \ref{rem:6-3}).

\smallskip
\noindent {\bf  2.\,(ii)} \ {\it The set
\begin{align}
 \label{eqn:ovh}
\Sigma :=  \mbox{  $\bigcup _{0\le s \le 1}\, l_1(s) \times l_2(s)$},
\end{align}
is a real, $3$-dimensional compact submanifold
of ${\bf  T}$, and $\bar { H}=\Sigma$. Given $0\le t\le 1$, if we define
$$
\Sigma (t) :=(t,0;0,0) + \bar { H}=(t,0;0,0)+\Sigma,
$$
a coset of $\bar { H}$, then $\Sigma(0)= \Sigma(1)=\Sigma= \bar { H}$ and otherwise $\Sigma(t') \cap \Sigma(t'')=\emptyset $ in ${\bf  T}$ if $t'\ne t''$}.

\smallskip
{\bf  3.}\  {\it We have
\begin{align}
 \label{eqn:tsq}
{\bf  T} = \mbox{$  \bigcup _{0\le t\le 1}\, \Sigma(t)$}.
\end{align}
The quotient space ${\bf  T}/ \bar {H} = {\mathbb{R}}/[1]=
S^1$ and
\begin{align}
 \label{eqn:abd}
D&=\mbox{  $\bigcup _{ t_1< t<t_2} \,  \Sigma  (t)$},
\end{align}
where $0\le  t_1<t_2<1$}.

\smallskip
\noindent We will also prove a converse statement:

\smallskip
{\bf  4.}\ {\it Given integers
$m,n,m'\in {\boldsymbol{Z}}; n',p, q \in {\boldsymbol{Z}}^+$ with $(m,n),$ $(m',n')=\pm 1, (p,q)= 1$, we can find $a, b\in {\mathbb{R}}$ satisfying
(\ref{eqn:a-b}) and (\ref{eqn:defMprime}) to construct a pseudoconvex domain $D \subset {\bf
T}$ with smooth boundary which is not Stein. This domain has the
property that $D( \tau\exp tX)^{-1}
=D$ for all $t\in {\mathbb{C}} $ and for all $\tau\in D$ where $X$ is a nonzero
holomorphic vector field with the property that the Lie subgroup
$H$ of ${\bf  T}$ corresponding to the Lie subalgebra $\mathfrak g _0=\{cX  \in \mathfrak g : c\in {\mathbb{C}}\}$ is equal to $\{w=(a+bi)z\}/[e_1,e_2,e_3,e_4]$. Moreover, every holomorphic function on $D$ is constant.}

\smallskip
We proceed with the proofs of items {\bf 1.} through {\bf 4.}
Following (\ref{eqn:fundam}) we write
\begin{equation}
\begin{split}\label{eqn:F}
 F&: (x_1,x_2) \mapsto (x_3,x_4)=(Ax_1+Bx_2,Cx_1-Ax_2), \\
F^{-1}&: (x_3,x_4) \mapsto (x_1,x_2)=(-Ax_3-Bx_4,-Cx_3+Ax_4)
\end{split}\end{equation}
 so that $F, \, F^{-1}$ are linear mappings and
\begin{align*}
 \label{anal}
H &= \bigl\{(x_1,x_2, F(x_1,x_2))\in {\mathbb{R}}^4 : (x_1, x_2)
 \in  {\mathbb{R}}_{x_1}\times {\mathbb{R}}_{x_2}\bigr\}/ [e_1,e_2,e_3,e_4]  \quad
 \mbox{  in } \ {\bf  T}.
\end{align*}

By Sard's theorem there exists a point $x^0\in D$ such that
$$\nabla
\Lambda(x^0)=(\frac{\partial \Lambda}{\partial x_1 }, \ldots ,
\frac{\partial \Lambda}{\partial x_4 } )(x^0)\ne {\mathbf{0}}.$$
The $c$-Robin function $\Lambda(x)$ is invariant under parallel
translation in ${\mathbb{R}}^4$ since the operator $\Delta -c =\Delta-1$ associated to the Euclidean metric $\sum_{ i=1  }^{ 4  }dx_i^2 $ is invariant under parallel
translation.
Thus if we write $D_1:=D-x^0$ in ${\bf  T}$ and $\Lambda_1(x)$ is the
$c$-Robin function for $D_1$,
then $\Lambda_1( x)= \Lambda(x +x^0)$. Therefore we may assume that $D$ contains the origin $ {\mathbf{0}}$ and $\nabla\,
\Lambda({\mathbf{0}})\ne {\mathbf{0}}$. Thus $H \Subset D$ and $\Lambda(x)\equiv
\Lambda({\mathbf{0}})$ on $H$.

 For  $m,n \in {\boldsymbol{Z}}$, we set
$$
 P(m,n):=(m,n, F(m,n))=(m,n ,A m +B n, C m - A n)\in H.
$$
Note that such a point is equal to
$(0,0, F(m,n))=(0,0, A m +B n, C m - A n)$
in  ${\bf  T} $. We set $ {\cal P}:=\{F(m,n)\in {\mathbb{R}}_{x_3} \times {\mathbb{R}}_{x_4}: m,n\in {\boldsymbol{Z}}\}$ and
$$
\widehat{ {\cal P}}=\{ F(m,n)+ (K,l)\in {\mathbb{R}}_{x_3}\times {\mathbb{R}}_{x_4}
: m,n,k,l \in {\boldsymbol{Z}}  \}
$$
where $K=k+l \xi$. Then $ {\cal  P}/[e_3', e_4']= \widehat{ {\cal
P}}/[e_3',e_4']$; i.e., ${\cal P}$ and $\widehat{ {\cal P}}$ define
the same set in $T_2$. We show:
\medskip

{\it The closure $Cl[ \widehat{ {\cal P}}]$ of  $\widehat{ {\cal P}}$ in ${\mathbb{R}}_{x_3}\times
 {\mathbb{R}}_{x_4} $ consists of an infinite
number of parallel lines the same distance apart, one of which passes through the origin $(0,0)$}.

\smallskip
Note that, by definition, $\widehat {\cal P}$ is an additive subgroup in
 ${\mathbb{R}}_{x_3}\times {\mathbb{R}}_{x_4} $. Moreover, $\widehat {\cal P}$ cannot lie on a single line so that $Cl[\widehat{ {\cal P}}]$ is not a single line containing the origin. To verify the italicized statement we must first rule out two other possibilities:
\begin{enumerate}
 \item [i.] $\widehat{ {\cal P}}$ is an isolated set in ${\mathbb{R}}_{x_3}\times
 {\mathbb{R}}_{x_4} $;
  \item [ii.] \  $Cl[ \widehat{ {\cal P}}]= {\mathbb{R}}_{x_3}\times
 {\mathbb{R}}_{x_4} $.
\end{enumerate}
To prove that  i. does
not occur, let $\Delta$ be a fundamental parallelogram  of $T_2$
with vertices $(0,0),(1,0),(\xi,1)$, and $(1+\xi,1)$. Then each
point $P(m,n)$ in ${\cal P}$ is equivalent to some point $P'(m,n)$ in
$\Delta$. Since $\Delta$ is bounded, it is enough to show that $\{P'(m,n)\}_{m,n\in {\boldsymbol{Z}}}$ consists of
infinitely many distinct points in $\Delta$. We first verify that for any three pairs $(m_i,n_i)\in {\boldsymbol{Z}} \times {\boldsymbol{Z}}, \
i=0,1,2$,
\begin{eqnarray}\label{1sts}
P'(m_0,n_0)=P'(m_1,n_1)=P'(m_2,n_2)  \   \hbox{implies}  \   \frac{ m_1-m_0}{n_1-n_0}=
\frac{ m_2-m_0}{n_2-n_0}.  \quad
\end{eqnarray}
For the sake of obtaining a contradiction, suppose $\frac{ m_1-m_0}{n_1-n_0}\ne
\frac{ m_2-m_0}{n_2-n_0}.$
The condition $P'(m_0,n_0)=P'(m_1,n_1)=P'(m_2,n_2)$ implies that  there exist $(p_i,q_i)\in {\boldsymbol{Z}} \times {\boldsymbol{Z}}, \ i=0,1,2$  such
that the three points $$\{(A m_i +B n_i -(p_i+q_i \xi), C  m_i - A
       n_i-q_i)\}_{i=0,1,2}$$ are the same in $\Delta$.
This means that
\begin{eqnarray} \nonumber
&{} A(m_1-m_0) +B (n_1-n_0)=(p_1-p_0) +(q_1-q_0) \xi,\\ \nonumber
&{} A (m_2-m_0) +B(n_2-n_0)=(p_2-p_0) +(q_2-q_0) \xi, \label{eqn:easy}\\
 &{}C (m_1-m_0) - A (n_1-n_0)=q_1-q_0,\\ \nonumber
 &{}C (m_2-m_0) - A(n_2-n_0)=q_2-q_0.\nonumber
\end{eqnarray}
From the last two equations it follows that $C$ and $A$ are rational. Using the relation $A^2+BC+1=0$, $B$ must be rational as well. Since $\xi $ is irrational, we have
$q_1-q_0=q_2-q_0=0$; hence $A=C=0$, yielding a contradiction. This proves (\ref{1sts}). Now suppose $\{P'(m,n)\}_{m,n\in {\boldsymbol{Z}}}$ is a finite set in $\Delta$. Then $\{\frac{ m_1-m_0}{n_1-n_0}\}_{m_0,m_1,n_0,n_1\in {\boldsymbol{Z}}}$ would have to be a finite set of rational numbers, which it is not.

We next rule out case ii. Assume, for the sake of obtaining a contradiction, that
$Cl[\widehat{ \cal P}]= {\mathbb{R}}_{x_3}\times {\mathbb{R}}_{x_4}$. Since
$\Lambda(x)\equiv \Lambda({\mathbf{0}})$ on $H$ and since $(0,0;{\cal P})\subset H$, from continuity it follows that $\Lambda(x)\equiv \Lambda({\mathbf{0}})$ on $(0,0;Cl[ \widehat{
{\cal P}}])$. Hence $\Lambda(0,0;
x_3,x_4)\equiv  \Lambda({\mathbf{0}})$ on $T_2$.
In particular,
$
\frac{\partial \Lambda}{\partial x_3}= \frac{\partial \Lambda}{\partial
x_4 }=0 $  at ${\mathbf{0}}$. Moreover, since $\Lambda \equiv \Lambda({\mathbf{0}})$ on
$H=\{w= (a+bi)z\}/[e_1,e_2,e_3,e_4]$, i.e.,
$\Lambda(z,(a+bi)z) \equiv \Lambda( {\mathbf{0}})$ in ${\mathbb{C}}_z$ where $z=x_1+ix_3$, we have
$$\frac{\partial \Lambda}{\partial z}+ (a+bi)\frac{\partial \Lambda}{\partial w}
\equiv 0$$  on $ H$. Equivalently,
$$
\frac{\partial \Lambda}{\partial x_1}- i\frac{\partial \Lambda}{\partial
x_3}
 +(a+bi)( \frac{\partial \Lambda}{\partial x_2}- i\frac{\partial \Lambda}{\partial
x_4})\equiv 0   \quad  \mbox{  at } \ {\mathbf{0}}.
$$
Since $\frac{\partial \Lambda}{\partial x_3}= \frac{\partial \Lambda}{\partial x_4 }=0 $ at ${\mathbf{0}}$ and $b\ne 0$ it follows that $
\frac{\partial \Lambda }{\partial x_1}= \frac{\partial \Lambda
}{\partial x_2}   =0$ at ${\mathbf{0}}$. This contradicts the hypothesis that $\nabla \Lambda({\mathbf{0}})\ne {\mathbf{0}}$.

We conclude that $ Cl[\widehat{ {\cal
P}}]=\cup _{\nu=-\infty}^\infty L_2^{\nu}$, where $L_2:=L_2^0$ passes through the origin, and $L_2^{\nu}=L_2+\nu
(d,0)$ for some $d>0$; i.e., $L_2^{\nu}$ is parallel to $L_2$
translated by $\nu(d,0)$ along the $x_3$-axis. By a similar calculation as in
(\ref{eqn:easy}) one can check that the line $L_2$ cannot be either the $x_3$-axis or
the $x_4$-axis.

We claim that the line $L_2$ defines a simple closed curve $l_2:=L_2/ [(1,0),
 (\xi,1)]$ in $T_2$. Indeed, if not, $L_2$ defines a non-closed curve which must be
dense in $T_2$. Since $\Lambda(0,0; L_2)\equiv \Lambda({\mathbf{0}})$, we have
$\Lambda(0,0;T_2 )= \Lambda({\mathbf{0}})$. As above, this contradicts $\nabla \,
\Lambda({\mathbf{0}})\ne {\mathbf{0}}$.

There thus exists a unique $m'\in {\boldsymbol{Z}}, \ n'\in {\boldsymbol{Z}}^+$ with  $(m',n')=\pm1$  such
that $L_2$ passes through
$$
m'(1,0)+n'(\xi, 1)=(m'+n' \xi,n')= : (M',n')
$$
which is equal to $(0,0)$ as a point in $T_2$. Following 2. in Remark \ref{rem:6-3}, we give $l_2$ a positive orientation determined from that of $L_2$ directed from $(0,0)$ to $(M',n')$. Considering the standard
fundamental  domain $(\Delta'')$ of $T_2$ associated to $L_2$ defined by
$$ \{(0,0),(M_1',n_1'), (M',n'), (M_1'+M', n_1'+n')\}$$
 with $m_1'n_1-m_1n_1'=1$, we see that the distance $d $ between successive curves $L_2^{\nu}$, e.g., between $L_2 $ and $L_2(1)$, must be of the form
$d=\frac l{kn'} $ where $k,l \in {\boldsymbol{Z}}$ with $(k,l)=\pm 1$. For otherwise $\cup_{ \nu=-\infty  }^{ \infty   }L_2^{\nu} $ is dense in
${\mathbb{R}}_{x_3}\times {\mathbb{R}}_{x_4}$. Since $(0,0;\widehat{ {\cal P}})\subset H$ and hence
 $(0,0; Cl[ \widehat{ {\cal P}}])= (0,0;\cup_{ \nu=-\infty
 }^{ \infty   }L_2^{\nu})\subset \bar H $, we again get a contradiction to $\nabla \,
\Lambda({\mathbf{0}})\ne {\mathbf{0}}$. We will determine $k$ and $l$ explicitly in Remark
\ref{rem:6-2}.

\smallskip
Define $L_1:=F^{-1}(L_2)\subset {\mathbb{R}}_{x_1}\times
{\mathbb{R}}_{x_2}$. We show:
\begin{eqnarray}
\begin{array}{llll}
 &1. & L_1 \times L_2 \subset \bar { H}  \ \ \  \mbox{  and }
  \ \ \    \ H+L_1
  \times L_2 \subset \bar { H}; \nonumber \\
  &2. & L_1 \mbox{  defines a simple closed curve
 $l_1$ in $T_1$;} \nonumber
\label{eqn:zyuu-2}\\
&3.  & \Lambda (x)\equiv \ \Lambda({\mathbf{0}})  \  \  \mbox{  on } \
  L_1  \times L_2. \nonumber
\end{array}
\end{eqnarray}
We first prove the following  inclusion:
$$
(*) \quad    H+(0,0;L_2) \subset \bar H.
$$
To prove this, fix $(x_1,x_2;F(x_1,x_2))\in H$
 and $ F(m,n) +[k(1,0)+l(\xi,1)] \in \widehat{
{\cal P}}$. Then, as points in ${\bf  T}$,
\begin{align*}
& (x_1,x_2; F (x_1,x_2)+ F (m,n)+[k(1,0)+l(\xi,1)] )\\ &=(x_1+m,x_2+n;F
 (x_1+m,x_2+n)).
\end{align*}
By definition, this point is in $H$. We conclude that for all $(x_1,x_2;F(x_1,x_2))\in H$ we have $(x_1,x_2;F (x_1,x_2)+ \widehat{ {\cal P}}) \subset H$.
 Since $L_2\subset Cl[\widehat{ {\cal P}}]$, it follows
that  $(*)$ holds.

We verify 1. and 3.: Fix $p:=(x_1^0,x_2^0;x_3',x_4')\in
 L_1\times L_2$. Then there exists $(x_3^0,x_4^0)\in L_2 $ with $F(x_1^0,x_2^0)=(x_3^0,x_4^0)$ so that
$$
p=(x_1^0,x_2^0;x_3',x_4')= (x_1^0,x_2^0;F(x_1^0,x_2^0))+(0,0; x_3'-x_3^0, x_4'-x_4^0).
$$
Since $(x_3'-x_3^0, x_4'-x_4^0)\in L_2-L_2=L_2$, we have $p\in
 H+(0,0;L_2)$. Then $(*)$ implies $p\in \bar { H}$,
so that  $L_1 \times L_2 \subset \bar { H}$. Since $\bar { H}$ is an additive group, we conclude that $H+L_1 \times L_2 \subset \bar {
 H}$ and 1. is proved. Since $\Lambda(x) \equiv
 \Lambda({\mathbf{0}})$  on $\bar { H}$, 1. implies 3.

We prove 2. by contradiction. If 2. is not true, then $L_1$ is
dense in ${\mathbb{R}}_{x_1}\times {\mathbb{R}}_{x_2}  $, or, equivalently,
$L_1/[e_1',e_2']$ is dense in $T_1$. It
follows from 1. that $T_1 \times L_2 \subset \bar { H}$, and hence
$\Lambda(x)\equiv \Lambda({\mathbf{0}})$ on $T_1
 \times L_2$. In particular, this holds on $(T_1;0,0)$.
Using an argument similar to that used in the proof that $Cl[ \widehat{ {\cal P}}]\not = {\mathbb{R}}_{x_3}\times
 {\mathbb{R}}_{x_4} $, this together with $\Lambda(x_1,x_2; F (x_1,x_2))\equiv
\Lambda({\mathbf{0}})$ for $(x_1,x_2)\in {\mathbb{R}}_{x_1}\times {\mathbb{R}}_{x_2} $
implies $\nabla \Lambda({\mathbf{0}})={\mathbf{0}}$, a contradiction.

\smallskip
Before proceeding, we make a remark about the orientation of the curves $l_1$ and
$l_2$. Recall we fix the direction of the line $L_2$ as that given by the directed line segment from $(0,0)$ to $(M',n')$ in
${\mathbb{R}}_{x_3}\times {\mathbb{R}}_{x_4} $; this determines the positive orientation of $l_2$ in
$T_2$. Similary, we fix the direction of $L_1$ as that given by the directed line segment from
$(0,0)$ to $F^{-1}(M',n')$ in ${\mathbb{R}}_{x_1}\times {\mathbb{R}}_{x_2}$,
which induces the positive orientation of $l_1$ in $T_1$. We say that $F$ maps the
positively oriented curve $l_1$
to the positively oriented curve $l_2$.
It follows from 1. in Remark  \ref{rem:6-3} and assertion 2. that
there exist unique $m,n\in {\boldsymbol{Z}}$ with the property that $(m,n)=\pm 1$ such that the closed
curve $l_1$ has a positive orientation associated with that of the directed line segment from $(0,0)$ to $(m,n)$.

\smallskip
We replace ${\cal P}=\{F(m,n) \in {\mathbb{R}}_{x_3}\times
{\mathbb{R}}_{x_4}  : m,n\in {\boldsymbol{Z}}\}$ in $T_1$ by ${\cal Q}:=\{F^{-1}(M',n')\in
{\mathbb{R}}_{x_1}\times {\mathbb{R}}_{x_2} : m',n'\in {\boldsymbol{Z}}\}$ in $T_2$, where $M'=m'+n'\xi$. Similarly, we replace
$$
\widehat{ {\cal P}}=\{ F(m,n)+ (K,l)\in {\mathbb{R}}_{x_3}\times {\mathbb{R}}_{x_4}
: m,n,k,l \in {\boldsymbol{Z}}  \}
$$
by
$$
\widehat{ {\cal Q}} =\{F^{-1}(M',n')+(k',l') \in {\mathbb{R}}_{x_1}\times
{\mathbb{R}}_{x_2}: m',n', k',l'\in {\boldsymbol{Z}}\}.
$$
As we found $L_2\subset  Cl[ \widehat{ {\cal P}}]$, there exists a line
$\widetilde { L}_1 \subset Cl[\widehat{ {\cal Q}}]$ in ${\mathbb{R}}_{x_1}\times {\mathbb{R}}_{x_2}  $ passing through the origin
$(0,0)$. This line is neither the $x_1$-axis nor the $x_2$-axis. We claim that
$$
\widetilde { L}_1=L_1    \quad  \mbox{  in } {\mathbb{R}}_{x_1}\times
{\mathbb{R}}_{x_2}.
 $$

 To verify this we set $\widetilde { L}_2=F(\widetilde {
L}_1)$ and we assume, for the sake of obtaining a contradiction, that $\widetilde L_2 \ne L_2$. From 3. applied to  $\widetilde L_1$ and $\widetilde L_2$, we have $\Lambda(x)\equiv  \Lambda({\mathbf{0}})$ on $\widetilde { L}_1 \times
 \widetilde { L}_2$. It follows that $\Lambda(x)\equiv  \Lambda({\mathbf{0}})$
on the following four distinct lines passing through the origin in ${\mathbb{R}}^4$:
$$(L_1;0,0),  \quad   (0,0;L_2),  \quad   (\widetilde { L}_1;0,0),  \quad
(0,0;\widetilde { L}_2).
$$
We conclude that $\nabla \, \Lambda({\mathbf{0}})={\mathbf{0}}$, a contradiction.

We define the lines
$ L_2(\nu):=L_2 +\nu (\frac{ 1}{n'},0),\  \nu=0,\pm 1,\pm 2,\ldots $,  and $
L_1(\nu):=F^{-1}(L_2(\nu))$. Then $L_2(0)=L_2$; \, $L_1(0)=L_1$, and the
lines $\{L_1(\nu)\}_{\nu=0,\pm 1, \pm 2,...}$ are parallel translates of $L_1$ in ${\mathbb{R}}_{x_1}\times {\mathbb{R}}_{x_2}   $ by integer multiples of the vector $(d^*,0)$ for some $d^*>0$. We note that each line $L_2(\nu)$ defines the same closed curve $l_2$ in $T_2$, and  each $L_1(\nu)$ defines a  closed curve $l_1(\nu)$ in $T_1$, but the curves $l_1(\nu)$ for distinct $\nu$ may be different. We note also that
$$
L_1(\nu)\times L_2(\nu ) \subset \bar{
H}.
$$

To verify this, fix  a point $(x_1^0, x_2^0)\in L_1(\nu)$ so that
$F(x_1^0,x_2^0)\in L_2(\nu)$. Since $L_1(\nu)= (x_1^0, x_2^0)+L_1$ and
$L_2(\nu)= F(x_1^0,x_2^0)+L_2$, it follows that $L_1(\nu) \times L_2
(\nu)= ((x_1^0, x_2^0)+L_1, F(x_1^0, x_2^0)+L_2)\in H+L_1 \times L_2$. From 1., $H+L_1 \times L_2
\subset \bar { H}$ and the result follows. We thus have that
$$\Lambda(z)\equiv \Lambda({\mathbf{0}})  \quad   \mbox{on} \  l_1(\nu) \times l_2(\nu),  \ \nu=0,\pm 1,\ldots .$$

Our next goal is to show:
 \begin{enumerate}
  \item [$(*)$] there exist unique $p,q\in {\boldsymbol{Z}}^+$ with  $(p,q)=1$  such
 that

$L_1(q)=
L_1+p\,(\frac{1}n,0 )$,
 or, equivalently, $L_2+q\,(\frac{ 1}{n'},0) = F(L_1+ p\,(\frac{ 1}{n},0))$.
\end{enumerate}
We prove $(*)$ by contradiction. If $(*)$ is false, then writing
$d^*=d_0\, \frac{ 1}n$, we see that $d_0$ must be irrational.
Hence the set of curves $l_{1}(\nu), \ \nu=\pm 1,\pm 2,\ldots $ is
dense in $T_1$. By the previous paragraph, $\Lambda \equiv \Lambda({\mathbf{0}})$
on $T_1 \times l_2$, which contradicts $\nabla \Lambda({\mathbf{0}})\ne {\mathbf{0}}.$
Thus $(*)$ is proved.

For $s\in (-\infty, \infty)$, we define
\begin{eqnarray*}
\begin{array}{llll}
{\cal L}_1(s)&:= L_1+ ps({ 1}/n,0) \quad  &\mbox{ in } \ {\mathbb{R}}_{x_1}\times
 {\mathbb{R}}_{x_2};\\
{\cal L}_2(s)&:=  L_2+ qs({ 1}/{n'},0) \quad  &\mbox{ in } \
 {\mathbb{R}}_{x_3}\times {\mathbb{R}}_{x_4}
\end{array}
\end{eqnarray*}
so that ${\cal L}_2(s)= F ( {\cal L} _1(s))$.
Each line ${\cal L}_1(s)$ defines a simple closed curve $l_1(s)$
in $T_1$, and each line ${\cal L}_2(s)$ defines a simple closed curve
$l_2(s)$ in $T_2$. Apriori, we define
$$
\Sigma:= \mbox{  $\bigcup_{-\infty < s < \infty}$}\, l_1(s)\times l_2(s) \subset {\bf  T}.
$$
However, ${\cal L}_i(s+k)=
{\cal L} _i(s), \ i=1,2$ as a set in $T_i$, for all $k\in {\boldsymbol{Z}}$, so that
$$
\Sigma= \mbox{  $\bigcup_{0\le s\le 1}$}\, l_1(s)\times l_2(s) \subset {\bf  T}.
$$
By the definition of $p,q\in {\boldsymbol{Z}}^+$
 we have $({\cal L} _1(s) \times {\cal L}_2(s)) \cap
({\cal L} _1(s') \times {\cal L}_2(s'))=\emptyset$ for $s\ne s', \ 0\le
s,s'\le 1$, so that
\begin{eqnarray}
 (l_1(s) \times l_2(s))\  \cap\ ( l_1(s') \times l_2(s'))=\emptyset,  \quad
s\ne s', \ 0\le s<s'\le 1. \label{eqn:remain}
\end{eqnarray}
Then, since $l_1(s+1)=l_1(s), \ l _2(s+1)=l_2(s)$,
and each $l_1(s), \ l_2(s)$ is a simple closed curve, it follows that  $\Sigma$ is a three-dimensional, compact manifold in ${\bf  T}$.

We show that
\begin{equation}
\label{barh}
 \Sigma = \bar { H} =H+ L_1 \times L_2  \quad  \mbox{  in } \ T_1 \times
 T_2.
 \end{equation}
We proved in  1. that $H+L_1 \times L_2 \subset \bar { H}$. We next show $\Sigma \subset H+L_1 \times L_2$. Since
$${\cal L} _2(s)=F({\cal L} _1(s))=F(L_1+ps({1}/{n},0 ))=
 L_2+F(ps({1}/{n},0 )),$$
it follows that
\begin{align*}
 {\cal L} _1(s) \times {\cal L}_2(s)
&= \bigl((ps({1}/{n},0 )+L_1;
F(ps({1}/{n},0 )+L_2 \bigr)\\
&= \bigl( \,ps({1}/{n},0);F(ps({1}/{n},0 ))\bigr)+(L_1, L_2 )
\in H+ L_1 \times L_2.
\end{align*}
Finally we show $H \subset \Sigma$; and hence $\bar { H} \subset \bar \Sigma =\Sigma$. Fix $((x_1,x_2),F(x_1,x_2))\in H$. There is a unique $k\in
{\boldsymbol{Z}}$ and $0\le s <1$
  such that $(x_1,x_2)\in {\cal L} _1(k+s)$, or, equivalently,
$ {\cal L}_1(k+s)= (x_1,x_2)+L_1 $. Then we have $F(x_1,x_2) \in {\cal
  L} _2(k+s)$,  and hence $ {\cal L}_2(k+s)= F(x_1,x_2)+L_2 $.
Since ${\cal L} _1(k+s)= {\cal L} _1(s)$
  in $T_1$ and ${\cal L} _2(k+s)= {\cal L}_2(s)$ in $T_2$, it follows
  that
$$
((x_1,x_2),F(x_1,x_2))\in ((x_1,x_2)+L_1, F(x_1,x_2)+L_2) =({\cal L}
_1(s), {\cal L}_2(s)),
$$
which proves the inclusion $H \subset \Sigma$.

\smallskip
Equation (\ref{barh}) will lead us to the concrete formulas for the quotient
space ${\bf  T}/ \bar { H}$ in ${\bf 3.}$ and the canonical projection
$\widetilde { \pi}:\ {\bf  T} \mapsto
{\bf  T}/ \bar { H}$.

\smallskip For $-\infty <t <\infty$, define
$$
\Sigma(t):= (pt ({ 1}/n,0); 0,0) + \Sigma  \quad  \mbox{  in } \
{\bf  T},
$$
which is a coset of $\Sigma$ and hence is a real three-dimensional compact manifold in ${\bf  T}$. To be precise,
\begin{align*}
 \Sigma (t)&= \mbox{  $\bigcup_{0\le s\le 1}$} \,  (pt(1/n,0)+ {\cal L}_1(s),
{\cal L}_2(s))/[e_1,e_2,e_3, e_4]\\
&= \mbox{  $
\bigcup_{0\le s\le 1}$} \,  ({\cal L}_1(t+s),
{\cal L}_2(s))/[e_1,e_2,e_3,e_4]\\
&=\mbox{  $\bigcup_{0\le s\le 1}$} \,  l_1(t+s) \times l_2(s).
\end{align*}
By (\ref{eqn:remain}), we have
$\Sigma(t)=\Sigma(t')$ for $t-t'\in {\boldsymbol{Z}}$ and $\Sigma(t)\cap
\Sigma(t')=\emptyset$ for $ t-t' \not \in {\boldsymbol{Z}}$.
We shall show
$$
(**) \qquad {\bf  T}= \mbox{ $\bigcup _{0\le t\le 1}\, \Sigma(t)$}.
$$
Note since $l_1(s)=l_1(s+1)$ in $T_1$, this is equivalent to (\ref{eqn:tsq}). Fix $p_0=(x_1,x_2;x_3,x_4)/[e_1,e_2,e_3, e_4] \in
{\bf  T} $.  Clearly we can choose real numbers $s_1,s_2$
so that $s_1p\,( { 1}/n,0)$ is the $x_1-$intercept of the line $(x_1,x_2)+L_1$ and $s_2q\,( { 1}/{n'},0)$ is the $x_3-$intercept of the line $(x_3,x_4)+L_2$. Since ${\cal L} _1(s_1)= {\cal L} _1(s_2) + t'p(1/n,0)$ where
$t'=s_1-s_2$, it follows that
\begin{align*}
&(x_1,x_2;x_3,x_4)\in ({\cal L} _1(s_1), {\cal L}_2(s_2) )= ({\cal L} _1(s_2), {\cal L} _2(s_2)) + (t'p(1/n,0); 0,0),
\end{align*}
showing that $p_0\in \Sigma + (l_2(t');0,0) = \Sigma(t')$.
We now simply choose $n_0\in {\boldsymbol{Z}}$ and $ t\in [0,1)$ so that $t'=n_0+t$;  then $\Sigma(t)= \Sigma(t')$, and hence $p_0\in \Sigma(t)$, proving $(**)$ and thus ${\bf 2.(ii)}$. Since $n_0,t$ are uniquely determined by the point
$p_0=(x_1,x_2,x_3,x_4)/[e_1,e_2,e_3, e_4 ]$ in ${\bf
T}$, it follows that
$$
{\bf T}/\bar { H}  \approx \mathbb{R}/[1] =S^1 \ \mbox{as real
manifolds};$$
$\widetilde
{\pi}:  {\bf T} \to {\bf T}/\bar { H}$ where
$\widetilde {\pi}(p_0)= pt( {1}/n,0)+ \Sigma,$
or equivalently, $\widetilde { \pi}(p_0)=t$.

\smallskip
Applying Remark \ref{re:5-2} to our domain $D \subset {\bf  T}$, there exist $0\le
t_1<t_2<1$  such that $D= \cup_{t_1<t<t_2} \Sigma(t)$.
Thus we have proved ${\bf  3.}$ The proofs of
 ${\bf  1.}$ and ${\bf  2.(i)}$ are elementary calculations. Note that
 $F(L_1)=L_2$; $L_1$ is the line in ${\mathbb{R}}_{x_1}\times {\mathbb{R}}_{x_2}  $ of
 slope $n/m \ne 0, \infty$ passing through the origin; and
$L_2$ is the line in ${\mathbb{R}}_{x_3}\times {\mathbb{R}}_{x_4}  $ of
slope $n'/M' \ne 0, \infty$ passing through the origin.
It follows from (\ref{eqn:fundam}) that
$$
\frac{ C-A(n/m)}{A+B(n/m)}=\frac{ n'}{M'}.
$$
Furthermore, since  $F({\cal L}_1)={\cal L}_2(1)$, the point $F(p/n,0)$ lies
on the line ${\cal L}_2(1): \ x_4= \frac{ n'}{M'}(x_3-\frac{
q}{n'}).$ Hence
$$
\frac{C p}n=\frac{ n'}{M'}\ (\frac{A p}{n}- \frac{ q}{n'}).
$$
From these two equations and the
Jacobian condition $A^2+BC+1=0$, a simple calculation yields
 $$
A=\frac{ p_1p_2+q_1q_2}{p_2q_1-p_1q_2},  \quad
B=-\frac{ p_1^2+q_1^2}{p_2q_1-p_1q_2},  \quad  C=\frac{ p_2^2+q_2^2}{p_2q_1-p_1q_2},
$$
where $p_1,p_2,q_1,q_2$ are defined by (\ref{eqn:defMprime}). Then we
have
\begin{align}
 \label{ab-no-hyouzi}
b=-1/B=\frac{p_2q_1-p_1q_2}{ p_1^2+q_1^2},  \qquad  a=bA= \frac{
p_1p_2+q_1q_2}{p_1^2+q_1^2},
\end{align}
which proves (\ref{eqn:a-b}) and hence ${\bf  1. (i)}$. By direct calculation
$$
F(m,n)=(Am+Bn,Cm-An)=\frac{ p}{q}(M',n'),$$ so that
$(\,q(m,n), p(M',n')\,)$ lies on the integral curve $S$ in (\ref{eqn:ell}).  We also have
$$ F(1/n,0)=(A/n, C/n)=\frac{
q}{p}(1/n',0)+ \frac{ C}{nn'}(M',n').$$
Since $((p/n,0),F(p/n,0))\in S$, we have $(\,
p(1/n,0),q(1/n',0)+\eta(M',n')\,)\in S$, where $\eta=\frac{
pC}{nn'}=\frac{ p}{nn'} \frac{ p_2^2+q^2_2}{p_2q_1-p_1q_2}$.
This proves ${\bf  1.(ii)}$.

Since $F(p/n,0)=(Ap/n,Cp/n)$, using $F(0,0)=(0,0)$ it follows from the
  linearity of $F$ that
 $(x_1,0; Ax_1,  Cx_1) \in
H$. This together with $$(L_1,0,0), (0,0,L_2) \subset \bar { H}$$ implies that the real Lie subalgebra $\mathfrak h _0$ of $\mathfrak h $ corresponding to
the real three-dimensional Lie subgroup $\bar { H}$  of ${\bf  T}$
is generated by $\{X_1, X_2, X_3\}$ in $\mathfrak g $  where
\begin{align*}
X_1&= m \frac{\partial  }{\partial x_1} +n \frac{\partial }{\partial
 x_2}, \ \   X_2=M'\frac{\partial  }{\partial x_3} +n' \frac{\partial
 }{\partial x_4},  \ \  X_3= \frac{\partial }{\partial x_1}+ A \frac{\partial }{\partial
 x_3}+C \frac{\partial }{\partial x_4 }.
\end{align*}
Substituting $A$ and $C$ in the above formulas, we have by simple
calculation that $\{q X_1, p X_2, (p_1^2+p_2^2)X_3\} $ is equal to
(\ref{eqn:hzero}). Thus ${\bf  2.}$ is proved.

It remains to prove ${\bf  4.}$ Given integers $m,n,m',n',p,q$ as in ${\bf  4}$, we define $a,b$ by (\ref{eqn:a-b}).
Using the same method as in the above proof, we construct $D$ of the form (\ref{eqn:abd}). Since each $\Sigma(t)$ is
 a smooth, closed, Levi flat surface in ${\bf  T}$ and $\Sigma(t)\Subset D$ for $t_1<t<t_2$, it follows that $D \subset {\bf  T}$
is a pseudoconvex domain with smooth boundary which is not Stein. Moreover,
there are no nonconstant holomorphic functions on $D$. To verify this, we may asume ${\mathbf{0}} \in D$. Suppose $f$ is a holomorphic function on
$D$. Since $H$ is an analytic curve which is relatively compact in $D$, $f\equiv c$ on $H$. Hence $f\equiv c$ on
$\bar { H}=\Sigma(0)=\Sigma$. Since $\Sigma$ is a real three-dimensional
manifold, $f\equiv c$ on $D$.
\begin{remark} \label{rem:6-2} \
  {\rm The explicit formulas
 (\ref{ab-no-hyouzi}) for $a$ and $b$ in terms of
 $m,n,m'\in {\boldsymbol{Z}}, \ n',p,q\in {\boldsymbol{Z}}^+ $ and hence the explicit formulas involving the mapping
 $F: \ (x_1,x_2) \mapsto (x_3,x_4)$ imply the following results.

\smallskip
  {1.} \  Recall that  $\widehat{ {\cal P}}:=\{F(i,j)+(K,l)\in
 {\mathbb{R}}_{x_3}\times
{\mathbb{R}}_{x_4} : i,j,k,l\in {\boldsymbol{Z}}\}$ where $K=k+l\xi$. We showed that
$Cl[ \widehat{ {\cal P}}]$ consists of the family of parallel lines $L_2^\nu, \ \nu=0,\pm 1, \pm 2, \ldots $, where $L_2^\nu=L_2+ \nu(d,0)$ and
$L_2=\{t(M',n')\in {\mathbb{R}}_{x_3}\times {\mathbb{R}}_{x_4} : t\in {\mathbb{R}}\}$. We
shall
show that $d=\frac{ 1}{pn'}$. Fix $F(i,j)+(K,l)\in
\widehat{ {\cal P}}$. From our previous calculations, we have
\begin{align}
 \label{eqn:explicit-ab}
F(m,n)&=\frac{ p}q (M',n') \ \ \ \mbox{and} \ \ \   F(1/n,0)= \frac{ q}p(1/n',0)+ \eta (M',n')
\end{align}
where $\eta$ is irrational. Since
\begin{align*}
 (i,j)&= \frac{ j}n(m,n)+ (ni-jm)(1/n,0),  \quad  (K,l)=\frac{
 l}{n'}(M',n')+
(n'k-m'l)(1/n',0),
\end{align*}
it follows that
\begin{align*}
F(i,j)+(K,l)&= \frac{ jp}{nq} (M',n') +(ni-jm)( \frac{q}{p}(1/n',0)+ \eta(M',n'))
\\ & \quad \qquad \qquad \quad  + \frac{
 l}{n'}(M',n')  +
(n'k-ml)(1/n',0)\\
&=t\,(M',n')+ \nu\, (\frac{ 1}{pn'},0),
\end{align*}
where $t= \frac{ jp}{nq}+ \frac{ l}{n'}+(ni-jm)\eta \in {\mathbb{R}}$ and  $\nu= (ni-jm)q+(n'k-m'l)p\in {\boldsymbol{Z}}$.
This shows that
$\widehat{\cal P} \subset L_1 +\cup_{\nu=-\infty}^\infty\nu(\frac{
 1}{pn'},0)$, so that $d/ \frac{ 1}{pn'}$ is a positive integer. On the other
 hand, in the above proof, from $(p,q)=1, (m,n)=\pm 1$ and $(m',n')=\pm 1$ we can choose $i,j,k,l\in {\boldsymbol{Z}}$  such that
$\nu=1$.  Hence $\frac{ 1}{pn'}/d$ is an integer, so that $d=\frac{
 1}{pn'}$.

\smallskip
 {2.} \ Let $$\Delta_1:
(0,0), (m_1,n_1), (m,n), (m_1+m,n_1+n)$$ and
$$\Delta_2:
(0,0), (M'_1,n'_1), (M',n'), (M'_1+M', n_1'+n')$$ be fundamental
domains associated to $L_1$ and
$L_2 $ in ${\mathbb{R}}_{x_1}\times
{\mathbb{R}}_{x_2}  $ and ${\mathbb{R}}_{x_3}\times {\mathbb{R}}_{x_4} $.
Since $(m,n)= (\frac{ 1}{n},0)+ \frac{ n_1}{n} (m,n)$, we obtain
from
(\ref{eqn:explicit-ab})
\begin{align*}
F(q(m,n)) &=p(M',n'),\\
 F(p(m_1,n_1))&=qF(\frac{ 1}{n},0)+ \frac{ pn_1}{n} F(m,n)
= q(\frac{ 1}{n'},0) + \eta''(M',n'),
\end{align*}
where $\eta''= \frac{ p^2n_1}{nq}+ pn$, so that $\eta''$ is irrational.
}
\end{remark}

\vspace{ 5mm}


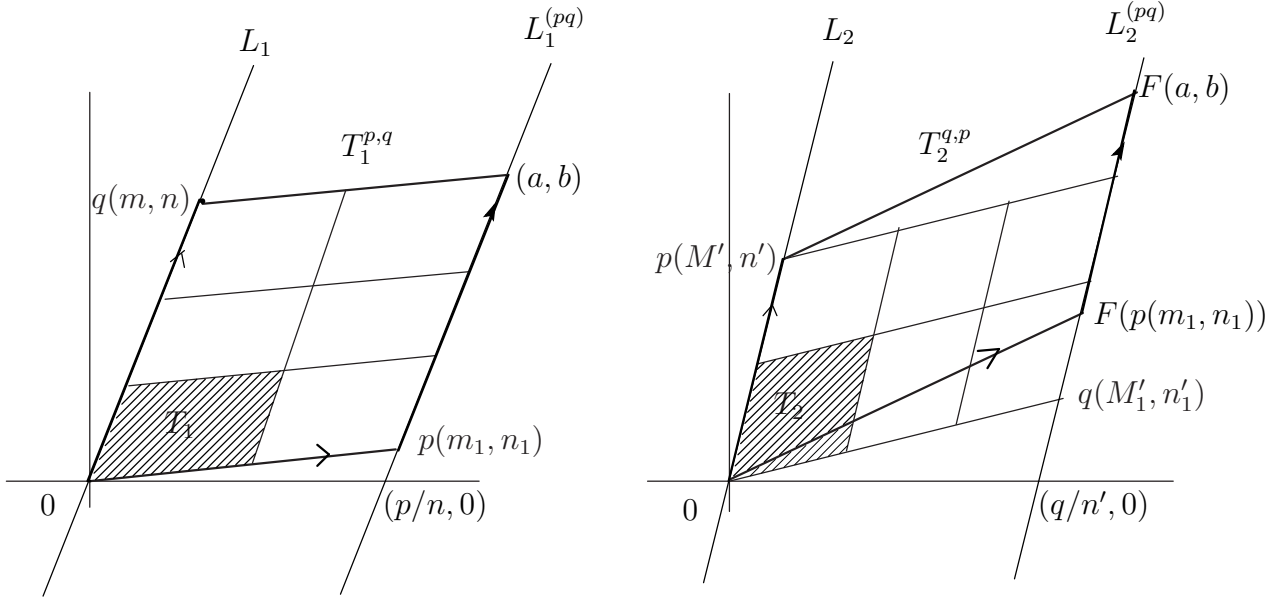
\unitlength 0.1in
\begin{picture}( 60.4100, 27.0500)( 15.8000, -33.5500)
%
\special{pn 8}%
\special{pa 1112 940}%
\special{pa 1112 3106}%
\special{fp}%
%
\special{pn 8}%
\special{pa 680 2970}%
\special{pa 3128 2970}%
\special{fp}%
%
\special{pn 13}%
\special{pa 1696 1522}%
\special{pa 3266 1372}%
\special{fp}%
%
\special{pn 13}%
\special{pa 1112 2970}%
\special{pa 2696 2804}%
\special{fp}%
%
\special{pn 8}%
\special{pa 1314 2472}%
\special{pa 2906 2314}%
\special{fp}%
%
\special{pn 8}%
\special{pa 1502 2020}%
\special{pa 3078 1876}%
\special{fp}%
%
\special{pn 8}%
\special{pa 2438 1450}%
\special{pa 1954 2876}%
\special{fp}%
%
\special{pn 8}%
\special{pa 4424 946}%
\special{pa 4424 3128}%
\special{fp}%
%
\special{pn 8}%
\special{pa 4000 2970}%
\special{pa 6722 2970}%
\special{fp}%
%
\special{pn 8}%
\special{pa 4424 2970}%
\special{pa 6152 2538}%
\special{fp}%
%
\special{pn 8}%
\special{pa 4706 1810}%
\special{pa 6434 1378}%
\special{fp}%
%
\special{pn 8}%
\special{pa 5030 2812}%
\special{pa 5302 1646}%
\special{fp}%
%
\special{pn 8}%
\special{pa 5598 2674}%
\special{pa 5872 1508}%
\special{fp}%
%
\special{pn 8}%
\special{pa 4568 2358}%
\special{pa 6296 1926}%
\special{fp}%
\put(15.8000,-26.7400){\makebox(0,0){$T_1$}}%
\put(47.3400,-25.7300){\makebox(0,0){$T_2$}}%
\put(13.8000,-15.0000){\makebox(0,0){$q(m,n)$}}%
\put(31.4000,-27.6000){\makebox(0,0){$p(m_1,n_1)$}}%
\put(43.6000,-18.1000){\makebox(0,0){$p(M',n')$}}%
\put(65.6000,-25.2000){\makebox(0,0){$q(M_1',n_1')$}}%
%
\special{pn 13}%
\special{pa 4424 2962}%
\special{pa 6254 2092}%
\special{fp}%
%
\special{pn 13}%
\special{pa 4706 1810}%
\special{pa 6534 940}%
\special{fp}%
%
\special{pn 13}%
\special{pa 3186 1580}%
\special{pa 3208 1530}%
\special{fp}%
\special{sh 1}%
\special{pa 3208 1530}%
\special{pa 3162 1584}%
\special{pa 3186 1580}%
\special{pa 3200 1600}%
\special{pa 3208 1530}%
\special{fp}%
%
\special{pn 13}%
\special{pa 6448 1242}%
\special{pa 6462 1214}%
\special{fp}%
\special{sh 1}%
\special{pa 6462 1214}%
\special{pa 6414 1264}%
\special{pa 6438 1260}%
\special{pa 6450 1282}%
\special{pa 6462 1214}%
\special{fp}%
%
\special{pn 8}%
\special{pa 1904 2414}%
\special{pa 1386 2934}%
\special{fp}%
\special{pa 1854 2424}%
\special{pa 1336 2942}%
\special{fp}%
\special{pa 1810 2424}%
\special{pa 1286 2948}%
\special{fp}%
\special{pa 1760 2430}%
\special{pa 1234 2956}%
\special{fp}%
\special{pa 1710 2438}%
\special{pa 1192 2956}%
\special{fp}%
\special{pa 1666 2438}%
\special{pa 1142 2962}%
\special{fp}%
\special{pa 1616 2444}%
\special{pa 1134 2928}%
\special{fp}%
\special{pa 1566 2452}%
\special{pa 1162 2856}%
\special{fp}%
\special{pa 1522 2452}%
\special{pa 1192 2784}%
\special{fp}%
\special{pa 1472 2458}%
\special{pa 1220 2712}%
\special{fp}%
\special{pa 1422 2466}%
\special{pa 1250 2640}%
\special{fp}%
\special{pa 1378 2466}%
\special{pa 1278 2568}%
\special{fp}%
\special{pa 1954 2408}%
\special{pa 1430 2934}%
\special{fp}%
\special{pa 1998 2408}%
\special{pa 1480 2928}%
\special{fp}%
\special{pa 2048 2400}%
\special{pa 1530 2918}%
\special{fp}%
\special{pa 2098 2394}%
\special{pa 1574 2918}%
\special{fp}%
\special{pa 2106 2430}%
\special{pa 1624 2912}%
\special{fp}%
\special{pa 2084 2496}%
\special{pa 1674 2904}%
\special{fp}%
\special{pa 2062 2558}%
\special{pa 1718 2904}%
\special{fp}%
\special{pa 2042 2624}%
\special{pa 1768 2898}%
\special{fp}%
\special{pa 2012 2696}%
\special{pa 1818 2890}%
\special{fp}%
\special{pa 1990 2760}%
\special{pa 1868 2884}%
\special{fp}%
\special{pa 1970 2826}%
\special{pa 1912 2884}%
\special{fp}%
%
\special{pn 8}%
\special{pa 5122 2222}%
\special{pa 4454 2890}%
\special{fp}%
\special{pa 5166 2222}%
\special{pa 4446 2942}%
\special{fp}%
\special{pa 5152 2280}%
\special{pa 4518 2912}%
\special{fp}%
\special{pa 5138 2336}%
\special{pa 4604 2870}%
\special{fp}%
\special{pa 5130 2386}%
\special{pa 4684 2832}%
\special{fp}%
\special{pa 5116 2444}%
\special{pa 4770 2790}%
\special{fp}%
\special{pa 5102 2502}%
\special{pa 4850 2754}%
\special{fp}%
\special{pa 5086 2558}%
\special{pa 4936 2712}%
\special{fp}%
\special{pa 5072 2616}%
\special{pa 5014 2674}%
\special{fp}%
\special{pa 5066 2236}%
\special{pa 4468 2832}%
\special{fp}%
\special{pa 5008 2250}%
\special{pa 4482 2774}%
\special{fp}%
\special{pa 4950 2264}%
\special{pa 4490 2726}%
\special{fp}%
\special{pa 4892 2280}%
\special{pa 4504 2668}%
\special{fp}%
\special{pa 4834 2294}%
\special{pa 4518 2610}%
\special{fp}%
\special{pa 4778 2308}%
\special{pa 4532 2552}%
\special{fp}%
\special{pa 4720 2322}%
\special{pa 4546 2496}%
\special{fp}%
\special{pa 4662 2336}%
\special{pa 4562 2438}%
\special{fp}%
\special{pa 4604 2352}%
\special{pa 4568 2386}%
\special{fp}%
%
\special{pn 4}%
\special{pa 5058 2674}%
\special{pa 4878 2856}%
\special{fp}%
\special{pa 4986 2702}%
\special{pa 4820 2870}%
\special{fp}%
\special{pa 4906 2740}%
\special{pa 4762 2884}%
\special{fp}%
\special{pa 4820 2784}%
\special{pa 4706 2898}%
\special{fp}%
\special{pa 4742 2818}%
\special{pa 4648 2912}%
\special{fp}%
\special{pa 4662 2856}%
\special{pa 4590 2928}%
\special{fp}%
\special{pa 4576 2898}%
\special{pa 4532 2942}%
\special{fp}%
\special{pa 5050 2726}%
\special{pa 4936 2840}%
\special{fp}%
\special{pa 5036 2784}%
\special{pa 4994 2826}%
\special{fp}%
%
\special{pn 8}%
\special{pa 1530 1824}%
\special{pa 1602 1766}%
\special{fp}%
\special{pa 1602 1854}%
\special{pa 1588 1782}%
\special{fp}%
%
\special{pn 8}%
\special{pa 4598 2106}%
\special{pa 4648 2048}%
\special{fp}%
%
\special{pn 8}%
\special{pa 4670 2106}%
\special{pa 4648 2054}%
\special{fp}%
%
\special{pn 8}%
\special{sh 1}%
\special{ar 1696 1508 10 10 0  6.28318530717959E+0000}%
\special{sh 1}%
\special{ar 1696 1508 10 10 0  6.28318530717959E+0000}%
\put(8.9600,-30.9100){\makebox(0,0){$0$}}%
\put(42.2200,-31.2000){\makebox(0,0){$0$}}%
\put(25.5000,-12.3000){\makebox(0,0){$T_1^{p,q}$}}%
\put(55.4000,-12.2000){\makebox(0,0){$T_2^{q,p}$}}%
\put(67.6000,-21.1000){\makebox(0,0){$F(p(m_1,n_1))$}}%
%
\special{pn 13}%
\special{pa 2290 2780}%
\special{pa 2350 2830}%
\special{fp}%
\special{pa 2290 2890}%
\special{pa 2350 2830}%
\special{fp}%
%
\special{pn 13}%
\special{pa 5760 2380}%
\special{pa 5820 2280}%
\special{fp}%
\special{pa 5820 2280}%
\special{pa 5710 2290}%
\special{fp}%
%
\special{pn 8}%
\special{pa 1960 800}%
\special{pa 870 3570}%
\special{fp}%
%
\special{pn 8}%
\special{pa 3500 790}%
\special{pa 2410 3560}%
\special{fp}%
\put(19.7000,-6.8000){\makebox(0,0){$L_1$}}%
\put(35.2000,-6.0000){\makebox(0,0){$L_1^{(pq)}$}}%
\put(29.0000,-31.0000){\makebox(0,0){$(p/n,0)$}}%
%
\special{pn 8}%
\special{pa 4960 780}%
\special{pa 4290 3500}%
\special{fp}%
%
\special{pn 8}%
\special{pa 6570 760}%
\special{pa 5900 3480}%
\special{fp}%
\put(49.9000,-6.2000){\makebox(0,0){$L_2$}}%
\put(65.3000,-5.7000){\makebox(0,0){$L_2^{(pq)}$}}%
\put(63.1000,-31.0000){\makebox(0,0){$(q/n',0)$}}%
\put(34.9000,-13.9000){\makebox(0,0){$(a,b)$}}%
\put(67.8000,-9.2000){\makebox(0,0){$F(a,b)$}}%
%
\special{pn 13}%
\special{pa 1680 1500}%
\special{pa 1100 2970}%
\special{fp}%
\special{pa 1110 2970}%
\special{pa 1110 2970}%
\special{fp}%
%
\special{pn 13}%
\special{pa 3280 1370}%
\special{pa 2710 2810}%
\special{fp}%
%
\special{pn 13}%
\special{pa 4700 1810}%
\special{pa 4420 2970}%
\special{fp}%
%
\special{pn 13}%
\special{pa 6520 930}%
\special{pa 6250 2090}%
\special{fp}%
\put(35.4000,-38.8000){\makebox(0,0){where $(a,b)=(qm+pm_1,qn+pn_1)$}}%
\end{picture}%

\vspace{ 20mm}
\noindent  Recall from Remark \ref{rem:6-2} that $l_1,l_1^*$ and $l_2,l_2^*$ are
simple closed curves in $T_1$ and $T_2$. The set $T_1^{p,q}$ is the $pq$-sheeted torus over $T_1$ winding $p $ times along $l_1^*$
and $q$ times along $l_1$, while
$T_2^{q,p}$ is the $pq$-sheeted torus over $T_2$ winding $q $ times along $l_2^*$
and $p$ times along $l_2$ (see the figure for the case $p=2,q=3$). This gives a visual  interpretation of the Lie subgroup $H$ of ${\bf  T}$ and its closure
$\bar { H}$.

The topic of Levi flat surfaces in two dimensional complex tori has
recently been studied by T. Ohsawa \cite{o-2}. See also O. Suzuki \cite{suzuki}.

\section{Complex homogeneous spaces.}
\setcounter{section}{6}
In this section, we let $M$ be an $n$-dimensional complex space with
the property that there exists a connected complex Lie group
$G\subset$Aut\,$M$ of complex dimension $m \ge n$ which
acts transitively on $M$. As prototypical examples, we can take
$M={\mathbb{P}}^n$ (complex projective space) and $G=GL(n+1,{\mathbb{C}})$, or more generally, we can take
$M=G(k,n)$
(complex Grassmannian manifold) or $M={\cal F}_n$ (complex flag space),
and $G=GL(n,{\mathbb{C}})$.
As
noted in the previous section, another  example is provided by $M={\mathbb{C}}$ in
which case $G=$\,Aut\,$M$ is a noncommutative complex Lie group of complex
dimension $2$.

We summarize some well-known results concerning homogeneous spaces which
will be used in this section. We write $e$ for the identity element
of $G$.

{\bf  1.} \ For a fixed $z\in M$, let
$$H_z:=\{g\in G:g(z)=z\}$$
be the isotropy subgroup of $G$ for $z$.
Then $H_z$ is a complex $(m-n)$-dimensional analytic set in $G$ without
singular points. In particular, $H_z$ is closed and hence is a Lie subgroup of $G$, but $H_z$ is not always connected. We let $G/H_z$ denote the set of all left cosets $gH_z$.
This quotient space $G/H_z$ has the structure of a complex
$n$-dimensional manifold; moreover, if we let $\pi_z:G \mapsto G/H_z$ be the
coset mapping
$\pi_z(g)=gH_z$ and we let $\psi_z: G \mapsto M$ be the
 mapping $\psi_z(g)=g(z)$, then there exists an isomorphism
 $\alpha _z: G/H_z \to  M$ such that
$ \alpha _z \circ \pi_z =\psi_z$ on $G$.
 This is equivalent to $\alpha_z(gH_z)=g(z)$ for $g \in G$. These three
 mappings are holomorphic; hence  $\psi_z$ gives a holomorphic fibering
of $G$
 over $M$ with fibers  $\psi _z^{-1}(\zeta) \approx
H_z$ as complex manifolds; and, given  $g \in G$,  $\psi _z^{-1}(g(z))=gH_z$ as sets
in $G$.

{\bf  2.} Fix $z \in M$ as in {\bf  1}. Let $g_0\in G$ and set $ z_0=g_0(z)\in M$. Then there exists a local
holomorphic section $\sigma$ of $G$ with respect to $\psi_z$ over  a
neighborhood $V
\subset M$ of $z_0$ passing through $g_0$; i.e.,
\begin{align}
\nonumber \sigma(z_0)= g_0 \quad \mbox{ and } \quad \psi_z \circ \sigma
 =  \mbox{identity on $V$}
\end{align}
so that $\sigma(\zeta)(z)=\zeta$  for $\zeta\in V$.
 If we set $\tau: g\in \psi_{z}^{-1}(V)\mapsto \sigma(g(z))^{-1}g \in H_z$,
 then $\tau(g)$ depends holomorphically on $g$ and, for a fixed $z'\in
 V$, the restriction of $\tau$ to $\psi_z^{-1}(z')$ is a bijection from
 $\psi^{-1}_z(z')$ onto $H_z$. It follows
that
\begin{align}
 \label{eqn:sigma2}
T_z: g \in \psi_z^{-1}(V) \to (g(z), \tau(g)) \in V \times H_z
\end{align}
provides a local trivialization of $G$ over $M$ with $T_z:  \sigma(\zeta) \to (\zeta, e), \ \zeta \in V$.
Similarly we have the following result. Fix
$z,w\in M$ and $g\in G$ with $g(z)=w$. There exists a
neighborhood  $V $ of $z$ in $M$ and for each $\zeta\in V$, there exists an element $g_\zeta\in G$
such that $g_z=g$; $g_\zeta(z)=w$ for $\zeta\in V$; and $g_\zeta$
depends holomorphically on $\zeta\in V$.

{\bf  3.} \ Fix  $z_0, \ z_1\in M$. Let $h(z_1)=z_0$ where $h \in G$,
so that
$ h^{-1}H_{z_0}h=H_{z_1}$.  We then have two fiberings of $G$ over $M$
with respect to $\psi_{z_0}$ and $\psi_{z_1}$.  We consider the right
translation $R_h: G \to G$ where $R_h(g)= gh, \ g \in G$. Then $R_h$
induces a holomorphic isomorphism ${\cal R}_h: G/H_{z_0} \to
G/H_{z_1}$ via ${\cal R}_h(gH_{z_0})=ghH_{z_1}, \ g \in G$.
In other words, for $z\in M$, take any $g,k\in  G$ with
$z=g(z_0)=k(z_1)$. Then ${\cal R}_h(gH_{z_0})=kH_{z_1}$.

\smallskip
Let $D \Subset   M$ be a domain with $C^\infty $ boundary in $M$.
For $z\in M$, we let
$$D(z):=\{g\in G:g(z)\in D\}=\psi_z^{-1}(D)=\pi_{z}^{-1}\circ
\alpha_z^{-1}(D)$$
be a possibly unbounded and possibly disconnected domain in $G$. Thus

\vspace{1cm}


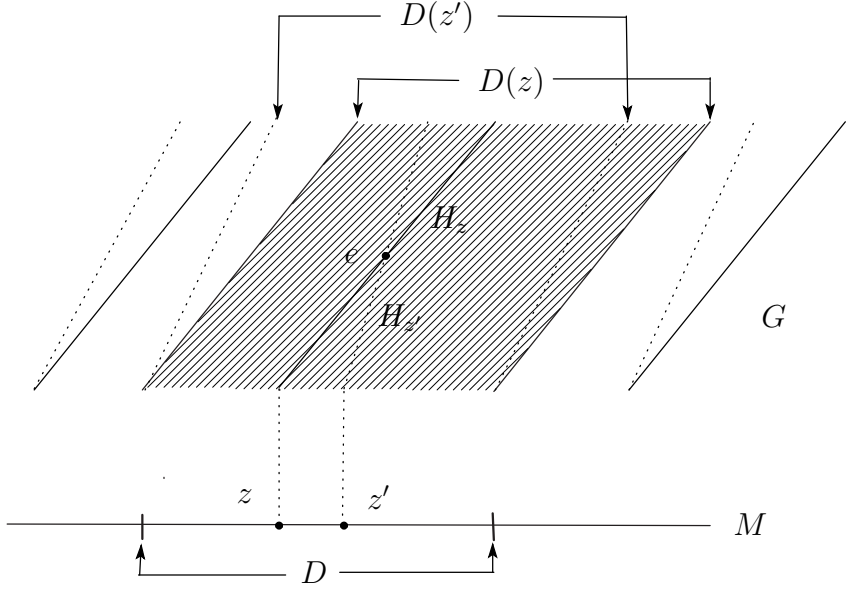
\unitlength 0.1in
\begin{picture}( 43.4000, 29.9700)( 24.0000,-43.0600)
%
\special{pn 8}%
\special{pa 2400 4040}%
\special{pa 6040 4048}%
\special{fp}%
%
\special{pn 8}%
\special{pa 3212 3802}%
\special{pa 3212 3796}%
\special{fp}%
%
\special{pn 13}%
\special{pa 4914 3984}%
\special{pa 4920 4124}%
\special{fp}%
%
\special{pn 13}%
\special{pa 3100 3992}%
\special{pa 3100 4118}%
\special{fp}%
%
\special{pn 8}%
\special{pa 3808 3340}%
\special{pa 3808 3340}%
\special{fp}%
%
\special{pn 8}%
\special{pa 4920 3348}%
\special{pa 6040 1940}%
\special{fp}%
%
\special{pn 8}%
\special{pa 3100 3340}%
\special{pa 4214 1940}%
\special{fp}%
%
\special{pn 20}%
\special{sh 1}%
\special{ar 4360 2640 10 10 0  6.28318530717959E+0000}%
\special{sh 1}%
\special{ar 4360 2640 10 10 0  6.28318530717959E+0000}%
\put(46.8900,-24.5100){\makebox(0,0){$H_z$}}%
%
\special{pn 8}%
\special{pa 5620 3340}%
\special{pa 6740 1940}%
\special{fp}%
%
\special{pn 8}%
\special{pa 2540 3340}%
\special{pa 3660 1940}%
\special{fp}%
\put(44.3700,-29.6200){\makebox(0,0){$H_{z'}$}}%
\put(39.8900,-42.9200){\makebox(0,0){$D$}}%
%
\special{pn 8}%
\special{pa 3808 3334}%
\special{pa 3808 4040}%
\special{dt 0.045}%
%
\special{pn 20}%
\special{sh 1}%
\special{ar 3808 4048 10 10 0  6.28318530717959E+0000}%
\special{sh 1}%
\special{ar 3808 4048 10 10 0  6.28318530717959E+0000}%
%
\special{pn 20}%
\special{sh 1}%
\special{ar 4144 4048 10 10 0  6.28318530717959E+0000}%
\special{sh 1}%
\special{ar 4144 4048 10 10 0  6.28318530717959E+0000}%
%
\put(39.7500,-42.2200){\makebox(0,0)[lb]{}}%
\put(36.3200,-38.9300){\makebox(0,0){$z$}}%
\put(43.2500,-39.0000){\makebox(0,0){$z'$}}%
%
\special{pn 8}%
\special{pa 4144 3348}%
\special{pa 4136 4040}%
\special{dt 0.045}%
%
\put(48.5000,-15.6200){\makebox(0,0)[lb]{}}%
\put(50.1100,-17.3700){\makebox(0,0){$D(z)$}}%
\put(63.6900,-29.5500){\makebox(0,0){$G$}}%
\put(62.5000,-40.4000){\makebox(0,0){$M$}}%
\put(41.8500,-26.4000){\makebox(0,0){$e$}}%
%
\special{pn 8}%
\special{pa 4220 1724}%
\special{pa 4214 1912}%
\special{fp}%
\special{sh 1}%
\special{pa 4214 1912}%
\special{pa 4236 1846}%
\special{pa 4216 1860}%
\special{pa 4196 1846}%
\special{pa 4214 1912}%
\special{fp}%
%
\special{pn 8}%
\special{pa 6040 1730}%
\special{pa 6040 1920}%
\special{fp}%
\special{sh 1}%
\special{pa 6040 1920}%
\special{pa 6060 1852}%
\special{pa 6040 1866}%
\special{pa 6020 1852}%
\special{pa 6040 1920}%
\special{fp}%
%
\special{pn 8}%
\special{pa 4214 1716}%
\special{pa 4738 1724}%
\special{fp}%
\special{pa 5228 1716}%
\special{pa 6040 1716}%
\special{fp}%
%
\special{pn 8}%
\special{pa 6040 1724}%
\special{pa 6040 1710}%
\special{fp}%
%
\special{pn 8}%
\special{pa 3086 4306}%
\special{pa 3094 4152}%
\special{fp}%
\special{sh 1}%
\special{pa 3094 4152}%
\special{pa 3070 4218}%
\special{pa 3092 4206}%
\special{pa 3110 4220}%
\special{pa 3094 4152}%
\special{fp}%
%
\special{pn 8}%
\special{pa 3080 4306}%
\special{pa 3842 4306}%
\special{fp}%
\special{pa 4122 4292}%
\special{pa 4920 4292}%
\special{fp}%
%
\special{pn 8}%
\special{pa 4914 4292}%
\special{pa 4914 4146}%
\special{fp}%
\special{sh 1}%
\special{pa 4914 4146}%
\special{pa 4894 4212}%
\special{pa 4914 4198}%
\special{pa 4934 4212}%
\special{pa 4914 4146}%
\special{fp}%
%
\special{pn 8}%
\special{pa 3808 1388}%
\special{pa 3800 1920}%
\special{fp}%
\special{sh 1}%
\special{pa 3800 1920}%
\special{pa 3822 1854}%
\special{pa 3802 1866}%
\special{pa 3782 1852}%
\special{pa 3800 1920}%
\special{fp}%
%
\special{pn 8}%
\special{pa 5606 1388}%
\special{pa 5614 1926}%
\special{fp}%
\special{sh 1}%
\special{pa 5614 1926}%
\special{pa 5632 1860}%
\special{pa 5612 1874}%
\special{pa 5592 1860}%
\special{pa 5614 1926}%
\special{fp}%
\put(46.2600,-13.9400){\makebox(0,0){$D(z')$}}%
%
\special{pn 8}%
\special{pa 3808 1388}%
\special{pa 4332 1388}%
\special{fp}%
\special{pa 4942 1388}%
\special{pa 5606 1394}%
\special{fp}%
%
\special{pn 8}%
\special{pa 3100 3348}%
\special{pa 4914 3340}%
\special{ip}%
\special{pa 4200 1954}%
\special{pa 6040 1948}%
\special{ip}%
%
\special{pn 4}%
\special{pa 5046 1954}%
\special{pa 4368 2634}%
\special{fp}%
\special{pa 5088 1954}%
\special{pa 3710 3334}%
\special{fp}%
\special{pa 4354 2648}%
\special{pa 3668 3334}%
\special{fp}%
\special{pa 5004 1954}%
\special{pa 3626 3334}%
\special{fp}%
\special{pa 4962 1954}%
\special{pa 3584 3334}%
\special{fp}%
\special{pa 4920 1954}%
\special{pa 3548 3326}%
\special{fp}%
\special{pa 4878 1954}%
\special{pa 3506 3326}%
\special{fp}%
\special{pa 4836 1954}%
\special{pa 3464 3326}%
\special{fp}%
\special{pa 4794 1954}%
\special{pa 3422 3326}%
\special{fp}%
\special{pa 4752 1954}%
\special{pa 3380 3326}%
\special{fp}%
\special{pa 4710 1954}%
\special{pa 3338 3326}%
\special{fp}%
\special{pa 4668 1954}%
\special{pa 3296 3326}%
\special{fp}%
\special{pa 4626 1954}%
\special{pa 3254 3326}%
\special{fp}%
\special{pa 4584 1954}%
\special{pa 3212 3326}%
\special{fp}%
\special{pa 4542 1954}%
\special{pa 3170 3326}%
\special{fp}%
\special{pa 4500 1954}%
\special{pa 3128 3326}%
\special{fp}%
\special{pa 4458 1954}%
\special{pa 3220 3194}%
\special{fp}%
\special{pa 4416 1954}%
\special{pa 3388 2984}%
\special{fp}%
\special{pa 4374 1954}%
\special{pa 3548 2780}%
\special{fp}%
\special{pa 4332 1954}%
\special{pa 3710 2578}%
\special{fp}%
\special{pa 4290 1954}%
\special{pa 3878 2368}%
\special{fp}%
\special{pa 4248 1954}%
\special{pa 4038 2164}%
\special{fp}%
\special{pa 5130 1954}%
\special{pa 3752 3334}%
\special{fp}%
\special{pa 5172 1954}%
\special{pa 3794 3334}%
\special{fp}%
\special{pa 5214 1954}%
\special{pa 3836 3334}%
\special{fp}%
\special{pa 5256 1954}%
\special{pa 3878 3334}%
\special{fp}%
\special{pa 5298 1954}%
\special{pa 3920 3334}%
\special{fp}%
\special{pa 5340 1954}%
\special{pa 3962 3334}%
\special{fp}%
\special{pa 5382 1954}%
\special{pa 4004 3334}%
\special{fp}%
\special{pa 5424 1954}%
\special{pa 4046 3334}%
\special{fp}%
%
\special{pn 4}%
\special{pa 5466 1954}%
\special{pa 4088 3334}%
\special{fp}%
\special{pa 5508 1954}%
\special{pa 4130 3334}%
\special{fp}%
\special{pa 5550 1954}%
\special{pa 4172 3334}%
\special{fp}%
\special{pa 5592 1954}%
\special{pa 4214 3334}%
\special{fp}%
\special{pa 5634 1954}%
\special{pa 4256 3334}%
\special{fp}%
\special{pa 5676 1954}%
\special{pa 4298 3334}%
\special{fp}%
\special{pa 5718 1954}%
\special{pa 4340 3334}%
\special{fp}%
\special{pa 5760 1954}%
\special{pa 4382 3334}%
\special{fp}%
\special{pa 5802 1954}%
\special{pa 4424 3334}%
\special{fp}%
\special{pa 5844 1954}%
\special{pa 4466 3334}%
\special{fp}%
\special{pa 5886 1954}%
\special{pa 4500 3340}%
\special{fp}%
\special{pa 5928 1954}%
\special{pa 4542 3340}%
\special{fp}%
\special{pa 5970 1954}%
\special{pa 4584 3340}%
\special{fp}%
\special{pa 6012 1954}%
\special{pa 4626 3340}%
\special{fp}%
\special{pa 5914 2094}%
\special{pa 4668 3340}%
\special{fp}%
\special{pa 5754 2298}%
\special{pa 4710 3340}%
\special{fp}%
\special{pa 5592 2500}%
\special{pa 4752 3340}%
\special{fp}%
\special{pa 5424 2710}%
\special{pa 4794 3340}%
\special{fp}%
\special{pa 5264 2914}%
\special{pa 4836 3340}%
\special{fp}%
\special{pa 5096 3124}%
\special{pa 4878 3340}%
\special{fp}%
%
\special{pn 8}%
\special{pa 3808 3326}%
\special{pa 4928 1940}%
\special{fp}%
%
\special{pn 8}%
\special{pa 3114 3340}%
\special{pa 3794 1920}%
\special{dt 0.045}%
%
\special{pn 8}%
\special{pa 4144 3334}%
\special{pa 4578 1940}%
\special{dt 0.045}%
%
\special{pn 8}%
\special{pa 4928 3340}%
\special{pa 5614 1934}%
\special{dt 0.045}%
%
\special{pn 8}%
\special{pa 2540 3334}%
\special{pa 3296 1940}%
\special{dt 0.045}%
%
\special{pn 8}%
\special{pa 5620 3340}%
\special{pa 6264 1948}%
\special{dt 0.045}%
\end{picture}%

\vspace{1cm}

\noindent  $D(z)$ is a set of cosets modulo $H_z$ in $G$,\, i.e.,
$D(z)$ may be considered as a subset of $G/H_z$ in such a way that
$\alpha_z$ maps $D(z)/H_z$ biholomorphically to $D$. We remark
that the connectedness or disconnectedness of $D(z)$ in $G$ is
independent of the point $z\in M$. Indeed, given $z'\in M$, take
$h\in G$ such that $h(z)=z'$. Since the equality $D(h(z))=
D(z)h^{-1}$ holds (this will be proved  as II. in Proposition
\ref{prop:useful}), we have $D(z')=D(z)h^{-1}$, so that $D(z)$ and
$D(z')$ are simultaneously connected or disconnected in $G$. Note
that by homogeneity we have, for any $z\in M$,
\begin{align} \label{eqn:triviality}
 D=\{g(z)\in M : g\in D(z)\}.
\end{align}
We also mention that in the case when the isotropy subgroup $H_z$ of $G$ is
 connected, $D(z)$ is connected in $G$.
For, we may assume $z\in D$, so that $e\in D(z)$ and $H_z \subset
 D(z)$. Let $g\in D(z)$, i.e., $g(z)\in D$. We take a continuous
curve $\gamma: t\in [0,1] \to \gamma(t)$ in  $D$  such that
$\gamma(0)=g(z)$ and $\gamma(1)=z$. We use property {\bf  2.} to piece together a continuous section $\sigma: t\in [0,1] \to \sigma(t)$ of $G $ over $\gamma $ via $\psi_z$, i.e., $\sigma(t)(z)= \gamma (t)$ for
$t\in [0,1]$, with $\sigma(0)=g$. Since $\gamma \subset D$ and
 $\sigma(0)\in D(z)$, we have $\sigma \subset D(z)$. Also, since $\sigma(1)(z)=
z$, we have $\sigma(1)\in H_z$. From the assumption that $H_z$ is connected, we can find a continuous
curve $\eta: t\in [0,1] \to \eta(t)$ in $H_z$ such that
$\eta(0)= \sigma(1)$ and $\eta(1)= e$. Putting together the curves $\sigma$ and $\eta$; i.e., defining the concatenation
$$\Gamma(t):= \sigma (2t) \ \hbox{for} \  t\in [0,1/2]; \ \ \hbox{and} \
 \ \Gamma(t)= \eta (2t-1) \ \hbox{for} \   t\in [1/2,1],$$
gives a continuous curve $\Gamma: t\in [0,1] \to \Gamma(t)$ in $G$ starting at $g$ and terminating at $e$ which lies entirely in $D(z)$. Hence $D(z)$ is a connected subset of $G$.

 We next consider the case when $H_z$ is not
 connected in $G$. We write $H_z=\cup_{k=1}^\infty H^{(k)}_z$ as the
 union of its connected components in $G$. We let
 $H_z'=H_z^{(1)}$  denote the
 connected component of $H_z$  which contains the identity $e$.
Note that $\{H_z^{(k)}\}_k$  are isolated sets in $G$ (see 1. in
 Proposition  \ref{prop:disc}). If
 we
 take  $h_k\in H_z^{(k)}, \, k=1,2,\ldots$, then
 \begin{align}
 \label{eqn:visible}
H_z^{(k)}=h_kH_{z}'=H_{z}'h_k , \quad k=1,2, \ldots
\end{align}
To see this, since $h_kH'_{z}$ is connected in $G$ and contains $h_k$, it
  follows from the inclusion $h_kH'_{z} \subset H_{z}$ that $h_kH'_{z} \subset
  H^{(k)}_{z}$. Conversely, since $h_k^{-1}H_{z}^{(k)} $ is
  connected in $G$, it follows from $e\in h_k^{-1}H_{z}^{(k)}
  \subset H_{z}$ that $h_k^{-1}H_{z}^{(k)}\subset H_{z}'$. Thus the
  first equality holds. The second one can be verified in a similar fashion; it also follows from the first by the normality
  of $H_z'$ in $H_z$.

  In the case when $H_z$ is not
 connected the connectivity of $D(z)$
 depends on the domain $D $ in $M$. In fact, if $D$
 is very small and $z\in D$, then each connected component of $D(z)$  contains
 only one component of $H_z$. For a general  domain $D \Subset M$ and
 $z\in D$, we decompose
$D(z)$ into a finite or a countably infinite collection of connected components:
\begin{align}
 \label{eqn:connectedcom}
D(z)= \mbox{  $\bigcup_{i=1}^\infty$}\,  D^{(i)}(z).
\end{align}
 We write
$D'(z):=D^{(1)}(z)$  for the connected
component which contains the identity $e$,
so that $H'_{z}\subset D'(z)$.
We set ${\cal H}'(z):=H_z \cap D'(z)$. Then ${\cal H}'(z)$ consists of
a finite or a countably infinite collection of components $\{H_z^{(k)}\}_{k=1}^\infty $, say
\begin{align}
 \label{eqn:H-prime}
{\cal H}'(z)= \{H_z', H_z^{(\alpha_2)}, H_z^{(\alpha_3)}, \ldots
\}= \mbox{ $ \bigcup_{j=1}^\infty $} \, h_j'H_z',
\end{align}
where $h_1'=e$ and $h_j'\in H_z ^{(\alpha_j)} , j=2,3,\ldots $ (using (\ref{eqn:visible})). We have that $g\in D'(z)$ implies that $g {\cal H}'(z) \subset D'(z)$. To see this, let $h\in {\cal H}'(z)$. We connect $e$ and $g$ by a continuous curve $\gamma: t\in
[0,1]\to \gamma(t)$ in
$D'(z)$. Then the continuous curve
$ \widetilde { \gamma} : t\in [0,1] \to \gamma(t)h $
in $G$ connects $h$ and $gh$. Since $\psi_z( \widetilde { \gamma} )=
\psi_z (\gamma) \subset D$ and $\widetilde { \gamma}(0)=h\in {\cal H}'(z) \subset
D'(z)$, we have $\widetilde { \gamma} \subset D'(z)$, and hence  $gh= \widetilde { \gamma}(1)\in
D'(z)$.

Next we show that for each $i\ge 2$ we have $D^{(i)}(z) \cap H_z \ne \emptyset$. We begin with $g_0\in D^{(i)}(z)$ and set $z_0:=g_0(z)\in D$. Take a
continuous curve $\gamma: t\in [0,1]\to \gamma(t)$ in $D$ such
that $\gamma (0) =z_0$ and $\gamma(1)=z$. Using
property {\bf  2.}, we can construct a continuous section $\sigma: t\in [0,1]\to
\sigma(t)$  of $G$ over $\gamma$ with respect to $\psi_z$ such that
$\sigma(0)=g_0$. Hence  $\sigma \subset D(z)$. Then
$\sigma(0)=g_0\in D^{(i)}(z)$ implies that $\sigma \subset D^{(i)}(z)$.
On the other hand,  since $\sigma(1)(z)=\gamma (1)=z$, we have
$\sigma(1)\in H_z$, so that $\sigma (1)\in H_z\cap D^{(i)}(z)$.

Thus for $i=2,\ldots$ we define the non-empty sets ${\cal H}^{(i)}(z):= D^{(i)}(z)\cap H_z$ and we can write
\begin{align}
 \label{qn:yasisugiru}
H_z&=\mbox{  $\bigcup_{i=1}^\infty$} \,  {\cal H}^{(i)}(z)\ \
\end{align}
where ${\cal H} ^{(1)}={\cal H}'(z) $ and the union is disjoint. If we fix ${\bf  h}^{(i)}\in
{\cal H}^{(i)}(z),\,i=2,3, \ldots $, then we have
\begin{align}
 \label{eqn:atode}
D^{(i)}(z)=  D'(z){\bf  h}^{(i)} ,   \quad   \quad
 {\cal H}^{(i)}(z)  = {\cal H}'(z)  {\bf  h}^{(i)}.
\end{align}
To verify these equalities, first observe that since $D'(z){\bf  h}^{(i)}$ is a connected subset of $D(z)$ which contains
${\bf  h}^{(i)}$,  we have $D'(z){\bf  h}^{(i)} \subset  D^{(i)}(z)$. A similar
argument implies that $D^{(i)}(z)({\bf  h}^{(i)})^{-1} \subset D'(z)$, so that
the first formula holds. An interpretation of this formula is that each  $D^{(i)}(z), \,i=2,3, \ldots$ is isomorphic to $D'(z)$ via right translation ${\cal R}_{{\bf  h}_i}:
 g\in D'(z)
 \to g {\bf  h}^{(i)} \in D^{(i)}(z)$. The second formula clearly follows from the first and the definitions of ${\cal H}^{(i)}(z)$ and ${\cal H}'(z)$.
 \begin{proposition}
     \label{eqn:elm} Let $z\in D$. Then
\begin{enumerate}
 \item [1.]  $D=\{g(z)\in M: g\in
       D'(z)\}$;
 \item [2.] $D'(z)h \subset D'(z)$ \ for $ h\in {\cal H}'(z)$;
  \item [3.]
let $g\in D'(z)$ and $h\in H_z$. If $gh\in D'(z)$, then $h\in
      {\cal H}'(z)$.
\end{enumerate}
\end{proposition}
\noindent {\bf  Proof}.
To prove 1., we have from (\ref{eqn:triviality}) that
$$\{g(z)\in M: g\in
       D'(z)\}\subset D.$$
       For the reverse inclusion, fix $\hat z \in D$ and consider a continuous curve $\gamma:  t\in [0,1]\to \gamma(t) $
in $D$
such that $\gamma(0)=z$ and $\gamma(1)=\hat z$. By property {\bf  2.}, we can find a continuous
section $\sigma: t\in [0,1] \to \sigma(t)$ of $G$ over $\gamma$ via $\psi_{z}$  such that
  $\sigma(0)=e$. Consequently, since $\sigma(0)=e\in D'(z)$, we have
$\sigma \subset D'(z)$. In particular,  $\sigma(1)\in D'(z)$ and
 $\sigma(1)(z)=\gamma(1)=\hat z$, so that taking $g:=\sigma(1)$
verifies the reverse inclusion. To prove 2.,
let $g\in D'(z)$ and $h\in {\cal H}'(z)$. We take
a continuous curve $\eta: t\in
[0,1]\to \eta(t)$ in $D'(z)$ such that $\eta(0)=e$ and
$\eta(1)=g$. Then the continuous curve $\xi:
t\in
[0,1]\to \xi(t):=\eta(t)h$ in $G$ satisfies
$\xi(t)(z)=\eta(t)h(z)=\eta(t)(z)\in D, \ t\in [0,1]$.
 It
follows from $\xi(0)=h\in D'(z)$ that $\xi(t)\in D'(z), t\in [0,1]$,
and hence  $gh=\xi(1)\in D'(z)$, which proves 2.
To prove (3), we set $g'=gh \in D'(z)$. Since $g\in D'(z)$ we can find a
continuous curve $\tau:t\in [0,1] \to g(t)$ in $D'(z)$  such that
$g(0)=e$ and $g(1)=g$. Consider the continuous curve $\tau': t\in
[0,1]\to \tau'(t):=g(t)h$ in $G$.  We have $\tau'(t)(z)=
g(t)h(z)=g(t)(z)\in D$. Since $\tau'(1)= gh \in D'(z)$, it follows
that $\tau' \subset  D'(z)$. In particular, $h=eh= \tau'(0)\in
D'(z)$, which proves (3).
\hfill $\Box$

\smallskip

\noindent We remark that we can rephrase 1. of Proposition \ref{eqn:elm} as $[D'(z_0)](z_0)=D$. In a similar fashion, it can be shown that
\begin{equation} \label{awkward}
[\partial D'(z_0)](z_0)=\partial D \end{equation}
where $\partial D'(z_0)$ is the boundary of $D'(z_0)$ in $G$ and $\partial D$ is the boundary of $D$ in $M$.
\smallskip

\begin{corollary}\label{eqn:happy} \
Let $z\in D $. Then ${\cal H}'(z)$ is a closed Lie subgroup of $H_z$.
\end{corollary}
\noindent {\bf Proof}. We first prove that ${\cal H}'(z)$ is a subgroup of
$H$. Since $e\in {\cal H}'(z)$, it suffices to prove that
$g:=h_1h_2^{-1}\in {\cal H}'(z)$ for $h_1,h_2 \in {\cal H}'(z)$. Since $g(z)=z$, we have $g\in H_z$. Then using
$h_1=gh_2$ and $h_1\in D'(z)$, it follows from 3. in Proposition
\ref{eqn:elm} that $g\in {\cal H}'(z)$. We next verify
that ${\cal H}'(z)$ is closed in $G$; since a closed subgroup of a Lie group is a Lie subgroup, this will complete the proof. Let $h_n\in {\cal H}'(z), \,n=1,2,\ldots $ and fix
$g\in G$ with $h_n \to g\ (n\to \infty)$. Since $H_z$ is closed,
we have $g\in H_z$. Since $e\in D'(z)$, $D'(z)$ is open in $G$, and $gh_n ^{-1}\to e$, it follows that $gh_n^{-1}\in D'(z)$ for
sufficiently large $n$. Thus $g\in D'(z)h_n$. Using 2. in Proposition
\ref{eqn:elm} we conclude that $g\in D'(z)$, and hence $g\in {\cal H}'(z)$.
     \hfill $\Box$

 \begin{proposition}\label{prop:useful} \  We have
  \begin{enumerate}
   \item [{I.}]\
$e \in D(z)\  ( \partial D(z), \ \overline{  D(z)}^c )$ if
     and only if $ z\in D\  ( \partial D, \ \overline{
     D}^c).
$
   \item [{II.}] \  $D(h(z))=
   D(z)h^{-1}$ for $z\in M$ and $h\in G$.
 \item [{III.}] \  $D'(h(z))=
   D'(z)h^{-1}$ for $z\in D$ and $h\in D'(z)$.
     \end{enumerate}
 \end{proposition}
\noindent {\bf Proof}. By definition, we have  $e\in D(z)$ if and only if $z \in
D$. Let $e \in \partial D(z)$. Then there exists a sequence $ \{g_n\}_n
\subset D(z)$ such
that $g_n \to e$. Hence, $g_n(z) \to z$ in $M$. Since $g_n(z)\in D$ and
$z \not\in D$, we have $z \in \partial D$.
Conversely, let $z \in
\partial D$. Then $e \not \in D(z)$. We fix a holomorphic section
$\sigma$ of
$G$ via  $\psi _z$ over a neighborhood $V\subset
M$  of $z$ with $\sigma(z)=e$; thus $\sigma(\zeta)(z)=\zeta$ for all
$\zeta\in V$. We take $\{\zeta_n\}_n \subset V \cap D$ such
that $\zeta_n \to z$ as $n \to \infty$. It follows that
$\sigma(\zeta_n)\in D(z)$ and
$\sigma(\zeta_n) \to e$ as $n \to \infty$, so that $ e \in \partial D(z)$.
We thus have $e\in \partial D(z) \Leftrightarrow z \in \partial
D$, which proves I.

Property II. is based on (\ref{eqn:triviality}). Indeed, $g \in D(h(z))$ if and only if $gh(z)\in D$  if and only if $gh \in D(z)$ if and only if $g \in D(z)h^{-1}$.

To prove III., we
fix $z\in D$ and $h\in D'(z)$.  Then II. implies $D(h(z))=D(z)h^{-1}$.
From (\ref{eqn:connectedcom}) we thus  obtain
 $D(h(z))= \cup_{ i=1  }^{\infty   }[ D^{(i)}(z)h^{-1}]$ where this is a disjoint union. Hence
 $D'(h(z))$ coincides with some $D^{(i)}(z)h^{-1}$ which contains
 the identity element $e$. On the other hand, since $h\in D'(z)$, we have $e\in
 D'(z)h^{-1}$; hence $i=1$ and $D'(h(z))= D'(z)h^{-1}$. \hfill
 $\Box$

\smallskip
Let $\partial D$ denote the
boundary of $D$ in $M$ and let $\partial D(z)$ denote the relative
boundary of
$D(z)$ in $G$. Note that $\partial D(z)=\psi_z^{-1}(\partial D)$ is
smooth in $G$, but it is not necessarily compact in $G$.
A similar argument to that used in proving I. shows that for
$g\in G$ and $z \in M$, $g \in \partial D(z)$ if and only if
$g(z) \in \partial D$.  Moreover, we note that, for $z\in D$ and $h\in D'(z)$,
\begin{equation}
\begin{split}
  \label{eqn:baoundary-prime}
\partial D'(h(z))&=(\partial D'(z))h  ^{-1} \quad \mbox{ and }  \quad
 \overline { D'(h(z))}^c =[ \,\overline { D'(z)}^c\,]\,h  ^{-1}.
\end{split}
\end{equation}

It suffices to prove the first formula. Let $g\in \partial D'(h(z))$. Then $g\not\in D'(h(z))$ and
there exists $g_n\in D'(h(z))$  such that $\lim_{ n\to \infty}g_n= g$
in $G$. From III. there exists $k_n\in D'(z)$  such that $g_n=k_n h  ^{-1}  $. Thus, $k:=\lim_{ n\to \infty}k_n= gh$. Since  $k=gh \not \in D'(h(z)) h
=D'(z)$ by III., it follows that $k \in \partial D'(z)$, and hence $g=k
h  ^{-1}  \in
(\partial D'(z))h  ^{-1}  $, so that $\partial D'(h(z))\subset (\partial D'(z))h  ^{-1}$. To prove the reverse inclusion, let $g\in \partial D'(z)$. Then $gh^{-1}\not \in D'(z) h^{-1}=D'(h(z))$. We take $\{g_n\}\subset   D'(z), \ n=1,2,...$ with $g_n\to g$ in $G$ as $n\to \infty$. Since $g_nh^{-1}\in D'(z)h^{-1}=D'(h(z))$, we have $gh^{-1}\in \partial D'(h(z))$.

\smallskip
From II. we see that for $z,
z'\in M$,  $D(z)$ and $D(z')$ are biholomorphically
equivalent. In case $z,\, z'\in D$, from III. and 1. in Proposition \ref{eqn:elm}, $D'(z)$ and
$D'(z')$ are also  biholomorphically equivalent.
\begin{lemme}\label{re:ebound} \   Let
$\zeta\in \partial D$ and let $\{z_n\}_n \subset D$ converge to $\zeta$ as $n \to \infty$. Then there exist $g_n \in \partial D'(z_n), \ n=1,2,
\ldots $, such that
$g_n \to e$ in $G$ as $n \to \infty$.
\end{lemme}
\noindent {\bf Proof}. We take a holomorphic section
$\sigma: V \to  G$ of the fiber space
$\psi_\zeta: \,  G \to M $  over a neighborhood $V\subset M$ of $\zeta$
with $\sigma(\zeta)=e$. For sufficiently large $n$ we have
$z_n\in V$ and hence $\sigma(z_n )(\zeta)=z_n$. It follows that $g_n: =
\sigma(z_n)^{-1}$ is close to $e$ in $G$; moreover, we have $g_n \in
\partial D(z_n)$, for $g_n(z_n)=\sigma(z_n)^{-1}\sigma(z_n)(\zeta)=\zeta
\in \partial D$. Thus, $g_n\in \partial D'(z_n)$.
 \hfill $\Box$

\begin{proposition}\label{prop:dz} \
 Assume that  $D$ is pseudoconvex in $M$, and fix $z \in M$. Then we
 have
\begin{enumerate}
 \item  [1.]  $D(z)$ is pseudoconvex in $G$; and
\item [2.] there exists a sequence of piecewise smooth pseudoconvex
 domains $\{D_n(z)\}_n$ in $G$ such that $D_n(z)
\subset D_{n+1}(z) \Subset D(z)$ and $D(z)=\bigcup_{n=1}^\infty D_n(z)$.
\end{enumerate}
\end{proposition}
\noindent {\bf Proof}.
Recall that $\partial D(z)=\psi_z^{-1}(\partial D)$ is smooth in
$G$. Let $g_0\in \partial D(z)$ in $G$ and set $z_0=g_0(z)$ in $M$.
Then
$z _0\in \partial D$. By (\ref{eqn:sigma2}), we can find  a neighborhood
$V$ of
$z_0$ in $M$ such that setting ${\cal V}=\psi^{-1}_z(V) \subset G$,
we have ${\cal V}\cap \partial D(z)$ is biholomorphically equivalent to
$(V \cap \partial D) \times
H_z$. Noting that $H_z$ is an $(m-n)$-dimensional complex
manifold, the pseudoconvexity of $\partial D$ at $z_0$
yields the same property  for $\partial D(z)$ at $g_0$, proving 1.

To show 2., we need the following result of Kazama~\cite{K}: a complex
Lie group $G$ is {\it weakly complete}: $G$ admits a $C^{\infty}$
plurisubharmonic
exhaustion function; i.e., there exists a $C^{\infty}$ plurisubharmonic
function $\varphi $ in $G$ with $\{g\in
G:\varphi (g) < \alpha\}\Subset G$ for each $\alpha \in {\mathbb{R}}$.
Given $n \gg 1$, we set $$ D_n(z):= D(z) \cap \{g \in G : \varphi  (g)< n\}.
$$ By 1., we see that $D_n(z)$ is a pseudoconvex domain with piecewise
smooth boundary in $G$; clearly $D_n(z) \subset D_{n+1}(z)\Subset D(z)$ and
$D(z)=\cup_{n=1}^\infty D_n(z)$.  \hfill $\Box$

\smallskip
Now we discuss two types of variations of domains in $G$. First,
we have a variation of domains $D(z)$ in $G$ with parameter
$z\in M$:
$$
{\cal D}_M:z\in M \to D(z) \subset G.
$$
As usual we identify the variation
$ {\cal D}_M$ with the $(m+n)$-dimensional domain ${\cal D}_M=\cup _{z\in
M}(z,
D(z)) \subset M \times G$.
Given a domain
$V\subset M$, we consider the variation ${\cal D}_V: z\in V \to D(z)$,
and identify ${\cal D}_V$ with $\cup _{z\in V}(z, D(z))
\subset V
\times G$.

\smallskip
Next we define ${\cal D}: \ = \cup_{z\in D} (z, D(z))$, a variation of domains
$D(z)\subset G$ with parameter $z \in D$:
$$
{\cal D}: z\in D \to D(z) \subset G.
$$
Furthermore we define ${\cal D}':=\cup_{z\in D}(z,D'(z))$, where recall that $D'(z)$ is
the connected component of $D(z)$ containing the identity $e$. This is  a variation of domains
$D'(z)\subset G$ with parameter $z\in D$:
$$
{\cal D}': z\in D \to D'(z) \subset G.
$$
\begin{lemme} \label{lem:locally} \
 ${\cal D}$ and ${\cal D}'$ are locally holomorphically  trivial variations.
\end{lemme}
\noindent {\bf  Proof}.  Since the proofs are similar, we only give the
proof of the lemma for ${\cal D}'$.  Let $z_0\in D$. By
taking  $z=z_0$ in property {\bf  2.}, we can find a
neighborhood $V\subset D$ of $z_0$ and a holomorphic section $\sigma$ over
$V$ via $\psi_{z_0}$ such that
$\sigma(\zeta)(z_0)=\zeta$ for $\zeta\in V$ with $\sigma(z_0)=e$.
Since
$\sigma(\zeta)(z_0)=\zeta \in V \subset D$, we have $\sigma(\zeta)\in
D(z_0)$. It follows from $\sigma(z_0)=e\in D'(z_0)$ that
$\sigma(\zeta)\in D'(z_0)$.  By III in Proposition \ref{prop:useful} we
have $D'(z_0)= D'(\sigma(\zeta) ^{-1} (\zeta))= D'(\zeta)
\sigma(\zeta)$. Hence,
$$
T: \ (\zeta, g) \in {\cal D}'_V:=\cup_{\zeta\in V}(\zeta, D'(\zeta)) \to (\zeta, g \sigma(\zeta))\in V \times
D'(z_0)
$$
provides a trivialization.   \hfill $\Box$

\smallskip
Fix a K\"ahler metric $ds^2$ on  $G$ (such a metric exists by~\cite{K2})
and let $c$ be a strictly positive
$C^{\infty}$ function on $G$. As shown in the previous section, for a
fixed  $z\in D$, since $e\in D(z)$, we can form the $c$-Robin
constant $\lambda(z)$ for $(D(z),e)$. We recall
that the $c$-Robin constant is defined by the usual exhaustion method
in the case of an unbounded connected domain $D(z)$
(see, for example, Chapter IV
in ~\cite{AS}).
Furthermore, if $D(z)$ is not connected, we consider
the $c$-Robin
constant $\lambda(z)$ for $(D'(z),e)$
and we call this the $c$-Robin constant for $(D(z),e)$.
Using standard  methods
from potential theory, from Lemma \ref{lem:locally}
 we see that  $\lambda(z)$ is smooth in
$D$. Furthermore, since $\partial D$ is smooth in $G$,
Lemma \ref{re:ebound} implies that $\lambda(z) \to -\infty$ as $z \to
\partial
D$; i.e., $-\lambda (z)$ is a smooth exhaustion function for $D$. We
call $\lambda (z)$ the $c$-{\bf Robin function } for $D$.
Under these circumstances we have the following result.
\begin{theorem}\label{thm:exh} \
If $D\Subset M$ is a smoothly bounded pseudoconvex domain, then
the $c$-Robin function   $-\lambda(z)$ for $D$ is a
plurisubharmonic
 exhaustion function for $D$.
\end{theorem}
\noindent{\bf Proof}. It remains to prove the plurisubharmonicity of
$-\lambda (z)$ on $D$. Using $D_n(z) \Subset G$ from 2. of
Proposition \ref{prop:dz}, we set
${\cal D}_n=\cup _{z\in D}(z, D_n(z))$ and
${\cal D}_n'=\cup _{z\in D}(z, D'_n(z))$  where $D'_n(z)$ is the
connected component of $D_n(z)$ containing $e$.
By Lemma \ref{lem:locally}, ${\cal D}_n'$ as well as
${\cal
D}_n$ is pseudoconvex in $D \times G$.
Let $V \Subset D$ be a domain. We take $n$ sufficiently large so that $e\in D_n(z)$ for $z\in V$. Then for $z \in V$ we may consider
the $c$-Robin constant $\lambda_n(z)$ for $(D_n(z),e)$, i.e., for
$(D'_n(z),e)$. It follows from
Theorem \ref{thm:psm} that $-\lambda_n(z)$ is a plurisubharmonic
function on $V$. (We remark that we assumed $\partial {\cal D}$ was  smooth in
this theorem; however, it is standard to see  that the result remains true in the case $\partial {\cal D}$ is only assumed to be
piecewise smooth.) Since
$-\lambda _n(z)$ decreases to $-\lambda(z)$ in $V$, $-\lambda (z)$ is
plurisubharmonic on $V$, and hence on $D$.
\hfill $\Box$

\smallskip
We next discuss conditions under which $-\lambda (z)$ is {\it strictly} plurisubharmonic on
$D$. Suppose not; i.e., suppose there exists a point $z_0\in D$ at which
the complex Hessian $[{\partial ^2(-\lambda)\over
\partial z_j \partial \overline z_k}(z_0)]$ has a zero eigenvalue so that
$[{\partial ^2(-\lambda(z_0+at))\over
\partial t \partial \overline t}]|_{t=0}$ $=0$ for some direction
$a\in {\mathbb{C}}^n$ ($a\not = 0$).  A standard result in the
theory of homogeneous spaces is that $M \approx G/H_{z_0}$ via the map $\psi_{z_0}(g)=g(z_0)$; using this it follows that there exists a left-invariant holomorphic
vector field $X$ on $G$ such that
  the tangent vector of $(\exp t X)(z_0) $ at $t=0$ in $M$
is equal to $a$; and thus
$[\frac{\partial ^2(-\lambda( \exp tX(z_0))}{\partial t \partial \overline
{ t} }]|_{t=0}=0.$

\smallskip
Note that $\exp tX \in D'(z_0)$ for $|t|<\rho$ if $\rho$ is sufficiently small. It  follows from III. in
Proposition \ref{prop:useful} that
$$D'(\exp tX(z_0))=D'(z_0)\cdot (\exp
tX)^{-1}.$$
We thus  apply
Lemma \ref{thm:lie}
to $G, D'(z_0),$ and $e$ corresponding to $M,D,$ and $\zeta$ in
the lemma -- note that the domain $D$ in the lemma is
bounded, but the argument is valid for unbounded $D$ -- to obtain
 the following facts: if we assume $g\in D'(z_0)$, then
\begin{enumerate}
 \item [a.]  the integral curve $\{g \exp tX: \  t\in {\mathbb{C}}\}$ for $X$ is contained in
       $  D'(z_0)$ and;
\item [b.] $D'(z_0)\cdot g^{-1}= D'(z_0)\cdot (g\exp tX)^{-1}$ for all
       $t\in {\mathbb{C}}$; i.e., $D'(g(z_0))=D'(g\exp tX(z_0))$ {for all} $t\in {\mathbb{C}}$.
\end{enumerate}
\noindent Moreover, from the above formula
and the definition of $\lambda(z)$,
\begin{enumerate}
 \item [c.] $\lambda (g(z_0))=\lambda (g\exp tX(z_0))$ for all
$ t\in {\mathbb{C}}.$
\end{enumerate}
The set $\{g\exp tX(z_0):t\in {\mathbb{C}}\}$ is a one-dimensional holomorphic curve in $D$ since $a\not =
0$. Since $-\lambda$ is an exhaustion
function for
$D$, we have
\begin{enumerate}
 \item [d.] $\{g \exp tX(z_0) : t\in {\mathbb{C}}\}$ is relatively compact in $D$.
\end{enumerate}
Using the same argument that we used for the last assertion in Corollary
\ref{oldremark5.1} for $A=\partial D'(z_0)$ or $A= \overline { D'(z_0)}^c$, we have
\begin{eqnarray} \label{eqn:bd-and-out}
\begin{array}{lll}
 & g \exp tX \in A,  \ \    &\mbox{  for} \  t \in {\mathbb{C}}, \vspace{1mm}\\
& Ag  ^{-1}  =A(g \exp tX) ^{-1}  , \ \  &\mbox{  for} \ t\in {\mathbb{C}},
\end{array}\end{eqnarray}
for $g\in A$.

We conclude that if $-\lambda$ is not strictly plurisubharmonic at
$z_0$, then for any $g \in D'(z_0)$
we obtain a
one-dimensional holomorphic curve $g(C)\Subset D$
(where $C:=\{ \exp tX(z_0): t\in {\mathbb{C}}$\}) on which $\lambda$ is
constant, with value $\lambda (g(z_0))$.
Now take any point $z'\in D$. Using 1. in Proposition
\ref{eqn:elm}, there exists $g'\in D'(z_0)$ such that
$g'(z_0)=z'$. Then, as
mentioned above, the one-dimensional curve $C':= g'(C)$ passes through $z'$ and is relatively compact in $D$. Thus {\sl if $-\lambda$ is not strictly plurisubharmonic at {\bf some} point $z_0\in D$, then for {\bf each} point $z'\in D$, there exists a one-dimensional holomorphic curve
$C'\Subset D$ passing through $z'$}. (We mention that as a curve in $D$, $C'$ is  conformally equivalent to
 ${\mathbb{C}}, \, {\mathbb{C}}^*:={\mathbb{C}}\setminus \{0\},\,
 \mathbb P^1$ or a torus.) Hence any
 exhaustion function
$s(z)$ which is plurisubharmonic on $D$, in particular, the function
$-\lambda(z)$, must be constant on $C'$, and hence
 the complex Hessian $[{\partial ^2s\over
\partial z_j \partial \overline z_k}(z')]$ has a zero eigenvalue.

 In this vein,
we mention a result of Michel~\cite{Mi} for compact $M$: if $D$
admits a strictly pseudoconvex
boundary point, then $D$ is Stein. This follows in the case of
a domain $D$ with $C^{\infty}$ boundary
from our results. For if $p\in \partial D$ is a strictly pseudoconvex
boundary point and if $-\lambda$ is not strictly plurisubharmonic on $D$, then
$[{\partial
^2(-\lambda)\over
\partial z_j \partial \overline z_k}(z_0)]$ has a zero eigenvalue at some
$z_0\in D$; i.e., $[{\partial ^2\lambda(z_0+at)\over
\partial t \partial \overline t}]|_{t=0}=0$ for some direction
$a\in {\mathbb{C}}^n$ ($a\not = 0$). This implies (a. above) that there exists a
nonvanishing
left-invariant holomorphic vector field $X$ on $G$ such that if $g\in
\partial D'(z_0)$, then $g \exp tX\in  \partial D'(z_0)$ for all $t\in
{\mathbb{C}}$. Take
$g'\in \partial D'(z_0)$ with $g'(z_0)=p$; such a $g'$ exists by homogeneity.
Since
$a\not = 0$, the projection of this integral curve $\{g \exp tX: t\in
{\mathbb{C}}\}$ under $\psi_{z_0}$ is
a one-dimensional holomorphic curve through $p$ lying in $\partial D$,
contradicting the
strict pseudoconvexity of $\partial D$ at
$p$.

We also mention that Diederich-Ohsawa \cite{D-O} show that if $D$ is a pseudoconvex domain with $C^\omega$ boundary in a complex $2$-dimensional
manifold and $D $ admits a strictly
pseudoconvex boundary point, then $D$ is Stein.

\smallskip
What are some sufficient conditions on the pair  $(M,G)$ which insure
that for {\it any} relatively compact, smoothly bounded pseudoconvex domain
$D$
in $M$, the $c$-Robin function $\lambda(z)$ for $(D(z),e)$ has the property that
$-\lambda (z)$ is strictly plurisubharmonic on $D$\,? We discuss
two such
conditions.  Recall that $H_{z_0}'$ denotes the connected component of
 the isotropy subroup $H_{z_0}$ of $G$ for $z_0$ which contains the
 identity  $e$. First, we say that the pair $(M,G)$ satisfies the {\bf three
point property} if for any triple of distinct points $z_0,z_1$ and $z_2$
in $M$, there exists an element
$g\in H'_{z_0}$ with $g(z_1)=z_2$.
 As an example, take $M={\mathbb{C}}$ and $G=\hbox{Aut}M$, or for, $n>1$,
$M={\mathbb{C}}^n$ and $G$ the subgroup of Aut$M$ generated by $GL(n,{\mathbb{C}})$ and translations. Another example is $M={\mathbb{P}}^n$ and
$G=GL(n+1,{\mathbb{C}})$.
However the complex Grassmannian manifold $M=G(k,n)$ with $G=GL(n,{\mathbb{C}})$ does
not
 satisfy
the three point condition if $2\le k\le n-2$ (see {\bf 2.} in Appendix A).
 Next, we say that the
pair $(M,G)$ satisfies the {\bf spanning property} if for any $z_0
\in M$ and for any one-dimensional complex curve $C$ in a
neighborhood of $z_0$ in $M$ passing through $z_0$, we can find $h_i, \ i=1,\ldots ,n$, in
 $H'_{z_0}$ with the following
 property: let $C_i:=h_i(C)$ be the image of $C$ under $h_i$; this is a
curve in $M$ passing through $z_0$; if $a_i$ denotes the complex tangent
 vector of $C_i$ at $z_0$, then we require that the $n$ vectors $\{a_1, \ldots , a_n\}$
span ${\mathbb{C}}^n$. The
Grassmannian manifolds $M=G(k,n)$ with $G=GL(n,{\mathbb{C}})$ satisfy the spanning condition for any $k=1,...,n$; this is {\bf 3.} in Appendix A.

The importance of these notions occurs in
the following result.
\begin{theorem}
\label{3ptspan}\ Let $(M,G)$ satisfy the three point property or the spanning property.
Then for any pseudoconvex domain $D\Subset M$ with $C^{\infty}$ boundary,
if we let $\lambda(z)$ be the $c$-Robin function for $D$, then
$-\lambda(z)$
is a strictly plurisubharmonic exhaustion function for $D$; i.e., $D$ is
Stein.
\end{theorem}
\noindent {\bf Proof}.   By Theorem \ref{thm:exh},  $-\lambda(z)$
is a plurisubharmonic exhaustion function
for $D$. We prove the strict plurisubharmonicity by contradiction.
Suppose that $-\lambda(z)$ is not strictly plurisubharmonic at $z_0\in
D$. Then
the complex Hessian $[{\partial ^2(-\lambda)\over
\partial z_j \partial \overline z_k}(z_0)]$ has a zero eigenvalue. From
the above discussion, we conclude that there exists a left-invariant
holomorphic vector field $X$ on $G$ such
that setting $C: =\{\exp tX (z_0) :  t \in {\mathbb{C}}\}$, the tangent vector
of $C$ at $t=0$ is not zero, and for any $g\in D'(z_0)$
 we
obtain a one-dimensional holomorphic curve $g(C): =\{g \exp tX (z_0) :
t\in {\mathbb{C}}\} \Subset D$ on
which $\lambda$ is constant with value $\lambda (g(z_0))$, and
$D'(g \exp tX(z_0))=D'(g(z_0))$.

 We first treat the case when $(M,G)$ satisfies the three
point property.  Fix $z_1\in C\setminus \{z_0\}$.
Given $z\in D$ with $z\not =z_0,z_1$, from the three-point property we can find $g\in
H'_{z_0}$ such that
 $g(z_1)=z$. Since $g\in H'_{z_0} \subset D'(z_0)$,
it follows that $\lambda$ is constant on the curve
$g(C)$ with value $\lambda(g(z_0))=\lambda(z_0)$.
From $z=g(z_1)\in g(C)$, we have, in particular, $\lambda(z)=\lambda(z_0)$,\ i.e., $\lambda \equiv\  \hbox{const}.$ on $D$,
a contradiction.

We next treat the case when $(M,G)$ satisfies the
spanning condition.
Then we can find $h_i, \ i=1,\ldots, n$ in $H'_{z_0} \subset D'({z_0})$   such that
 $\{a_1, \ldots , a_n\}$ span ${\mathbb{C}}^n$, where $a_i$ is the complex
 tangent vector of the curve $h_i(C):=h_i \exp tX (z_0), t\in {\mathbb{C}}$  at
 $z_0$.
It follows from c. that
$$
\lambda(h_i \exp tX(z_0))= \lambda (h_i(z_0))=\lambda(z_0), \ \  \
\hbox{for all} \ t \in {\mathbb{C}}.
$$
In other words, $\lambda$ assumes the same value $\lambda(z_0)$ on each
curve $h_i(C), \ i=1,\ldots , n$ with different tangent directions at
$z_0$ which span ${\mathbb{C}}^n$. Therefore all eigenvalues of the complex Hessian
$[{\partial ^2(-\lambda)\over
\partial z_j \partial \overline z_k}(z_0)]$ are $0$. By the positive
semi-definiteness of this matrix, we conclude that $[{\partial ^2(-\lambda)\over
\partial z_j \partial \overline z_k}(z_0)]$ is the zero matrix.

Now take $z'\in D$ with $z'\not = z_0$. By 1. in Proposition
\ref{eqn:elm}
 we can find $g'\in D'(z_0) $ such that
$g'(z_0)=z'$. Let $C'=g'(C)$. Then $\lambda(z)$ is constant on $C'$ by
property c., so that
$[{\partial ^2(-\lambda)\over
\partial z_j \partial \overline z_k}(z')]$ has a zero eigenvalue.
Then we can find $h_i'\in H'_{z_0}, \ i=1,\ldots, n$ such that
 $\{b_1, \ldots , b_n\}$ span ${\mathbb{C}}^n$, where $b_i$ is the complex tangent vector of the curve $h_i'(C')$ at
$z'$. By the previous argument, all eigenvalues of the complex Hessian
$[{\partial ^2(-\lambda)\over
\partial z_j \partial \overline z_k}(z')]$ are $0$; hence this matrix is
the zero matrix. Since this argument is valid for any $z'\in D$, we
conclude that
$-\lambda$ is pluriharmonic on
$D$. Since $-\lambda$ is an exhaustion function on $D$, $-\lambda$ is
constant on $D$ by the minimum principle. This is a contradiction.
\hfill $\Box$

\smallskip
As mentioned, the complex Grassmannian manifold $G(k,n)$ satisfies the spanning
condition. Therefore, the $c$-Robin function
$-\lambda(z)$ for $ D$ is
a strictly plurisubharmonic exhaustion function for any pseudoconvex
domain $D$ with
$C^{\infty}$ boundary in $G(k,n)$. This result has been obtained by T.
Ueda~\cite{U} (also see
 A. Hirschowitz~\cite{hirs}) by a different method in the more general
situation when
 $D$ is any finite or infinitely sheeted unramified pseudoconvex
 domain over $G(k,n)$.
However, as will be shown in {\bf 4.} of Appendix A, the
flag space $M={\cal F}_n$ with $G=GL(n,{\mathbb{C}})$ for $ n\ge 3$ does not satisfy the
spanning property.
To analyze the pseudoconvex domains in ${\cal F}_n$ with
smooth boundary which are not Stein, we shall establish Theorem
\ref{thm:mainth} which addresses this issue for general homogeneous spaces.

We use the same notation $M,G, e $ with $\dim G=m$ and $\dim M=n \le
m$. Fix $z\in M$ and let $H_{z}$ be
 the isotropy subgroup of $G$ for $z$.
We let $\psi_{z}$ denote the
projection from $G$ onto $M\approx G/H_{z}$.
We write $H_z'$ for the connected component of $H_z$ containing
$e$.  Then $H_z'$  is a connected, closed, normal Lie subgroup of $H_z$
  with $\dim
H_z'=m_0:=m-n$.
Let $\mathfrak X   $ denote the (complex) Lie algebra consisting of all left-invariant holomorphic vector fields $X$ on $G$. Finally, $\mathfrak h _{z}$ denotes the Lie  subalgebra of $\mathfrak X $ corresponding to
the Lie subgroup
$H'_{z}$ of $G$ with $\dim \,\mathfrak h _{z}=m_0$; thus
$$
{\textstyle     H'_{z}=\{\,\prod _{j=1}^{\nu} \exp t_jX_j \ :\,\nu\in {\boldsymbol{Z}}^+, \  t_j\in {\mathbb{C}}, \
 X_j\in \mathfrak h_z\,\}}.
$$

Let  $D \Subset
 M$ be a pseudoconvex
domain with smooth boundary which is not Stein.
We form
the $c$-Robin function $\lambda(z)$ for $D$. Recall that we defined, for a given $z\in D$,
$$
D(z):=\{ g\in G \mid g(z) \in D\}=\psi_z^{-1}(D) \subset G,
$$
and considered the variation ${\cal D} $ of domains $D(z)$ in $G$ with
parameter $z\in D$:
$$
{\cal D}: z\in D \to D(z) \subset G.
$$
The $c$-Robin function $\lambda(z)$ was defined as the $c$-Robin constant
for $(D(z), e)$; precisely, for the connected component $D'(z)$ of $D(z)$
containing $e$.
Then $-\lambda (z)$ is  a plurisubharmonic
exhaustion function for $D$. Under the assumption that $D$ is not Stein,
there exists a point $z_0\in
D$  such
that $-\lambda(z)$ is not strictly  plurisubharmonic at $z_0$, i.e.,
there exists $a\in {\mathbb{C}}^n \ (a \ne 0) $ such that $[ \frac{\partial^2
\lambda (z_0+ at)}{\partial t \partial \overline
{ t}}]|_{t=0}=0.$ If  we take an arbitrary vector field $X\in \mathfrak X $   such that
$[\frac{d \exp tX (z_0)}{d t }]|_{t=0}=a$, then we have:
\begin{enumerate} {\it
 \item [o.]  \   $[\frac{\partial ^2 \lambda(\exp tX(z_0)) }{\partial t \partial \overline {
t}}]|_{t=0} =0$;
 \item [i.]  \ $\{\exp tX(z_0): \  t\in {\mathbb{C}}\}\subset D$ \,and \,  $\lambda
 ( \exp
 tX(z_0))= \lambda (z_0), \  t\in {\mathbb{C}}; $ indeed
 \item [ii.] \ $ \{\exp tX(z_0)  :\  t \in {\mathbb{C}} \} \Subset
         D$; moreover
 \item [iii.] \   $\{\exp tX: \ t\in {\mathbb{C}}\} \subset  D'(z_0)$ \,and \,$D'(\exp tX(z_0))= D'(z_0),   \ t\in {\mathbb{C}}$.}
\end{enumerate}
From properties of the Robin function $\lambda(z)$, these four
conditions are equivalent.  We write $o.\sim iii.$ and we consider
\begin{align}
 \label{eqn:X-0}
\mathfrak X _{z_0}:=
\{\, X\in \mathfrak X  : X  \ \mbox{satisfies  one (all) conditions(s) } \, o. \sim iii. \, \}.
\end{align}
We note that
\begin{equation} \label{strictinclusion} \mathfrak h _{z_0}\subsetneq  \mathfrak X _{z_0} \subsetneq \mathfrak   X. \end{equation} For suppose $X\in
 \mathfrak h _{z_0}$. Then $\exp t X\in H'_{z_0}$ for $t\in {\mathbb{C}}$.
 Hence $\exp tX(z_{0})=z_0$,  so that $X \in
 \mathfrak X _{z_0}$. Under the assumption that $D$ is not Stein, we have already shown that $\dim \mathfrak X _{z_0} > \dim \mathfrak h_{z_0}=m_0$. Next, since $-\lambda$ is a plurisubharmonic exhaustion function for $D$, $-\lambda$ is not pluriharmonic. Following the argument at the end of the proof of Theorem \ref{3ptspan}, we see that $\dim \mathfrak X> \dim \mathfrak X _{z_0}$.

We want to prove that $\mathfrak X _{z_0}$ is a complex Lie subalgebra of $\mathfrak
 X $, and verify some properties of $\mathfrak X _{z_0}$. To this end, we begin with a technical lemma.
\begin{lemme}\label{lem:fundamen} \
 Let $ \nu\in {\boldsymbol{Z}}^+$ and let $ g_i, h_i \in H'_{z_0}; \,  X_i\in
 \mathfrak X _{z_0}; \ t_i\in {\mathbb{C}}, \ i=1,2,\ldots,\nu
 $ and $g\in D'(z_0)$. Then
\begin{align*}
 (1) &  \quad
g \bigl[ {\textstyle  \prod_{i=1}^\nu g_i} (\exp t_i X_i)
 h^{-1}_i\bigr]\in D'(z_0);  \nonumber \\
(2) &  \quad    \lambda \bigl(g\bigl[{\textstyle\prod _{i=1}^\nu } (g_i (\exp t_i X_i)
 h^{-1}_i\bigr](z_0)\bigr)
= \lambda(g(z_0)).
\end{align*}
\end{lemme}
\noindent {\bf Proof}.
By the definition of the $c$-Robin constant
$\lambda(z)$ for $(D(z),e)$, to verify (2) it suffices to prove
\begin{align}
 \label{eqn:simple-but-good}
(2')\ & {\textstyle  \quad D'\bigl(g\bigl[\prod _{i=1}^\nu g_i (\exp t_i X_i) h^{-1}_i\bigr](z_0)\bigr)
= D'(g(z_0)).}\nonumber
\end{align}
We prove (1) and $(2')$ by induction. We will utilize some of the
properties
$a.-d.$ mentioned after the proof of Theorem \ref{thm:exh}. Let $\nu=1$. By 2. in Proposition
  \ref{eqn:elm}, we have $gg_1\in D'(z_0)$. This together with
property a. implies $C:=\{gg_1\exp tX_1 \in G:t\in {\mathbb{C}}\}\subset D'(z_0)$.
We see from property III. of Proposition \ref{prop:useful} and $h_1^{-1}\in H_{z_0}'$ that
$
Ch_1^{-1}\subset D'(z_0)h_1^{-1}=D'(h_1(z_0))=D'(z_0),
$
so that $(1)$ for $\nu=1$ is  proved.
Furthermore, since $gg_1\in D'(z_0)$ and $h_1^{-1}(z_0)=z_0$, it follows from property b. that
$$
  D'\bigl(gg_1 (\exp \,t_1 X_1 \bigr)h^{-1}(z_0))= D'\bigl(gg_1 (\exp
  \,t_1 X_1 )(z_0)\bigr)
=D'\bigl(gg_1(z_0)\bigr)=D'(g(z_0)),
$$
which proves $(2')$ for $\nu=1$.
 Next we assume that
(1) and $(2')$
are true for $\nu\le k$ and we shall prove them  for
 $\nu=k+1$.
We define, for $t\in {\mathbb{C}}$,
\begin{align*}
 k(t)&:= {\textstyle  g\bigl[ \prod _{i=1}^k g_i (\exp t_i X_i)
 h^{-1}_i \bigr]\,  ( g_{k+1}(\exp tX_{k+1})h_{k+1}^{-1})\in G, }
\end{align*}
and we consider the curve $C:=\{k(t)\in G :  t\in {\mathbb{C}}\}$ in $G$.
It suffices to prove that $C \subset D'(z_0)$ and
$D'(k(t)(z_0))=D'(g(z_0)), \
  t\in {\mathbb{C}}$. Note that
$$
k(t)= k_0\cdot  (\exp tX_{k+1}) h_{k+1}^{-1},$$ where
$$
k_0:={\textstyle   g \bigl[\prod _{i=1}^{k-1} g_i (\exp t_i X_i)
 h^{-1}_i \bigr]\,  (\,g_{k}(\exp t_kX_{k})(g_{k+1}^{-1}h_{k})^{-1}\,)}.
$$
Then the induction hypothesis for (1) in the case $\nu=k$ implies that $k_0\in
D'(z_0)$ since $g^{-1}_{k+1}h_k\in H_{z_0}'$; thus property a. yields that
$k_0 \exp tX_{k+1}\in D'(z_0), \ t\in {\mathbb{C}}$.
Since $h_{k+1}^{-1}\in H'_{z_0}$, it follows from
2. in Proposition \ref{eqn:elm}  that $C \subset D'(z_0)$, so that (1) for
the case $\nu=k+1$ is valid. Moreover, since
$g \prod _{i=1}^{k} g_i (\exp t_i X_i)
 h^{-1}_i \in
D'(z_0)$ (by the induction hypothesis (1) in the
 case $\nu=k$) and $ g_{k+1}\exp tX_{k+1}h_{k+1}^{-1}\in D'(z_0) $ (which follows by (1) in the case $\nu=1$ for $g=e$),  we again utilize property III of Proposition \ref{prop:useful}
 (twice) to obtain
\begin{align*}
 D'(k(t)(z_0))&
  =   D'
\bigl(\,[g {\textstyle\prod _{i=1}^{k}}g_i(\exp
 tX_i)h_i^{-1}] (g_{k+1}\exp tX_{k+1}h_{k+1}^{-1})(z_0)\,\bigr)\\
&\  = D'((g_{k+1}\exp tX_{k+1}h_{k+1}^{-1})(z_0))
[g {\textstyle \prod _{i=1}^{k}}g_i(\exp
 tX_i)h_i^{-1}]^{-1}\\
&\   = D'(z_0)[g {\textstyle\prod _{i=1}^{k}}g_i(\exp
 tX_i)h_i^{-1}]^{-1} \ \
 \\
&\  =  D'\bigl(\, \bigl[g {\textstyle \prod _{i=1}^{k}} g_i (\exp t X_i)
 h^{-1}_i\bigr](z_0)\, \bigr)    \quad  \\
&\ =D'(g(z_0)).
\end{align*}
The last equality follows from $(2')$ in the case $\nu=k$. This proves $(2')$ for $\nu=k+1$. \hfill $\Box$

\smallskip
Arguing as in (1) but using (\ref{eqn:bd-and-out}) and (\ref{eqn:baoundary-prime}) for $A=\partial D'(z_0), \  \overline { D'(z_0)}^c$ instead of a. and III. in
Proposition (\ref{prop:useful}) for $D'(z_0)$, we obtain
\begin{eqnarray}
\begin{split}
  \label{eqn:fund-boundary}
& {\textstyle  g \bigl[ \prod_{i=1}^\nu g_i (\exp t_i X_i)
 h^{-1}_i\bigr]\in A, \qquad  g\in
 A.}
\end{split}\end{eqnarray}
We can now prove the following fundamental result.
\begin{lemme} \label{le:2-is-new} \ {\empty}
\begin{enumerate}
 \item $\mathfrak X _{z_0}$ is a complex Lie subalgebra of $\mathfrak
 X $.
 \item Let $X\in \mathfrak X _{z_0}$ and $g\in H'_{z_0}$. Then any $Y\in
 \mathfrak X $ which satisfies
\begin{eqnarray} \label{eqn:kanari-imp}
 [\frac{d\exp tY(z_0) }{d t
 }]|_{t=0}\ &=& \
[\frac{dg \exp tX(z_0) }{dt }]|_{t=0}
\end{eqnarray}
belongs to $\mathfrak X _{z_0}$.
 \item For $z_1\in D$ we have $\dim\,\mathfrak X _{z_1}=\dim\, \mathfrak X
 _{z_0}$.
\end{enumerate}
 \end{lemme}

 \smallskip
 \noindent Note that (\ref{eqn:kanari-imp}) says that the tangent vector of the curve
       $C_Y:=\{\exp tY(z_0)\in M: t \in {\mathbb{C}}\}$ at $z_0$ in $M$
 coincides with that of the curve $g(C_X)$ at $z_0$.
 \smallskip

\noindent {\bf Proof}.
To verify 1., let $X,Y\in \mathfrak X _0$. Fix $\alpha, \beta\in {\mathbb{C}}$. If we set
$t_1=\alpha t, \ t_2= \beta t; \ X_1=X, \  X_2=Y; $ and
$g_1=g_2=h_1=h_2=e$ in Lemma \ref{lem:fundamen} (2), we obtain
$$
\lambda ( (\exp \alpha tX)(\exp \beta tY)(z_0))= \lambda(z_0),  \quad
\mbox{  for $t\in {\mathbb{C}}$}.
$$
Then as we deduced formula (\ref{eqn:tri-but}) from (\ref{eqn:well-known-1}), we obtain
$$
[\frac{\partial^2 \lambda (\exp (\alpha tX+ \beta tY)(z_0)) }{\partial t
\partial \overline { t}}]|_{t=0} = [\frac{\partial^2 \lambda (z_0)
}{\partial  t \partial \overline { t} }]|_{t=0}=0.
$$
Hence $\alpha X+ \beta Y\in \mathfrak X _{z_0}$. Next, if we set
$t_4=t_3=-\sqrt t, t_2= t_1=\sqrt t; \ X_4=X_2=Y, X_3=X_1= X; g_i=h_i=e \
(i=1,\ldots , 4)$ in Lemma \ref{lem:fundamen} (2), we obtain
$$
\lambda ([( \exp (- \sqrt t)Y)(\exp (- \sqrt t)X)(\exp \sqrt tY)(\exp
\sqrt tX)(z_0)])=
\lambda(z_0),  \ \ \mbox{  for $t\in {\mathbb{C}}$}.
$$
Then as we deduced  (\ref{eqn:tri-but2}) from (\ref{eqn:well-known-2}),  we obtain
$$
[\frac{\partial^2 \lambda ((\exp t[X,Y])(z_0) )}
{\partial t \partial \overline { t} }]|_{t=0}=[
 \frac{\partial^2 \lambda (z_0)}{\partial t \partial \overline { t} }]
|_{t=0}=0.
$$
This shows that $[X,Y]\in \mathfrak X _{z_0}$, and  $1.$ is proved.

To prove 2., we let $X\in \mathfrak X _{z_0}, g\in H'_{z_0}$ and assume $Y\in
\mathfrak X $ satisfies (\ref{eqn:kanari-imp}). By (2) of Lemma
\ref{lem:fundamen}, for $g\in H'_{z_0}$ we have
$$
\lambda (g \exp tX(z_0))= \lambda (g(z_0))= \lambda (z_0), \quad \mbox{
for $\in
{\mathbb{C}}$},
$$
so that  $[\frac{\partial^2 \lambda (g \exp tX(z_0)) }{\partial t
\partial \overline{  t} }]|_{t=0}=0 $.  It follows from condition
(\ref{eqn:kanari-imp}) and Proposition \ref{prop:often} that
$$ [\frac{\partial^2 \lambda (\exp tY(z_0)) }{\partial t
\partial \overline{  t} }]|_{t=0}=0,$$ and hence $Y\in \mathfrak X
_{z_0}$.

Assertion 3. is intuitively clear from the property that, for $g\in D'(z_0)$,  $$\lambda (g(z_0))=\lambda (g\exp tX(z_0))$$ for all
$ t\in {\mathbb{C}}$ (this is c. listed after the proof of Theorem \ref{thm:exh}). The
rigorous proof is fairly technical. We begin by setting $\dim \mathfrak X _{z_0}= m_0+m_1$ and $m_2=m-(m_0+m_1)$ where $m= \dim \, G$ and $m_0=\dim \mathfrak h
_{z_0}$. Note that $m_0$
is independent of the point $z_0$ but, apriori, $m_1 $  may depend
on $z_0$. We  consider the following direct sum decompositions:
\begin{equation}
\begin{split}
\label{eqn:decomposition}
\mathfrak X _{z_0}= \mathfrak h _{z_0} \  \dotplus\  \mathfrak k_{z_0}
 \,;  \qquad  \mathfrak X &= \mathfrak h _{z_0} \ \dotplus \  \mathfrak k_{z_0} \
 \dotplus\  \mathfrak
 l_{z_0}.
\end{split}
\end{equation}
Here $\mathfrak k_{z_0}$ and $\mathfrak l_{z_0}$ are defined by these relations. Let $$\{X_1,\ldots , X_{m_2}\}, \  \{ Y_1,\ldots , Y_{m_1}\} \ \hbox{and} \  \{Z_1,\ldots , Z_{m_0}\}$$ be bases of $\mathfrak l _{z_0}, \ \mathfrak
 k_{z_0}$ and $\mathfrak h_{z_0}$. We fix canonical local
 coordinates on $V\subset {\mathbb{C}}^m$ for a neighborhood ${\cal V} \subset G$ of $e$  such that $g\in {\cal V}$ corresponds to $(\xi;\eta;\zeta)\in V$; precisely,
\begin{eqnarray*}
 &   g&= {\textstyle  \prod_{i=1}^{m_2}\exp a_i X_i \prod_{j=1}^{m_1}
\exp b_jY_j \prod _{k=1}^{m_0} \exp c_k Z_k \in {\cal V}} \\
&{}& \longleftrightarrow \boldmath (\xi;\eta;\zeta)=( a_1, \ldots , a_{m_2};b_1, \ldots ,b_{m_1};c_1,
 \ldots , c_{m_0})\in V,
\end{eqnarray*}
where $|a_i|,|b_j|,|c_k|<\rho$ (for $\rho >0$ sufficiently small), $e$ corresponds to $0=(0;0;0)$, and
\begin{align*}
 & \psi_{z_0}(g) \longleftrightarrow \boldmath (\xi;\eta)= (a_1, \ldots ,
 a_{m_2};b_1, \ldots ,b_{m_1})\in V'\subset {\mathbb{C}}^{m_2+m_1}
\end{align*}
are local coordinates of the  neighborhood  ${\cal V}'= \psi_{z_0}({\cal
 V}) \subset M$ of
 $z_0$ where $z_0$ corresponds to $0=(0;0)$.

Fix $ z_1\in D$. Using 1. in Proposition \ref{eqn:elm} we can find
$g_1\in D'(z_0)$ with $g_1(z_0)=z_1$. As in
(\ref{eqn:decomposition})  we
consider the direct sums
 $$ \mathfrak X _{z_1}= \mathfrak h _{z_1} \, \dotplus\,  \mathfrak
 k_{z_1};  \qquad
 \mathfrak X = \mathfrak h _{z_1} \, \dotplus \,  \mathfrak k_{z_1} \,
 \dotplus \,  \mathfrak
 l_{z_1}.$$ We set $\dim \, \mathfrak X _{z_1}=m_0+ m_1'$ where
 $m_0={\rm dim} \, \mathfrak h_{z_1}$.
To verify 3. it suffices to show that $m_1=m_1'$. Similar to before, we consider
 bases $$\{X'_1,\ldots , X'_{m_2'}\}, \ \{ Y'_1,\ldots , Y'_{m_1'}\} \ \hbox{and} \ \{Z'_1,\ldots , Z'_{m_0}\}$$ of $\mathfrak l _{z_1}, \ \mathfrak
 k_{z_1}$ and $\mathfrak h_{z_1}$. Here, $m_2'=m-(m_0+m_1')$. We fix canonical local
 coordinates on $V_1\subset {\mathbb{C}}^m$ for a  neighborhood ${\cal V}_1\subset
 G$ of $e$ with
\begin{eqnarray*}
 &g&={\textstyle   \prod_{i=1}^{m'_2}\exp a_i X'_i\prod_{j=1}^{m'_1}\exp b_iY'_j
  \prod_{k=1}^
{m_0} \exp c_k Z'_k
 \in {\cal V}_1}\\
&{}& \longleftrightarrow \boldmath (\xi';\eta';\zeta')=( a_1, \ldots , a_{m_2'};b_1, \ldots ,b_{m_1'};c_1,
 \ldots , c_{m_0})\in V_1,
\end{eqnarray*}
where $|a_i|,|b_j|,|c_k|<\rho$ for $\rho >0$ sufficiently small; $e$ corresponds to $0=(0;0;0)$; and
\begin{align*}
 & \psi_{z_1}(g) \longleftrightarrow \boldmath (\xi';\eta')= (a_1, \ldots ,
 a_{m_2'};b_1, \ldots ,b_{m_1'})\in V_1'\subset {\mathbb{C}}^{m_2'+m_1'}
\end{align*}
are local coordinates for the neighborhood  ${\cal V}_1'=\psi_{z_1}
 ({\cal V}) \subset M$ of
 $z_1$ where $z_1$ corresponds to $0=(0;0)$.

Given $X+Y= \sum_{ i=1  }^{m_2   }a_iX_i + \sum_{ j=1  }^{ m_1  } b_jY_j
\in \mathfrak l_{z_0}+ \mathfrak k _{z_0} $, we consider the
curve $C: =\{\exp t(X+Y)(z_0) :  t\in
{\mathbb{C}}\}$ in $M$ with initial value $z_0$ and tangent
vector
$\tau :=[\frac{d \exp t(X+Y)(z_0) }{d t }]|_{t=0} \in {\mathbb{C}}^{m_2+m_1}={\mathbb{C}}^{m-m_0}$ at $z_0$ in $M$. Since
$$
\exp [t(X+Y)+O(t^2)]=
{\textstyle  \prod_{i=1}^{m_2} \exp ta_iX_i \prod_{j=1}^{m_1} \exp tb_jY_j}$$
near $t=0$, we can write $\tau = (a_1, \ldots, a_{m_2};b_1,\ldots, b_{m_1}
):=({\bf  a}, {\bf  b})$
in terms of the canonical local coordinates of $V'$.
We note that this correspondence $X+Y\in \mathfrak k _{z_0}+ \mathfrak l_{z_0}
\longleftrightarrow \tau=({\bf  a};{\bf  b}) \in {\mathbb{C}}^{m_2+m_1}$ is
linear and one-to-one; moreover, $X+0$ ($0+Y$) corresponds to $({\bf
a};0)$ ($(0;{\bf  b})$). We consider the image curve $g_1 (C)$ with initial value at $z_1$ in
$M$ and tangent vector $\tau':=[\frac{d g_1\exp t(X+Y)(z_0) }{\partial t
}]|_{t=0}
\in {\mathbb{C}}^{m'_2+m'_1}={\mathbb{C}}^{m-m_0}$ at $z_1$ in $M$  (in terms of the local coordinates $(\xi';\eta')$ for $V_1'$). Then the correspondence $L:\tau \in
{\mathbb{C}}^{m-m_0}\to \tau'\in {\mathbb{C}}^{m-m_0}$ is linear and one-to-one
since $\tau \ne 0 $ implies $\tau_1\ne 0$.
 We consider the case when $X=0$, that is, ${\bf  a}=0$, and
$Y=\sum_{ j=1  }^{ m_1  } b_jY_j\ne 0$, so that
$\tau=(0;b_1,\ldots,b_{m_1})\ne (0;0)$.
Since $Y\in \mathfrak X _{z_0}$
and $g_1\in D'(z_0)$, it follows from property c. that
$\lambda (z)$ is constant on the curve $g_1(C)$.
  Now by a standard result in the
theory of homogeneous spaces there exists $\widetilde { Y}\in
\mathfrak X $  such that the tangent vector  $[\frac{d \exp t \widetilde
{ Y} (z_1) }{dt}]|_{t=0}$
of the curve $\{\exp \,t \widetilde { Y}(z_1): \ t\in {\mathbb{C}}\}$ at $z_1$ in $M
$ equals $\tau'$
(in the local coordinates $(\xi';\eta')$ for $V_1'$); thus $[\frac{\partial ^2
\lambda(\exp t \widetilde { Y}(z_1))}{\partial t
\partial \overline { t}}]|_{t=0}=0 $, so that $\widetilde { Y}\in
\mathfrak X _{z_1}$. We conclude that $\tau'$ is
of the form $\tau'=(0;b_1', \ldots ,b_{m_1'}')\in
{\mathbb{C}}^{m_2'+m_1'}$. Since $L$ is linear and $\tau' \ne
(0;0)$, we have $m_1\le m_1'$. Reversing the roles of $z_0$ and $z_1$ gives the reverse inequality; hence $m_1'=m_1$.   \hfill $\Box$

\smallskip
Using 1. in  Lemma \ref{le:2-is-new} we can construct the connected Lie subgroup $\Sigma_{z_0}$ of
$G$ which corresponds to $\mathfrak X _{z_0}$, and
the subset
$\sigma_{z_0}$ of $M$:
\begin{equation}
\begin{split}
\label{eqn:Sigma-and-sigma}
\Sigma_{z_0}&=\{\,{\textstyle  \prod _{i=1}^\nu \exp t_i X_i \in G:
\nu\in {\boldsymbol{Z}}^+, \ t_i\in {\mathbb{C}},
\ X_i\in \mathfrak X_{z_0}}\,\};\\
\sigma_{z_0}&:= \Sigma_{z_0}(z_0):= \{\,g(z_0)\in M: g\in\Sigma_{z_0}\,\}.
\end{split}
\end{equation}
We have the following.
\begin{corollary} \label{cor:sigma-0-es}${}$
\begin{enumerate}
 \item If $g\in D'(z_0)$, then $g \Sigma_{z_0} \subset D'(z_0)$;
 $\lambda(z)=\,const.=\lambda (g(z_0))$ on $g \sigma_{z_0}$; and
 $g \sigma_{z_0} \Subset D$.
 \item For $z_1\in D$, we define $ \Sigma_{z_1} \subset G$
 and $\sigma_{z_1}= \Sigma_{z_1}(z_1)
 \subset M$
as in (\ref{eqn:Sigma-and-sigma}).
Then we have, for $g\in D'(z_1)$, $\lambda =\,const.=
 \lambda(g(z_1))$ on $g\sigma_{z_1}$; and $g \sigma_{z_1} \Subset D$.
 \item If $g\in \partial D'(z_0)$, then $g \Sigma_{z_0} \subset \partial
 D'(z_0)$ and $g \sigma_{z_0} \subset \partial D$.
\end{enumerate}
\end{corollary}
\noindent {\bf  Proof}. The first assertion follows from Lemma \ref{lem:fundamen} and
(\ref{eqn:Sigma-and-sigma}). To verify 2., let $z_1\in D$ and take
$g\in D'(z_1)$  with $z_1=g(z_0)$. Using 1., $z_0$ and $z_1$ play the same role as in 3. of
 Lemma \ref{le:2-is-new}, showing that 2. for $z_1$ as well as
 2. for $z_0$ are valid. The first formula in 3. follows from
 (\ref{eqn:fund-boundary}). The second formula follows from the first and the fact that $[\partial D'(z_0)](z_0)= \partial D$ (cf., (\ref{awkward})).
 \hfill $\Box$

\smallskip By 1. in Proposition \ref{eqn:elm} we thus have
\begin{equation}
 \label{eqn:bounda}
\begin{split}
 D&= \mbox{  $\bigcup_{g\in D'(z_0)}$}\, g\sigma_{z_0}\ \ \ \mbox{  with } \ g\sigma_{z_0}
 \Subset D; \\
 \partial D&= \mbox{  $\bigcup _{g\in \partial D'(z_0)}$}\, g\sigma_{z_0}.
\end{split}\end{equation}
Note that we have not proved that (\ref{eqn:bounda}) provide foliations of $D$ and $\partial D$. Moreover, we remark that if $g\in \overline { D'(z_0)}^c$, then $g \Sigma_{z_0}
\subset \overline { D'(z_0)}^c $ (as can be seen using an argument   similar to that used to verify the first assertions
in 1. or 3.); however, it is not necessarily true that $g \sigma_{z_0}
\subset \overline { D}^c$. In Theorem
\ref{thm:main-nonconnected} we will need a fairly technical argument  to
find  a domain  $K$ with $D \Subset K \subset M$  such that $g
\sigma_{z_0}\Subset  \overline { D}^c$ for $g\in K'(z) \cap \overline {
D'(z)}^c$.

\smallskip
We introduce the following terminology. Let $\sigma$ be a subset of
$M$. We call
$\sigma$ a {\it $\mu$-dimensional
non-singular, generalized ($f$-generalized)
analytic set in $M$} if, for any point $z\in \sigma$ ($z \in M$), there exists a neighborhood  $V$
of $z$ in $ M$ such that each connected component of $V \cap \sigma$ is a $\mu$-dimensional non-singular analytic set in $V$. Moreover, if $\sigma$ admits no non-constant bounded
 plurisubharmonic functions, then $\sigma $ is said to be
 {\it parabolic.} As a simple example to illustrate the difference between generalized and $f$-generalized analytic sets, let
 $$\Sigma=\{z=(z_1,z_2)\in {\mathbb{C}}^2: z_2=z_1\}$$
 and take $\sigma:=\Sigma \cap B$ where $B$ is the open unit ball
 $$B=\{z=(z_1,z_2)\in {\mathbb{C}}^2: |z_1|^2 + |z_2|^2 <1\}.$$
 Then $\sigma$ is a generalized analytic set in ${\mathbb{C}}^2$ but it is not an $f$-generalized analytic set in ${\mathbb{C}}^2$ because for $z=(z_1,z_1)\in \partial B$, there is no neighborhood $V$ of $z$ in ${\mathbb{C}}^2$ such that $V \cap \sigma$ is an analytic set in $V$.

\smallskip We first prove the main theorem in case the isotropy subgroup $H_z$
 is connected.
\begin{theorem}\label{thm:mainth}\ Let $M$ be a complex homogeneous
 space of finite dimension $n$ and
let $G$ be a  connected complex Lie transformation group of finite dimension $m\ge
 n$ which acts transitively on $M$. Assume that the isotropy subgroup $H_z$ of $G$ for $z\in M$ is
 connected. Let $D $ be a pseudoconvex domain in $M$ with smooth boundary. Assume that  $D $ is not
 Stein and fix $z_0\in D$. Then:
\begin{enumerate}
 \item [1-a.] The subset $\sigma_{z_0}=\Sigma_{z_0}(z_0)$ of $M$ defined in (\ref{eqn:Sigma-and-sigma})
is  a parabolic irreducible $m_1$-dimensional non-singular
$f$-generalized analytic set in   $M$ passing through $z_0$
satisfying the following properties:
\begin{enumerate}
\item [(o)]\   $\sigma_{z_0} \Subset D$;  $\Sigma_{z_0}= G \cap Aut\,\sigma_{z_0}$;  $\Sigma_{z_0}$ acts transitively on
 $\sigma_{z_0}$; and $\sigma_{z_0} \approx  \Sigma_{z_0} /H_{z_0}$.
 \item [(i)]\ $M$ ($D$) is foliated by cosets $g\sigma_{z_0}, \ g\in G$ ($g\in D(z_0)$) and each set $g\sigma_{z_0}$ is relatively compact  in $M$ ($D$):
$$
M=\mbox{  $\bigcup _{g\in G}$} \, g\sigma_{z_0}  \quad   and\  \quad   D
= \mbox{  $\bigcup_{g\in D(z_0)}$}\,
 g\sigma_{z_0}.\
$$
 \item [(ii)] \  Any parabolic non-singular  generalized
analytic set $\sigma$ in $M$  which is relatively
      compact  in $D$ is contained in a set $g\sigma_{z_0}$ for some $g\in D(z_0)$.
\end{enumerate}
 \item [1-b.] Assume that $\sigma_{z_0}$ in 1-a. is closed in
       $M$. Then
\begin{enumerate}
 \item [(o)]  $\sigma_{z_0}$ is an irreducible, $m_1$-dimensional compact
       non-singular analytic set in $D$ and $\Sigma_{z_0}$ is a closed
       connected Lie
       subgroup of  $G$.
   \item [(i)] Let $M_0:=G/\Sigma_{z_0}$ and let  $
 \pi_0:\ M=G/H_{z_0}
 \mapsto M_0$ be the canonical projection. Then ${\pi}_0^{-1}(\zeta) \approx \sigma_{z_0}$ (as complex
       manifolds) for each $\zeta\in
       M_0$; moreover, there exists a Stein domain $D_0 \Subset M_0$ with
       smooth boundary  such that $D= { \pi}_0  ^{-1}(D_0)$.
\end{enumerate}
\end{enumerate}
 \end{theorem}
\noindent {\bf Proof}. From the hypotheses of the theorem and 2. of Corollary
\ref{cor:sigma-0-es} the function $-\lambda(z)$ is not strictly
plurisubharmonic at $z_0$.
We maintain the notation $\mathfrak h _{z_0}$ and $\mathfrak X _{z_0}$ utilized in the proof of Lemma
 \ref{le:2-is-new}. Since $H_{z_0}$ is assumed to be connected,
 i.e., $H_{z_0}=H_{z_0}'$,  the Lie subalgebra  $\mathfrak h _{z_0}$ of
 $\mathfrak X  $ corresponds to the
 Lie subgroup $H_{z_0}$  in $G$; $D(z_0)$ is connected in $G$, i.e.,
 $D(z_0)=D'(z_0)$  and $H_{z_0} \subset D(z_0)$; moreover, each $D(z)\subset G,  \
 z\in D$, is connected. We let $\Sigma_{z_0}$ denote the connected Lie subgroup of $G$ corresponding to the
 Lie subalgebra $\mathfrak X _{z_0}$. We consider the direct sum decompositions from (\ref{eqn:decomposition}):
  $$
 \mathfrak X _{z_0} = \mathfrak h _{z_0} \dotplus  \mathfrak k_{z_0}; \qquad
\mathfrak X = \mathfrak h _{z_0}  \dotplus  \mathfrak k_{z_0} \dotplus
 \mathfrak l_{z_0}.
$$
Here $\dim \mathfrak h_{z_0}=m_0$ and we set $\dim \mathfrak X
 _{z_0}:=m_0+m_1$
  and $\dim \mathfrak l_{z_0}:=m_2$ so that
 $m=m_0+m_1+m_2$.
Let $\boldmath X=\{X_1, \ldots , X_{m_2}\},  \boldmath Y=\{Y_1, \ldots
 , Y_{m_1}\}$ and $\boldmath Z=\{Z_1, \ldots , Z_{m_0}\}$ be
bases of $\mathfrak l_{z_0}, \mathfrak k_{z_0}$ and $ \mathfrak
 h_{z_0}$. Then we have
\begin{equation}
 \begin{split} \label{eqn:imphyo}
   H_{z_0}&=\{\,{\textstyle   \prod _{j=1}^{\nu} \exp t_jA_j :\,\nu\in {\boldsymbol{Z}}^+, \  t_j\in {\mathbb{C}}, \
 A_j\in \boldmath Z}\,\};\\
\Sigma_{z_0}&= \{\,{\textstyle   \prod _{j=1}^{\nu} \exp t_jA_j  :\,\nu\in {\boldsymbol{Z}}^+, \  t_j\in {\mathbb{C}}, \
 A_j\in \boldmath Y \cup \boldmath Z} \,\}.
\end{split}\end{equation}
Note first that $H_{z_0}$ is an $m_0$-dimensional, non-singular closed
analytic set in $G$. From the (holomorphic) Frobenius theorem stated in the beginning of section 5, $\Sigma_{z_0}$ is an  $(m_0+m_1)$-dimensional non-singular
$f$-generalized anaytic set in $G$ with the property that
$G$ is foliated by the sets $\{gH_{z_0}\}_{g\in G}$ as well as by the sets $\{
g\Sigma_{z_0}\}_{g\in G}$. Thus,
if $g,\ g'\in G$,
\begin{equation}
 \begin{split} \label{eqn:fol-HS}
& g'\in gH_{z_0} \Longleftrightarrow  gH_{z_0}\cap g'H_{z_0}\ne \emptyset  \  \Longleftrightarrow \
gH_{z_0}=g'H_{z_0};\\
&g'\in g\Sigma_{z_0} \Longleftrightarrow g \Sigma_{z_0}\cap g'
 \Sigma_{z_0}\ne \emptyset \ \  \Longleftrightarrow \ \
g\Sigma_{z_0}=g' \Sigma_{z_0}.
\end{split}
\end{equation}
Given $g\in G$, we will call $gH_{z_0}$ ($g\Sigma_{z_0}$) an integral
manifold for $\mathfrak h_{z_0}$ ($\mathfrak X _{z_0}$) in $G$
passing through $g$.  Since $H_{z_0}$ is closed and connected in $G$, we
have the following result on the integral
manifolds for $\mathfrak h_{z_0}$: if we fix $g_0\in G$, then there exists a
neighborhood  $V$ of $g_0$ in $G$  such that for $g,\ g'\in V$,
 \begin{align*}
 [(gH_{z_0}) \cap V] \cap [(g'H_{z_0}) \cap V]=\emptyset \ \
  \Longrightarrow \ \
(gH_{z_0})\cap (g'H_{z_0})=\emptyset.
\end{align*}
The analogous property does not necessarily hold for $\Sigma_{z_0}$ since
$\Sigma_{z_0}$ is not necessarily closed in $G$.

Fix $g^0\in G$.
From item 3. in the discussion of the Frobenius theorem in section 5, we have the following canonical local
coordinates $(\Phi_0, V_0, \Delta_0)$ of $g^0$ in $G$, which we call
$F$-local
coordinates of $g^0$:
\begin{enumerate}
 \item [i)]\  $V_0$ is a  neighborhood  of $g^0$ in $G$;
 \item [ii)]\ $\Delta_0=\Delta'_0 \times \Delta''_0 \times
  \Delta'''_0 \subset {\mathbb{C}}^{m_2} \times {\mathbb{C}}^{m_1} \times {\mathbb{C}}^{m_0}$ where
$$
\Delta'_0 ={\textstyle  \prod}_{i=1}^{m_2}\{|\xi_i|<r_i\},  \quad
 \Delta''_0 ={\textstyle  \prod}_{i=1}^{m_1}\{|\eta_i|<r_i\},
 \quad  \Delta'''_0 ={\textstyle  \prod}_{i=1}^{m_0}\{|\zeta_i|<r_i\};
$$
 \item [iii)]\ $\Phi_0:g\in V_0 \mapsto (\xi(g),\eta(g),\zeta(g))\in\Delta_0$ is a holomorphic
 isomorphism  with $\Phi_0(g^0)=0\in \Delta_0$
\end{enumerate}
which satisfies the properties:
\begin{enumerate}
 \item [i')]\ given $\xi^0=(\xi_1^0, \ldots , \xi_{m_2}^0)\in
 \Delta_0'$, the $(m_0+m_1)$-dimensional polydisk in
 $\Delta_0$ defined by the equations
$$
S_{\xi_0}:  \ \  \xi_i=\xi_i^0 , \qquad  i=1,\ldots, m_2,
$$
corresponds to  an
integral manifold for $\mathfrak X _0$ from the map  $\Phi_0$, i.e.,
if $g\in \Phi_0^{-1}(S_{\xi_0})$, then $\Phi_0^{-1}(S_{\xi_0})=$
  the connected component of $[g
 \Sigma_{z_0}] \cap V_0 $ containing $g$;

\item [ii')]\ given $(\xi^0, \eta^0)=(\xi_1^0, \ldots , \xi_{m_2}^0;
\eta_1^0,\ldots , \eta_{m_1}^0)\in
 \Delta_0' \times \Delta_0''$, the $m_0$-dimensional polydisk in
 $\Delta_0$ defined by the equations
$$
S_{(\xi^0,\eta^0)}:  \ \  \xi_i=\xi_i^0 , \ \eta_j=\eta_j^0,   \qquad
i=1,\ldots, m_2; \ j=1,\ldots , m_1 ,
$$
corresponds to an
integral manifold for $\mathfrak h _0$ from the map  $\Phi_0$, i.e.,
if $g\in \Phi_0^{-1}(S_{(\xi^0, \eta^0)})$, then
       $\Phi_0^{-1}(S_{(\xi^0,\eta^0)}) =[g
 H_{z_0}] \cap V_0 $.
\end{enumerate}
From property ii')  we have:
\begin{enumerate}
 \item [iii')]\ fix $(\xi^0,\eta^0), \
(\xi^1,\eta^1) \in \Delta_0' \times \Delta_0''$, and let   $g^0\in
\Phi_0^{-1}(S_{(\xi^0, \eta^0)})$ and $g^1\in
\Phi_0^{-1}(S_{(\xi^1, \eta^1)})$. If
$
(\xi^0,\eta^0)\not=(\xi^1,\eta^1)$ then $g^0H_{z_0}\cap
        g^1H_{z_0}=\emptyset$.
\end{enumerate}
Using the local coordinate system $F = \{(\Phi_0, V_0, \Delta_0)\}_{g^0\in G}$ for $G$ at $g^0$, we canonically obtain a local coordinate system $
 {\bf   f} = \{(\phi_0, v_0, \Delta_0' \times \Delta_0'')\}_{z^0\in M}
$
for $M$ at a point $z^0$ as follows: fix $z^0\in M$ and take $g^0\in G$ with $\psi_{z_0}(g^0)=g^0(z_0)=z^0$. Then we have $F$-local coordinates
$(\Phi_0,V_0,\Delta_0)$ at $g^0$ for $G$. We set
$v_0=\psi_{z_0}(V_0 ) \subset M$; this is a neighborhood of $z^0$ in
$M$. Take $z'\in v_0$ and choose $g'\in V_0$ such that
$\psi_{z_0}(g')=z'$. Let $\Phi_0(g')=(\xi',\eta',\zeta')\in
\Delta_0$ and consider the correspondence
$$
\phi_0: z' \mapsto (\xi', \eta')\in \Delta_0' \times \Delta_0''.
$$
Since $M=G/H_{z_0}$,  we see from condition ii')  that $\phi_0(z')$ does not depend on the
choice of $g'\in V_0$ with $\psi_{z_0}(g')=z'$. Thus $\phi_0$ defines
a mapping from $v_0$ onto $\Delta_0' \times \Delta_0''$.  We also see from
iii') that  $\phi_0$ is one-to-one. Moreover, $\phi_0$ is
holomorphic. To see this, we
 consider a  holomorphic section $\tau$ of $G$ over $v_0$ via
 $\psi_{z_0}$  such that
$$
\tau:\ z'\in v_0 \mapsto \tau(z')\in V_0
$$
and $\tau(z^0)=g^0$.  We let $\Phi_0(\tau(z'))=(\xi(z'), \eta(z'),
\zeta(z')) \in \Delta_0$, so that $z'\to(\xi(z'), \eta(z'),
\zeta(z'))$ is a  holomorphic mapping from $v_0$ into $\Delta_0$.
 On the
other hand, by the definition of $\phi_0$ we have
$
\phi_0(z')= (\xi(z'), \eta(z'))\in \Delta'_0 \times
\Delta''_0.
$
Hence $\phi_0$ is  holomorphic on $v_0$.
Thus the family
$
 {\bf  f}  :=\{(\phi_0, v_0, \Delta'_0 \times \Delta''_0)\}_{z^0\in M}
$
gives a local coordinate system for $M$. We will
call $(\phi_0, v_0,
\Delta_0' \times \Delta_0'')$ $F$-local coordinates at $z^0$ for $M$ via $g^0$ (i.e.,
via $F$-local coordinates $(\Phi_0,V_0, \Delta_0)$ at $g^0$).

We recall from (\ref{eqn:Sigma-and-sigma}) the subset $\sigma_{z_0}:=\psi_{z_0}(\Sigma_{z_0})=\Sigma_{z_0}(z_0)$ of $M$. Then:
\begin{enumerate} {\it
 \item [(1')]\ $\sigma_{z_0}$ is a parabolic, irreducible, $m_1$-dimensional
       non-singular f-generalized analytic set in $D $.
Moreover, $\Sigma_{z_0}=G \cap \,Aut\,\sigma_{z_0} $;
$\Sigma_{z_0}$ acts transitively on $\sigma_{z_0}$; and
       $\sigma_{z_0} \approx \Sigma_{z_0} /H_{z_0}$.
 \item [(2')]\ $g \sigma_{z_0} =\psi_{z_0}(g \Sigma_{z_0})$ for $g\in
 G$. Moreover, $M=\cup_{g\in G} \,g \sigma_{z_0} $ satisfies the property that if
$g_1 \sigma_{z_0} \cap g_2 \sigma_{z_0} \ne \emptyset$ for $g_1,g_2 \in G $ then $g_1 \sigma_{z_0} =g_2 \sigma_{z_0} $.
 \item [(3')]\ fix $z^0\in M$ and take $g^0\in G$ with
 $\psi_{z_0}(g^0)=z^0$. Let $(\phi_0, v_0, \Delta_0 '\times
 \Delta _0'')$ be $F$-local coordinates at $z^0$ for $M$ via $g^0$. Take $z'\in v_0$ with $\phi_0(z')=(\xi',\eta')\in \Delta_0'
 \times \Delta_0''$ and let $g'\in G$ with $z'\in
 g'\sigma_{z_0}$. Then there exists a connected component
 $\sigma'$ of $[g'\sigma_{z_0}]\cap \boldmath v_0$ such that $\sigma'=\phi_0^{-1}(
\{\xi'\} \times \Delta_0'')$ in $M$.  }
\end{enumerate}

{\it  Proof.} The first assertion in (2') follows from the definition
of $\sigma_{z_0}$.  We next verify (3'). Let $z'\in v_0$ with $\phi_0(z')=(\xi',\eta')$ and let $g'\in
G$ with $z'\in
g'\sigma_{z_0}$. We take
$h'\in V_0$ with $\psi_{z_0}(h')=h'(z_0)=z' $; then $
\Phi_0(h')=(\xi ', \eta', \zeta')$ for some $\zeta'\in \Delta_0'''$. Since  $z'\in
g'\sigma_{z_0}$, we can find $s\in \Sigma_{z_0}$ such that
$z'=g's(z_0)$, and $(h')^{-1}g' s \in H_{z_0}$.
 Since $\Sigma_{z_0}$ is a connected Lie subgroup of $G$ which contains $H_{z_0}$,  it follows that
$$
g'\Sigma_{z_0}=g's\Sigma_{z_0} \subset
h'H_{z_0}\cdot\Sigma_{z_0}= h'\Sigma_{z_0}  \quad  \mbox{so that }
\quad
g'\Sigma_{z_0}=h' \Sigma_{z_0}.
$$
 Since $
\Phi_0(h')=(\xi ', \eta', \zeta')$, we see  from i') that
the connected component of
$[g'\Sigma_{z_0}]\cap V_0=[h'
\Sigma_{z_0}]\cap V_0$ containing $h'$ is equal to $\Phi_0^{-1}(\{\xi '\}\times
\Delta_0'' \times \Delta_0''')$. Together with the observation that $\Phi_0^{-1}(\xi, \eta,
\zeta)(z_0)= \phi^{-1}_0(\xi,\eta)$ for $(\xi,\eta, \zeta)\in \Delta_0$, this  yields (3').

We now prove the first assertion in (1'). Fix $z^0\in \sigma_{z_0}$ and let $g^0\in
G$ with $g^0(z_0)=z^0$. As before, we consider $F$-local coordinates $(\Phi_0, V_0,
\Delta_0)$    at $g^0$  for $G$ and $(\phi_0, v_0, \Delta_0' \times \Delta_0'')$ at
$z^0$ for $M$ via $g^0$. Since $\phi_0(z^0)=(0,0 )\in \Delta_0' \times
\Delta_0''$ and $z^0\in \sigma_{z_0}$, it follows from (3') that
$\phi_0^{-1}(\{0\} \times \Delta_0'')$ coincides with some connected
component $\sigma'$ of $[\sigma_{z_0}]\cap v_0$. Thus, $\sigma'$ is an $m_1$-dimensional non-singular analytic set in
$v_0$. Now $\Sigma_{z_0}$ is connected;
this  implies that
$\sigma_{z_0}$ is an irreducible $m_1$-dimensional
non-singular $f$-generalized analytic set in $M$. Furthermore, using the representation for $\Sigma_{z_0}$ in (\ref{eqn:imphyo}), one can show that
  $\sigma_{z_0}=\Sigma_{z_0}(z_0)$ is parabolic. Now take $g\in D(z_0)$.
 Since $\lambda(z)$ is  continuous on $D$, it
follows from Lemma \ref{lem:fundamen} and (\ref{eqn:imphyo})
that
\begin{eqnarray}
 \label{eqn:kanarii}
\lambda(z)= \ {\rm  const.} \ = \lambda(g(z_0))  \qquad \mbox{ for  } \,
z\in  g\sigma_{z_0}.
\end{eqnarray}
Since $-\lambda(z)$ is an exhaustion function for $D$, we must have $g
\sigma_{z_0} \Subset D$, proving the first assertion in (1').

To prove the second assertion in (2'), assume that $g'\sigma_{z_0}\cap
 g''\sigma _{z_0}\ne \emptyset$ and take
 $z^0\in g'\sigma_{z_0}\cap g''\sigma _{z_0}$. Then there exist
 $s',s''\in \Sigma_{z_0}$ with $z^0=g's'(z_0)=g''s''(z_0)$, so that we
 find $h\in H_{z_0}$  such that $g's'=g''s''h$. Since $\Sigma_{z_0} $ is
 a group containing $H_{z_0}$, it follows that $ g'=g''s''h (s')
 ^{-1}  \in g'' \Sigma_{z_0} $, and hence $g' \Sigma_{z_0}  \subset g''
 \Sigma_{z_0}$. We thus have $g' \sigma_{z_0} \subset g''
 \sigma_{z_0} $, and hence $g' \sigma_{z_0} = g'' \sigma_{z_0} $.

For the second assertion in (1'),
 since $G \subset Aut\,M$,
it suffices to prove that $\{g\in G: g \sigma_{z_0} =\sigma_{z_0}
 \}=\Sigma_{z_0} $ and to verify that $\Sigma_{z_0}$ acts transitively on $\sigma_{z_0}$.
Let $g\in \Sigma_{z_0} $. We have $g
\sigma_{z_0} = \sigma_{z_0} $ by (2') which shows that $\Sigma_{z_0} \subset \{g\in G: g \sigma_{z_0} =\sigma_{z_0}
 \}$. To
prove the reverse inclusion, fix $g\in G$ with $g \sigma_{z_0} =\sigma_{z_0}
 $. There exist $s\in \Sigma_{z_0}$ and $h\in
H_{z_0}$  such that $g= sh$. Since $H_{z_0 }\subset \Sigma_{z_0} $,  it follows that $g\in \Sigma_{z_0} $, which proves $\{g\in G: g \sigma_{z_0} =\sigma_{z_0}
 \} \subset \Sigma_{z_0}$. To verify the
transitivity, let $\zeta\in \sigma_{z_0}$. From (\ref{eqn:triviality}) we can find $g\in D(z_0)$
with $g(z_0)=\zeta$. Using (2') we have $g \sigma_{z_0} = \sigma_{z_0}$ so that $g\in \Sigma_{z_0}$ as required.
\hfill $\Box$

\smallskip
Using
(2'), (3') and (\ref{eqn:triviality})
we have
\begin{equation}
\begin{split}
 \label{eqn:sliced}
 M&=\mbox{ $\bigcup _{g\in G}$} \, g\sigma_{z_0} \quad \mbox{  and } \quad
 D=\mbox{  $\bigcup _{g\in
 D(z_0)}$} \, g\sigma_{z_0}, \end{split}
\end{equation}
which  provide foliations of $M$ and $D$.
 We see from (1')  that each leaf $g \sigma_{z_0}$ in $D $ for $g\in
 D(z_0)$ is relatively
 compact in $D$, and hence each leaf $g \sigma_{z_0}$ in $M$ for $g\in G$ is relatively
 compact in $M$. Thus (i) in
 $1-a.$ is proved. Assertions (1') and (2')
 verify (o) in $1-a.$

\smallskip
To prove (ii) in $1-a.$ and $1-b.$, we begin with the following two
items:
\begin{enumerate}
 \item [$\alpha.$] {\it Let $z'\in D$ and let $g'\in G$ with
       $g'(z_0)=z'$. Let $l$ be a one-dimensional complex analytic curve passing through $z'$ and contained
in a  neighborhood $U$ of $z'$ in $M$ such that $\lambda(z)$ is constant on $l$. Then $l\subset g'\sigma_{z_0}$. }
\item [$\beta.$] {\it $\Sigma_{z_0}$ is closed in $G$ if and only if
$\sigma_{z_0}$  is closed in $M$.}
\end{enumerate}

\smallskip
\noindent {\it Proof.}
We prove  $\alpha.$ by contradiction. It suffices to consider the case when
 $l$ intersects the $m_1$-dimensional generalized analytic
 set $g'\sigma_{z_0}$ transversally at $z'$. Indeed, assume that $\alpha.$ is true for any $g' \sigma_{z_0}$
 which $l$ intersects transversally. Now
 assume that there exists some set $g'\sigma_{z_0}$ such that $l$ does not intersect $g'\sigma_{z_0}$
transversally at $z'$. To verify $\alpha. $ in such a case, proceeding by contradiction, we
assume that $l \not\subset g'\sigma_{z_0}$. Since we have a foliation $D =\cup _{g\in D(z_0)}\,g\sigma_{z_0}$ of $D$ in (\ref{eqn:sliced}), we can find $g'' \in D(z_0)$ sufficiently
close to $g'$  such that $g''\sigma_{z_0}\cap l\ne \emptyset $;
$g''\sigma_{z_0} \cap g'\sigma_{z_0}=\emptyset$; and
 $l$ intersects $g''\sigma_{z_0}$ transversally. But then by assumption we have $l \subset
g''\sigma_{z_0}$, which is a contradiction.

Thus we assume that $l \not\subset g'\sigma_{z_0}$ and also that $l$  intersects $g'\sigma_{z_0}$ transversally at $z'$.
Using standard results for homgeneous spaces, we can choose $X\in \mathfrak X$  such
 that  the curve $C:= \{\exp
tX(z'):\,t\in {\mathbb{C}}\}$ in $M$ has the same complex tangent direction at $z'$ as the line $l$. Call this direction $\boldmath a$. Since $\lambda(z)= const.$ on $l$, we have
 $[\frac{\partial^2
\lambda (z' +at)}{\partial t \partial \overline { t} }] |_{t=0}=0 $.
It follows from the equivalent conditions $o.\sim iii.$ before Lemma
\ref{lem:fundamen} that $\lambda(\exp
tX(z'))= const. $ for $t\in {\mathbb{C}}$.
Since $(g')^{-1}(z')=z_0\in  D$, we have $(g')^{-1}\in D(z')=D'(z')$. Thus, from property c., we see that $\lambda((g')^{-1} \exp tX(z'))= const. $ for $t\in
{\mathbb{C}}$.  However the curve $C': = (g')^{-1} C$
intersects $\sigma_{z_0}$ transversally at $z_0$ since the curve $C$ intersects $g'\sigma_{z_0}$ transversally at $z'$. We can choose  $Y\in  \mathfrak X$
such that the curve $\widehat{ C}:= \{\exp tY(z_0): t\in {\mathbb{C}}\}$ has the same complex tangent direction at $z_0$ as $ C'$. This shows that  $\widehat{ C} \not\subset \sigma_{z_0} $.  We again use the conditions $o.\sim
iii.$ to conclude that $\lambda(z)=const.$ on $\widehat{ C}$, so that $Y\in \mathfrak X
_{z_0}$. Since $\sigma_{z_0}=\Sigma_{z_0}(z_0)$ and $\{\exp tY: \ t\in {\mathbb{C}}\} \subset \Sigma_{z_0}$, we have $\widehat{ C}\subset \sigma_{z_0} $. This is
a contradiction.

To prove the sufficiency in $\beta$., let $z^\nu\in \sigma_{z_0},
\ \nu=1,2, \ldots $,  and let $z^0\in D$ with $z^\nu\to z^0$ as $\nu\to
\infty$. We show that $z^0\in \sigma_{z_0}$. Fix $g^0\in
G$ with $g^0(z_0)=z^0$. Consider $F$-local coordinates
$(\Phi _0,V_0,\Delta_0)$ at $g^0$ for $G$ and $(\phi_0,v_0,\Delta_0'
\times \Delta_0'')$ at $z^0$ for $M$ via $g^0$. Since $z^\nu \in v_0$
for sufficiently large $\nu$, we can find $g^\nu\in V_0$ with
$g^\nu(z_0)=z^\nu$. By taking a subsequence of $\{g^\nu\}_\nu$ if necessary, we may
assume that there exists $g^*\in \overline{  V_0}\subset G$ such that
$g^\nu \to g^*$ as $\nu\to \infty$, so that $g^*(z_0)=z^0$. On the other
hand, since $z^\nu\in \sigma_{z_0}$, we can find $s^\nu\in \Sigma_{z_0} $
such that $z^\nu=s^\nu(z_0), \ \nu=1,2,\ldots $. Hence
$g^\nu(z_0)=s^\nu(z_0)$, so that $(s^\nu)^{-1}g^\nu\in H_{z_0} \subset
\Sigma_{z_0}$. Thus $g^\nu\in s^\nu\Sigma_{z_0}=
\Sigma_{z_0}$.
 Since $g^\nu\to g^*$ as $\nu\to \infty$ and since $\Sigma_{z_0}$ is assumed
 to be closed
 in $G$, it follows that $g^*\in \Sigma_{z_0}$, so that $z^0=g^*(z_0)\in
 \sigma_{z_0}$.
To prove the necessity in $\beta$., let $g^\nu\in \Sigma_{z_0},\
\nu=1,2,\ldots   $ and
$g^0\in G$ with $g^\nu \to g^0$ as $\nu\to \infty$. We show that
$g^0 \in \Sigma_{z_0}$. We set $z^\nu=g^\nu(z_0)\in \sigma_{z_0}$ and
$z^0=g^0(z_0)\in M$. Since $z^\nu\to z^0$ as $\nu\to \infty$ and since
$\sigma_{z_0}$ is assumed to be closed in $M$, it follows that
$z^0\in \sigma_{z_0}$. Thus we can find $s^0\in \Sigma_{z_0}$ with $s^0(z_0)=z^0$. Hence $g^0(z_0)=s^0(z_0)$, so that $g^0\in
s^0H_{z_0} \subset s^0\Sigma_{z_0}=\Sigma_{z_0}$
and $\beta.$ is proved.
\hfill $\Box$

\smallskip Since $-\lambda(z)$ is a plurisubharmonic exhaustion function
for $D$, assertion $\alpha.$ implies (ii) of $1-a.$, completing the proof
of $1-a$.

\smallskip
To prove $1-b$. assume that $\sigma_{z_0}$ is closed in $M$. From
$1-a$.,  $\sigma_{z_0}$ is an irreducible,
 $m_1$-dimensional compact non-singular analytic set in $D$. From
 $\beta.$, $\Sigma_{z_0}$ is closed, which
 proves (o) in $1-b.$

To prove (i) in $1-b.$ we note  that  $\Sigma_{z_0}$ is closed in $G$, so that
 $\Sigma_{z_0}$ is an $(m_0+m_1)$-dimensional connected  closed
Lie subgroup of $G$ which corresponds to the Lie subalgebra $\mathfrak
 X_{z_0} $ of $\mathfrak X $. We consider the quotient space $M_0:= G/
 \Sigma_{z_0}$; this is an $m_2$-dimensional complex homogeneous space
 with Lie transformation group $G$ which acts transitively on $M_0$.
Since $H_{z_0}\subset \Sigma_{z_0}$ as Lie groups, we have the
canonical holomorphic projection
$$
\pi_0: z=gH_{z_0} \in M\mapsto w=\pi_0(z)= g\Sigma_{z_0}\in M_0.
$$
Just as we identified $H_{z_0} \subset G$ with the point $z_0$ in $M$, we may consider $w_0:=\Sigma_{z_0}$ as a point in $M_0$, i.e.,
$\Sigma_{z_0}$ is the isotropy subgroup of $G$ for the point $w_0$ in
$M_0$. Then
\begin{equation}  \label{eqn:notessentialbut2}
 \begin{array}{llll}
\pi_0(g(z_0))&=g(w_0),  \quad  \ &\mbox{  for all } \ g\in G; \vspace{ 1mm}\\
\pi_0^{-1}(g(w_0))&=g\sigma_{z_0},   \quad  &\mbox{  for all } \ g\in G.
 \end{array}
\end{equation}

Indeed, the first formula comes from the definition of $\pi_0$.
For the second one, we first note that
$$
\pi_0(g \sigma_{z_0} )= \pi_0(g\Sigma_{z_0}(z_0))= g\Sigma_{z_0}
(w_0)=g(w_0),
$$
which proves the inclusion $g\sigma_{z_0}\subset \pi_0^{-1}(g(w_0))$. To prove the inclusion $$\pi_0^{-1}(g(w_0))\subset g\sigma_{z_0},$$   take
$z\in M$ with $\pi_{0}(z)=g(w_0)$. We can find $g'\in G$  with
$g'(z_0)=z$. Thus, $g'(w_0)=\pi_0(g'(z_0))=g(w_0)$ so that
$g'\in g \Sigma_{z_0}$. It follows that $g'(z_0)\in g \sigma_{z_0} $,
and $z\in g\sigma_{z_0}$. Thus, (\ref{eqn:notessentialbut2}) is proved.

Set $D_0:=\pi_0(D) \subset M_0$. We have $D_0
\Subset M_0$ since $D \Subset M$; moreover,
 let $w\in D_0 \ (\partial D_0)$ and take $g\in G$  such that
$w=g(w_0)$. Then
\begin{align}
 \label{eqn:d0d}
\pi_0  ^{-1}  (w)= g \sigma_{z_0} \Subset D \   (\partial D).
\end{align}
Indeed, since the proof is similar, we only give the proof in the case $w\in D_0$. Since $g(z_0)\in D$, we find
$g'\in D(z_0)$  such that $g(z_0)=g'(z_0)$, which implies $g=g'h \in
D(z_0)$. This together with (\ref{eqn:notessentialbut2}) and (i) in
1-a.  implies
(\ref{eqn:d0d}).

Therefore, we have
\begin{align}
 \label{eqn:inversed0}
D=\pi_0^{-1}(D_0) \qquad \mbox{  and } \qquad \partial D= \pi_0  ^{-1}
(\partial D_0).
\end{align}

To prove that $D_0$ is a Stein domain in $M_0$, we
define a real-valued function $\lambda_0(w)$ on $D_0$ as follows.
Fix $w\in D_0$ and take $z\in D$ with $\pi_0(z)=w$.
We define
$\lambda_0(w) := \lambda ( z)$, which is well-defined by (\ref{eqn:d0d}).
 In order to prove that $\lambda_0(w)$ is a strictly plurisubharmonic exhaustion
function on $D_0$ we shall utilize the usual coordinates for $M$
and $M_0$. Fixing $g^0\in G$, by the Frobenius theorem we have  $F$-local coordinates $(\Phi_0,
 V_0, \Delta_0)$ at $g^0$ for $G$ with the following properties:
$$
\Phi_0: g\in V_0 \mapsto (\xi(g),\eta(g), \zeta(g))\in \Delta_0:=\Delta_0'\times
 \Delta_0'' \times \Delta_0^{'''};
$$
for a given $ p=(\xi',\eta',\zeta')\in \Delta_0$, if we set
$g':=\Phi_0^{-1}(p)\in V_0$, then
\begin{equation}
 \begin{split}
  [g'H_{z_0}] \cap V_0&= \Phi_0^{-1}( \{(\xi', \eta')\}\times
  \Delta_0^{'''});\\
[g'\Sigma_{z_0}] \cap V_0&= \Phi_{0}^{-1}(\{ \xi' \}\times
  \Delta_0''\times \Delta_0^{'''});
\end{split} \end{equation}
and
\begin{equation}
 \begin{split}
&([g'H_{z_0}]\cap V_0)\cap
([g''H_{z_0}]\cap V_0)=\emptyset \Longrightarrow
[g'H_{z_0}]\cap [g''H_{z_0}]=\emptyset;\\
&([g'\Sigma_{z_0}]\cap V_0)\cap
([g''\Sigma_{z_0}]\cap V_0)=\emptyset \Longrightarrow
[g'\Sigma_{z_0}]\cap [g''\Sigma_{z_0}]=\emptyset.
\end{split} \end{equation}
Note that $\Sigma_{z_0} $ and $H_{z_0}$ are closed in $G$.  We set $U_1:=\{
gH_{z_0}\in G: g\in V_0\} \subset G$, and we regard $U_1$ as a neighborhood of
$z^0=g^0H_{z_0}$ in $M$. Similarly, we regard $U_2:=\{
g\Sigma_{z_0}\in G: g\in V_0\} \subset G$ as a neighborhood of
$w^0=g^0\Sigma_{z_0}$ in $M_0$. For $g\in V_0$, we consider the correspondences
\begin{eqnarray} \label{eqn:locoo}
\left.
\begin{array}{llll}
  \varphi _1=\phi_0&: z=gH_{z_0}\in U_1  \mapsto (\xi, \eta)\in
\Delta_0' \times \Delta_0''; \\
\varphi _2&: w=g\Sigma_{z_0}\in U_2 \mapsto \xi \in
 \Delta_0'.
\end{array}
\right.
\end{eqnarray}
Then $(\varphi _1, U_1, \Delta'_{0} \times \Delta_{0}'')$ and
 $(\varphi _2, U_2, \Delta_{0}')$
define  local coordinates at $z^0$ for $M$ and
local coordinates at $w^0$ for $M_0$.
 In particular,
 \begin{align}
  \label{eqn:varphi12}
\varphi_1( [g^0 \sigma_{z_0}] \cap U_1)=\{0\}\times \Delta_0'';  \qquad
\varphi _2([g^0 \Sigma_{z_0}] \cap U_2)= 0.
\end{align}
In terms of  the local coordinates (\ref{eqn:locoo}), if we
restrict $\pi_0$ from the neighborhood $U_1$ of
$z^0$ in $M$ to the  neighborhood $U_2$ of $w^0$ in $M_0$, we have
$$
\pi_0: (\xi, \eta )\in \Delta_0' \times \Delta_0'' \mapsto  \xi\in \Delta_0'.
$$

From the local  coordinate expression for $\pi_0$ and
 (\ref{eqn:inversed0}), we have that $\partial D_0 $ is smooth in $M_0$. Moreover, using the local coordinates (\ref{eqn:locoo}),  $\lambda_0(w)$
is of the form $\lambda_0(\xi)=\lambda(\xi,\eta)$ for $(\xi,\eta)\in \Delta_0'
\times \Delta_0''$. In particular, $\lambda_0(\xi)=\lambda(\xi,0)$ in
$\Delta_0'$. It follows from the plurisubharmonicity of $-\lambda (z)$ on
$D$ that $-\lambda_0(w)$ is plurisubharmonic on $D_0$. Since $\partial D_0=
\pi_0(\partial D)$,  the fact that $-\lambda(z)$ is an exhaustion function for $D$ implies that $-\lambda_0(w)$ is an exhaustion function for  $D_0$.

It remains to verify that
$\lambda_0(w)$ is strictly plurisubharmonic on $D_0$. If not, there
exists $w^0\in D_0$  and $a _0 \in {\mathbb{C}}^{m_2}\ (a_0\ne 0)$ such that
$
[\frac{\partial^2 \lambda_0(w^0+ a_0t)}{\partial t \partial \overline { t
}}]|_{t=0}=0 $. Take $z^0\in D$ and $g^0\in D(z_0)$ with $g^0(z_0)=z^0$ and
$\pi_0(z^0)=w^0$. Using local coordinates
$(\varphi _1, U_1, \Delta'_{0} \times \Delta_{0}'')$ at $z^0$ for $M$ and
$(\varphi _2, U_2, \Delta_{0}')$ at $w^0$ for $M_0$, we may assume that
 $w^0+a_0t, \ |t|\ll 1$ in $U_2$ is of the form $\xi=a_0t,\ |t|  \ll 1$ in
 $\Delta_0'$.  Since $\lambda_0(a_0t)= \lambda(a_0t, 0), \ |t|\ll  1$, we
 have  $\bigl[\frac{\partial
^2\lambda(a_0t,0)}{\partial t \partial \overline { t}}\bigr]|_{t=0}=0 $.
Take $X \in \mathfrak X $ with the property that the integral
curve
$C^0:= \{\exp t X(z^0): \ t \in {\mathbb{C}}\}$ in $M $ has direction $(a_0;0)\in {\mathbb{C}}^{m_1}
\times {\mathbb{C}}^{m_2}$ at
$z^0$.  From the equivalent conditions $o.\sim iii.$ we see
that $\lambda (z)=const.$ on $C^0$.  Since $z^0\in C^0\cap [g^0 \sigma_{z_0}] $, it follows from  $\alpha.$ that  $C^0 \subset
g^0\sigma_{z_0}$. On the other hand,
by (\ref{eqn:varphi12}), $[g^0 \sigma _{z_0}]\cap U_1$ corresponds to
$\{0\} \times \Delta''_0$ using the local coordinates $\varphi _1$. Since the direction of $C^0$ at $z^0$ is $(a_0;0)$
with $a_0\ne 0$, this is a contradiction. Condition $1-b$. is
proved.
\hfill $\Box$

\smallskip
  We would like to prove a result, Theorem
\ref{thm:main-nonconnected}, corresponding to
Theorem
\ref{thm:mainth} in the case when the isotropy
subgroup $H_z$ is not connected for use in studying special Hopf spaces. We begin with a discussion below which
makes clear the difference between the connected case and the non-connnected
case. The bulk of the proof of Theorem
\ref{thm:main-nonconnected} can be found in Appendix B; we commence with some background and motivation. The next two results, Propositions \ref{prop:disc} and \ref{prop:normalcover},
are standard. To keep the article self-contained, and since we will use some notation from the proofs
later on, we include the brief proofs. One preliminary notational remark: if $A,B$ are subsets (not necessarily subgroups) of a group $G$, we write $AB:=\{ab: a\in A, \ b\in B\}$.
 \begin{proposition}\label{prop:disc} \
Let $z_0\in M$ and suppose that $H_{z_0}$ is not connected. Consider the decomposition of
 $H_{z_0}$ into connected components $H_{z_0}^{(k)}, \ k=1,2,\ldots, $
 where $H_{z_0}^{(1)}=H_{z_0}'$ (see (\ref{eqn:visible})). Then there exists  a  neighborhood  $\boldmath v_0$ of $e$ in $G$  such that
\begin{enumerate}
 \item \ $\boldmath v_0 H_{z_0}^{(k)}\cap \boldmath v_0 H_{z_0}^{(l)} = \emptyset,
 \quad  k\ne l$;
 \item \ if $a,\, b\in \boldmath v _0$ and $a^{-1}b\in
  H_{z_0}$, then
  $a^{-1} b \in H_{z_0}'$.
\end{enumerate}
\end{proposition}
\noindent {\bf Proof}.
Since $H_{z_0}$ is closed in $G$, 2. automatically holds provided $\boldmath v_0$ is a sufficiently small neighborhood of $e$ in $G$.
To verify 1. we first prove that
 there exists a  neighborhood  $\boldmath v_0$ with property 2. such that
\begin{eqnarray}
 \label{eqn:locally1}
 {\boldmath v_0} H_{z_0}' \cap  (\mbox{  $\bigcup_{k=2}^\infty $} \, \boldmath  v_0
 H_{z_0}^{(k)}) =  \emptyset.
\end{eqnarray}
If no such neighborhood exists, then for any neighborhood $v_0$ of $e$, for $j=1,2,\ldots  $,\,  we can find
 $g_j\in \boldmath v_0, h_j\in H_{z_0}'$ and $g_{jk_j} \in
 \boldmath v_0, h_{jk_j}\in H_{z_0}^{(k_j)}$,    $k_j \ge 2$ such that
 $g_j,  g_{jk_j}
 \to e$ as $j\to \infty$ and $
g_jh_j = g_{jk_j}h_{jk_j}, \ j=1,2 \ldots
$.
Hence $ g_{jk_j}^{-1}g_j = h_{jk_j}h_j^{-1}$. The left-hand side $ g_{jk_j}^{-1}g_j$ tends
 to $e$ as $j\to \infty$ and the right-hand side $h_{jk_j}h_j^{-1}$ always belongs to $H_{z_0} $.
 Hence $h_{jk_j}h_j^{-1}\in H_{z_0}'$ for sufficiently large $j$, so
 that $h_{jk_j}\in H_{z_0}'$.  This contradicts $h_{jk_j}\in H_{z_0}^{(k_j)}$ for all $k_j \ge 2$, verifying (\ref{eqn:locally1}). We now prove that this neighborhood $v_0$
 satisfies 1. If not,
 $\boldmath v_0 h_k H_{z_0}' \cap
\boldmath v_0 h_l H_{z_0}' \ne \emptyset$ for some $k\ne l$. Then
 $\boldmath v_0( h_k H_{z_0}'h_k^{-1}) \cap
\boldmath v_0(h_l h_k^{-1})(h_kH_{z_0}'h_k ^{-1}) \ne \emptyset $. Since
 $H_{z_0}' $ is  a normal subgroup of $H_{z_0} $, we have
$\boldmath v_0 H_{z_0}' \cap
\boldmath v_0(h_l h_k^{-1})H_{z_0}' \ne \emptyset $. Since $k\ne l$ we
 have $(h_l h_k^{-1})H_{z_0}'= H_{z_0} ^{(m)}$ for some $m\ge 2$ and
$H_{z_0}^ {(m)} \cap H_{z_0}' =\emptyset$ by (\ref{eqn:visible}). Hence
$ \boldmath v_0H_{z_0}'\cap \boldmath {  v}_0 H_{z_0}^{(m)}\ne
 \emptyset$, which   contradicts (\ref{eqn:locally1}). Condition 1. is
 proved.  \hfill $\Box$

\smallskip
We consider the quotient space $M'=G/H_{z_0}'$ which is a homogeneous
space with Lie transformation group $G$ which acts transitively on
$M'$. We note that $M'$ depends on the point $z_0$ in $M$.
 The space  $M'$ has the same dimension, $m-n$, as $M$.
We write $w_0$ for the point in $M'$ which corresponds to
$H_{z_0}'$ and we introduce the notation $\widehat{ H}_{w_0}:=H_{z_0}' \subset G$.
Then we have the projections
$$
\widehat{ \psi}_{w_0}: \ g\in  G \mapsto
 g(w_0)\in M',  \qquad  \widehat{ \pi}_{w_0}:g\in G \mapsto g\widehat{ H}_{w_0} \in G/ \widehat{ H}_{w_0}.
$$
Recall we defined $\mathfrak h _{z_0}$ as the Lie subalgebra of $\mathfrak X $
which corresponds to $H_{z_0}'$; we now write
$ \widehat{ \mathfrak h} _{w_0}:= \mathfrak h_{z_0} \subset \mathfrak X $.
Thus,  $ \widehat{
 H}_{w_0}$ is the connected isotropy
 subgroup of $G$ for $w_0\in M'$ whose Lie subalgebra is  $\widehat{ \mathfrak h
 }_{w_0}$, and a point $w$ in
 $M'$ is identified with $g \widehat{ H}_{w_0}$ where $g(w_0)= w$.

We consider the mapping
$$
 \widehat{ \pi}:\  gH_{z_0}' \in M' \mapsto gH_{z_0}\in  M.
$$
We remark that $\widehat{
\pi} $  depends on $z_0$ in $M$.
This mapping is well-defined,
holomorphic and surjective.
 Note that for any $w\in M'$ and $g\in G$
\begin{eqnarray} \label{eqn:commutative-wg}
\widehat{ \pi} (g w)= g \widehat{\pi}( w), \ \   \qquad  \widehat{
\pi}^{-1}(gH_{z_0})= \mbox{  $\bigcup _{k=1}^\infty$ } \, g h_k \widehat{ H}_{w_0}
\end{eqnarray}
with $h_k\in H_{z_0}^{(k)}$ as in (\ref{eqn:visible})
since $H_{z_0}^{(k)}= h_k \widehat{ H}_{w_0}, k=1, 2,\ldots $ by
(\ref{eqn:visible}).

To verify (\ref{eqn:commutative-wg}), since $w= g' H_{z_0}' $ for some $g'\in G$, we have, by definition,
$ \widehat{ \pi}(gw)= \widehat{\pi}( gg' H_{z_0}' )= gg' H_{z_0} = g
\widehat{\pi}(w)$. This proves the first equality.
By (\ref{eqn:visible}) we  have
$\widehat{\pi}(gH_{z_0}^{(k)})=gH_{z_0}$, so that $\cup _{k=1}^\infty gH_{z_0}^{(k)} \subset \widehat{
\pi}^{-1}(gH_{z_0})$. Conversely, if we take $h\in G$
such that $\widehat{ \pi}(hH_{z_0}')=gH_{z_0}$, then we have
$g^{-1}h\in H_{z_0}= \cup_{k=1}^\infty H_{z_0}^{(k)}$, and hence $h\in
gH_{z_0}^{(k)}$ for some $k$. Thus, $hH_{z_0}'\subset
gH_{z_0}^{(k)} H_{z_0}'=gH_{z_0}^{(k)}$. Since both sets $hH_{z_0}'$ and $gH_{z_0}^{(k)}$ are irreducible,
analytic sets in $G$ of the same dimension, we have
$hH_{z_0}'=gH_{z_0}^{(k)}$. Thus the reverse inclusion holds, and the second equality is proved.
\begin{proposition}
 \label{prop:normalcover} \
The pair $(\widehat{ \pi}, M')$ gives a normal, unramified cover of
 $M$; i.e., if $z\in M$, $w\in \widehat{
 \pi}^{-1}(z)$, and $\gamma$ is a real, one-dimensional curve in $M$
 starting at $z$, then there exists a unique curve $\gamma'$ in $M'$
 which starts at $w$ and satisfies
 $\widehat{ \pi}(\gamma')= \gamma$.
\end{proposition}

\noindent{\bf Proof}. We use the notation in Proposition
\ref{prop:disc}; in particular, $v_0$ is a neighborhood of $e$ in $G$
satisfying the conditions 1 and 2. Let $z\in M$; thus $z=gH_{z_0}$ for
some $g\in G$, or equivalently $z=g(z_0)$. Then we have by
(\ref{eqn:commutative-wg})
\begin{align}
 \label{eqn:defofwk}\widehat{
 \pi}^{-1} (z)= \{w_k\}_{k=1,2,\ldots }  \quad  \mbox{   where } \
 w_k=gH_{z_0}^{(k)}=gh_kH_{z_0}'\in M'.
\end{align}
Note that $\{w_k\}_{k=1,2,\ldots }$ are isolated in $M'$. We define
\begin{align}
 \label{eqn:sectionw}
V:=g \boldmath v_0 H_{z_0} \ \    \mbox{  and }\ \ U_k:=g\boldmath v_0
 H_{z_0}^{(k)}=g \boldmath v_0 h_k \widehat{ H}_{w_0}=g \boldmath v_0
 \widehat{ H}_{w_0} h_k,
\end{align}
so that $V$ is a  neighborhood  of
 ${ z}$ in $M$ and $U_k$
  is a  neighborhood  of ${ w}_k$ in $M'$. It suffices to prove that
\begin{equation}
\begin{split}
 \label{eqn:bijective}
&\widehat{ \pi}^{-1}(V)= \mbox{  $\bigcup _{k=1}^\infty $} \, U_k\ \, \mbox{a disjoint
 union;}\\
&\widehat{ \pi}: \ U_k\mapsto V   \mbox{  is
 bijective for each $k=1,2,\ldots$.}
\end{split}
\end{equation}
The first equality, and the disjointness, follow from 1. in Proposition \ref{prop:disc}.
The surjectivity of $\widehat{ \pi}$ in the second assertion is clear by its definition. To verify  injectivity, it suffices to show that
$ga H_{z_0}=gb H_{z_0} \Longrightarrow gaH_{z_0}^{(k)}=gbH_{z_0}^{(k)}$
 for $a,\, b\in \boldmath v_0$, or, equivalently, since $H_{z_0}^{(k)}=h_k H_{z_0} '$,
$a^{-1}b \in H_{z_0} \Longrightarrow h_k^{-1}(a^{-1}b)h_k\in H_{z_0}'.$
Since $H_{z_0}'$ is a normal subgroup of $H_{z_0}$, this last implication follows from 2. in
 Proposition \ref{prop:disc}. \hfill $\Box$

\smallskip
We set up some more notation which will be useful in discussing the
case when $H_{z_0}$ is not connected. For our relatively compact domain $D$ in $M$, we take $z_0\in D$ and set
$\widehat{ D}:= \widehat{
  \pi}^{-1}(D)$ in $M'$.  Note that
$ \widehat{ D} $ need not be relatively compact in $M'$; moreover, $
\widehat{ D} $ may be disconnected. However, from  Proposition
\ref{prop:normalcover},  $\widehat{D}$ has no relative boundary over $D$;
 and from (\ref{eqn:bijective}), for $w \in
\partial \widehat{ D}$ we have $\widehat{ \pi}(w)\in \partial D$. We decompose $\widehat{ D}$ into its connected components $\widehat{
D}_j$ in $M'$: $\widehat{ D}= \mbox{  $\bigcup_{j=1}^\infty $} \,
\widehat{ D}_j$.
  We use the convention that $\widehat{ D}_1$ is the connected component containing $w_0$. Recalling the notation $D'(z_0), \ D^{(i)}(z_0)$ in
(\ref{eqn:connectedcom}) for the connected components of $D(z_0)$ in
$G$,
we will show in Appendix B that we have, after possibly relabeling indices,
\begin{align}
\nonumber
 \widehat{ D}_k=
\{g(w_0)\in M': g\in D^{(k)}(z_0)\},  \quad  k=1,2,\ldots.
\end{align}

\smallskip We introduce the notation
\begin{eqnarray}
 \label{defofwhl}
 \widehat{ \lambda }(w): = \lambda (\widehat{\pi}(w)),  \quad  w\in
\widehat{ D}_1.
\end{eqnarray}
 Note that this notation implies that $\widehat{ \lambda}(w)$
 is related to $\widehat{ \pi}:M'=G/\widehat{ H}_{w_0} \mapsto
 M=G/H_{z_0}$,  but we caution the reader that $\widehat{
 \lambda}(w) $ is not necessarily the $c$-Robin function for
 $\widehat{ D}_1$ in $M'$. However, since $\widehat{ D}_1$ is an unramified cover over $D$ without relative boundary, it follows that $-\widehat{ \lambda}(w)$ is, indeed, a plurisubharmonic
 function on the domain $\widehat{ D}_1 \,\subset M'$  such that
{for $w'\in \partial \widehat{
 D}_1$ we have $\lim_{ w \to w'} \widehat{ \lambda}(w)=-\infty$}. Moreover, if $K \Subset D$ and we set $\widehat{ K}_1:=\widehat{
 \pi}^{-1}(K)\cap \widehat{ D}_1$, then   $\widehat{ \lambda}(w)$ is bounded in
 $\widehat{ K}_1$.

We define
\begin{align}  \label{eqn:z-w-relation}
\widehat{ \mathfrak X }_{w_0}&:= \{ X\in \mathfrak X :
[\frac{\partial^2  \widehat{ \lambda} (\exp tX(w_0))}{\partial t \partial
 \overline {
    t }}] |_{t=0}=0\}.
\end{align}
With the notation above, we will show in Appendix B that
$$ \widehat{ \mathfrak X }_{w_0} = \mathfrak X _{z_0}=\{ X\in \mathfrak X : \exp tX(w_0) \in
 \widehat{ D}_1, \  \widehat{\lambda}(
\exp tX(w_0))= \widehat{ \lambda }(w_0), \  t\in {\mathbb{C}} \}$$
and for $X \in \widehat{ \mathfrak X }_{w_0}$ and $g\in  \widehat{ D}_1(w_0)$,
 $$g\exp tX(w_0)\in  \widehat{ D}_1,  \quad   \widehat{ \lambda}( g\exp tX(w_0))= \widehat{
 \lambda}(g(w_0))$$
 for all $t\in {\mathbb{C}}$. Associated to $ \widehat{ \mathfrak X }_{w_0}$ is the Lie subgroup
\begin{align}\label{eqn:whSigam}
\widehat{\Sigma}_{w_0}= \{\,{\textstyle  \prod}_{i=1}^\nu \exp t_iX_i\in G:\nu\in {\boldsymbol{Z}}^+,
 t_i\in {\mathbb{C}}, X_i\in \widehat{ \mathfrak X }_{w_0}\,\}.
\end{align}

 We now state the main theorem in case $H_z$ is not connected, and then
 we make a remark on the roles of the Lie subalgebra $ \widehat{
 \mathfrak X }_{w_0}$ of ${ \mathfrak X }$ and Lie subgroup
 $\widehat{\Sigma}_{w_0}$. In assertion (o) in $2-a.$ in the theorem
 below,
  ${\cal
 H}'(z_0)$ is the Lie subgroup of  $H_{z_0}$ studied in
Corollary \ref{eqn:happy}\,;
$\Sigma_{z_0} $ is the connected Lie subgroup of $G$ corresponding to
the Lie subalgebra $\mathfrak X_0$ defined in (\ref{eqn:X-0}), and
 \begin{align}
  \label{eqn:ahprime}
\Sigma_{z_0} {\cal H}'(z_0) =\{sh\in G: s\in \Sigma_{z_0} , \ h\in {\cal
  H}'(z_0) \}.
\end{align}
As mentioned, the proof of the theorem will be deferred to Appendix B.
\begin{theorem}\label{thm:main-nonconnected}\ Under the
 same notation as in Theorem \ref{thm:mainth}; i.e., $M$ is a homogeneous space of dimension $n$ and the Lie
 transformation group $G$ of dimension $m$ acts transitively on $M$,   assume that the isotropy subgroup $H_z$ of $G$ for $z\in M$ is not connected. Let $D $ be a pseudoconvex domain in $M$ with smooth boundary. Assume that  $D $ is not
 Stein. Fix $z_0\in D$. Then:
\begin{enumerate}
 \item [2-a.]  The subset $\sigma_{z_0}=\Sigma_{z_0} (z_0)$ of $M$ defined in
       (\ref{eqn:Sigma-and-sigma}) is
 a parabolic $m_1$-dimensional non-singular $f$-generalized analytic set in
       $M$ passing through $z_0$
 with the following properties:
\begin{enumerate}
\item [(o)]   $\sigma_{z_0} \Subset D$; $\Sigma_{z_0}{\cal
      H}'(z_0)= D'(z_0) \cap\,Aut\,\sigma_{z_0}$;  $\Sigma_{z_0}{\cal
      H}'(z_0)$ is a Lie
      subgroup of $G$ which acts transitively on $\sigma_{z_0}$;  and
      $$\sigma_{z_0}\approx \Sigma_{z_0}{\cal H}'(z_0)  / {\cal
      H}'(z_0).$$
   \item [(i)] There exists a domain $K$ with  $D \Subset K
 \subset M      $  such that $D (\,K\,)$ is foliated by the sets $g\sigma_{z_0}, \
 g\in K'(z_0)$ ($g\in D'(z_0)$) and each such set $g\sigma_{z_0}$ is
 relatively compact  in
$K (\,D\,)$:
$$
K= \mbox{  $\bigcup_{g\in K'(z_0)}$} \,
g\sigma_{z_0} \qquad and \qquad  D\ =\mbox{  $\bigcup_{g\in D'(z_0)}$}
 \,
 g\sigma_{z_0}.
$$
 \item [(ii)] {The same result as (ii) in 1-a. of Theorem
 \ref{thm:mainth} holds. }
\end{enumerate}
 \item [2-b.]
Assume that $\sigma_{z_0}$ is closed in
       $M$. Then
\begin{enumerate}
 \item [(o)]  The same result as (o) in 1-a. of Theorem
\ref{thm:mainth} holds.
   \item [(i)]
There exists  a complex manifold
 ${ K}_0$    and a holomorphic  map  $
 \pi_0:\ K
 \mapsto K_0$ such that $ { \pi}_0^{-1}(\zeta) \approx \sigma_{z_0}$ (as
       complex manifolds) for each $\zeta  \in K_0$;  moreover, there exists a Stein domain $D_0 \Subset K_0$ with
       smooth boundary  such that $D= { \pi}_0  ^{-1}(D_0)$.
\end{enumerate}
\end{enumerate}
 \end{theorem}

\begin{remark}\label{rem:lastremark} \
{\rm We first remark that, as in the connected case, we have
$$\mathfrak h _{z_0}\subsetneq  \mathfrak X _{z_0} \subsetneq \mathfrak   X$$
where $\mathfrak h _{z_0}$ is the Lie subalgebra of $\mathfrak   X$ corresponding
to $H'_{z_0}$ and $\mathfrak X _{z_0} = \widehat {\mathfrak X} _{w_0}$ from (\ref{eqn:z-w-relation}). Also, in the proof of Theorem \ref{thm:main-nonconnected}
 the connected Lie subalgebra $\mathfrak X _{z_0} =\widehat{ \mathfrak X }_{w_0}$ of $\mathfrak X $
defined by (\ref{eqn:z-w-relation}) and the connected Lie subgroup $\Sigma_{z_0}=
 \widehat{ \Sigma}_{w_0}$ of $G$ which corresponds to $\mathfrak X
 _{z_0}$ will play essential roles. In the special cases when $H_{z_0} \subset \Sigma _{z_0}$ and $\Sigma_{z_0}$ is closed in $G$, we will see from the
 proof that $2-b.$ in Theorem \ref{thm:main-nonconnected} implies
 that  $K=M; \ K_0=G/ \Sigma_{z_0}$ (hence $K_0$ is a homogeneous
 space); \ $\pi_0: \ gH_0\in M \mapsto g \Sigma_{z_0}\in K_0$; and $D_0$ is
 a Stein domain in $K_0$ with smooth boundary. That is, in this case,  the result is the
 same as in $1-b.$ of Theorem \ref{thm:mainth}}.
 \end{remark}

 This special case described in Remark \ref{rem:lastremark} occurs when $M$ is a special Hopf manifold; we now turn our attention to this case. We let  $
({\mathbb{C}}^n)^*:={\mathbb{C}}^n \setminus \{0\}
$ and we fix $\alpha \in {\mathbb{C}}$ with $|\alpha|\ne 1$.  For $z,\ z'\in ({\mathbb{C}}^n)^*$,
we define the following equivalence relation in $({\mathbb{C}}^n)^*$:
$$
 z\sim w  \quad   {\rm  iff}  \quad   w= \alpha^k\,z \      \mbox{  for
 some integer} \ k.
$$
We consider the following space:
\begin{eqnarray}
 \label{eqn:special-hoph}
\mathbb{H}_n=({\mathbb{C}}^n)^* / \sim,
\end{eqnarray}
this is  an $n$-dimensional compact manifold.  The space $\mathbb{H}_n$
 is called a special Hopf manifold.  We write $[z]\in \mathbb{H}_n$ to
 denote the equivalence class of a point $z=(z_1,\ldots , z_n)\in
 ({\mathbb{C}}^n)^*$.
In case $n=1$,  $\mathbb H_1$ is the usual one-dimensional complex
 torus $T_{\alpha}$.
 The space
$\mathbb{H}_n$ clearly has the following property.
\begin{proposition}\label{prop:glnc}  ${}$
\begin{enumerate}
 \item  The group $GL(n,{\mathbb{C}})$ is a Lie transformation group of
 \
 $\mathbb{H}_n$ and acts transitively on $\mathbb{H}_n$; i.e.,
 $\mathbb{H}_n$, equipped with the Lie group $GL(n,{\mathbb{C}})$,  is a
 homogeneous space.
 \item We have the canonical projection $\pi_0: [z_1, \ldots , z_n]\in \mathbb H_n \mapsto  [z_1:\ldots :z_n]
 \in \mathbb P^{n-1}$  such that $\pi_0  ^{-1}  (\zeta) \approx T_{\alpha}$ for
 each $\zeta \in \mathbb P^{n-1}$.
\end{enumerate}
\end{proposition}

We write $
O := [(1,0,\ldots ,0)] \in \mathbb{H}_n
$
and call $O$ the base point of $\mathbb{H}_n$. Then the isotropy subgroup $H_0$ of
$GL(n,{\mathbb{C}})$ for $O$ is
\begin{align*}
 H_0=&\bigl\{\
\left(
\begin{array}{c|cccc}
\alpha^k& &   (*) &  \\
\hline 0 &   \\
\vdots &   & \large A\\
0 &  &
\end{array}
\right) \in GL(n,{\mathbb{C}})\ :
  k\in {\boldsymbol{Z}}, \ (*)\in {\mathbb{C}}^{n-1},  \det \, A\ne 0
\bigr\},
\end{align*}
which is closed but not connnected in $GL(n,{\mathbb{C}})$.
We let $H'_0$ denote the connected component of $H_0$ which contains the identity
$I_n$ in $GL(n,{\mathbb{C}})$, i.e.,
\begin{align*}
 H'_0=&\bigl\{ \
\left(
\begin{array}{c|cccc}
1& &   (*) &  \\
\hline 0 &   \\
\vdots &   & \large A\\
0 &  &
\end{array}
\right) \in GL(n,{\mathbb{C}})\ :
  (*)\in {\mathbb{C}}^{n-1}, \ \det \, A\ne 0
\bigr\}.
\end{align*}
Therefore,
$$
\mathbb{H}_n \approx  GL(n,{\mathbb{C}})/ H_0= \{g H_0: g\in GL(n,{\mathbb{C}})\}.
$$
To be precise, let
\begin{align}
 \label{eqn:matrix-g}
g=\left(
\begin{array}{cccc}
 g_{11}& \ldots & g_{1n} \\
\vdots &\ddots& \vdots  \\
g_{n1}& \ldots &g_{nn}
\end{array}
\right)\in GL(n,{\mathbb{C}}).
\end{align}
Let
$ {\bf  g}:=(g_{11}, \ldots ,g_{n1})\in ({\mathbb{C}}^n)^*$ denote
the first column vector of $g$. Then the mapping
$$
\alpha_0: \ gH_0\in GL(n,{\mathbb{C}})/H_0 \to g(O)= [{\bf  g}]\in \mathbb{H}_n
$$
is bijective.

For local coordinates in a  neighborhood  $V$ of the base point $O$ we can
take
$$
\left(
\begin{array}{clccll}
1+t_1& &  &  \\
t_2 & 1&  &  \\
 \vdots    &  &\ddots  \\
t_n &      &    &1
\end{array}
\right),  \quad |t_i|<\rho, \ i=1,\ldots ,n,
$$
where the missing entries are all $0$. Equivalently, let $U: =\{{\bf  t}=(t_1,\ldots , t_n)\in
{\mathbb{C}}^n: |t_i|< \rho\}  $ where the base point $O$ of
$\mathbb{H}_n$ corresponds
to the origin $0$  of $ {\mathbb{C}}^n$. That is, let
$g\in GL(n,{\mathbb{C}})$ in
(\ref{eqn:matrix-g}) be
close to the identity $I_n$. Corresponding to
$gH_0\in \mathbb{H}_n$ the point
${\bf  t}(g):=(g_{11}-1, g_{21}, \ldots , g_{n1})\in {\mathbb{C}}^n,
$
is close to $(0,0,\ldots,0 )$. We have (i)\, if $g_1H_0\ne g_2H_0$ for
$g_1,g_2 \in V$, then
${\bf  t}(g_1)\ne {\bf  t}(g_2)$; (ii)\, given ${\bf  t}'\in U$, we can find $g\in GL(n,{\mathbb{C}})$ close to $I_n$ with ${\bf  t}(g)={\bf
t}'$.
We call the local coordinates ${\bf  t}$ at $O$ the standard local
coordinates at $O$ in $\mathbb{H}_n$.

We consider the Lie algebra
$\mathfrak X $ consisting of all left-invariant holomorphic vector
fields $X$ on $GL(n,{\mathbb{C}}) $. We identify $\mathfrak X $ with
$M_n({\mathbb{C}})$, the set of all $n\times n$ square matrices, as follows:  to $X=(\lambda_{ij}) \in M_n({\mathbb{C}})$ we associate a left-invariant holomorphic vector field $v_X$ on $GL(n,{\mathbb{C}})$ via, for $g=(x_{ij})\in GL(n,{\mathbb{C}})$,
\begin{align}
 \label{eqn:co-v-x}
v_X (g):= \sum_{i,j=1}^n \lambda_{ij} X_{ij}(g),  \quad  \mbox{  where
 } \ \
 X_{ij}(g)= \sum_{ k=1  }^{ n  } x_{ki}\frac{\partial }{\partial
x_{kj } }.
\end{align}
Hence we identify the vector
   field $v_X$ on $GL(n,{\mathbb{C}})$ with the matrix $X=(\lambda_{ij})$ in
 $M_n({\mathbb{C}})$ as additive groups. The integral curve ${C}_X$ for $v_X$ with initial value $I_n$ is given by
$$
{C}_X=\{\exp\, tX\in GL(n,{\mathbb{C}}) : t \in {\mathbb{C}}\},
$$
and  the integral curve of $v_X$ with initial value $g\in GL(n,{\mathbb{C}})$
 is given by  $g\,{C}_X \in GL(n,{\mathbb{C}})$.

 We let $\mathfrak
h_0$ denote the corresponding Lie subalgebra for $H'_0$,
so that
\begin{align}
 \label{eqn:g0m}
\mathfrak h_0
 =\bigl\{
\left(
\begin{array}{cccc}
 0&     \\
\vdots   & \ \  A \ \ \\
0       &
\end{array}
\right)\in M_n({\mathbb{C}}): \ A \in M_{n,\,n-1}({\mathbb{C}})\bigr\},
\end{align}
where $M_{n,\,n-1}({\mathbb{C}}) $ denotes the set of all
$n \times (n-1)$-matrices. We have
$$
H_0'=\{{\textstyle  \prod}_{i=1}^\nu \exp t_iX_i \in GL(n,{\mathbb{C}}): \nu\in {\boldsymbol{Z}}^+, \ t_i\in
{\mathbb{C}}, \ X_i\in \mathfrak h _0\}.
$$

\begin{theorem}\label{thm:special-hoph}\
 Let $D \Subset \mathbb{H}_n$ be a pseudoconvex domain with smooth boundary
 which is not Stein. Then there exists a Stein domain $D_0$ in $\mathbb P^{n-1}$ with
 smooth boundary such that  $D=\pi^{-1} _0(D_0)$.
\end{theorem}
\noindent {\bf Proof}.  We may assume $D$ contains the base point $O$.
Following Theorem
\ref{thm:main-nonconnected} for such domains $D$
in a homogeneous space, we fix a K\"ahler metric $ds^2$ on $G=GL(n,{\mathbb{C}})$ and a strictly positive
$C^{\infty}$ function $c$ on $G$ and we consider the $c$-Robin function
$\lambda ([z])$ for  $D$. Define the following subset of the Lie algebra
$\mathfrak X =  M_n({\mathbb{C}})$:
$$
\mathfrak X _0=\{X\in \mathfrak X : [\frac{\partial^2\lambda (\exp tX
(O))}{\partial t \partial \overline { t }}]|_{t=0}=0 \}.
$$
Under our assumptions for $D$ we showed  that $\mathfrak X _0$ is a Lie subalgebra of $\mathfrak X $ with
\begin{align}
\label{eqn:g0x0x}
\mathfrak h _0 \subsetneqq \mathfrak X _0 \subsetneqq \mathfrak X
\end{align}
(recall (\ref{strictinclusion}) and Remark \ref{rem:lastremark}). Fix $X\in \mathfrak X _0 \setminus \mathfrak h _0$. We may assume
that $X$ is of the form
$$
X=
\left(
\begin{array}{cccc}
 c_1&     \\
\vdots   & \ \ (*) \ \ \\
c_n       &
\end{array}
\right)\in M_n({\mathbb{C}}),
$$
where some $c_i\ne 0, \ i=1,\ldots,n$, and $(*)$ is a matrix in $M_{n,\,n-1}({\mathbb{C}})$.
 We show that
\begin{align}
 \label{eqn:j20}
c_j=0,   \quad   j=2,\ldots ,n.
\end{align}

\smallskip
We verify (\ref{eqn:j20}) by
 contradiction. Thus we assume that some $c_j \ne 0, \ j=2,\ldots ,n$.
Take any $Y\in \mathfrak h _0$ of the form
$
Y=({\bf  0}, \ldots ,{\bf  0}, {\bf  y}, {\bf  0}, \ldots , {\bf
0})
$,
where  ${\bf  0}$  is the zero $n$-colmun vector and ${\bf  y}$ is an $n$-column vector in the $j$-th column. An elementary calculation yields
$$
[X,Y]=
\left(
\begin{array}{cccc}
 -y_1c_j&     \\
\vdots &   (0)\\
-y_nc_j
\end{array}
\right) + Y',
$$
where $(0)$ is the zero matrix in $M_{n,n-1}({\mathbb{C}})$ and $Y'\in  \mathfrak
h_0$. Recall that $\mathfrak X_0$ is a Lie subalgebra of $\mathfrak X
$ with $\mathfrak h _0 \subset \mathfrak X _0$. Since the points $y_j, \ j=1,\ldots, n$, can be arbitrarily chosen in ${\mathbb{C}}$, if $c_j\ne 0$ it follows that
$\mathfrak X_0= \mathfrak X$. This contradicts (\ref{eqn:g0x0x}),
which proves (\ref{eqn:j20}).

\smallskip
We thus see that if $\mathfrak X_0 \ne  \mathfrak h _0$, then, $\dim \mathfrak X_0= 1+n(n-1)$ and
$$
\mathfrak X_0=
\bigl\{
\left(
\begin{array}{ccccccc}
 x &  \quad    &\\
 0 & \\
\vdots &  \quad (*)&\\
0 &
\end{array}\right):\  x \in {\mathbb{C}}, \  (*) \in M_{n,n-1}({\mathbb{C}})
\bigr\}.
$$
To complete the proof of the theorem we let $\Sigma_0$ denote the integral manifold for $\mathfrak X _0$ in
$GL(n,{\mathbb{C}})$:
$$
\Sigma_0=\{{\textstyle  \prod}_{i=1}^\nu \exp t_iX_i \in
GL(n,{\mathbb{C}}): \nu\in {\boldsymbol{Z}}^+,  \ t_i\in
{\mathbb{C}}, \  \ X_i\in \mathfrak X _0\},
$$
so that
$$
\Sigma_0=
\bigl\{
\left(
\begin{array}{c|cccccc}
 a &  \quad  (*)  &\\
\hline 0 & \\
\vdots &  \quad  A&\\
0 &
\end{array}\right):\  a \in {\mathbb{C}}^*; \ (*)\in \mathbb C^{n-1}, \  A\in
GL(n-1, {\mathbb{C}})
\bigr\}.
$$
Hence $\Sigma_0$ is a closed Lie subgroup of $GL(n,{\mathbb{C}})$ with
$ H_0\subset \Sigma_0$. Moreover, writing $^tO$ for the transpose of the row vector $O = [(1,0,\ldots ,0)]$,
$$
\sigma_0:= \psi_0(\Sigma_0)= \Sigma_0(^tO)= \{\,[(a,0,\ldots ,0)]\in \mathbb{H}_n: a\in {\mathbb{C}}^*\,\},
$$
which coincides with the torus $T_{\alpha}$.
 Since $GL(n,{\mathbb{C}})/
\Sigma_0= {\mathbb P}^{n-1}  $ and the projection
$gH_0 \mapsto g \Sigma_{z_0}$ from $\mathbb H_n$ to  $\mathbb P^{n-1}$ mentioned in Remark \ref{rem:lastremark} coincides with $\pi_0$ in 2. of
Proposition \ref{prop:glnc},
Theorem  \ref{thm:special-hoph} follows from Remark \ref{rem:lastremark} (taking $K=\mathbb H_n$ and $K_0=
{\mathbb{P}}^{n-1} $). \hfill $\Box$

\section{Flag space}\setcounter{section}{7}
In this section, we apply $1-b$. in Theorem \ref{thm:mainth} to
the  flag space ${\cal F}_n$ to determine all pseudoconvex domains $D
\Subset {\cal F}_n$ with smooth boundary which are not Stein (see
Theorem \ref{lasttheorem}). There are few known results on the Levi
problem in ${\cal F}_n $ (cf., Y-T. Siu \cite{siu}, K. Adachi
\cite{Adachi}). By definition, the flag space ${\cal F}_n$ is the set of all nested
sequences
\begin{align}
 \label{eqn:pointof -function}
z: \ \{0\} \subset F_1 \subset \ldots \subset F_{n-1} \subset {\mathbb{C}}^n,
\end{align}
where $F_i, \ i=1,\ldots ,n-1$ is an $i$-dimensional vector subspace of
${\mathbb{C}}^n$. We describe the structure of ${\cal F}_n$ as a homogeneous space. Given  $A=(a_{ij})\in GL(n,{\mathbb{C}}) $ we shall define an isomorphism
${\cal A}$ of ${\cal F}_n$. Consider the linear transformation
of ${\mathbb{C}}^n$ given by
$$
 {{} A}:
\begin{pmatrix}
 Z_1 \\
\vdots \\
Z_n
\end{pmatrix}=
\begin{pmatrix}
 a_{11} & \cdots &a_{1n} \\
  \vdots  &\ddots &\vdots \\
 a_{n1} & \cdots & a_{nn} \\
\end{pmatrix}
\begin{pmatrix}
 z_1 \\
 \vdots \\
z_n
\end{pmatrix}.
$$
For $z \in {\cal F}_n$ as in (\ref{eqn:pointof -function}), we then
define ${\cal A}(z)\in {\cal F}_n$ via
\begin{align}
 \label{gl-A-z}
 {{\cal A}}(z): \ \ \{0\} \subset   {{} A}(F_1) \subset   {{} A}(F_2)\subset
{{} A}(F_{n-1}) \subset  {\mathbb{C}}^n.
\end{align}
In this way $GL(n,{\mathbb{C}})$ acts transitively on ${\cal F}_n$; i.e., ${\cal
F}_n$ is a homogeneous space with  Lie transformation group $GL(n,{\mathbb{C}})$.

 We fix the following point $O$ in ${\cal F}_n$:
\begin{align*}
O: \ \  \
& \{0\}\subset F_1^0
\subset F_2^0
\subset \cdots \subset
F_{n-1}^0\subset {\mathbb{C}}^n,
\end{align*}
where
$$
F_i^0: z_{i+1}=\cdots = z_n=0, \  \quad  i=1,\ldots ,n-1.
$$
We call $O$ the {\it  base point} of ${\cal F}_n$. The isotropy subgroup $H_0$
of $GL(n,{\mathbb{C}}) $ for the point $O$ is the set of all upper triangular
non-singular  matrices:
$$
H_0= \left\{\begin{pmatrix}
 a_{11}& a_{12}& \cdots &a_{1n}\\
0&       a_{22}& \cdots & a_{2n}\\
  & \ddots & \ddots &  \\
0& 0& 0& a_{nn}
\end{pmatrix} \in GL(n,{\mathbb{C}})
\right\}.
$$
In particular,  $H_0$ is  connected in $GL(n,{\mathbb{C}})$ and $\dim H_0=
{{n(n+1)}\over 2}$.
 Since ${\cal F}_n \approx GL(n,{\mathbb{C}})/H_0$, the flag space ${\cal F}_n$  is
 a compact manifold with
$\dim {\cal F}_n=N:= {n(n-1)\over 2}$.

For local coordinates of a neighborhood of the base point $O$ in ${\cal
F}_n$ we can take
 $$
\begin{pmatrix}
 1&  &  &\\
&  \ddots    & 0 & \\
 & t_{ij}  & \ddots &  \\
& & & 1
\end{pmatrix}, \quad t_{ij} \in {\mathbb{C}}, \ \ \  1\le j<i \le n-1,
$$
where the point $O$ corresponds to the identity $I$ in $GL(n,{\mathbb{C}})$.
Equivalently
\begin{align*}
 {\bf  t}&=(t_{21}, t_{31}, \ldots, t_{n1}; t_{32}, \ldots,
 t_{n2};\ldots ;  t_{n\, n-1})\\
& \equiv  \ ({\bf  t}_1;\ {\bf  t}_2;\  \ldots ;\  {\bf   t}_{n-1})\in {\mathbb{C}}^N
 \end{align*}
where ${\bf  t}_j, \ j=1,\ldots, n-1$ is an $(n-j)$-column vector and the
base point
$O$ corresponds to the origin ${\bf  0}$ in ${\mathbb{C}}^{N}$.
We call these  local coordinates the {\it  standard local
coordinates}  at $O$.

For any $A=(a_{ij})\in GL(n,{\mathbb{C}}) $ which is sufficiently close  to
$I$,
each subspace ${A}(F_k^0),\ k=1,\ldots , n-1,$ in (\ref{gl-A-z})
can be written as
\begin{align} \label{eqn:can-rep}
 A(F_k^0): \
\left\{
\begin{array}{ccccccc}
 Z_{k+1}&=\alpha _{k+1\,1}^{(k)}Z_1+ \cdots &+\  \alpha
 _{k+1\,k}^{(k)}Z_k\\
\vdots & \vdots &\vdots \\
 Z_{n}&=\alpha _{n\,1}^{(k)}\,Z_1+ \cdots &+ \ \alpha _{n\,k}^{(k)}\, Z_{k},
\end{array}
\right.
\end{align}
where
\begin{eqnarray*}
 \left(
\begin{array}{cccccc}
 \alpha^{(k)}_{k+1\,1}& \cdots  & \alpha^{(k)}_{k+1\,k} \\
\vdots & \ddots & \vdots\\
\alpha^{(k)}_{n\,1}& \cdots  & \alpha^{(k)}_{n\,k} \\
\end{array}
\right)
&=&
\left(\begin{array}{cccc}
 a_{k+1\,1}& \ldots &a_{k+1\,k}\\
\vdots &  \ddots& \vdots\\
a_{n\,1}& \ldots  & a_{n\,k} \\
\end{array}
\right) \
\left(
\begin{array}{cccc}
 a_{11}& \cdots  & a_{1\,k} \\
\vdots &  \ddots& \vdots\\
a_{k\,1}& \cdots  & a_{k\,k} \\
\end{array}
\right)^{-1}  \quad  \nonumber \\
&\equiv& A_{n-k}^k\  A_{k}  ^{-1}.
\end{eqnarray*}
We call (\ref{eqn:can-rep}) the {   canonical representation} of the
$k$-dimensional vector space $A(F_k^0)$.
Using  the standard local coordinates ${\bf  t}$, we let
$$
 \gamma := (\gamma_{21}, \ldots ,\gamma_{n1}; \gamma_{32}, \ldots
 ,\gamma_{n2}; \ldots ; \gamma_{n\,n-1})
$$
denote the point $
 {\cal A}(O)$  in ${\cal F}_n$.
 Then
\begin{eqnarray} \label{eqn:gamma1}
 \left(
\begin{array}{ccccc}
 \gamma_{k+1\,k} \\
\vdots \\
\gamma_{n\,k}
\end{array}
\right) =
 \left(
\begin{array}{ccccc}
 \alpha_{k+1\,k}^{(k)} \\
\vdots \\
\alpha_{n\,k}^{(k)}
\end{array}
\right),  \quad k=1,\ldots ,n-1,
\end{eqnarray}
which is equal to the $k$-column vector determined by the coefficients of the variable
$Z_k$ in the canonical representation (\ref{eqn:can-rep}) of $A(F_k^0)$.

Recall that $M_n({\mathbb{C}})$ denotes the set of all $n \times n$-matrices with coefficients in ${\mathbb{C}}$ and we associate $X=(\lambda_{ij}) \in M_n({\mathbb{C}})$ to a left-invariant holomorphic vector field $v_X$ on $GL(n,{\mathbb{C}})$ as in (\ref{eqn:co-v-x}): for $g=(x_{ij})\in GL(n,{\mathbb{C}})$,
$$
v_X (g):= \sum_{i,j=1}^n \lambda_{ij} X_{ij}(g),  \quad  \mbox{  where
 } \ \
 X_{ij}(g)= \sum_{ k=1  }^{ n  } x_{ki}\frac{\partial }{\partial
x_{kj } }.$$
Then  ${C}_X=\{\exp\, tX\in GL(n,{\mathbb{C}}) : t \in {\mathbb{C}}\}$
is the integral curve of $v_X$ with initial value $I$ and $A\,{C}_X \in GL(n,{\mathbb{C}})$ is the integral curve of $v_X$ with initial value $A\in GL(n,{\mathbb{C}})$.

 This identification of left-invariant holomorphic vector fields and matrices will be heavily utilized. Also, unless otherwise stated, in the matrices occurring in this section all missing entries are $0$.

\begin{lemme}\label{prop:elp} ${}$
\begin{enumerate}
 \item Let $X= (\lambda_{ij})\in M_n({\mathbb{C}})$ and let
$ { C}_X(O): =\{(\exp\, tX)(O)\in {\cal F}_n \,: t \in {\mathbb{C}}\}$ be an
 analytic curve passing through the base point $O$ at $t=0$. Then
 the tangent vector of $C_X(O)$ at
$O$, in terms of the standard
 local coordinates ${\bf t}$  at $O$, has the form
$$
(\lambda_{21}, \lambda_{31,}\ldots ,\lambda_{n1};\lambda_{32}, \ldots
 ,\lambda_{n2};\ldots ;\lambda_{n\,n-1}).
$$
 \item Let $X=(\lambda_{ij})\in M_n({\mathbb{C}})$ and let $X'\in M_n({\mathbb{C}})$ be the matrix
$$
X': =
\left(\begin{array}{llllllllllllllll}
 0& &   & & \\
  \lambda_{21}     & 0&  & &\\
   \lambda_{31} & \lambda_{32}  & 0&& \\
  \vdots   &  \vdots & \ddots& \ddots & \\
\lambda_{n1}&\lambda_{n2} & \cdots & \lambda_{n\,n-1}&0
\end{array}
\right).
$$
For any $A\in H_0$, the direction of the analytic curve $A \exp tX(O)$ at
       $O$ in ${\cal F}_n$ is equal to that of the curve $A\exp tX'(O)$.
\end{enumerate}
\end{lemme}
{\bf Proof}. We have
$$
\exp tX =
 \left(
\begin{array}{ccccccccllllllll}
 1+  \lambda_{11}t &  \cdots
  &
 \lambda_{1n}t \\
\vdots &                 \ddots &    \vdots\\
 \lambda _{n1}t   & \cdots & 1+
\lambda_{nn}t
\end{array}
\right)
+ O(t^2).
$$
In terms of the standard local coordinates at $O$, we write
$$
{\bf  r}(t):= ({\bf  r}_1(t); {\bf  r}_2(t);\ldots ; {\bf  r}_{n-1}(t))
$$
for the point $\exp \, tX(O)$ in ${\cal F}_n$.
Fix $1\le k \le n-1$. Following (\ref{eqn:gamma1}) we consider the  $(n-k,k)$-matrix
$A_{n-k}^k(t)$ and $(k,k)$-matrix
$A_k(t )  ^{-1} $ for $A(t)=\exp tX$:
\begin{eqnarray*}
 A_{n-k}^k(t) =\left(
\begin{array}{ccccccccll}
   \lambda_{k+1\,1}t   & \cdots
   &   \lambda_{k+1\, k}t \\
 \vdots &                 \ddots  & \vdots    \\
 \lambda_{n1}t  & \cdots &
    \lambda_{nk}t
\end{array}
\right),
\end{eqnarray*}
\begin{eqnarray*}
 A_k(t)  ^{-1}&=& \left(
\begin{array}{ccccccccc}
 1+ \lambda_{11}t &   \cdots  &
   \lambda_{1k}t \\
\vdots &            \ddots &    \vdots\\
 \lambda_{k1}t &   \cdots &
  1+ \lambda_{kk} t
\end{array}
\right)^{-1} \\
&\equiv &( {\bf  d}^{(k)}_1(t), \ldots ,{\bf  d}_{k}^{(k)}(t) ),
\end{eqnarray*}
where each ${\bf  d}_j^{(k)}(t) \ j=1,\ldots ,k$ is a
$k$-column vector. We then have
$$
{}^t {\bf  r}_k(t)= A_{n-k}^k(t)\ {\bf  d}_k^{(k)}(t) +O(t^2),
$$
so that
\begin{eqnarray*}
\bigl[ \frac{ d}{dt}\,{}^t {\bf  r}_k(t)  \bigr]|_{t=0}=(A_{n-k}^k)'(0)\ {\bf  d}_k^{(k)}(0)=
 (A_{n-k}^k)'(0) \
\left(\begin{array}{ccccc}
 0\\   \vdots\\ 0\\ 1 \end{array} \right)=\left(\begin{array}{ccccc }
 \lambda_{k+1\,k}\\  \vdots \\ \lambda_{{n-1}k}\\ \lambda _{nk}
 \end{array}
\right),
\end{eqnarray*}
which proves $1$.

To prove 2. we write $X= (\lambda_{ij})\in M_n({\mathbb{C}})$ and $A=(a_{ij})\in
 H_0$.
 We denote by
$
{\bf  v}:= ( {\bf  v}_1\ ; {\bf  v}_2\ ; \ \cdots \ ; {\bf  v}_{n-1})
$
and $
{\bf  v}':= ( {\bf  v}'_1\ ; {\bf  v}'_2\ ; \ \cdots \ ; {\bf
v}'_{n-1})
$
the directions of the curves $A\exp tX(O)$ and   $A\exp tX'(O)$ at
$O$ in ${\cal F}_n$.  Following (\ref{eqn:gamma1}), for $k=1,\ldots , n-1$ we have
\begin{align}
 \label{eqn:d-use}
{}^t {\bf  v}_k&= \left(
\begin{array}{ccccccccccccc}
 \sum_{i=k+1   }^{ n  } a_{k+1\,i}\lambda _{i1},& \ldots ,&
\sum_{i=k+1   }^{ n  } a_{k+1\,i}\lambda _{ik}\\
 \sum_{i=k+2   }^{ n  } a_{k+2\,i}\lambda _{i1},& \ldots ,&
\sum_{i=k+2   }^{ n  } a_{k+2\,i}\lambda _{ik}\\
\vdots & \ddots & \vdots \\
a_{nn}\lambda_{n1}, & \ldots & a_{nn}\lambda_{nk}
\end{array}\right)
\left(
\begin{array}{cccc}
 a_{1k}^{(k)}\\
a_{2k}^{(k)}\\
\vdots\\
a_{kk}^{(k)}
\end{array}
\right),
\end{align}
where
the $k$-column vector in the right-hand side is the
$k^{th}$ column vector of the inverse matrix of $(a_{ij})_{i,j=1,\ldots k}$.
Hence, ${\bf  v}_k$ does not depend on $\lambda_{ij} $ for $(i,j)$ satisfying $1\le
i\le j\le n$, so that ${\bf  v}_k= {\bf  v}'_k$. Thus  2. is proved.
\hfill $\Box$

\smallskip
We define generalized flag spaces ${\cal
F}_n^\mathfrak M  $.
Let
\begin{eqnarray}
 \mathfrak M :=( m_1, \cdots ,m_\mu) \label{eqn:mfM}
\end{eqnarray}
be a fixed, finite sequence of positive integers with $1\le m_j\le n$ and
 $m_1+\cdots+ m_\mu=n$. Set
$
 n_j:=m_1+m_2+ \ldots +m_j$ for $j=1,\ldots , \mu$, and consider
a nested sequence $\zeta$ of vector spaces in ${\mathbb{C}}^n$
$$
\zeta:  \quad
\{0\}\subset S_{n_1} \subset S_{n_2} \subset \cdots \subset S_{n_{\mu-1}}
\subset {\mathbb{C}}^n,
$$
where $S_{n_j}, \ j=1,\ldots , \mu-1,$ is an $n_j$-dimensional vector
 space. Let ${\cal F}_n ^\mathfrak M  $ denote the set of all such
 sequences $\zeta$. We call ${\cal F}_n ^{\mathfrak m } $  the $\mathfrak M$-flag
 space in $ {\mathbb{C}}^n $. In particular, ${\cal F}_n ^\mathfrak M  $ coincides
 with ${\cal F}_n$ if  $\mu=n$. Clearly $GL(n,{\mathbb{C}})$ acts transitively on $ {\cal F}_n^{\mathfrak M} $. We fix the
 following point in $ {\cal F}_n^{\mathfrak M}  $ as the base point:
$$
O^{\mathfrak M}  :  \quad
  S_{n_j}= \{z_{n_j+1}= \ldots =z_n=0\},  \qquad  j=1,\ldots,\mu-1.
$$
Thus the isotropy subgroup $H ^{\mathfrak M}_0 $ of $GL(n,{\mathbb{C}})$ for $O
 ^{\mathfrak M} $ is the set of all  matrices
\begin{eqnarray*}
 \left(
\begin{array}{cccccc|ccc|cccc}
 &  & & & & & & & & & \\
 & &h_1&& & &&(*)&& &(*)  \\
 &  & & & & & & & & & \\
 &  & & & & & & & & & \\  \cline{1-12}
 &  & & & & & & & & &\\
 & &0 & & & &&h_j&& &(*)\\
 & &  & & & & & & & &\\  \cline{1-12}
 & &  & & & & & & & & \\
 & &0 & & & & &0 && & h_\mu&\\
 & &  & & & & & & & &
\end{array}
\right)
\end{eqnarray*}
where  $h_j \in GL(m_j,{\mathbb{C}}), \ j=1,\ldots, \mu$. Hence
\begin{enumerate}
\item[(i)] $ {\cal F}_n^{\mathfrak M} \ \approx \ GL(n,{\mathbb{C}})/ H _0^{\mathfrak
M}$;
 \item [(ii)] $H_0 \subset H_0 ^{\mathfrak M}  $ \ \ and\ \
$ H _0 ^{\mathfrak M}  / H_0 \ \approx \ {\cal F}_{m_1} \times \cdots \times
{\cal F}_{m_\mu}$,
  where ${\cal F}_{m_j}$ is the usual flag space in ${\mathbb{C}}^{m_j}$;
 \item [(iii)]  there exists a holomorphic projection
$
\pi_{\mathfrak M}  :  \ gH_0 \in {\cal F}_n \mapsto g H_0 ^{\mathfrak
       M}\in  {\cal F}_n ^{\mathfrak
       M}
$, such that
$$
(\pi _{\mathfrak M}  )^{-1}(\zeta)\  \approx  \  {\cal F}_{m_1} \times
       \cdots \times\ {\cal F}_{m_\mu},  \ \
\zeta\in {\cal F}_n ^{\mathfrak M}.
$$
\end{enumerate}
Our main result is the following.
\begin{theorem}
 \label{lasttheorem} \
Let $D \Subset  {\cal F}_n$ be a pseudoconvex domain with smooth
 boundary which is not Stein. Then there exists a unique sequence
$ \mathfrak M=(m_1,\ldots ,m_\mu)$ with $1< \mu<n$ and a Stein domain
 $D_0 \Subset  {\cal F}_n^{\mathfrak M} $ with smooth boundary such that
$
D =(\pi _{\mathfrak M}  )^{-1}(D_0)$.
\end{theorem}

To prove the theorem, we isolate the main ingredient, Lemma \ref{lem:vil}. To specialists in Lie theory, this is the statement that any  parabolic Lie subgroup of $GL(n,\mathbb C)$ containing the Borel Lie subgroup $H_0$ must be of the form $H_0 ^{(\mathfrak M)}
 \cap GL(n,\mathbb C)$ where $\mathfrak M$ is defined in
 (\ref{eqn:mfM}). We give an elementary proof accessible to non-specialists. We first introduce some useful
notation. For $1\le m\le n-1$ and $q=n-m\ge 1$,
\begin{align*}
O_n \ \mbox{or} \ O \ (\mbox{resp.}\ O_q)&= \mbox{ the base point of ${\cal F}_n$
 (resp. ${\cal F}_q$)};\\
M_q({\mathbb{C}})&= \mbox{ the set of all $q\times q$ square matrices}; \\
I_q&= \mbox{the identity in} M_q;
\\
M_{m,q}({\mathbb{C}}) &=\mbox{ the set of all $m \times q$ matrices} ; \\
 H_0^{(q)}&=\mbox{ the set of all upper triangular
matrices
 in $M_q({\mathbb{C}})$};
\\
M^{(m,n)}&=\left\{
\left(
\begin{array}{c|cc}
 X&Y  \\ \hline
0&0
\end{array}\right)\in M_n({\mathbb{C}})\ : \ X\in M_{m}({\mathbb{C}}), \ Y\in M_{m,q}({\mathbb{C}})\right\},
\end{align*}
and, for  an element $X_q\in M_q({\mathbb{C}})$ and
a subset $\mathfrak X^{(q)} \subset
M_q({\mathbb{C}})$,
\begin{align*}
\omega_n(X_q)&=
\left(
\begin{array}{c|cc}
             0&0  \\ \hline
         0& X_q
\end{array}
\right)\in M_n({\mathbb{C}});\\
\omega_n(\mathfrak X ^{(q)})&=\{ \omega_n(X_q)\in M_n({\mathbb{C}})\ : X_q \in
             \mathfrak X ^{(q)}\}.
\end{align*}

We note that from
 (\ref{eqn:co-v-x}),  the Lie subalgebra $\mathfrak H _0$ of $\mathfrak X$
which corresponds to the  isotopy group $H_0$ of $GL(n, {\mathbb{C}})$ for the identity
 $I$ can be written as  $
\mathfrak H_0=\{ v_X\in \mathfrak X : X\in H_0^{(n)}\}$. We identify
 $\mathfrak H_0$ with $H_0  ^{(n)}  $ just as we
 identify $ \mathfrak X $ with $M_n ({\mathbb{C}})$. Then $H_0=H_0 ^{(n)}\cap GL(n,{\mathbb{C}}) $.

\begin{lemme} \label{lem:vil} \
  Let $\mathfrak X_0 \subset M_n({\mathbb{C}})$ have  the following two
 properties:
\begin{enumerate}
 \item [(a)]\ $\mathfrak X_0$ is a Lie subalgebra of $M_n({\mathbb{C}})$ with
 $H_0^{(n)} \subset \mathfrak X_0 $;
 \item [(b)]\ let $X\in \mathfrak X_0$ and $A \in H_0$. Then  any
 $Y\in M_n({\mathbb{C}})$ which satisfies
\begin{align}
 \label{eqn:direc-con}
\bigl[\frac{d A \exp tX(O_n  )}{
dt}\bigr]|_{t=0}  =
\bigl[\frac{d \exp tY(O_n  )}{
dt}\bigr]\bigr|_{t=0}
\end{align}
 belongs to $\mathfrak X_0$.
\end{enumerate}
Assume  $H_0^{(n)} \subsetneq  \mathfrak X _0 \subsetneq M_n({\mathbb{C}})
 $. Then
there exist a unique integer $\mathfrak m, \ 1\le \mathfrak m\le n-1$,  and a unique subset
 $\mathfrak{X}_0^{(q)} \subset M_q({\mathbb{C}})$ where $q=n-\mathfrak m\ge 1$ such that
\begin{enumerate}
 \item [(1)]\ $\mathfrak{X} _0^{(q)}$ satisfies propeties $(a)$ and $(b)$ with $n$ replaced by $q$;
 \item [(2)]\ $\mathfrak X _0= M^{(\mathfrak m,n)}+ \omega_n( \mathfrak{X}_0^{(q)}).$
\end{enumerate}
\end{lemme}

 \begin{remark} \ {\rm  According to T. Morimoto (private communication), condition $(a)$ implies condition $(b)$. We include $(b)$ as a hypothesized property to maintain the elementary character of the proof of the lemma. Note that $(b)$ is very similar in content with 2. of Lemma \ref{le:2-is-new}. }
\end{remark}

\medskip
\noindent {\bf  Proof}. The uniqueness is clear. To verify the
existence of $\mathfrak X_0^{(q)}$ satisfying (1) and (2) we  divide
the proof into seven short steps.
For $1\le k \le \nu \le n$  we write $X_{\nu 1}\in M_n({\mathbb{C}})$
for the matrix whose $(\nu,1)$-entry is $1$ and whose other entries are
all $0$.

\smallskip
$1^{st} \ step.$  \ Fix $\nu$ with $2\le \nu \le n$.
 If $X_{\nu 1}\in \mathfrak X _0$, then $\mathfrak X _0$ contains
all matrices  of the form
$$
 Y_{(a_2, \ldots ,a_\nu)}= \left(
\begin{array}{ccccc|ccc}
 0& & \ \ \  & &  &&& \\
 a_2& & \ \ \  & &  & &&\\
 \vdots &&&&  &        & & \\
a_\nu &&&&   &      &&\\ \hline
&&&&      &        && \\
&&&&&&&
\end{array}
\right),  \quad a_j\in \mathbb{C}, \ j=2,\ldots , \nu,
$$
where  the missing  entries are all $0$.

\smallskip \noindent {\it Proof.}
Let $a_j \in {\mathbb{C}},  \ j=1,\ldots ,\nu$ with $a_\nu \ne 0$.
 We consider
$$
A=
\left(
\begin{array}{cccc|cc}
 1 &{}      &{}         &  a_1  & {} &{}\\
{} &\ddots  &{}         & {\vdots}         & {} &{}          \\
{} &{}      & 1         & a_{\nu-1}        & {}  &{}          \\
{} &{}      &           & a_{\nu}  & {} &{}          \\ \hline
{} &{}      &            &{}         & {}  & {}  \\
{} &{}      &            &{}         & {I_{n-\nu}}  & {} \\
{} &{}      &            &{}         & {}  & {}
\end{array}
\right)\in H_0.
$$
We let
$$
 {\bf  a(t)}:=({\bf  a}_1(t); {\bf  a}_2(t);\ldots ; {\bf  a}_{n-1}(t))
$$
represent the point $A\exp tX_{\nu1}(O)$ in ${\cal F}_n$
in terms of the standard local coordinates at $O$.  Then we have by
(\ref{eqn:gamma1})
\begin{align*}
 {\bf  a}_1(t)&=\frac{ 1}{1+a_1 t}\ (a_2 t, \ldots ,a_{\nu-1}t, a_\nu t, 0, \ldots
 , 0);\\
 {\bf  a}_j(t)&=(0,\ldots , 0),  \quad j=2,\ldots,n-1.
\end{align*}
It follows that the direction of the  curve $A\exp tX_{\nu1}(O)$ at $O$
in ${\cal F}_n$ is
$$
{\bf  v}=(a_2,\ldots ,a_\nu,0,\ldots ,0;0,\ldots ,0;\ldots ;0).
$$
From 1. in  Lemma \ref{prop:elp}  the curve $\exp tY_{(a_2, \ldots ,a_\nu
)}(O)$ at $O$ in ${\cal F}_n$
has the same direction ${\bf  v}$. Since $X_{\nu\,1}\in \mathfrak X _0$,
it follows from  property (b) in
Lemma \ref{lem:vil} that  $Y_{(a_2,
\ldots ,a_\nu)}\in \mathfrak X _0$.  This proves the $1^{st}$ step in the case where
$a_\nu\ne 0$. Suppose now $a_\nu=0$. From the previous case, we have $Y_{(a_2,\ldots ,a_{\nu-1},
 \varepsilon)}\in \mathfrak X _0$ for any $\varepsilon \in {\mathbb{C}}\setminus \{0\}$. Since $\mathfrak X
_0$ is a vector subspace of $M_n({\mathbb{C}})$, letting $\varepsilon \to 0$ we see that
the $1^{st}$ step is true for $a_\nu=0$.
  \hfill
$\Box$

\smallskip
$2^{nd} \ step.$  \  Fix $\nu$ with $2\le \nu\le n$.
 If $X_{\nu 1}\in \mathfrak X_0$, then $\mathfrak X _0$ contains all matrices
of the form
$$
Z_{(a_1, \ldots , a_{\nu-1})}=\left(
\begin{array}{cccccccccccc}
 0&    & & & && \\
 &0 & & &  &&\\
&&\ddots&&&&\\
  & & &0 & && \\
 a_1& a_{2}& \cdots &a_{ \nu-1} &0 &\cdots &0  \\
      &&  &&  & \ddots &\\
      &&&&&&  0
\end{array}
\right).
$$

\smallskip
\noindent  {\it Proof.}  It suffices to prove $Z_{(1,a_2, \ldots ,
a_{\nu-1})}\in \mathfrak X_0$ for any $a_2,\ldots ,a_\nu\in {\mathbb{C}}$. We consider
$$
A=
\left(
\begin{array}{cccc|cc}
 1 & -{a_{2}}      &{\ldots }         &  -a_\nu  & {} &{}\\
{} & 1  &{}         &          & {} &{}          \\
{} &{}      & \ddots         &          & {}  &{}          \\
{} &{}      &           & 1          & {} &{}          \\ \hline
{} &{}      &            &{}         & {}  & {}  \\
{} &{}      &            &{}         & {I_{n-\nu}}  & {} \\
{} &{}      &            &{}         & {}  & {}
\end{array}
\right)\in H_0.
$$
We let
$$
 {\bf  a(t)}:=({\bf  a}_1(t); {\bf  a}_2(t);\ldots ; {\bf  a}_{n-1}(t))
$$
denote the point $A\exp tX_{\nu1}(O)$ in ${\cal F}_n$
in terms of the standard local coordinates at $O$.  Then we have by
(\ref{eqn:gamma1})
\begin{align*}
 {\bf  a}_j(t)&=\frac{ 1}{1-a_\nu t}\ (0,\ldots , 0, a_jt,0, \ldots , 0),
 \quad \mbox{ $j=1,2,\ldots,\nu-1$},
\end{align*}
where we set $a_1=1$ and  $a_jt$ is  located in the $(\nu-j)$-th slot.
 For $\nu\le j\le
n$ we have ${\bf  a}_j(t)=(0,\ldots , 0)$.
It follows that the direction ${\bf  v}$ of the analytic curve $A\exp tX_{\nu 1}(O)$ at $O$
in ${\cal F}_n$ is
\begin{align*}
 {\bf  v} =(&0,\ldots ,0,a_1,0, \ldots , 0;\ \cdots\  ;
          0,\ldots ,0,a_{\nu-1},0, \ldots , 0; 0, \ldots,0 ;\ \cdots \ ;\ 0).
\end{align*}
From
1. in  Lemma \ref{prop:elp} and (b) in Lemma
\ref{lem:vil} we have $Z_{(1,a_2,\ldots ,a_{\nu-1})} \in \mathfrak X_0
$. \hfill $\Box$

\smallskip
Let $1\le k<\nu \le n$ and write $X_{\nu\, k}\in M_n({\mathbb{C}})$ for the matrix
whose $(\nu, k)$-entry is $1$ and whose other entries are all $0$. By
the
 $1^{st}$ and $2^{nd}$ steps we have:
\begin{align}
X_{\nu 1}\in \mathfrak X _0  \   \mbox{ {\it     implies} }   \
 X_{ij} \in \mathfrak X_0, \ \  1\le j <i\le \nu.
 \label{eqn:xnu}
\end{align}

\smallskip
$3^{rd}\ step$. \  Fix $\nu$ with $2\le \nu\le n$. Assume that $\mathfrak X _0$
contains a matrix of the form
\begin{eqnarray} \label{eqn:Y}
Y_\nu=
\left(\begin{array}{c|ccccccccccc}
 0 & & & & & \\
a_2& 0& &&&\\
a_3& &0&&&  & \mbox{all} \ \ (0) \\
 \vdots  & & &  &\ddots & &&\\
a_{\nu-1}& & & &&0\\
1 &  &  &(*)&&     & 0 \\ \cline{1-1}
0&   &&&&&&0\\
\vdots && &&&&&&\ddots  \\
0&&&&&&&&&0
 \end{array}
\right)
\end{eqnarray}
where the $\frac{ (n-2)(n-1)}2$\  entries $(*)$ are complex
numbers. Then $\mathfrak X _0$ contains a matrix of the form
\begin{eqnarray}\label{eqn:ystar}
 Y^*_\nu&=&
\left(\begin{array}{c|ccccccccccc}
 0 & & & & & \\
0& 0& &&&\\
0& &0&&&  & \mbox{all} \ \ (0) \\
\vdots  & & &  &\ddots & &&\\
0& & & &&0\\
1 &  &  &(\star)&&     & 0 \\ \cline{1-1}
0&   &&&&&&0\\
\vdots && &&&&&&\ddots  \\
0&&&&&&&&&0
 \end{array}
\right)
\end{eqnarray}
where the $(\nu,1)$-entry is $1$ and the entries $(\star)$ are complex numbers (i.e., $Y^*_\nu$ is a matrix of the same form as $Y_\nu$ but with the entries $a_2=\cdots = a_{\nu-1}=0$).

\smallskip \noindent {\it Proof.}  We consider
\begin{eqnarray*}
A= \left(\begin{array}{cccccc|cccccc}
 1 & & &        & &0 & & & \\
   & 1 &        & & & -a_{2} & & & & \\
& & 1 &        & &  -a_3 & & & & \\
& &  &1        & &  \vdots  & & & & \\
&  &    & & \ddots  &-a_{\nu-1}& & & & \\
&  &    &        & & 1 &  & & & \\ \hline
& &    &        & &  & 1& & & \\
&  &    &        & &  & &\ddots  & & \\
&  &    &        & &  & & & 1 \\
 \end{array}
\right)
\in H_0.
 \end{eqnarray*}
Then  the direction  ${\bf  v}$  of the curve $A\exp tY_\nu(O)$ at $O$
 is of the form
$$
{\bf  v}=(0,\ldots ,0,1,0 \ldots ,0; b_{32}, \ldots b_{n2};\ldots ; b_{n\, n-1})
$$
where $b_{ij}, \  2\le j<i\le n$ are complex numbers.
From
1. in  Lemma \ref{prop:elp} and (b) in Lemma
\ref{lem:vil},  $Y^*_\nu$ in (\ref{eqn:ystar}) for $(\star)=(b_{ij})$ belongs to $\mathfrak X _0$.
  \hfill $\Box$

\smallskip
$4^{th}\ step.$  \  Fix $\nu$ with $2\le \nu \le n$.
If $\mathfrak X _0$ contains at least one  matrix $Y_\nu^*\in M_n({\mathbb{C}})$  of
the form (\ref{eqn:ystar}),
then  $X_{\nu\,1}\in \mathfrak X _0$.

\smallskip
\noindent {\it Proof.} \ We set $(\star)=(a_{ij}),\,  2\le j<i \le n-1$ in $Y_\nu^*$ in (\ref{eqn:ystar}).
 For a given sequence
$$
{\{\bf  M};\varepsilon\}:\ \ \  M_2 \gg \ldots \gg M_{\nu-1} \gg 1 \gg \varepsilon _1 \gg \ldots \gg\varepsilon
_{n-\nu}>0,
$$
we consider the diagonal matrix
\begin{align}
 \label{eqn:ame}
A= A( {\bf  M}; \varepsilon )=
\left(
\begin{array}{cccccccccc}
1& & & & & & & & & \\
&M_2 & & & & & & & & \\
& & &\ddots  & & & & & & \\
& & & &M_{\nu-1} & & & & & \\
& & & & & 1     & & & & \\
& & & & &       &\varepsilon _1 & & & \\
& & & & &       & & \ddots  & & \\
& & & & &       & & & \varepsilon _{n-\nu}& \\
\end{array}
\right)
\in H_0.
\end{align}
 By direct calculation,  $A\exp tY_\nu^*$ is equal to
\begin{eqnarray*}
 \left(\begin{array}{c|ccccccccccc}
 1 & & & & & \\
0& M_2& &&&\\
0&M_3a_{32}t &M_3&&&  & \mbox{all} \ \  (0) \\
\vdots   & \vdots&\vdots  &  &\ddots & &&\\
0&M_{\nu-1}a_{\nu-1\,2}t &M_{\nu-1}a_{\nu-1\,3}t & \cdots  &&M_{\nu-1}\\
t &a_{\nu2}t  & a_{\nu3}t  &\cdots & &    a_{\nu\nu-1}t & 1 \\ \hline
0&\varepsilon _1 a_{\nu+12}t   & \varepsilon _1 a_{\nu+13}t & \cdots && \varepsilon
 _1 a_{\nu+1\nu-1}t&\varepsilon _1 a_{\nu+1\nu}t& \varepsilon _1\\
\vdots && \vdots &&&&&&\ddots & &  \mbox{all}\  \ (0) \\
\vdots &&\vdots &&&&&& &\varepsilon _{n-\nu-1}  \\
0& \varepsilon _{n-\nu}a_{n2}t   & \varepsilon _{n-\nu} a_{n3}t& \cdots
 & &
 &\cdots &&
 &\varepsilon _{n-\nu} a _{n n-1}t & \varepsilon  _{n-\nu}
 \end{array}
\right)
\end{eqnarray*}
up to terms of order $O(t^2)$. We write
$$
{\bf  v}=({\bf  v}_1 ; \ {\bf  v}_2; \ldots :{\bf  v}_\nu; \ldots ,
:{\bf  v}_{n-1})
$$
for the direction of $A\exp t Y_\nu^*(O)$ at $O$. Using
(\ref{eqn:gamma1}) we have  $
{\bf  v}_1=(0,\ldots ,0,1,0,
\ldots ,0),
$
where $1$  is in the $(\nu-1)$-st slot. We also have
$$
{\bf  v}_2= \frac{ 1}{M_2}(M_3a_{32},\,\ldots , M_{\nu-1}a_{\nu-1\,2},\,a_{\nu2},\,\varepsilon _1a_{\nu+1,\,2}, \ldots , \varepsilon _{n-\nu}
a_{n2}).
$$
Thus by taking
$$
M_2 \gg M_3, \ldots ,M_{\nu-1} , \varepsilon _1, \ldots , \varepsilon
       _{n-\nu}>0
$$
we can make $ {\bf  v}_2$ as close to the $(n-2)$-row vector ${\bf  0}$
as we like.

Similar results hold for $ {\bf  v}_j, \ j=3,\ldots ,\nu-1$; i.e.,
by taking
$$
M_j \gg M_{j+1}, \ldots ,M_{\nu-1} , \varepsilon _1, \ldots , \varepsilon
       _{n-\nu}>0,
$$
we can make $ {\bf  v}_j$ as close to the $(n-j)$-row vector ${\bf  0}$
 as we like. We have
\begin{eqnarray*}
{\bf  v}_{\nu}=
(\varepsilon _{1}a_{\nu+1\,\nu},\, \ldots ,\,\varepsilon
_{n-\nu}a_{n\nu});
 \end{eqnarray*}
thus by taking
$$
 1\gg \varepsilon _1, \ldots , \varepsilon
       _{n-\nu}>0,
$$
we can make $ {\bf  v}_{\nu}$  as close to the $(n-\nu)$-row vector ${\bf  0}$
 as we like. We have
$$
{\bf  v}_{\nu+1}=\frac{ 1}{\varepsilon _1} (
 \varepsilon _{2}a_{\nu+2\,\nu+1}, \ldots , \, \varepsilon
 _{n-\nu}a_{n\nu+1})
$$
Thus by taking
$$
 \varepsilon _1\gg \varepsilon _2, \ldots , \varepsilon
       _{n-\nu}>0,
$$
we can make $ {\bf  v}_{\nu+1}$ as close to the $(n-\nu-1)$-row vector ${\bf  0}$
 as we like.

Similar results hold for $ {\bf  v}_j, \ j=\nu+2,\ldots,
 n-1$. For example, for $ {\bf  v}_{n-1}$, we have $ {\bf  v}_{n-1}=
 \frac1{  \varepsilon _{n-\nu-1}} \ { \varepsilon _{n-\nu}}a_{n\, n-1}
$, so that, by  taking  $ \varepsilon _{n-\nu-1}\gg \varepsilon_{n-\nu}>0$
we make $ {\bf  v}_{n-1}$ as close to the complex number $0$
 as we like.

Therefore by taking
$$
\{{\bf  M};\varepsilon \}: \  \ \ M_2 \gg \ldots \gg M_{\nu-1} \gg 1 \gg \varepsilon _1 \gg \ldots \gg\varepsilon
_{n-\nu}>0
$$
and considering $A=A(
{\bf  M};\varepsilon   ) \in H_0 $ as in (\ref{eqn:ame}),  the direction
${\bf  v}$ of  $A \exp
tY_\nu^* (O)$ at $O$ in ${\cal F}_n$ can be made as close to
$$
{\bf  U}:=(0,\ldots ,0,1,0,\ldots 0;0,\ldots ,0;\ \cdots ;0).
$$
as we like. Since ${\bf  U}$ is
the direction of $\exp tX_{\nu\,1}(0)$ at $t=0$, it follows from
1. in  Lemma \ref{prop:elp} and (b) in Lemma
\ref{lem:vil} that $X_{\nu\,1}\in \mathfrak X _0$.
\hfill $\Box$

\smallskip As  the integer $\mathfrak m$ in Lemma \ref{lem:vil} we take
$$
\mathfrak m := \mathop {{\rm  max }}_{X\in \mathfrak X_0}\ \left[\,\mathop{{\rm  max}} \left\{i:\,\mbox{each $(j,1)$-entry of $X$ with $j>i$  is  $0$}\,\right\} \,\right].
$$
Since $H_0 ^{(n)} \subsetneq \mathfrak X _0 $ we have $\mathfrak m \geq 2$; on the other hand, from the $4^{th}$ step, we have
$\mathfrak m \le n-1$.

\smallskip
 $5^{th}\ step.$ \ The number $\mathfrak m $ has the following
properties:
\begin{enumerate}
 \item [(i)]\
$
 \bigl\{ \left(
\begin{array}{c|c}
 X_\mathfrak m  & 0\\ \hline
0 & 0 \\
\end{array}
\right)\in M_n({\mathbb{C}}) \ : \ X_ \mathfrak m \in M_ \mathfrak m
 ({\mathbb{C}})\bigr\}\ +\ H_0 ^{(n)}  \subset \mathfrak X _0;
$
 \item [(ii)]\ for any $X\in \mathfrak
       X_0
 $, each $(i,j)$-entry  with $\mathfrak m +1 \le i\le n; \  1 \le
       j\le \mathfrak m$ is $0$.
\end{enumerate}

\smallskip
\noindent {\it Proof.}    By the $3^{rd}$  and the $4^{th}$   steps we have
 $X_{\mathfrak m \,1} \in \mathfrak X _0$; then (\ref{eqn:xnu}) implies
assertion (i).
 We shall prove (ii) by contradiction.
   Assume that there
exists $A=(a_{ij})\in \mathfrak X _0$ with $a_{i_0j_0}\ne 0$ for
some $ \mathfrak m+1 \le i_0\le n; \ 1\le j_0 \le \mathfrak m $.
 Since $1\le j_0\le  \mathfrak m$, it follows from (i)  that
$X_{j_0\,1}\in \mathfrak X _0$.
Now $\mathfrak X _0$ is a Lie
subalgebra of $M_n({\mathbb{C}})$; thus $\mathfrak X _0$ contains the Lie bracket $[A,X_{j_0}]$. A calculation shows that the $(i_0,1)$-entry of $[A, X_{j_0}]$ is
equal to $a_{i_0j_0}\ne 0$. Since $i_0>\mathfrak m$, this contradicts the definition of
$\mathfrak m$. \hfill $\Box$

\smallskip
Note that we only used the fact that $\mathfrak X _0$ is closed under Lie brackets and the structure of $H_0$ in  the proof of Lemma \ref{lem:vil}. We set $q=n- \mathfrak m \ge 1$ and we consider ${\mathbb{C}}^n= {\mathbb{C}}^\mathfrak m
\times {\mathbb{C}}^q$.

\smallskip
$6^{th} \ step.$ \ We define
\begin{align}
 \label{eqn:x-0-q}
\mathfrak X _0^{(q)}:= \{Y_q\in M_q({\mathbb{C}}): \omega_n (Y_q)\in \mathfrak
X _0\}.
\end{align}
Then $\mathfrak X _0^{(q)}$ satisfies (1) in Lemma \ref{lem:vil}; i.e., $\mathfrak
X_0^{(q)}$ satisfies propeties $(a)$ and $(b)$ with $n$ replaced by $q$.

\smallskip
\noindent {\it Proof.}   Property ($a$) with $n$ replaced by $q$ clearly holds.
To prove (b), let $X_q \in \mathfrak X _0^{(q)}, \ A_q
       \in H_0^{(q)}\cap GL(q, {\mathbb{C}})$ and let $Y_q\in M_q({\mathbb{C}})$ satisfy
\begin{align}
 \label{eqn:q-derivative} \bigl[\frac{d A_q \exp tX_q(O_q)}{
dt}\bigr]|_{t=0}  =
\bigl[\frac{d \exp tY_q(O_q )}{
dt}\bigr]|_{t=0}.
\end{align}
To prove our claim that $Y_q\in \mathfrak X_{0} ^{(q)}  $,
we first write
$$ {\bf V}:=({\bf  v}_1\ ; \ {\bf  v}_2\ ; \ \cdots ; \ {\bf
v}_{q-1})  \quad   \mbox{and}  \quad
{\bf  v}_k=( v_{k+1\,k}, \ldots , v_{q\,k}), \ 1\le k\le q-1,
$$
 for the direction of the curve $A_q \exp tX_q(O_q)$ at $O_q$
in ${\cal F}_q$ in terms of the standard local coordinates ${\bf  t}
\in {\mathbb{C}}^Q$. Here $Q=\frac{ q(q-1)}2$. If we set
$$
Y'_q: =
\left(\begin{array}{llllllllllllllll}
 0& &   & & \\
  v_{21}     & 0&  & &\\
   v_{31} & v_{32}  & 0&& \\
  \vdots   &  \vdots & \ddots& \ddots & \\
v_{q1}&v_{q2} & \cdots & v_{q\,q-1}&0
\end{array}
\right) \in M_q({\mathbb{C}}),
$$
then  1. in Lemma \ref{prop:elp} with $n$ replaced by $q$ together with (\ref{eqn:q-derivative})
implies  $Y_q=  Y_q'+ B_q^*$ for
some $B_q^*\in H_0^{(q)}$.

We next consider
\begin{eqnarray*}
 X:= \omega_n (X_q) \in M_n({\mathbb{C}}),  \quad A:=
\left(
\begin{array}{c|c}
 I_\mathfrak m  & 0\\ \hline
0 & A_q \\
\end{array}\right)
\in H_0,
\end{eqnarray*}
so that  $X \in \mathfrak
 X _0$. By  (\ref{eqn:d-use}),
the direction of the curve $A \exp \, tX(O)$ at $O$ in ${\cal
F}_n$ is
$$
{\bf  U}=(
{\bf  0}_{n-1}; \  \ldots \ ;\ {\bf  0}_{n- \mathfrak m }; \ {\bf  v}_1;
 \
{\bf  v}_2; \ \ldots\ ;  {\bf  v}_{q-1}),
$$
where ${\bf  0}_k$ is  the zero $k$-row vector.
 It follows from 1. in Lemma \ref{prop:elp}
 and  (b) in Lemma
\ref{lem:vil} that $ \omega_n(Y_q')\in
\mathfrak X _0$.  Consequently,
 $
\omega_n(Y_q)=\omega_n(Y_q'+B_q^*)=Y'+ \omega_n(B^*_q)\in \mathfrak X
_0$, which is what we needed to show.
  \hfill $\Box$

\smallskip
\smallskip $7^{th} \ step.$ \ Property (2) in Lemma
\ref{lem:vil} holds for $\mathfrak X _0^{(q)}$ defined in (\ref{eqn:x-0-q}).

\smallskip
\noindent {\it Proof.} Let $X\in \mathfrak X _0$. By (ii) in the
 $5^{th}$  step,   $X$ is   of the form
$$
X= \left(
\begin{array}{c|cc}
X_{ \mathfrak m} & X_{\mathfrak m\,q}  \\ \hline
0 & X_q
\end{array}
\right),  \quad  X_{\mathfrak m}\in M_{\mathfrak m}({\mathbb{C}}), \  X_{\mathfrak m\,q}\in
M_{\mathfrak m,q}({\mathbb{C}}), \ X_q\in M_q({\mathbb{C}}).
$$
Since $\mathfrak X _0$ is a vector subspace of $M_n({\mathbb{C}})$,  it follows from (i) in the
 $5^{th}$  step that  $\omega_n (X_q)\in
\mathfrak X _0$; thus $X_q\in \mathfrak X _0^{(q)}$ and $X\in
M^{(\mathfrak m,n)}+
\omega_n(X_q)$, so that
$\mathfrak X _0 \subset M^{(\mathfrak m,n)}+
\omega_n(\mathfrak X_0 ^{(q)} ).$  The reverse inclusion follows from
(i) of the $5^{th}$ step and the
definition of $\mathfrak X _0^{(q)}$. The
$7^{th}$ step, and Lemma \ref{lem:vil}, are proved. \hfill $\Box$

\smallskip
{\bf Proof of Theorem  \ref{lasttheorem}}. It suffices to prove the
theorem under the assumption that $D$
contains the base point $O$ in ${\cal F}_n$. We consider the $c$-Robin
function $\Lambda(z)$ for $D$ and define
$$
\mathfrak X _0:=\{X\in M_n({\mathbb{C}}): \bigl[\frac{\partial^2  \Lambda(\exp
tX(O))}{\partial t \partial \overline{  t} } \bigl]|_{t=0}=0 \}.
$$
We see from the assumption for $D$ and Remark \ref{rem:lastremark} that $H_0 ^{(n)}  \subsetneq \mathfrak X _0
\subsetneq M_n({\mathbb{C}})$. Using 1. and 2. in Lemma \ref{le:2-is-new} for general
homogeneous spaces we see that
$\mathfrak X_0 $ satisfies properties (a) and (b) in Lemma
\ref{lem:vil}. Using Lemma \ref{lem:vil} inductively, we can find
$\mathfrak M=(m_1,\cdots, m_\mu)$ with $1<\mu<n$ such that
$\mathfrak X _0$ is the subset $H_{0}^{(\mathfrak M)}$ of $M_n({\mathbb{C}})$ which
consists of all matrices  of the form
\begin{eqnarray}\nonumber
 \left(
\begin{array}{cccccc|ccc|cccc}
 &  & & & & & & & & & \\
 & &h_{m_1}&& & &&(*)&& &(*)  \\
 &  & & & & & & & & & \\
 &  & & & & & & & & & \\  \cline{1-12}
 &  & & & & & & & & &\\
 & &0 & & & &&h_{m_j}&& &(*)\\
 & &  & & & & & & & &\\  \cline{1-12}
 & &  & & & & & & & & \\
 & &0 & & & & &0 && & h_{m_\mu}&\\
 & &  & & & & & & & &
\end{array}
\right),
\end{eqnarray}
where $h_{m_j}\in M_{m_j}({\mathbb{C}}), \ j=1,\ldots , \mu,$ and each $(*)$ is
an arbitrary element in the corresponding space $M_{m_j,m_k}({\mathbb{C}})$ (here $m_j<m_k$).
 It follows that the integral manifold $\Sigma_0$ of the Lie subalgebra $\mathfrak X _0$ passing
through $O$ in $GL(n,{\mathbb{C}})$,\,i.e., the connected Lie
subgroup of $GL(n,{\mathbb{C}})$ corresponding to $\mathfrak X _0$,    is
$H_0^{(\mathfrak M)} \cap GL(n,{\mathbb{C}})$ and hence it is equal to the
isotropy subgroup $H_0 ^{\mathfrak M}  $  of $GL(n,{\mathbb{C}})$ at the identity $I$ for the
generalized flag
space ${\cal F}_n ^{\mathfrak M} $.
   Thus, for  the flag space ${\cal F}_n$,  the space $M_0:=GL(n,{\mathbb{C}})/
\Sigma_0$, which was considered in the proof of Theorem \ref{thm:mainth}
for general homogeneous spaces, coincides with the space ${\cal F}_n ^\mathfrak
M$. Consequently,
 the projection $\pi_0$ and the analytic set $\sigma \Subset D \subset M$ defined in $1-b$. in
Theorem \ref{thm:mainth} coincide with the projection
$\pi_ \mathfrak M: AH_0\in {\cal F}_n \to A\Sigma_0=AH_0 ^\mathfrak M  \in {\cal F}_n
^\mathfrak M$ (where $A\in GL(n,{\mathbb{C}})$) and the analytic set $ H_0 ^ \mathfrak M/H_0
\Subset D \subset {\cal F}_n$.
Therefore,  $\pi_\mathfrak M
^{-1}  (\zeta)\approx {\pi_ \mathfrak M}  ^{-1}  (O ^{ \mathfrak M})= H_0 ^\mathfrak M/ H_0\approx \prod_{j=1}^\mu {\cal F}_{m_j}$ for $\zeta\in {\cal F}_n
^{\mathfrak M} $. Using $1-b$. in Theorem \ref{thm:mainth}, there exists a Stein domain $D_0 \Subset {\cal
F}_n ^\mathfrak M $ with smooth boundary such that  $D=
{ \pi}_\mathfrak M^{-1}(D_0)$.  Theorem \ref{lasttheorem} is completely
proved.
  \hfill $\Box$

The following remark is from T. Ueda.
  Consider two generalized flag spaces
${\cal F}_n^ \mathfrak M $ and
${\cal F}_n^ \mathfrak L$ in ${\mathbb{C}}^n$, where  $\mathfrak M=(m_1, \ldots
, m_\mu)$, $ \mathfrak L=(l_1,\ldots ,l_\nu)$, $\mu >\nu$, and
 $$
l_1=m_1+ m_2+ \ldots +m_{j_1}, \   \ldots ,   \
l_\nu= m_{j_{\nu-1}+1}  +      \ldots +m_{\mu}.
$$
We introduce the notation $\mathfrak M \prec \mathfrak L $ for this situation. Then we have the canonical projection
$$
\pi^\mathfrak M_ \mathfrak L:\ gH_0 ^\mathfrak M \in {\cal F}_n^
\mathfrak M \ \mapsto \
gH_0 ^\mathfrak L \in {\cal F}_n^
\mathfrak L,
$$
where $H_0 ^\mathfrak M$ is the isotropy subgroup of $G$ for ${\cal
F}_n^\mathfrak M$ at the base point $O^\mathfrak M$.  Thus,
for each $z\in {\cal F}_n^ \mathfrak L$,
$$
(\pi_{\mathfrak L}^\mathfrak M)  ^{-1}  (z) \approx {\cal F}_{ l_1}^
{\mathfrak M_1} \times \ldots \times
{\cal F}_{ l_\nu}^ {\mathfrak M_\nu}  \quad  \mbox{as complex manifolds},
$$
where
$\mathfrak
M_k=(m_{j_k+1}, \ldots , m_{j_k}),\ k=1,2,\ldots, \nu$.
If $\mathfrak M=(1,\ldots ,1)$, i.e., ${\cal F}_n^ \mathfrak M=
{\cal F}_n$, we simply write $\pi_ \mathfrak L^\mathfrak M= \pi
_\mathfrak L$.

We have the following result.
\begin{corollary}\label{due-ueda}\ Let $D$ be a pseudoconvex domain with smooth boundary in ${\cal F}_n^
 \mathfrak M$
 which is not Stein. Then there exists a unique $\mathfrak L$  such that
 $\mathfrak M \prec \mathfrak L $ and a Stein domain $D_0  $ in
 ${\cal F}_n ^ \mathfrak L$ with smooth boundary  such that $D =(\pi_
 \mathfrak L^ \mathfrak M)^{-1}(D_0)=D$.
\end{corollary}
\noindent {\bf Proof}. We assume that the base point $O ^\mathfrak M$ of
${\cal F}_n ^\mathfrak M$ is contained in $D$. Define $\widetilde { D}:=
(\pi _ \mathfrak M)  ^{-1}  (D) \ \subset {\cal F}_n$. Then
$\widetilde { D}$ is a pseudoconvex domain with smooth boundary in
${\cal F}_n$. We consider the $c$-Robin function $\lambda(z)$ for
$\widetilde { D}$ and define
$$
\mathfrak X_0= \{X \in \mathfrak X:\ \bigl[\frac{\partial^2 \lambda (\exp tX(
O)) }{\partial t \partial \overline { t}}\bigr]|_{t=0} \}
$$
where $O$ is the base point of ${\cal F}_n$ and $\mathfrak X$ is the Lie algebra consisting of all left-invariant
holomorphic vector fields on ${\cal F}_n$. Let $\Sigma$ be the Lie
subgroup of $GL(n, {\mathbb{C}})$ which corresponds to $\mathfrak X_0$. Finally, let $\mathfrak g_0^ \mathfrak M$ be the Lie subalgebra
corresponding to $H_0 ^ \mathfrak M \subset GL(n,{\mathbb{C}})$. Since $\lambda= \ const. $ on
the fiber $(\pi_ \mathfrak M)  ^{-1}  (O^ \mathfrak M)$,
we have $\mathfrak g_0 ^\mathfrak M\subset \mathfrak X_0
$ (we identify both sets as subsets of
$M_n({\mathbb{C}})$).  Hence
$H_0 ^\mathfrak M \subset \Sigma(O)$. Fix $\zeta\in D$ and take  a point  $z\in
(\pi_\mathfrak M)  ^{-1}  (\zeta)$. Then
$$
\widehat{ \lambda}( \zeta ):=  \lambda (z)
$$
is well-defined and $-\widehat{ \lambda}$ is a plurisubharmonic exhaustion
function on $ D$. Since $D$ is not Stein, there exists $\zeta_0\in D$ and
a direction ${\bf  a}\in
{\mathbb{C}}^{N} \setminus \{0\}$ (here $N= \dim \,{\cal F}_n^ \mathfrak M$) such that $ \bigl[\frac{\partial^2 \widehat{
\lambda}(\zeta+ {\bf  a}t) }{\partial t \partial \overline {
t}}\bigr]|_{t=0}=0 $. Thus if we take $z_0\in (\pi _\mathfrak
M) ^{-1}  ( \zeta_0)$, then  $\bigl[\frac{\partial \lambda (\zeta_0+ {\bf
a}'t)}{\partial t \partial \overline { t}}\bigr]|_{t=0}=0 $ where
$\pi _\mathfrak M(\zeta_0+ {\bf
a}'t)= z_0+ {\bf  a}t$. If we take $X\in \mathfrak X$  such that
$
\bigl[\frac{ d \exp tX(O)}{dt}\bigr]|_{t=0}= {\bf  a}'$, then
$X \in \mathfrak X_0  \setminus \mathfrak g_0^\mathfrak M$.
It follows from the proof of Theorem \ref{lasttheorem} that there
exists $\mathfrak L$ with $\mathfrak M \prec \mathfrak L$  and
a Stein domain $D_0$ with smooth boundary in ${\cal F}_n ^\mathfrak L$
such that $\widetilde { D}= (\pi _\mathfrak L)  ^{-1}  (D_0)$. Since
$\mathfrak M \prec \mathfrak L$, it follows from the definition of
$\widetilde { D}$ that $\pi_ \mathfrak L^\mathfrak M= \pi_
\mathfrak L| _{\widetilde { D}} \circ (\pi
_\mathfrak M)  ^{-1}|_D $ is well-defined;  $
\pi_\mathfrak L^ \mathfrak M:\,D\,\mapsto\,D_0$ is
 surjective; and $D= (\pi _\mathfrak L^ \mathfrak M) ^ {-1}(D_0)$,
as claimed.
 \hfill $\Box$

\smallskip  T. Ueda has another proof of this corollary
following ideas in the paper \cite{Adachi} (which is based on \cite{U}).

\section{Appendix A}
We discuss the three point property and the spanning property for projective space, Grassmannian manifolds, and flag spaces. Our first observation is elementary:

\smallskip
\noindent{\bf  1.} Projective space $M={\mathbb P}^n$ with Lie transformation group $G=GL(n,{\mathbb{C}})$ satisfies the
three point property.
\smallskip

On the other hand, we have:

\smallskip
\noindent{\bf  2.} The Grassmannian manifold $M=G(k,n)$ with
Lie transformation group $G=GL(n,{\mathbb{C}})$ and $n\ge 4, \ n-2\ge k \ge 2$ does not satisfy the three point property.

\smallskip
Recall that $M$ is the set of all $k$-dimensional subspaces of  ${\mathbb{C}}^n$; we use coordinates $x=(x_1,\ldots
,x_n)$. Given $ z\in M$, we may write $z$ as
$$
z:  \ \   \begin{pmatrix}
        x_1\\
\vdots \\
x_n
       \end{pmatrix} =
\begin{pmatrix}
 \alpha_{11} &\cdots &\alpha_{1k}\\
{}&\ddots& {}\\
\alpha_{n1}& \cdots & \alpha_{nk}
\end{pmatrix}
 \begin{pmatrix}
        \zeta_1\\
\vdots \\
\zeta_k
       \end{pmatrix}, \quad \mbox{i.e.,}  \quad  {}^t
       x=\alpha\
       {}^t\zeta,
$$
where $\zeta= (\zeta_1, \ldots , \zeta_k)\in {\mathbb{C}}^k$.
Then $g=(g_{ij})\in G$ acts on $z$ as follows:
$$
g(z): \ \
         {}^tx =   (g\cdot \alpha) \ {}^t \zeta \in M.
$$
 We let $O$ denote the point in $M$ defined as
$
 x_{k+1}=\ldots =x_n=0$.
Then the isotopy subgroup $H_{0}$ of $G$ for the point $O$ consists of
all elements in $G$ of the form $\left(
\begin{array}{c|ccccc}
 A_{k}& *  \\ \hline
   O &B_{n- k}
\end{array}
  \right)
$,
where $A_{k}$ and $B_{n-k}$ are nonsingular  square matrices of order $k$ and
$n- k$.

We prove that $(G,M)$ does not satisfy the three point property by contradiction. Thus we take
the following three points:
$$z_0=O; \  z_1: \, x_{k}=x_{k+2}=\ldots =0; \ z_2: \, x_{k-1}=x_k=\ldots =0.$$
Then we have
$$z_0\cap z_1 \ :\, x_k=x_{k+1}= \ldots =0 \  \hbox{and} \ z_0\cap z_2:\,
x_{k-1}=x_{k} =\ldots =0.$$
Assuming the three point property holds,
we can find  $ g \in H_{0}'$ so that $g(O)=O$ and
 $g(z_1)=z_2$.
Since $g$
 is one-to-one, we have
$ g(z_0\cap z_1)= g(z_0)\cap g(z_1)=z_0\cap z_2$, which is
$(k-2)$-dimensional.
 On the other hand, $g(z_0\cap z_1)$ is  $(k-1)$-dimensional,
 since $z_0\cap z_1$ is $(k-1)$-dimensional and $g$ is one-to-one. This is a contradiction.

\smallskip
\noindent{\bf  3.}  The Grassmannian manifold $M=G(k,n)$ with Lie transformation group $G=GL(n,{\mathbb{C}})$ satisfies the spanning property.

\smallskip
For simplicity we set $n=p+q, k=p$; $ M=G(p+q,p)$ and
$G=GL(p+q,{\mathbb{C}})$. We set $O:\ x_{p+1}=\ldots =x_{p+q}=0$ to be our base point of $M$. Then the isotropy sugroup $H_0$ of $G$ for  $O$ is
$$
H_0= \left\{
\begin{pmatrix}
    A_{p} & * \\
O & B_{q}
   \end{pmatrix} \Bigm| \det A_p, \det B_q \ne 0 \right\}.
$$
To prove {\bf 3.} it suffices to prove that $M$ satisfies the spanning
property for the point $O$ and $H_0$.
We identify $M$ as the space $G/H_0$ of all  cosets $\{gH_0:
g\in G\}$. Note that $\dim M= pq$ and $\dim\, H_0=(p+q)^2-pq$.
As local coordinates in a  neighborhood  of $O$ in $M$
we can take
 \begin{align*}
T(t)&=
\left(
\begin{array}{c|ccccc}
I_p  &     0 \\ \hline
(t)   & I_q
\end{array}
\right)
\end{align*}
where $(t)=(t_{ij})$ is a $q\times p$-matrix and the base point $O$ corresponds to
the identity matrix $I_{p+q}$. We identify $T(t)$ with
\begin{eqnarray}
  \label{coor}
   \mathbf {t} =  (t_{11}, t_{12}, \ldots ,
 t_{1p};  \ldots ;t_{q1}, \ldots ,  t_{qp})\in {\mathbb{C}}^{pq}
\end{eqnarray}
and we call these local coordinates the standard local coordinates at  $O$. Let $M_{p+q}({\mathbb{C}})$ be the set of all
$(p+q)$-square matrices $X$.

Let $  X\in M_{p+q}({\mathbb{C}}) $ satisfy condition
$$(\star) :\ \ \ \ \ \ \ \
 { \lim_{ t\to 0} \frac{d\exp tX(O) }{d t}|_{t=0}\ne {\mathbf{0}}} \in
 {\mathbb{C}}^{pq}.$$
 We decompose $X=A+B$ in
 $ M_{p+q}({\mathbb{C}}) $ where
$$
A= \left( \begin{array}{c|ccc}
   0 & 0 \\ \hline
   (a)&0
\end{array}\right) , \qquad
B= \left( \begin{array}{c|ccc}
   * & * \\ \hline
   0&    *
\end{array}\right),
$$
$(a)=(a_{ij})$ in $A$ is a $q\times p$-matrix and $0$ in $B$ is the
$q\times p$-zero matrix. Then we have
$\exp tX = \exp t(A+B)=(I_{p+q}+tA)(I_{p+q}+tB) + O(t^2)$. Let $h\in H_0$ and let $|t| \ll  1$. Since
$$
h(\exp tX)(O)=h(\exp tX)H_0 =h(I_{p+q}+tA)h^{-1}[h (I_{p+q}+tB)H_0]+O(t^2),
$$
it follows from $h (I_{p+q}+tB)H_0\in H_0$ that
$$
h(\exp tX)(O) = h(I_{p+q}+tA)h^{-1}(O) +O(t^2)  \quad  \mbox{  as points in $M$.}
$$
In particular, if we take $h=I_{p+q}\in H_0$, condition $(\star)$ for $X$ implies that
$(a)\not \equiv 0$, i.e., there exists some $a_{\lambda\nu}\ne 0$.
We show that we can choose a finite number of $h\in H_0$ with
$$
h(I_{p+q}+tA)h^{-1}= \left(\begin{array}{c|cc}
           I_p  & 0 \\ \hline
         (h(t))&I_q
    \end{array}\right),
$$
where $(h(t))=(h_{ij}(t))$ is a $p\times q$-holomorphic matrix in $|t|  \ll 1$
 and   the set of  tangent vectors
 \begin{align}
 \label{eqn:vh}
V[h]: =\left(
\begin{array}{ccccc}
 h'_{11}(0), &\ldots , &{ h}_{1p}'(0)\\
  \vdots & \ddots  & \vdots \\
{ h}'_{q1}(0), &\ldots , &{
h}'_{qp}(0))
\end{array}
\right)
\end{align}
form $pq$-linearly independent vectors in ${\mathbb{C}}^{pq}$. Note that $h(I_{p+q}+tA)h^{-1}\in H_0'$.
\smallskip

We first show that we may assume $a_{11}\ne 0$. We have $a_{\lambda\nu}\ne 0$ for some  $1\le  \lambda \le p;\ 1\le
\nu\le q$; thus we consider
the following matrix $\mathbf h_{ \nu \lambda}\in H_{0}$:
$$
\mathbf  h_{\nu\lambda}=\left(
\begin{array}{c|c|c|c}
 I^*_\nu   & {}       &       \multicolumn{2}{c}{}        \\ \cline{1-2}
 {}        & I_{p-\nu}  &     \multicolumn{2}{c}{}          \\ \cline{1-4}
  \multicolumn{2}{c|}{}          &  I_\lambda^*  & {}         \\ \cline{3-4}
 \multicolumn{2}{c|}{}           & {}          & I_{q-\lambda}
\end{array}
\right).
$$
Here $I_\nu^*$ is the anti-diagonal identity  matrix of degree $\nu$ and  the empty blocks are
$0$ matrices.
An elementary calculation yields
\begin{align*}
 {\cal A}(t):=\mathbf h_{\nu\lambda}\ (I_{p+q}+tA)\  {\mathbf h_{\nu\lambda}}^{-1}&=
\left(
\begin{array}{c|c}
I_p & 0  \\ \hline
(*) & I_q
\end{array}
\right),
\end{align*}
where
 \begin{align*}
(*)=
\left(
\begin{array}{ccccccccllllllll}
a_{\lambda \nu}t& a_{\lambda \nu-1}t&\cdots  & a_{\lambda 1}t & \cdots &
a_{\lambda p}t \\
\vdots         & \vdots           & \vdots & \vdots         & \vdots  &
\vdots \\
a_{1\nu}t       & a_{1\nu-1}t       & \cdots & a_{11}t        & \cdots &
a_{1p}t
 \\
\vdots       & \vdots             & \vdots                & \vdots  &
\vdots &
\vdots \\
a_{q\nu}t     & a_{q\nu-1}t         & \cdots & a_{q1}t      & \cdots & a_{qp}t
\end{array}
\right).
\end{align*}
Thus we see that if our claim were  true for ${\cal A}(t)$, then it is also true for $I_{p+q}+tA$. Hence we may assume $a_{11}\ne 0$.

\smallskip
We next  show that we can select a finite number of $h\in H_0$ so that  the tangent vectors $V[h]$ in (\ref{eqn:vh}) span
${\mathbb{C}}^{pq}$ under the condition that $a_{11}\ne 0$ in $I_{p+q}+tA$. For let $1\le i\le p$ and $1\le j \le q$ be fixed, and let $K\gg 1$. Consider
$$
\mathbf  h_{ij}(K)=\left(
\begin{array}{c|c|c|c}
 \mathbf m    & {}       &       \multicolumn{2}{c}{}        \\ \cline{1-2}
 {}        & I_{p-i}  &     \multicolumn{2}{c}{}          \\ \cline{1-4}
  \multicolumn{2}{c|}{}          &  \mathbf n  & {}         \\ \cline{3-4}
 \multicolumn{2}{c|}{}           & {}          & I_{q-j}
\end{array}
\right) \in H_0
$$
where  $\mathbf m  $ is an $i$-square matrix and $\mathbf n$ is a
$j$-square matrix with
$$
\mathbf m=\left(
\begin{array}{lllll}
 & & & &  1/K  \\
&  & & 1& \\
& &  \cdot  \\
&\cdot \\
1
\end{array}
\right), \qquad
\mathbf n=\left(
\begin{array}{lllll}
 & & & &  1  \\
&  & & \cdot& \\
& &  \cdot  \\
&1 \\
K
\end{array}
\right).
$$
A calculation gives
\begin{align*}
\mathbf h_{ij}(K)\ (I_{p+q}+tA) \mathbf h_{ij}(K)^{-1}
&=
\left(
\begin{array}{ccc|ccc}
1  &       &         &    &         &         \\
   & \ddots &         &    &         &         \\
   &       & 1       &    &         &          \\ \cline{1-6}
  &       &        &  1 &         &          \\
   &   D_{ji}(t)  &         &    &  \ddots &           \\
   &       &       &    &         & 1
\end{array}
\right),
\end{align*}
where
 \begin{align*}
D_{ji}(t)=
\left(
\begin{array}{ccccccccllllllll}
a_{ji}t& a_{j\, i-1}t&\cdots  &K a_{j1}t& a_{j\, i+1}t & \cdots &
a_{jp}t \\
\vdots         & \vdots      & \vdots     & \vdots & \vdots         & \vdots  &
\vdots \\
Ka_{1i}t       & Ka_{1\,i-1}t       & \cdots & \boxed{K^2 a_{11}t}  & Ka_{1\,i+1}t      & \cdots &
Ka_{1p}t
 \\
\vdots       & \vdots     &\vdots        & \vdots                & \vdots  &
\vdots \\
a_{qi}t     & a_{q\,i-1}t         & \cdots & Ka_{q1}t &a_{q\, i+1}t     & \cdots & a_{qp}t
\end{array}
\right).
\end{align*}
Thus the tangent vector $\mathbf v_{ij}:= V[ \mathbf h_{ij}(K)]$ of the curve $$\mathbf h_{ij}(K)\ (I+tA)
\mathbf h_{ij}(K)^{-1}(O)$$ at $O$ in $K$ is
\begin{align*}
\mathbf v_{ij}=
\left(
\begin{array}{ccccccccllllllll}
a_{ji}& a_{j\, i-1}&\cdots  &K a_{j1}& a_{j\, i+1} & \cdots &
a_{jp} \\
\vdots         & \vdots      & \vdots     & \vdots & \vdots         & \vdots  &
\vdots \\
Ka_{1i}       & Ka_{1\,i-1}       & \cdots & \boxed{K^2 a_{11}}  & Ka_{1\,i+1}      & \cdots &
Ka_{1p}
 \\
\vdots       & \vdots     &\vdots        & \vdots                & \vdots  &
\vdots \\
a_{qi}     & a_{q\,i-1}         & \cdots & Ka_{q1} &a_{q\, i+1}     & \cdots & a_{qp}
\end{array}
\right)\in {\mathbb{C}}^{pq},
\end{align*}
where  $K^2$ occurs only in the $(j,i)$-entry. Thus if we take $K\gg 1$ sufficiently large, then $a_{11}\ne 0 $ implies
that $\{\mathbf v_{ij}\in {\mathbb{C}}^{pq}: i=1,\ldots ,p; j=1,\ldots , q\}$ are linearly
 independent in ${\mathbb{C}}^{pq}$ (under the
 identification (\ref{coor})) and {\bf 3.} is proved.

 \bigskip
\noindent{\bf 4.}  The flag space $M={\cal F}_n$ for $n\ge 3$ with Lie transformation group $G=GL(n,{\mathbb{C}})$
does not satisfy the spanning property.

\smallskip
We use the notation from section 7: $O:\, \{0\}\subset F_1^0
\subset F_2^0 \subset \ldots \subset F_{n-1}^0$ (the base point of ${\cal F}_n$);
$H_0$ (the isotropy subgroup of $G$ for the point $O$); ${\bf
t}=(t_{21}, \ldots , t_{n1};t_{32},  \ldots, t_{n2} ; \ldots
;t_{n\,n-1})\in {\mathbb{C}}^{n(n-1)/2}$ where $n(n-1)/2=\dim {\cal F}_n$
(the standard coordinates of a  neighborhood of $O$).
 To show that ${\cal F}_n$ does not satisfy the spanning property, we
take $z_0=O\in {\cal F}_n$ and $X=X_{21}\in \mathfrak X \,\equiv M_n({\mathbb{C}})$  such that the $(2,1)$-entry is $1$ and all other entries are
$0$. A calculation gives, for $t\in {\mathbb{C}}$,
\begin{align*}
\exp \, tX (O):  \quad
& \{0\}\subset F_1(t)
\subset F_2(t)
\subset \cdots \subset
F_{n-1}(t)\subset {\mathbb{C}}^n,
\end{align*}
where
\begin{align*}
 &F_1(t): \ z_2=tz_1, z_3= \ldots =z_n=0;\\
& F_j(t): \ z_{j+1}= \ldots
=z_n=0,  \quad  j=2,\ldots ,n-1.
\end{align*}
Fix $h\in H_0$. Since $h(F_i^0)=F_i^0, \ i=1,\ldots , n-1$, it
follows that
\begin{align*}
h \exp tX(O):   \quad
& \{0\}\subset  \widetilde {  F}_{1}(t)
\subset F_2^0 \subset
 \ldots \subset F_{n-1}^0\subset {\mathbb{C}}^n,
\end{align*}
where $  \widetilde {  F}_1(t)$ is a holomorphic
function (depending on $h$) in $|t|\ll 1$. Using (\ref{eqn:gamma1}),
 this point $h \exp tX(O)$ may be written,
in terms of  the standard local
coordinates $\mathbf t$, as
$$
(\widetilde { f}_{21}(t), \ldots ,\widetilde { f}_{n1}(t);0,\ldots ,0; \ldots ;0),
$$
where $\widetilde { f}_{j1}(t), \ j=2,\ldots ,n,$  are holomorphic functions in
$|t| \ll   1$.
Therefore,
$$
\bigl\{\,\bigl[\frac{ d(h \exp tX(O))}{dt}\bigl]|_{t=0} \in {\mathbb{C}}^{n(n-1)/2 }
 \ :
 \ h \in H_0 \,\bigr\}
$$
contains at most $n-1$ linearly independent vectors in
${\mathbb{C}}^{n(n-1)/2 } $, so that ${\cal F}_n$ does not satisfy the
spanning property.

\section{Appendix B}
In this appendix, we complete the proof of Theorem
\ref{thm:main-nonconnected}.  First we recall our notation.  We write $H_{z_0}'$ for the connected component of $H_{z_0}$ in
$G$ containing the identity element $e$.
We defined  the homogeneous space  $M':=G/H_{z_0}'$ with Lie
transformation group $G$. We  write $w_0$ for the point in $M'$ which
corresponds to $H_{z_0}'$ in $G/H_{z_0}'$, so that the isotropy subgroup
$ \widehat{ H}_{w_0}$ of $G $ for $w_0$ is equal to $H_{z_0}'$. Since $\widehat{ H}_{w_0} $ is a closed normal
subgroup of $H_{z_0}$, we can consider the canonical projection
$\widehat{ \pi}: w=g \widehat{ H}_{w_0}\in M' \to z=g H_{z_0} \in M $,
so that $(M', \widehat{ \pi})$ is a normal covering space over $M$.
We defined $\widehat{ D}=\widehat{
 \pi}^{-1}(D)=\cup_{j=1}^\infty \widehat{ D}_j \subset M'$.
Since $z_0\in D$,  we have $w_0\in \widehat{
\pi}^{-1} (D)$. We then focused on $\widehat{ D}_1$, the connected
component of $\widehat{ \pi}^{-1}(D)$ containing $w_0$. We defined
$\widehat{ \lambda}(w):=\lambda( \widehat{\pi} (w))$ for $w\in \widehat{
D}_1$, which is a smooth plurisubharmonic exhaustion function on $\widehat{ D}_1$.

We begin with the following.
\begin{lemme}
 \label{eqn:saigono} \
 Let $D \Subset M$ be a domain with piecewise smooth boundary in $M$
and let ${\cal V}$ be a
 neighborhood  of  $e$ in $G$. Then there exist a finite
 number of balls $V^{(i)}, \ i=1,2,\ldots ,N$ centered at some point
 $z_i\in D$  such that
\begin{enumerate}
 \item[(1)] $ \overline { D} \Subset \mbox{  $\bigcup_{i=1}^N$} \,
        V^{(i)} \Subset  {\cal V}( { D})$
where ${\cal V}(  { D})=\{g(z)\in M: g\in {\cal V}, z\in
        \overline { D} \};$
\item [(2)] $\widehat{ \pi}^{-1}(V^{(i)})= \mbox{  $\bigcup
 _{k=1}^\infty$} \,  U_k
 ^{(i)}\ \  (disjoint \ union)$   in $ M'$ and   $\widehat{ \pi}: \ U_k ^{(i)} \to V
 ^{(i)} $ is  bijective;
 \item [(3)] there exists a
 holomorphic section $\sigma_k ^{(i)} :w\in U_k ^{(i)} \to \sigma_k
 ^{(i)} (w) $ of $G$ over $U_k ^{(i)} $ via $\widehat{ \psi}_{w_0}$
 such that
$\sigma_i ^{(k)} (w')( \sigma _i ^{(k)} (w''))^{-1} \in {\cal V}
$
 for\ any
 $w', \ w''\in U_i ^{(k)}$,  $k=1,\ldots , N$ and
 $i=1,2,\ldots $.
\end{enumerate}
\end{lemme}
\noindent {\bf Proof}. From the Borel-Lebesgue theorem it  suffices to verify the following. Fix $z\in M$ and a neighborhood ${\cal V}$ of $e$ in $G$. We show that there exists a neighborhood $V$ of $z$ in $M$  satisfying
\begin{enumerate}
\item [(1)] $\widehat{ \pi}^{-1}(V)= \mbox{  $\bigcup _{k=1}^\infty$} \,
  U_k
 \ (\mbox{disjoint  union})$
  in $ M'$ and $\widehat{ \pi}: \ U_k  \to V
  $ is  bijective;
 \item [(2)] there exists a
 holomorphic section $\sigma_k  :w\in U_k  \to \sigma_k
(w) $ of $G$ over $U_k $ via $\widehat{ \psi}_{w_0}$
 such that
$\sigma_k (w')\,( \sigma _k (w''))^{-1} \in {\cal V}$
 for\ any $w', w''\in U_k$,  $k=1,2, \ldots $.
\end{enumerate}

To prove this, we first take $g\in G$ with $g(z_0)=z$.
 Let $\boldmath v_0$ be the  neighborhood  of $e$ in $G$ stated in
 Proposition \ref{prop:disc}. Let $V \subset M, U_k \subset
 M', \ k=1,2,\ldots ,
 w_k=gh_k(w_0)$\,(the center of $U_k$) as in
 (\ref{eqn:sectionw}) and  (\ref{eqn:defofwk}), where $\widehat{ \pi}
 ^{-1}  (z)=\{w_k\}_{k=1, 2,\ldots} $.
We also want to insure that $\boldmath v_0$ also satisfies
$ (g \boldmath v_0)( g \boldmath v_0)^{-1} \subset {\cal V}$.
There exists a  neighborhood
$V_0$ of $z_0$ in $M$ and
 a holomorphic section
$\sigma_0: \zeta \in V_0\to \sigma_0(\zeta)$ of $G$
over $V_0$ via $\psi_{z_0}$  such that $\sigma_0(z_0)=e$ and
 $\sigma_0(V_0) \subset \boldmath v_0$.  By redefining $
 \boldmath v_0:= \boldmath v_0 \cap \pi_{z_0}^{-1}(V_0)$, we may
 assume that $ \psi_{z_0}(\boldmath
 v_0)=V_0$, so that
$
\sigma_0(V_0) \subset \boldmath v_0 \subset \sigma_0 (V_0) H_{z_0}'$,
 and hence  $\boldmath v_0 H_{z_0}'=\sigma_0( V_0)H_{z_0}' $.
In this situation,  (\ref{eqn:sectionw}) implies that $V= g
\sigma_0(V_0)(z_0)=g(V_0)$,  and hence we have a bijection between $V$ and $V_0$ defined by $\zeta\in V_0\to \xi = g(\zeta) \in V$. Then\begin{align}
 \label{eqn:v0H}
\sigma: \ \xi\in V \to g\sigma_0(\zeta)= g \sigma_{0} (g  ^{-1}  (\xi))\in G,
\end{align}
is a holomorphic section of $G$ over $V$ via
$\psi_{z_0}$ with $\sigma(z)= g\sigma_0(z_0)=g$ and $\sigma(V)= g
\sigma_0 (V_0)$. It follows that
\begin{align}
 \label{eqn:ukgs}
 \nonumber U_k&= g \boldmath v_0 H_{z_0}' h_k(w_0)=
 g\sigma_0(V_0)H_{z_0}'h_k(w_0)=
 \sigma(V) H'_{z_0}h_k(w_0) \\
&=\sigma(V) h_k H_{z_0}'(w_0)
= \sigma_{0}(V)h_k(w_0)  , \  \hbox{for} \ k=1,2,\ldots,
\end{align}
and hence we have a bijection between $V$ and $U_k$ defined by $\zeta\in V\to w= \sigma(\xi)h_k(w_0) \in U_k$ with $w_k=\sigma(z)
h_k(w_0)=gh_k(w_0)$. Thus
$$
\sigma_k:\ w\in U_k \to \sigma(\xi)h_k \in G,
$$
is a holomorphic section of $G$ over $U_k$ via
$\widehat{ \psi}_{w_0}$ with $\sigma_k(w_k)=gh_k$. Since
$\sigma_0(V) \subset \boldmath v_0$, it follows from
(\ref{eqn:v0H}) that for $w',w''\in U_k$ there exist $\xi', \xi''\in
V$ and $\zeta', \zeta''\in V_0$ with
\begin{align*}
 \sigma_k(w') \sigma_k(w'')^{-1}&= (\sigma(\xi')h_k)(
 \sigma(\xi'')h_k)^{-1}
=\sigma(\xi')\sigma(\xi'')\\
&= g\sigma_0(\zeta') (g \sigma_0( \zeta''))
 ^{-1} \in   (g \boldmath v_0 )(g \boldmath v_0 ) ^{-1}  \subset  {\cal
 V},
\end{align*}
as required. \hfill $\Box$

\smallskip
With the decomposition of $\widehat{ D}$ into its connected components $
\widehat{ D}= \cup_{j=1}^\infty \widehat{ D}_j$, where $\widehat{ D}_1$
is the connected component containing $w_0$,  and $D'(z_0), \ D^{(k)}(z_0)$ in
(\ref{eqn:connectedcom}) are the connected components of $D(z_0)$ in
$G$, we
set
${\bf  h} ^{(k)}\in {\cal H}^{(k)} (z_0)\subset
D^{(k)}(z_0)$, $k=2,3,\ldots $ as in equation (\ref{eqn:atode}) for
$z=z_0$, so that $D ^{(k)} (z_0)= D'(z_0){\bf  h} ^{(k)}$. We next show that, after possibly relabeling indices,
\begin{align}
\label{eqn:indeces-arrang}
 \widehat{ D}_k=
\{g(w_0)\in M': g\in D^{(k)}(z_0)\},  \quad  k=1,2,\ldots.
\end{align}
To prove this, fix $k$ and the set ${ D}^{(k)}(z_0)$. Then ${\bf  h}
^{(k)} (w_0)\in \widehat{
\pi}^{-1}(z_0) \subset \widehat{ \pi}^{-1}(D)$ and there exists a unique
connected component of $\widehat{ \pi}^{-1}(D)=\sum_{j=1   }^{\infty
} \widehat{ D}_j$ which contains ${\bf
h}^{(k)} (w_0)$. We call this component $\widehat{ D}_k$ and we shall first prove equality of the two sides of
(\ref{eqn:indeces-arrang}) for this $k$.

We let $E_k$ denote the right-hand-side of
(\ref{eqn:indeces-arrang}). Using (\ref{eqn:atode}), we can write
$$E_k=D'(z_0)\, {\bf  h}^{(k)}  (w_0)=\{g'{\bf  h}^{(k)} (w_0) \in M': g'\in
D'(z_0)\}.
$$
This is a connected set in $M'$ which contains ${\bf
h} ^{(k)} (w_0)$ and satisfies $\widehat{ \pi}(E_k)\subset D$
 since $h ^{(k)} (w_0)=z_0$. Hence $E_k \subset \widehat{D}_k$.
 To prove the reverse inclusion we set $w_k= {\bf  h} ^{(k)} (w_0)\in
 \widehat{ D}_k$ and  let $w\in
 \widehat{ D}_k$. We take a continuous curve $\gamma:t\in [0,1] \to
 w(t)$ in $\widehat{ D}_k$  with $w(0)=w_k$ and $w(1)=w$. Then
 $\widehat{ \pi}(\gamma):t\in [0,1]\to z(t):= \widehat{ \pi}(w(t))$ is a
 continuous curve in $D$  with $z(0)= {\bf  h}^{(k)} (z_0)=z_0$. Using  property {\bf  2.} of complex homogeneous spaces from the beginning of section 6, we can find a continuous
 section $\sigma: t\in [0,1]\to \sigma(t)$ in $D'(z_0)$ over
 $\widehat{\pi}(\gamma) $ via $\psi_{z_0}$, i.e., $\sigma(t)(z_0)=z(t),
 \ t\in [0,1]$, with  $\sigma(0)=e$. We set $g_0:=\sigma(1)\in D'(z_0)$. We consider the continuous curve
$\sigma {\bf  h} ^{(k)} (w_0): t\in [0,1]\to \sigma(t){\bf  h}^{(k)} (w_0)$ in $M'$. Then $\widehat{\pi}(\gamma)= \widehat{ \pi}({\sigma(w_0)}) \subset D$. The two curves
 $\gamma$ and $\sigma {\bf  h}^{(k)}  (w_0)$ in $M'$ start at the same point
 $ w_k={\bf  h} ^{(k)} (w_0)\in
\widehat{ D}_1$. It follows that $ \gamma=\sigma {\bf  h}^{(k)} (w_0)$ as curves in
$M'$, and hence $w=w(1)=\sigma(1){\bf  h}^{(k)}  (w_0)=g_0 {\bf
h}^{(k)} (w_0)$. Thus
$\widehat{ D}_k \subset E_k$.

To finish the proof of (\ref{eqn:indeces-arrang}) it remains to show that
the sets $\{E_k\}_{k=1,2,\ldots }$ exhaust each set $\widehat{ D}_l$.
If not, by the previous argument, there exists a set $\widehat{ D}_l$ which is not equal to
any of the sets $E_k, \, k=1,2,\ldots $. From (\ref{qn:yasisugiru}) we have  $\cup_{k=1}^\infty {\cal
H}_k(z_0)=H_{z_0} $. Since ${\cal H}_k(z_0 ) \subset D^{(k)}(z_0)$, we
have $({\cal H}_k(z_0))(w_0) \subset E_k=\widehat{D}_k$. Take a point $w^*\in
\widehat{
\pi}^{-1} (z_0)\cap \widehat{ D}_l$. Then there exists  $h^*\in H_{z_0}$ with  $w^*=h^*(w_0)$. Since  $h^*\in {\cal H}_k(z_0)$ for some
$k$, we have  $h^*(w_0) \in \widehat{ D}_k$. Thus $h^*(w_0)\in \widehat{
D}_k \cap \widehat{ D}_l=\emptyset$, a contradiction. This proves (\ref{eqn:indeces-arrang}).

Since $\widehat{ D}_1$ is a domain in $M'$ which contains $w_0$, as with our discussion of $D(z_0)$ we will
consider the set:
$$
\widehat{ D}_1(w_0):=\{g\in G : g(w_0)\in \widehat{
D}_1\} \subset G.
$$
We have  the following equalities:
\begin{equation} \label{eqn:d-1-prime}
 \widehat{ D}_1=\{g(w_0)\in M' : g\in D'(z_0) \};  \quad
D'(z_0)= \widehat{ D}_1(w_0).
\end{equation}
The first equality is the case $k=1$ in (\ref{eqn:indeces-arrang}). The
inclusion $D'(z_0)\subset  \widehat{ D}_1(w_0)$ in the second equality follows from the first equality.
 To prove $\widehat{ D}_1(w_0)\subset D'(z_0)$, fix $g\in G$ with $g(w_0)\in \widehat{
D}_1 $. Take
a continuous curve $\gamma: t\in [0,1]\to w(t)\in \widehat{ D}_1$ with $w(0)=w_0$ and $w(1)=g(w_0)$. Since $\widehat{ D}_1 \subset M'=G/
\widehat{ H}_{w_0}$, we can find a continuous section $\sigma: t\in [0,1]
\to \sigma(t)$ in $G$ over $\gamma$ via $\widehat{ \psi}_{w_0}$ with  $\sigma(0)=e$. Thus $\sigma(t)(w_0)=w(t), \ t\in [0,1]$. We set
$g_1:=\sigma(1)\in G$. Since $g_1(w_0)= g(w_0)$, there exists $h\in \widehat{ H}_{w_0}$  such that $g=g_1h$. Since $\widehat{
H}_{w_0}=H_z'$ is connected in $G$, we can find a continuous curve
$\tau:t\in [1,2]\to \tau(t)$ in $\widehat{ H}_{w_0}$  such that
$\tau(1)=e$ and $\tau(2)= h$. Then $g_1\tau: t\in [1,2]\to g_1\tau(t)$ is
a continuous curve in $G$ starting at $g_1$ at $t=1$ and satisfying
$g_1 \tau(t)(w_0)= g_1(w_0)=g(w_0)\in \widehat{ D}_1, \ t\in
[1,2]  $
(since $\tau(t)\in \widehat{ H}_{w_0}$ and $\tau(t)(w_0)=w_0$). We form the continuous curve $\Gamma:t\in
 [0,2]\to \Gamma(t)$ in $G$ via the concatenation
$$\Gamma(t):= \sigma(t), \ t\in [0,1]; \ \hbox{and} \ \Gamma(t):=g_1\tau(t), \ t\in [1,2].$$
This curve $\Gamma$ starts  at
 $\Gamma(0)=e$ and terminates at $\Gamma(2)=g_1\tau(2)= g$. One can check that
 $\Gamma(t)(w_0)\in \widehat{ D}_1, \ t\in [0,2]$, so that
 $\Gamma(t)(z_0)\in D, \ t\in [0,2]$. Since
 $\Gamma(0)=e \in D'(z_0)$, we have $\Gamma(t)\in D'(z_0), \ t\in [0,2]$. In
 particular, $g\in D'(z_0)$, which verifies the inclusion $\widehat{ D}_1(w_0)\subset D'(z_0)$ and
 hence the second equality in (\ref{eqn:d-1-prime}).

\smallskip
Using (\ref{eqn:H-prime}) we set
$$ {\cal H}'(z_0)=
D'(z_0) \cap H_{z_0}= \mbox{  $\bigcup_{j=0}^\infty $} \, h'_{j} \widehat{ H}_{w_0}
 \ \ (\mbox{disjoint union}),
$$
where $h'_1=e$ and $h'_j\in H_{z_0}\cap D'(z_0), j=2,3,\ldots $.
 Thus, by
(\ref{eqn:d-1-prime})  we have
\begin{align}
 \label{eqn:wue3}
\widehat{\pi}^{-1}(z_0) \cap \widehat{ D}_1 =
\{w_0, h'_2(w_0),h'_3(w_0), \ldots \} \equiv:  \{w_0, w_0^{(2)}, w_0^{(3)}, \ldots \},
\end{align}
and these points are distinct. Moreover, if  $z\in D$, $w\in\widehat{ \pi}^{-1}(z)\cap  \widehat{
 D}_1 $ and $g\in D'(z_0)$ with $g(w_0)=w$ in $M'$, then
\begin{align}
 \label{eqn:tatto-deki}
\widehat{ \pi}^{-1}(z)\cap \widehat{ D}_1&=
\{g(w_0), g(w_0^{(2)}),
 g(w_0^{(3)}), \ldots  \}.
\end{align}
We first verify that the right-hand side is contained in the left-hand side. Take any element
 $g(w_0^{(j)})=gh'_{j}(w_0)$ in the right-hand side. Then we have  $\widehat{
 \pi}(g(w_0^{(j)})) = g h_{j}(z_0)=g(z_0)=\widehat{\pi}(g(w_0))=
 \widehat{ \pi}(w)=z$. Since $g\in D'(z_0)$ and $ h_j'\in
 D'(z_0)\cap H_{z_0}$, we see from 2. in Proposition \ref{eqn:elm} that
 $gh'_{j}\in D'(z_0)$.  Using (\ref{eqn:d-1-prime}) we have $g(w_0^{(j)})=gh'_j(w_0)\in \widehat{ D}_1$, which
 proves this inclusion. To prove  the reverse inclusion,
let $w' \in
 \widehat{ D}_1 \cap \widehat{ \pi}^{-1}(z)$. Using (\ref{eqn:d-1-prime}) we can find $g'\in D'(z_0)$ with $g'(w_0)= w'$. Taking $\widehat{ \pi}$ we have $g'(z_0)=z=
 g(z_0)$ in $M$,
 so that there exists an $h\in H_{z_0}$ with $g'=g h$.
From 3. in Proposition  \ref{eqn:elm} we have $ h\in D'(z_0)\cap H_z$.
We thus have $h\in h'_j \widehat{ H}_{w_0}$ for some $j$, and
 hence $w'= gh(w_0)=gh'_j(w_0)=g(w_0  ^{(j)}  )$, as claimed.

\smallskip
Let $E$ be a domain  with piecewise smooth boundary in $M$
which contains $D$. We write
  $E'(z_0)$  for the connected
  component of $E(z_0)$ containing $e$ in $G$.
 If
$E'(z_0)\cap H_{z_0}=D'(z_0)\cap H_{z_0}= {\cal H}'(z_0)$,
then we say that $E$ and $D$ have the {\it   same isotropy class}  at
$z_0$.
As for $\widehat{ D}_1$, we write $\widehat{ E}_1$  for
 the connected component of
 $\widehat{ \pi}^{-1}(E)$
 containing $w_0$ in $M'$.

Let $E$ be a domain containing $D$ with the same isotropy class as $D$ at $z_0$. Fix $z\in E$
and let ${\cal V}$ be a  neighborhood  of $e$ in $G$. Choose $g\in
E'(z_0)$ with $g(z_0)= z$. We can find a
neighborhood  $V$ of $z$ in $E$ which satisfies (1) and (2) in
the proof of Lemma  \ref{eqn:saigono}. For simplicity, if $h_k\in H_{z_0}$ is
contained in ${\cal H}'(z_0)$, say $h_k=h_j'$, then we write
$U_j'=\sigma(V)h_j'(w_0)$ for the set $U_k=\sigma(V)h_k(w_0)=g \sigma_0(V_0) h_k(w_0)$ in
(\ref{eqn:ukgs}).  Since $gh'_j\in E'(z_0) $ by 2. in Proposition
\ref{eqn:elm}, it follows from (\ref{eqn:d-1-prime})  that the center
$gh'_j(w_0)$  of $
U_j'$ is contained in $\widehat{ E}_1$. This together with $V \subset E$
 implies that $U_j' \subset \widehat{ E}_1 \cap \widehat{ \pi}  ^{-1}  (V)$.
Moreover we  show
\begin{align}
 \label{eqn:tatto-deki-2}
\widehat{ \pi}^{-1}(V)\cap \widehat{ E}_1= \{U_1'=U_1, U_2',U'_3,  \ldots \}.
\end{align}

It remains to prove that the left-hand side is contained in the
 right-hand side. Let  $w'\in \widehat{ \pi} ^{-1}  (V) \cap \widehat{ E}_1$.
 We can find $z'\in V_0$ and $h_k \in H_{z_0}$  with $w'=g
 \sigma_0(z')h_k(w_0)$. Take a continuous curve $\gamma:t\in [0,1] \to f(t)$ in
 $ \sigma_0(V_0)$  such that $f(0)=\sigma_0(z')$ and $f(1)=e$.
Then $\widetilde { \gamma}:t\in [0,1]\to \widetilde {
 \gamma}(t): = gf(t)h_k(w_0)$ is a continuous curve in $M'$ with
$\widehat{ \pi}( \widetilde { \gamma}(t))\in g
\sigma_0(V_0)(z_0)=g(V_0)=V \subset E$.
Since
$\widetilde { \gamma}(0)(w_0)=g \sigma_0(z')h_k(w_0)=w'\in
 \widehat{ E}_1 $, it follows that $\widetilde{\gamma} \subset
\widehat{  E}_1
$. In particular, we have $\widetilde { \gamma}(1)= gh_k(w_0)\in
\widehat{ E}_1$, so that $gh_k\in \widehat{ E}_1(w_0)= E'(z_0)$ by
(\ref{eqn:d-1-prime}). Since $g\in E'(z_0)$, this implies that $h_k\in {\cal H}'(z_0)$ by 3. in
Proposition      \ref{eqn:elm}. Thus (\ref{eqn:tatto-deki-2}) is proved.

\smallskip
We prove the following.
\begin{lemme} \label{lem:d1d2} \
 Let $D \Subset E \Subset M $ and $z_0\in D$. Let  $t\in [0,1] \to D(t)
 \Subset M$ be a one-parameter family of domains $D(t)$ satisfying $D \subset D(t) \subset E;\ D(0)=D; \  D(1)=E$;
and each boundary $\partial D(t)$ is
 piecewise smooth in $M$ with $\partial D(t)$ varying continuously with $t\in [0,1]$. Then
 \begin{enumerate}
 \item [(i)]  $D$ and $E$
 have the same isotropy class at $z_0$;
 \item [(ii)]    $\widehat{ D}_1 = \widehat{ E}_1 \cap \widehat{
 \pi}^{-1}(D)$ in $M'$;
 \item [(iii)] $D'(z_0)= E'(z_0)\cap D(z_0)$ in $G$.
\end{enumerate}
  \end{lemme}

\noindent {\bf  Proof}.  Let $D \Subset M$ be
a domain with piecewise smooth boundary $\partial D$. Let $z_0\in D$. Fix $z\in
\partial D$. We choose a neighborhood
$V= g (V_0)$ of $z$  in $M$
as constructed in the proof of
Lemma  \ref{eqn:saigono}
and we set $E:=D\cup V \subset M$.
To prove Lemma \ref{lem:d1d2}  it suffices to verify (i), (ii) and (iii) for such $D$ and $E$.
 To prove  (i), i. e.,
\begin{align}
 \label{eqn:restar}
  E'(z_0)\cap H_{z_0}= D'(z_0) \cap H_{z_0},
\end{align}
we consider  $\widetilde { D}_1 := \widehat{ \pi}^{-1}(D) \cap
\widehat{E}_1 \subset M'$.  Thus, if $g\in G $ with $g(w_0)\in  \widetilde { D}_1$, then $g(z_0)\in D$. Using (\ref{eqn:bijective}),  we have
$$
\widehat{ E}_1 \cap \widehat{ \pi}^{-1}(V)= \mbox{  $\bigcup
_{j=1}^\infty$} \,
U_{k_j},
$$
where the disjoint union $\{U_{k_j}\}_j$ is a subsequence of $\{U_k\}_k$ in
(\ref{eqn:bijective}).
To verify (\ref{eqn:restar}) it suffices to show that $E'(z_0)\cap
H_z\subset D'(z_0)$. Fix $h\in E'(z_0)\cap H_z$. To show that $h\in D'(z_0)$, fix a continuous curve
$\tau: t\in [0,1]\to \tau(t)$ in $E'(z_0)$ with $\tau(0)=e$ and
$\tau(1)=h$. Using  (\ref{eqn:d-1-prime}) for the set $E$, the curve $\gamma:=\tau(w_0)$
defined by $t\in
[0,1] \to w(t):=\tau(t)(w_0)$ is a continuous curve in $\widehat{ E}_1$
with $w(0)=w_0$ and $w(1)= h(w_0)$. We can take a partition
$0=t_1<t_2 < \ldots <t_{2\mu-1}<t_{2\mu}=1$ of $[0,1]$ such that
\begin{eqnarray*}
\left\{
 \begin{array}{llll}
w(t)\in \widetilde { D}_1   & {\rm  on} & [t_{2i-1}, t_{2i}],  \quad
 i=1,\ldots ,\mu,\\
w(t)\in U_{k_{j_i}} & {\rm  on} & [t_{2i},t_{2i+1}],  \quad  i=1,\ldots ,\mu-1,
 \end{array}
\right.
\end{eqnarray*}
where each set $U_{k_{j_i}}$ is one of the sets in the disjoint union
$\{U_{k_j}\}_{j=1,2,\ldots} $.
For simplicity in notation we write $I_i:=[t_{2i}, t_{2i+1}]$ and $U_i:=U_{k_{j_i}}$ for
$i=1,\ldots ,\mu-1$. For $ i=1,\ldots , \mu-1$,  if we
  set $\gamma_{2i}:\, t\in I_i\to
w(t)$, then $\gamma_{2i}$ is a continuous curve in $U_{i}$
 with $ w(t_{2i}), w(t_{2i+1})\in \widetilde { D}_1$.
 For each $i$, we
will construct a continuous curve $\widetilde { \tau}_{2i}:t\in I_{i}\to
\widetilde { \tau}_{2i}(t)$ in $G$  such that
$\widetilde { \tau}_{2i}(t)(w_0)\in \widetilde { D}_1$ in $M'$ and  $ \widetilde { \tau}_{2i}(t_{2i})= \tau(t_{2i}), \
\widetilde { \tau}_{2i}(t_{2i+1} )= \tau(t_{2i+1})$.

We first consider a continuous curve $ \widetilde { \gamma}_{2i}:
t\in I_i\to \widetilde { w}(t)$ in $\widetilde { D}_1 \cap U_i$ with $ \widetilde { w}(t_{2i})= w(t_{2i}), \
\widetilde {w}(t_{2i+1} )= w(t_{2i+1})$. From the proof of Lemma
\ref{eqn:saigono},  there is a holomorphic section $\sigma_i: w\in U_i \to
 \sigma_i(w) $ of $G$ over $U_i$ via $\widehat{ \psi}_{w_0}$, i.e.,
$\sigma_i(w)(w_0)= w, \ w\in U_i$. Thus, for $t\in I_i$,
$
\sigma_i(w(t))(w_0)= w(t)= \tau(t)(w_0)$  in $ M'$.
Hence we can find a continuous curve $t\in I_i\to h(t)$ in
$\widehat{ H}_{w_0}$ with
$\tau(t)= \sigma_i(w(t))h(t), \ t\in I_i$.
Next we consider the continuous curve
$$
\widetilde { \tau}_{2i}: t\in I_i\to \widetilde { \tau}_{2i}(t):=
\sigma_i(\widetilde { w}(t)) h(t)
$$
 in $G$. Then we have
$\widetilde {\tau}_{2i}(t)(w_0)=\sigma_i(\widetilde {
w}(t))(w_0)=\widetilde { w}(t) \in \widetilde { D}_1,
$
 and
$$
\widetilde { \tau}_{2i}(t_{2i})= \sigma_i(\widetilde { w}(t_{2i}))
h(t_{2i})=
 \sigma_i({ w}(t_{2i}))
h(t_{2i})=\tau(t_{2i});
$$
similarly $\widetilde { \tau}_{2i} (t_{2i+1})=\tau(t_{2i+1})$. Thus
$\widetilde { \tau}_{2i}$ satisfies the desired properties.

Now for $i=1,\ldots , \mu$, let $\tau_{2i-1} $ denote the restriction of $\tau$ to the subinterval
 $[t_{2i-1}, t_{2i}]$. Consider the concatenized curve
$$
T: = \tau_1+\widetilde {\tau}_2+ \tau_3 + \cdots + \widetilde {
 \tau}_{2\mu-2} +\tau _{2\mu-1}
$$
in $G$. Clearly $T$ is continuous on $[0,1]$ and $T(t)(w_0)\in  \widetilde { D}_1
 $.
 Hence  $T(t)(z_0)\in D, \ t\in [0,1] $.  Since $ T(0)=e$, it follows
 from the
definition of $D'(z_0)$ that
 $T(t)\in D'(z_0)$. In particular, $h=\tau(1)=T(1)\in D'(z_0)$, which proves
 (i).

To prove (ii), note first that it is clear that the left-hand-side is contained in the right-hand-side. To prove the reverse inclusion, let
$w\in \widehat{ \pi}^{-1}(D) \cap \widehat{ E}_1$, so that
$z:=\widehat{\pi}(w)\in D$, and take $g\in E'(z_0)= \widehat{ E}_1(w_0)$ with
$g(w_0)=w$. Thus $\widehat{ \pi}(g(w_0))=\widehat{ \pi}(w)$; i.e., we have $g(z_0)=z$ in $M$. Using 1. in Proposition
\ref{eqn:elm} we can find $g'\in D'(z_0)$  such that $g'(z_0)=z$. Thus there exists $h\in H_{z_0}$ with $g=g'h$ in $G$. Since $g\in E'(z_0)$ and
$g'\in D'(z_0) \subset E'(z_0)$, we apply 3. in Proposition
\ref{eqn:elm} to obtain $h\in E'(z_0)\cap H_{z_0}$, so that $h\in
D'(z_0)\cap H_{z_0}$ by (i). It follows from 2. in Proposition
\ref{eqn:elm} that $g=g'h\in D'(z_0)$. Hence $w=g(w_0)\in \widehat{
D}_1$ by (\ref{eqn:d-1-prime}), which proves the reverse inclusion in (ii).

For (iii), again, that the left-hand-side is contained in the right-hand-side  is clear. To prove the reverse inclusion, let $g\in E'(z_0)\cap D(z_0)$. By (\ref{eqn:d-1-prime}) we
have $w:=g(w_0)\in
\widehat{ E}_1$. We  set $z:=g(z_0)\in D$. Since $ \widehat{ \pi}(w)
=
g(z_0)=z$,
it follows that $w\in \widehat{
E}_1 \cap \widehat{ \pi}^{-1}(z)$, so that $w\in \widehat{ D}_1$ by
(ii). Again using (\ref{eqn:d-1-prime}) we can find $g'\in D'(z_0)$  such that
$w=g'(w_0)$, and hence there exists $h\in H'_{z_0}$ with
$g=g'h$. Therefore,
2. in Proposition \ref{eqn:elm} implies that $g\in D'(z_0)$, which proves
the reverse inclusion in (iii).
\hfill $\Box$

\smallskip
For the proof of Theorem \ref{thm:main-nonconnected} we next study
 the subset $\Sigma_{z_0} {\cal H}'(z_0)$ of $G$ defined in
 (\ref{eqn:ahprime}).

\begin{lemme}${}$
 \label{liegroup-hz}
 \begin{enumerate}
  \item If $h\in {\cal H}'(z_0)$, then $h  ^{-1} \Sigma_{z_0} h =
  \Sigma_{z_0}$.
 \item    The set $\Sigma_{z_0} {\cal  H}'(z_0) $ is a Lie subgroup
of $G$  such that
 \begin{enumerate}
  \item [(i)] $\Sigma_{z_0}$ is a normal Lie subgroup of $\Sigma_{z_0}
  {\cal H}'(z_0)$;
  \item [(ii)] $ \Sigma_{z_0} {\cal H}'(z_0)= {\cal H}'(z_0)
 \Sigma_{z_0}$ \ and \ $D'(z_0)\Sigma_{z_0} {\cal H}'(z_0) = D'(z_0) $;
\item [(iii)]
 using the
 notation $h_j', j=1,2,\ldots $ in
(\ref{eqn:H-prime}), we have
\begin{eqnarray*}
 \Sigma_{z_0} {\cal H}'(z_0) =\mbox{  $ \bigcup _{j=1}^\infty$} \,  h_j' \Sigma_{z_0}=
 \mbox{  $\bigcup _{j=1}^\infty $}\,  \Sigma_{z_0}h_j'
\end{eqnarray*}
and these are disjoint unions;
\item [(iv)]  $\Sigma_{z_0} {\cal H}'(z_0) = D'(z_0)\cap
    \,Aut\,\sigma_{z_0} $, which acts transitively on
    $\sigma_{z_0}$; and $\sigma_{z_0}$ is isomorphic to $\Sigma_{z_0} {\cal
    H}'(z_0)/
 {\cal H}'(z_0)$ (as complex manifolds).
\end{enumerate}
\item   If $\sigma_{z_0}$ is closed in $M$, then $\Sigma_{z_0} {\cal H}'(z_0)$ is a closed
 Lie subgroup of $G$.
 \end{enumerate}
\end{lemme}
\noindent{\bf  Proof}.
We prove 1. and  2. in the lemma in steps $(1) \sim (6)$.

\smallskip
 (1) \  {\it Each of $\Sigma_{z_0} {\cal H}'({z_0})$ and ${\cal
H}'({z_0})\Sigma_{z_0}$ are contained in $ D'({z_0})$.}

\smallskip
To see this, let $s\in \Sigma_{z_0}$ and $h\in {\cal H}'({z_0})$. To prove that $sh\in
D'({z_0})$ we take a continuous curve $\gamma: t\in [0,1] \to s(t)$ in
$\Sigma_{z_0}$ with $s(0)=e$ and $s(1)=s$ (note that $\Sigma_{z_0}$ is connected
and contains $e$). Then $\widetilde { \gamma}: t\in [0,1] \to \widetilde
{ \gamma}(t)=s(t) h$ is a continuous curve in $G$. We have
$$
\widetilde { \gamma}(t)({z_0})=s(t)h({z_0})=s(t)({z_0})\in \Sigma_{z_0}({z_0})=\sigma_{z_0} \Subset
D,  \quad  t\in [0,1].
$$
Since $\widetilde { \gamma}(0)=h \in {\cal H}'({z_0}) \subset D'({z_0})$, it follows that
$\widetilde { \gamma}({z_0}) \subset D'({z_0})$. Hence $sh=\widetilde {
\gamma}(1)\in D'({z_0})$, so that $\Sigma_{z_0} {\cal H}'({z_0}) \subset D'({z_0})$.

To verify the inclusion ${\cal H}'({z_0}) \Sigma_{z_0}
\subset D'({z_0})$, let $h\in {\cal H}'({z_0})$ and $ s\in \Sigma_{z_0}$. We
take the same continuous curve $\gamma: t\in [0,1] \to s(t)$ in
$\Sigma_{z_0}$ with $s(0)=e$ and $s(1)=s$; then $\widehat{  \gamma}: t\in [0,1] \to \widehat
{ \gamma}(t)=hs(t) $ is a continuous curve in $G$.We then have
$$
\widehat { \gamma}(t)({z_0})=hs(t)({z_0})\in h\sigma_{z_0} \Subset D,
\quad  t\in [0,1].
$$
The last
equality comes from $h\in D'(z_0)$ and formula (\ref{eqn:bounda}).
 Since $\widehat { \gamma}(0)= h\in D'({z_0})$, we
have $hs=\widehat { \gamma}(1)\in D'({z_0})$. Thus (1) is proved.

\smallskip
 (2) \  {\it  If  $h\in {\cal H}'({z_0}) $, then
$ h  ^{-1}  \Sigma_{z_0} h= \Sigma_{z_0} $. In particular,
        ${\cal H}'({z_0}) \Sigma_{z_0} =\Sigma_{z_0} {\cal H}'({z_0})$.}

\smallskip
Let $\Sigma':=h  ^{-1}  \Sigma_{z_0} h$. Then $\Sigma_{z_0} $
and $\Sigma'$ are Lie subgroups of $G$ of the same dimension.
 Since $H_{z_0}'$ is normal in $H_{z_0}$, $\Sigma'$ and $\Sigma_{z_0} $
 each contain $H_{z_0}'$, which is of dimension $m_0$. We write
 $$\mathfrak h _0=[X_1,
 \ldots ,X_{m_0}] \ \  \mbox{  and} \ \  \mathfrak X'= [X_1,\ldots ,X_{m_0}, X_{m_0+1},
 \ldots, X_{m_0+m_1} ]$$
for the Lie algebras of $H_{z_0}'$ and  $\Sigma'$. Thus
$$
\Sigma'=\{\, {\textstyle  \prod}_{i=1}^\nu \exp t_iX_i\in G: \nu\in {\boldsymbol{Z}}^+,\ t_i\in {\mathbb{C}},\
X_i\in \mathfrak X '\,\}.
$$
To verify (2), it suffices to prove that $\mathfrak X ' \subset \mathfrak X
 _{{z_0}}$, or equivalently,
$$
X_j \in \mathfrak X_{z_0}, \ \ j=m_0+1,\ldots , m_0+m_1.
$$
Fixing such a $j$, we consider the one-dimensional curve
$C:=\{\exp tX_j ({z_0})\in
M: t\in{\mathbb{C}}\}$
in $M$. Since $\exp tX_j \in \Sigma'$ and $h  ^{-1}  \in {\cal H}({z_0}')\subset D'({z_0})$,  we have
$$
C \subset \Sigma'({z_0})= h  ^{-1}  \Sigma_{z_0} h ({z_0}) = h  ^{-1}  \Sigma_{z_0} ({z_0})=h  ^{-1}
\sigma_{z_0} \Subset D.
$$
It follows that $\lambda= \  const. = \lambda({z_0})$ on $C$, and hence  $X_j\in
 \mathfrak X _{z_0}$ from the definition of $\mathfrak X_{z_0}$. This proves  (2).

\smallskip
(3) \  {\it   $D'(z_0)\Sigma_{z_0} {\cal H}'({z_0}) =
 D'(z_0)$.}

\smallskip
 The inclusion $D'(z_0)\subset D'(z_0)\Sigma_{z_0} {\cal H}'({z_0})$ is clear; we prove the reverse
 inclusion.  Let $g\in D'(z_0), s\in \Sigma_{z_0} , h\in {\cal H}'(z_0)$.
 Since $\Sigma_{z_0}$ is connected and contains $e$ and $s$, we
can find a continuous curve $l:t\in [0,1]\to l(t):
=gs(t)$ in $G$ with $s(t)\in \Sigma_{z_0} , \ t\in [0,1]$; $l(0)= g$
and $l(1)=gs$. Since $l(0)=g\in D'(z_0)$, we have  $l \subset D'(z_0)$,
and hence  $gs=l(1)\in D'(z_0)$. Since $h\in {\cal H}'(z_0)=D'(z_0)\cap
H_{z_0}$,  it follows from  2. in Proposition \ref{eqn:elm} that $gsh\in
D'(z_0)$, as required.

\smallskip
 (4) \  {\it   We have 2. (iii)  in the lemma.}

\smallskip
From (\ref{eqn:H-prime}) we have
\begin{eqnarray*}
\Sigma_{z_0}  {\cal H}'({z_0}) &=& \Sigma_{z_0}\,\mbox{
 $\bigcup_{j=1}^\infty$}\,
 h_j' H_{z_0}'\\
&=&\mbox{  $\bigcup_{j=1}^\infty$}\,  \Sigma_{z_0} h_j' H_{z_0}'\\
&=& \mbox{  $\bigcup_{j=1}^\infty$} \,  h_j'\Sigma_{z_0} H_{z_0}'  \quad  \mbox{    (by \ (2))} \\
 &=& \mbox{  $\bigcup_{j=1}^\infty$}\,  h_j'\Sigma_{z_0}   \quad  \quad   \mbox{ \ \  (since   $
  H_{z_0}' \subset \Sigma_{z_0} $)}.
\end{eqnarray*}
To prove that this is a disjoint union, assume that $h_j' \Sigma_{z_0} \cap h_k' \Sigma_{z_0}
 \ne \emptyset$. Then $h_j's_j= h_k's_k$ for some
 $s_j,s_k\in \Sigma_{z_0} $. If
 we define $\alpha:=( h_k')  ^{-1}  h_j'= s_k (s_j) ^{-1}  $, then
$\alpha\in {\cal H}'({z_0}) \cap \Sigma_{z_0} $, since both ${\cal H}'({z_0}) $ and
$\Sigma_{z_0} $
 are subgroups in $G$. It follows that
$
h_j' \Sigma_{z_0} =h_k' \alpha \Sigma_{z_0}\in h_k' \Sigma_{z_0},
$
and hence $h_j' \Sigma_{z_0} = h_k' \Sigma_{z_0} $, and (4) is proved.

\smallskip
 (5) \   {\it   The set $\Sigma_{z_0} {\cal H}'({z_0})$ is a Lie
 subgroup
  of $G$, and
$\Sigma_{z_0}$ is a normal Lie subgroup of $\Sigma_{z_0} {\cal
       H}'(z_0) $.}

\smallskip
 By (4) the set $\Sigma_{z_0} {\cal H}'({z_0}) $ is a disjoint union of
 $(m_0+m_1)$-dimensional nonsingular $f$-generalized analytic sets which are equivalent
 to $\Sigma_{z_0}$. Thus we need only verify that $\Sigma_{z_0} {\cal H}'({z_0})$
is a subgroup of $G$. Since $e\in \Sigma_{z_0} {\cal H}'({z_0})$, it suffices
to prove that
for any $s_1,s_2\in \Sigma_{z_0}$ and $h_1,h_2\in {\cal H}'({z_0}) $, we have
$x:= (s_1h_1)(s_2h_2)  ^{-1}\in \Sigma_{z_0} {\cal H}'({z_0})    $. To see this,
\begin{eqnarray*}
 x&=& s_1h_1 h_2  ^{-1}  s_2  ^{-1}  \\
&\in& \Sigma_{z_0} {\cal H}'({z_0}) \Sigma_{z_0}  \quad  \mbox{ (since ${\cal H}'({z_0}) $
 is a group)} \\
&\subset& \Sigma _{z_0} \Sigma_{z_0} {\cal H}'({z_0})  \quad \   \mbox{(by (2))}\\
& =& \Sigma_{z_0} {\cal H}'({z_0})  \quad  \mbox{ \  \ \ \  (since $\Sigma_{z_0} $ is a group)},
\end{eqnarray*}
so that $\Sigma_{z_0} {\cal H}'({z_0}) $ is a Lie subgroup of $G$. This
with (2) implies that $\Sigma_{z_0} $ is a normal Lie subgroup of
$\Sigma_{z_0} {\cal H}'({z_0}) $.

\smallskip
 (6) \   {\it We have 2. (iv) in the lemma.}

\smallskip
 For the first equality in (iv) it suffices to show
that
$$\Sigma_{z_0} {\cal H}'(z_0) = D'(z_0)
       \cap\,Aut\,\sigma_{z_0} .$$
 Let $g=sh\in \Sigma_{z_0} {\cal H}'(z_0)
       $. Then $g \sigma_{z_0}= sh \sigma_{z_0} = sh \Sigma_{z_0} (z_0)=
       s\Sigma_{z_0}h(z_0)= \Sigma_{z_0} (z_0)=\sigma_{z_0}  $. This together with (3) proves $\Sigma_{z_0} {\cal H}'(z_0) =\subset D'(z_0)  \cap\,Aut\,\sigma_{z_0} $.  For the reverse inclusion, let $g\in
       D'(z_0)$  with $g \sigma_{z_0} =\sigma_{z_0} $. Then $g(z_0)= s(z_0)$ for some $s\in \Sigma_{z_0}$, so that $g=sh$ for
       some $h\in H_{z_0}$. Since $g,s\in D'(z_0)$, we have $h\in {\cal
       H}'(z_0) $ by  3. in Proposition \ref{eqn:elm}. Hence $g\in
       \Sigma_{z_0} {\cal H}'(z_0) $, as claimed.
       To prove the transitivity of $\Sigma_{z_0} {\cal H}'(z_0) $ on $\sigma_{z_0}$, let
       $\zeta\in \sigma_{z_0}$; there exists $s\in \Sigma_{z_0} $ with $\zeta=s(z_0)$. Since $\sigma _{z_0} \subset D$, we have $g\in D'(z_0)$ with $g(z_0)=\zeta $, hence we can find $h\in H_{z_0} $ with
       $g=sh$; again, this implies  $h\in {\cal H}'(z_0)$.

       Finally to show $ \Sigma_{z_0} {\cal H}'(z_0)
       /{\cal H}'(z_0) \approx \sigma_{z_0}$, it suffices to prove that the isotropy
       subgroup $\{g\in \Sigma_{z_0} {\cal H}'(z_0): g(z_0)=z_0\}$ of
       $\Sigma_{z_0} {\cal H}'(z_0) $ for $z_0\in \sigma_{z_0} $ is
       equal to ${\cal H}'(z_0)$. This follows from the definition
${\cal H}'(z_0)=H_{z_0} \cap D'(z_0)$ and (3).

\smallskip
This completes the proofs of $(1) \sim (6)$, and hence
1. and 2. in the lemma are proved. To prove 3. we assume that $\sigma_z$ is closed in $M$.
 Let $s_n h_n \in \Sigma_z {\cal H}'(z), \ n=1,2,\ldots  $, and take $g\in
G$ with $s_n h_n \to g $ as $n \to \infty$. Since
$
s_nh_n(z)= s_n(z)\in \Sigma_z (z)= \sigma_z \Subset D,
$
and since $\sigma_z$ is assumed to be closed in $G$, the limit
$g(z)=\lim_{ n\to \infty} s_nh_n(z)$ is contained in $\sigma_z$. Here, if
necessary, we take a subsequence of $\{s_nh_n(z)\}_n$. Thus there exists $s^*\in \Sigma_z$ with $g(z)= s^*(z)$. By definition of the isotropy subgroup $H_z$, there exsists $h^*\in H_z$  such that
$g= s^*h^*$. Since $g(z)\in D$, we have  $g\not\in
\partial[D(z)]$. On the other hand, by (1) we have $s_nh_n\in D'(z)$, so that $g\in D'(z)\cup \partial
 [D'(z)]$. Using $\partial [D'(z)] \subset \partial [D(z)]$, we have
 $g\in D'(z)$. Moreover, we have $s^*\in \Sigma_z\subset D'(z)$. Applying 3. in Proposition 6.1 to $g= s^*h^*$, we have $h^*\in {\cal H}'(z) $. Hence
$g \in \Sigma_z {\cal H}'(z)$, and 3. in the lemma is proved.
 \hfill $\Box$

\smallskip
From 2. (ii) in Lemma \ref{liegroup-hz} we can show that the union  stated  in  (\ref{eqn:bounda}):
\begin{align}
 \label{eqn:lami}
\mbox{ $ D={\textstyle  \bigcup_{g\in {
D}}\, g\sigma_{z_0}} $  is a disjoint union.}
\end{align}
\noindent To see this, take
$s_1,s_2\in \Sigma_{z_0}$ with $g_1s_1(z_0)= g_2s_2(z_0)$. Then  we can find $h\in H_{z_0}$ with $g_1 s_1=g_2 s_2 h$. From $2$. (ii) in Lemma  \ref{eqn:elm} we have $g_1s_1, g_2s_2\in D'(z_0)$; then from 3. in Proposition \ref{eqn:elm} we have
 $h \in {\cal H}'(z_0) $. Next, from  $1$. in
 Lemma \ref{liegroup-hz} we have $g_1s_1= g_2hs_2'$ for some $s_2'\in
 \Sigma_{z_0}$. Since $\Sigma_{z_0}$ is a group, it follows from
$1$. in
 Lemma \ref{liegroup-hz} that
\begin{align*}
 g_1\sigma_{z_0}&= g_1 \Sigma_{z_0} (z_0)= g_1 s_1 \Sigma_{z_0} (z_0)
=
g_2hs_2' \Sigma_{z_0} (z_0)\\&= g_2 h \Sigma_{z_0} (z_0)=g_2 \Sigma_{z_0}
 h(z_0)
= g_2 \Sigma_{z_0} (z_0)= g_2 \sigma_{z_0},
\end{align*}
as claimed.

\smallskip Recall the notation $\widehat{ \lambda}(w)$ and $ \widehat{
\Sigma} _{w_0} $ in  (\ref{defofwhl}) and (\ref{eqn:z-w-relation}):
\begin{align*}
 \widehat{ \lambda }(w) &=\lambda (\widehat{\pi}(w)),  \quad  w\in
\widehat{ D}_1;\\
\widehat{ \mathfrak X }_{w_0}&= \{ X\in \mathfrak X :
[\frac{\partial^2  \widehat{ \lambda} (\exp tX(w_0))}{\partial t \partial
 \overline {
    t }}] |_{t=0}=0\};\\
\widehat{\Sigma}_{w_0}&= \{\,{\textstyle  \prod}_{i=1}^\nu \exp t_iX_i\in G:\nu\in {\boldsymbol{Z}}^+,
 t_i\in {\mathbb{C}}, X_i\in \widehat{ \mathfrak X }_{w_0}\,\}.
\end{align*}

\begin{lemme}\label{lem:widehatXw} \
With the above notation,
 ${}$
\begin{enumerate}
 \item [(1)]\ $ \widehat{ \mathfrak X }_{w_0} = \mathfrak X _{z_0}$ and
 hence $\Sigma_{z_0}=\widehat\Sigma_{w_0}$;
 \item [(2)]\
$\widehat{ \mathfrak X }_{w_0}=\{ X\in \mathfrak X : \exp tX(w_0) \in
 \widehat{ D}_1, \  \widehat{\lambda}(
\exp tX(w_0))= \widehat{ \lambda }(w_0), \  t\in {\mathbb{C}} \}$.
\end{enumerate}
\end{lemme}
\noindent {\bf  Proof}.  If $|t| \ll 1$, then $\exp tX(w_0) \in
\widehat{ D}_1 $ and $\exp tX(z_0)\in D$. From the definition of
$\widehat{\lambda}$ and $\widehat{\pi}$ and (\ref{eqn:commutative-wg}), we have
$\widehat{\lambda}(\exp tX(w_0))= \lambda(\exp t X(z_0))$. Thus (1) follows
from the definition of the sets $\widehat{ \mathfrak X }_{w_0}$ and $\mathfrak X _{z_0}$. We write $\widetilde
{ \mathfrak X }_{w_0}$ for the right-hand side of (2). Let $X \in \widehat{  \mathfrak X }_{w_0}$. By (1) we have $X\in \mathfrak X
 _{z_0}$, so that, in particular, $X$ satisfies condition i. in the definition of
 $\mathfrak X _{z_0}$ in (\ref{eqn:X-0}). It follows from the definitions of
 $\widehat{ D}_1 $ and $\widehat{ \lambda}$ that $X\in \widetilde {
 \mathfrak X }_{w_0}$, so that $\widehat
 { \mathfrak X }_{w_0}\subset \widetilde{ \mathfrak X }_{w_0}$. The reverse inclusion is clear, and (2) is
 proved.
 \hfill $\Box$

\smallskip
Using standard arguments this lemma implies the following results.
 From (1)
 in Lemma \ref{lem:widehatXw},  $\widehat{ \mathfrak X }_{w_0}$ is a Lie subalgebra
of $\mathfrak X $  with  $\widehat{ H}_{w_0}=H_{z_0}' \subsetneq  \mathfrak X _{z_0} =\widehat{ \mathfrak X }_{w_0}$.
 Now we let $\widehat{ \mathfrak h} _{w_0}$ denote the
Lie subalgebra of $\mathfrak X $ which corresponds to the closed Lie
subgroup $ \widehat{H} _{w_0}$ of $G$. Note that  $\dim
\widehat{ \mathfrak h}_{w_0}=m_0= \dim H_{z_0}$, and let
$\dim \widehat{ \mathfrak X
}_{w_0}=m_0+m_1 = \dim \mathfrak X _{z_0}$.
Recall we used the notation $\widehat \Sigma_{w_0}$ for the  connected
Lie subgroup of $G$  which
 corresponds to $\widehat{ \mathfrak X} _{w_0}$, even though this is the same as
 $\Sigma_{z_0}$.
From the Frobenius theorem, $\widehat{
 \Sigma}_{w_0}$  is an irreducible, $(m_0+m_1)$-dimensional non-singular
 $f$-generalized analytic set in $G$. Define
\begin{align}
 \label{eqn:wsw0}
\widehat{ \sigma}_{w_0}:=
\widehat{ \psi}_{w_0}( \widehat{ \Sigma}_{w_0}) = \widehat{ \Sigma}_{w_0}(w_0)=
\{g(w_0)\in M' : \ g\in \widehat{ \Sigma}_{w_0}\} \subset M'.
\end{align}
 Since $\widehat{ H}_{w_0}$ is a closed connected Lie subgroup of $G$,
 it follows from the inclusion $\widehat{H}_{w_0} \subset \widehat{ \Sigma}_{w_0}$
 that $\widehat{
 \sigma}_{w_0}$ is an irreducible,
$m_1$-dimensional non-singular $f$-generalized analytic set in $M'$ passing
through $w_0$.
Now to each tuple
$$\bigl\{
G;\ M=G/H_{z_0};\ D \Subset M; \ z_0\in D; \ \mathfrak
X_{z_0} \subset \mathfrak X;\ \Sigma_{z_0} \subset D(z_0) \subset  G;\ \sigma_{z_0}
\Subset D\bigr\}$$ from Theorem \ref{thm:mainth}  where $H_{z_0}$ is connected in $G$ (so that $H_{z_0} \subset \Sigma_{z_0})$ we assign a corresponding tuple
$$\bigl\{G;\ M'=G/ \widehat{ H}_{w_0};\
\widehat{ D}_1 \subset M'; \ w_0\in \widehat{ D}_1; \ \widehat {\mathfrak
X}_{w_0} \subset \mathfrak X; \ \widehat{ \Sigma}_{w_0}\subset D'(z_0) \subset G; \ \widehat\sigma_{w_0} \subset \widehat{
D}_1\bigr\}$$
from Theorem \ref{thm:main-nonconnected}
where $\widehat{ H}_{w_0}$ is connected in $G$ (so  that $\widehat{ H}_{w_0} \subset \widehat{ \Sigma}_{w_0})$. Note that, in general, $\widehat\sigma_{w_0} \not\Subset \widehat{D}_1$. Then by the same  argument used to obtain
 the foliations of $M$ and $D$ in (\ref{eqn:sliced}), we  obtain the
 following foliations of $M'$ and $\widehat{ D}_1$
\begin{align}
 \label{sliced-2}
M'=\mbox{  $\bigcup_{g\in G}$} \,  g\, \widehat{ \sigma}_{w_0} \quad
 \mbox{  and } \quad
 \widehat{ D}_1= \mbox{  $\bigcup _{g\in
 D'(z_0)}$} \,  g\, \widehat{ \sigma}_{w_0},
\end{align}
despite the fact that $\widehat{ D}_1$
 is not necessarily relatively compact in $M'$. Therefore, given $w\in M'(\widehat{ D}_1)$, we have a unique leaf $g \widehat{
\sigma}_{w_0} \subset M'\,(\widehat{ D}_1)$
which passes through $w$,. Here $g\in G(D'(z_0))$ is uniquely determined up to
$\widehat{ H}_{w_0}$. i.e., if another leaf $g' \widehat{ \sigma}_{w_0}$ passes through
$w$, then $g'=gh$ for some $h\in H'(z_0)=\widehat{ H}_{w_0}$.

\smallskip
With these preliminaries, we commence with the following lemma related
to the representation of $D$ and $\partial D$ in (\ref{eqn:bounda}):
$D= \bigcup_{g\in D'(z_0)} g\sigma_{z_0}$ and $\partial D= \bigcup _{g\in \partial
 D'(z_0)}g\sigma_{z_0}$.
\begin{lemme}
 \label{lem:s0nonsing} \
\begin{enumerate}
 \item The set  $\sigma_{z_0}$ in $M$ is an $m_1$-dimensional non-singular
$f$-generalized analytic set in $M$ with $\sigma_{z_0} \Subset D$.
 \item The union $ { \bigcup} _{g\in D'(z_0)} \, g
       \sigma_{z_0}$ is a  foliation.
\end{enumerate}
\end{lemme}
\noindent {\bf  Proof}. We verified that $\sigma_{z_0} \Subset D$ in
(\ref{eqn:bounda}). Since  1. follows from 2., we prove 2.
First we verify the following elementary equality:
\begin{align}
 \label{eqn:hkprimehat}
\hbox{for any} \ g\in D'(z_0), \ \  \widehat{ \pi} ^{-1}  (g\sigma_{z_0})\cap \widehat{ D}_1=
 \mbox{  $\bigcup_{k=1}^\infty$} \, g h_k'
 \widehat{ \sigma}_{w_0}.
\end{align}

Indeed, letting  $h_k'\in {\cal H}'(z_0)$, we have $h_k' \widehat{ \Sigma}_{w_0}=h_k'
\Sigma_{z_0} \subset D'(z_0)= \widehat{
D}_1(w_0)$ from $2$ (ii). in Lemma
\ref{liegroup-hz} so that $h_k' \widehat{ \sigma}_{w_0}=h_k'\widehat{ \Sigma}_{w_0}(w_0) \subset \widehat{
D}_1$. Using  1.
 in Lemma \ref{liegroup-hz} we obtain
\begin{align*}
 \widehat{\pi}( gh_k' \widehat{ \sigma}_{w_0})=\widehat{ \pi}(g h_k'
\widehat{  \Sigma}_{w_0} (w_0))=g h_k' \Sigma_{z_0}(z_0)
= g \Sigma_{z_0} h_k'(z_0)=g \sigma_{z_0}.
\end{align*}
This proves that $\mbox{  $\bigcup_{k=1}^\infty$} \, g h_k'
 \widehat{ \sigma}_{w_0}\subset \widehat{ \pi} ^{-1}  (g\sigma_{z_0})\cap \widehat{ D}_1$. Conversely, let $w'\in \widehat{
\pi}  ^{-1}  (g\sigma_{z_0}) \cap \widehat{ D}_1$. Using
(\ref{eqn:d-1-prime}) we can find $g'\in D'(z_0) $ with
 $g'(w_0) =w'$. We thus have $g'(z_0)= \widehat{
 \pi}(g'(w_0))=\widehat{\pi}( w')\in
g\sigma_{z_0}=g \Sigma_{z_0}(z_0)  $; hence
we can find $s\in \Sigma_{z_0}$ and $h\in H_{z_0}$
with $g'=gsh$ in $G$. From $2.$ (ii) in Lemma \ref{liegroup-hz}, $gs$ and $g'$ are contained in
$D'(z_0)$. It follows from 3. in Propositon \ref{eqn:elm} that $h\in {\cal
H}'(z_0)$; hence there exists $h_k'$ with $h\in h_k'\widehat{
H}_{w_0}$.  Using 1. in Lemma \ref{liegroup-hz}, we have
\begin{align*}
 w'=gsh(w_0)\in g \Sigma_{z_0} h(w_0)= g \Sigma_{z_0}h_k'(w_0)= g
h_k'\Sigma_{z_0}(w_0) =g h_k' \widehat{ \sigma}_{w_0},
\end{align*}
proving the reverse inclusion.

Finally, we note that $\widehat{ \pi} : M \to M'$ is a normal covering
and
$\widehat{ D}_1=\bigcup_{g\in D'(z_0)}g \widehat{ \sigma}_{w_0}$ is
a foliation of leaves of  $m_1$-dimensional nonsingular
generalized analytic sets $g \widehat{ \sigma}_{w_0}$.
 It follows from formula (\ref {eqn:hkprimehat})   that $D= \bigcup _{g\in D'(z_0)}g \sigma_{z_0}$
is a foliation of leaves of  $m_1$-dimensional non-singular
$f$-generalized analytic sets $g{ \sigma}_{z_0}$; hence
Lemma \ref{lem:s0nonsing} is proved.  \hfill $\Box$

\smallskip
We construct a domain $K$ as in $2-a$ (i) of Theorem
\ref{thm:main-nonconnected}. in the following lemma.
\begin{lemme} \label{lem:saigo}
 \ Let $D \Subset M $ and fix $z_0 \in D $ as above.
 Let $E$ be a domain with smooth boundary in $M$  such that $D \Subset E
 \subset M$  and  $D$ and $E$ are of the
 same isotropy class at $z_0$ in $D$ (such $E$ exist by Lemma
 \ref{lem:d1d2}).  Then there exists a domain
$K$ in $M$
with $D \Subset K \subset E$
 such that $K$ is foliated by $m_1$-dimensional non-singular
 $f$-generalized analytic sets $g\sigma_{z_0}$ in $M$:
\begin{align}
 \quad  {K}&= \mbox{  $\bigcup_{g\in K'(z_0)} $} \,   g \sigma_{z_0}  \quad    \mbox{  with }
 \ g\sigma_{z_0} \Subset K. \label{saig-saigo-2}
\end{align}
 \end{lemme}
\noindent
{\bf  Proof}.
Fix a neighborhood  ${\cal V}$ of $e$
in $G$  such that $
f( \overline { D}) \Subset E  $ for any $ f\in {\cal V}$.
 We apply Lemma \ref{eqn:saigono} for $  D$ and ${\cal V}$ to
 obtain balls $V ^{(i)}, \,i=1,\ldots ,N,$ centered at
$z_i\in  { D}$ in $M$ and $U_k^{(i)}$ in $M'$ which satisfy conditions
(1) $\sim$ (3)
in the lemma as well as the additional condition that $V ^{(i)}
\Subset E, \ i=1,\ldots N$.
We set $\widehat{ \pi}^{-1}(z_i)= \{w_k
^{(i)}\}_{k=1,2,\ldots } \in \widehat{\pi}^{-1}({ D}) \subset
M'$; i.e., $w_k ^{(i)}$ is the center of $U_k^{(i)}$.
 We may assume that $
A:= \cup_{i=1}^N V ^{(i)}
$
 is a connected domain in $M$ with $ D \Subset A \Subset E$.
Since $E$ and $D$
are of the same isotropy class at $z_0$, so are $A$ and  $D$; i.e.,
 $D'(z_0)\cap
 H_{z_0}= A'(z_0)\cap H_{z_0}=E'(z_0)\cap H_{z_0}$.

For each $i=1,\ldots ,N$, we  consider the sets $U_k ^{(i)} \subset
M' $ defined in (2) in Lemma \ref{eqn:saigono}. We only consider the sets $U_k ^{(i)}$  in $\{U_k ^{(i)}
 \}_{k=1,2,\ldots }$ which are
contained in $\widehat{ E}_1$, the connected component of
$\widehat{ \pi}^{-1}(E)$ containing $w_0$ in $M'$. By convention,
we use the same notation $\{ U_k ^{(i)} \}_{i,k}$ for these sets. From (ii) in Lemma \ref{lem:d1d2},
$\widehat{ A}_1= \widehat{ E}_1 \cap \widehat{\pi} ^{-1}(A)$; thus the connected
 component $\widehat{ A}_1$ of $\widehat{  \pi}^{-1}(A)$ in $M'$ containing
 $w_0$ may be decomposed as
$$
\widehat{ A}_1= \mbox{  $\bigcup_{i,k} $} \, U_k^{(i)}, \ \  \
i=1,\ldots,N ;\,
 k=1,2,\ldots.
$$
We claim that
\begin{align}
 \label{eqn:defofK}
K:= \mbox{  $\bigcup _{g\in A'(z_0)}$} \,  g \sigma_{z_0},
\end{align}
satisfies the conclusion of Lemma \ref{lem:saigo}.
Note that $A'(z_0)$ and $A$ in $M$ are domains in $G$ and
$\sigma_{z_0} $ is connected in $M$; thus $K$ is a domain in $M$.

We begin by showing that
\begin{align}
 \label{eqn:beginning}
K= \mbox{  $\bigcup _{g\in K'(z_0)}$} \,   g \sigma_{z_0}  \qquad   \mbox{  with }  \ \   g
 \sigma_{z_0} \Subset K.
\end{align}
We first show $A \subset K \subset E$.
Indeed, $A \subset K$ follows from $A=\{g(z_0)\in M: g\in A'(z_0)\}$ and
$z_0\in \sigma_{z_0}$. To prove $K \subset E$, we first show that
$$
\sigma_k^{(i)}  (w) \sigma_{z_0}  \Subset E,  \quad w\in U_k ^{(i)} ;\
i=1,\ldots ,N;
k=1,2,\ldots ,
$$
where $\sigma_k ^{(i)} $ is the holomorphic section of $G $ over $U_k
^{(i)} $ via $\widehat{ \psi}_{w_0}$ defined in (3) in Lemma
\ref{eqn:saigono}.

Since $\widehat{ \pi}( w_k  ^{(i)}  )=z_i \in D$, we see from (ii) in Lemma
\ref{lem:d1d2} that $\sigma_k^{(i)}(w_k ^{(i)}  )(w_0)= w_k ^{(i)}   \in  U_k^{(i)}
\subset \widehat{ E}_1 \cap \widehat{ \pi}^{-1}(D)=\widehat{ D}_1$,  so that $\sigma_k^{(i)}  (w_k ^{(i)}  ) \in
\widehat{D}_1(w_0)=D'(z_0)$ by (\ref{eqn:d-1-prime}).
It follows from 2. in Lemma \ref{lem:s0nonsing}  that
$\sigma_k^{(i)}  (w_k ^{(i)} )\sigma_{z_0}  \Subset { D}$.
 Now let $w\in U_k^{(i)}$,
and set $f:= [\sigma_k ^{(i)} (w)][\sigma_k ^{(i)} (w_k) ]^{-1}\in G$. Then
by (3) in Lemma \ref{eqn:saigono} we have  $f\in {\cal V}$.
It follows that $\sigma_k ^{(i)} (w) \sigma_{z_0}=f [\sigma_k ^{(i)} (w_k
^{(i)} )]
\sigma_{z_0} \subset f(D)  \Subset E$, as required.

Next we show that $
g \sigma_{z_0}  \Subset E$  for any $ g\in A'(z_0)$.
Let $g\in A'(z_0)$. Since  $g(w_0) \in \widehat{ A}_1$, we have
 $w:=g(w_0)\in U_k^{(i)}  $ for some
$i,k$.  Since $\sigma_k^{(i)}  (w)(w_0)= w= g(w_0)$, we can find an $h\in
\widehat{ H}_{w_0}$  such that  $g= \sigma_k^{(i)}(w) h $. Since
 $h \sigma_{z_0} = \sigma_{z_0} $, we have  $
g \sigma_{z_0} = (\sigma_k
^{(i)}(w)  h) \sigma_{z_0}= \sigma_k ^{(i)}  (w)\sigma_{z_0}  \Subset E,
$
as required.

Now we show that $K= \bigcup_{g\in K'(z_0)}g \sigma_{z_0} $. Since
$A \subset K$, we have the inclusion $K\subset \bigcup_{g\in K'(z_0)}g
\sigma_{z_0} $.  To prove the reverse inclusion,  let  $g'\in K'(z_0)$. Then $g'(z_0)\in
K$, so that we can find $g\in A'(z_0)$ and $s\in \Sigma_{z_0}$ with $g'(z_0)=
g s(z_0)$. Thus there exists $h\in H_{z_0}$ with
$g'= g s h$ in $G$. We note that $gs \in K'(z_0)$. To see this,
since $\Sigma_{z_0}$ is connected and contains $e$ in $G$, there exists a
continuous curve $\gamma: t\in [0,1] \to s(t)$ in $\Sigma_{z_0} $ with $s(0)=e$ and $s(1)=s$. The continuous curve $\widetilde {\gamma}: t\in [0,1]
\to \widetilde { \gamma}(t):= gs(t)$ in $G$ satisfies $\widetilde {
\gamma}( t)(z_0) \in g \sigma_{z_0} \subset K$. Since $\widetilde
{\gamma}(0)= g \in A'(z_0) \subset K'(z_0) $, we have $\widetilde {
\gamma} \subset K'(z_0)$, and hence $gs= \widetilde { \gamma}(1)\in
K'(z_0)$.

It follows from  2. in
Proposition \ref{eqn:elm} applied to $g'=gsh\in K'(z_0)$ that $h\in K'(z_0)\cap
H_{z_0}= E'(z_0)\cap H_{z_0}= {\cal H}'(z_0)$.
 Consequently,
\begin{align*}
g' \sigma_{z_0}& = gs h \sigma_{z_0} = gs h \Sigma_{z_0}(z_0)\\
&= gs \Sigma_{z_0} h(z_0)
  \quad  \mbox{  (by 1. in Lemma \ref{liegroup-hz})} \\
&= g \Sigma_{z_0}(z_0)  \quad \quad  \mbox{ (since $\Sigma_{z_0} $ is a
 group and $h(z_0)=z_0$})\\
  &=g \sigma_{z_0} \subset K,
\end{align*}
as claimed.

To prove that $g \sigma _{z_0} \Subset K$ for a given $g\in K(z')$, since $K$ is open it
suffices to show that for a given $\zeta\in \partial \sigma_{z_0}
\Subset{ D}$ we have
$g(\zeta) \in K$.
 Since $K'(z_0) $ is an open
set in $G$ containing $e$,
there is a ball $B$ centered at $e$ with $
gB \Subset K'(z_0) $. From the previous assertion we have $gB\sigma_{z_0}  \subset K$. It thus  suffices to show
$g(\zeta) \in gB \sigma_{z_0} $,  or
equivalently,  $\zeta \in  B\sigma_{z_0}$. This is also equivalent
to $ B  ^{-1}  (\zeta)\cap \sigma_{z_0}\ne \emptyset$. This is clear,
since $\zeta\in \partial \sigma_{z_0}$ and $B  ^{-1}  (\zeta)  $ is a
neighborhood  of $\zeta$ in $M$.

\smallskip

Finally, we shall  show  that the
union $K= \bigcup_{g\in K'(z_0)}g \sigma_{z_0}$ is a foliation. Indeed,  using a similar argument as in the proof of (\ref{eqn:hkprimehat}), under the hypothesis $K'(z_0) \cap
H_{z_0}={\cal H}'(z_0) $,
 we have the following fact:
\begin{align*}
  \widehat{ \pi} ^{-1}  (g\sigma_{z_0})\cap \widehat{ K}_1=
 {\textstyle  \bigcup_{k=1}^\infty }\,  g h_k'
 \widehat{ \sigma}_{w_0}\  \ \ \ \mbox{  for } \ g\in K'(z_0).
\end{align*}
Moreover, since  $\widehat{ \pi} : M \to M'$ is a normal covering
and
$\widehat{ K}_1=\cup_{g\in K'(z_0)}g \widehat{ \sigma}_{w_0}$ is
a foliation of leaves of  $m_1$-dimensional nonsingular
generalized analytic sets $g \widehat{ \sigma}_{w_0}$, it follows from the
 above fact  that
 $K= \cup _{g\in K'(z_0)}g \sigma_{z_0}$
is a foliation of leaves of   $g{ \sigma}_{z_0}$, and hence
Lemma \ref{lem:saigo} is proved.  \hfill $\Box$

\smallskip
The domain $K$ in $M$ as defined in
(\ref{eqn:defofK}) completes the
 proof  of items (o) and (i) in $2-a$ of Theorem
\ref{thm:main-nonconnected}.
 For item (ii),  let $\sigma \Subset D$ be a parabolic nonsingular
 generalized analytic set in $D$. Then we have $\lambda(z)=\ const.$ on
 $\sigma$.
 Since ${ D}= \bigcup_{g\in D'(z_0)} g \sigma_{z_0} $
 is a foliation, we use the same argument as in the proof of (ii) in $1-a$  of Theorem \ref{thm:mainth} to obtain ${ \sigma} \subset g {
 \sigma}_{z_0} $ for some $g\in D'(z_0)$,
 as required.

\smallskip
To prove $2-b.$
 we assume that $\sigma_{z_0}$ is closed in $M$. Since  $\sigma_{z_0}
 \Subset D$, $\sigma_{z_0}$
 is an $m_1$-dimensional, irreducible, nonsingular, compact analytic set in
 $D$. The outline of the proof of $2-b.$ is as follows. First, using (i) in $2-a$, we may consider the quotient space
 $K_0:=K/\sigma_{z_0}$ which is an $(n-m_1)$-dimensional connected complex manifold. Moreover, $D_0:= D/\sigma_{z_0}$ is a subdomain of $K_0$, and $\partial D_0$ is smooth in $K_0$. This last fact follows  from (\ref{eqn:bounda}):
 $\partial D= \cup_{g\in \partial D'(z_0)} g \sigma_{z_0}$, together with
 the smoothness of $\partial D$.

To be precise, we work in the Lie group $G$. We use the notation $G_1 \approx
G_2$ to mean that $G_1$ and $G_2$
are isomorphic as complex manifolds.
Since ${\cal H}'(z_0)$ and
$\Sigma_{z_0} {\cal H}'(z_0)$ are closed Lie subgroups of $G$ by 2. in
Lemma \ref{liegroup-hz}, we may consider the quotient spaces
$$
 \mathfrak A:= G/{\cal H}'(z_0)   \quad  \mbox{  and }  \quad  \mathfrak
 B:= G /\Sigma_{z_0}
{\cal H'}(z_0);
$$
then $\mathfrak A$ is  an $n$-dimensional connected complex manifold
and $\mathfrak B$ is an
$(n-m_1)$-dimensional connected complex manifold. Since ${\cal H}'(z_0)\subset
\Sigma_{z_0} {\cal H}'(z_0)$ we have the  canonical projection:
$$
\pi_b^a: \ g {\cal H}'(z_0) \in \mathfrak A \ \mapsto \ g \Sigma_{z_0} {\cal H}'(z_0)
\in \mathfrak B,
$$
and the  quotient space
$$
\mathfrak C := \Sigma_{z_0} {\cal H}'(z_0) / {\cal H}'(z_0) \ \
   \quad  \mbox{ with   $(\pi_b^a)  ^{-1}  (\zeta) \approx
\mathfrak C$ \  for \  $\zeta\in B$.}
$$
 From 2. (iv) in Lemma \ref{liegroup-hz} we have $\mathfrak C
 \approx \sigma _{z_0}$.

\smallskip
Since $D'(z_0){\cal H}'(z_0) = D'(z_0)$ and $ K'(z_0){\cal H}'(z_0) = K'(z_0)
 $, the quotients   $D'(z_0)/ {\cal H}'(z_0)$ and $K'(z_0)/
{\cal H}'(z_0)  $ are well-defined and define domains, say  $D_a$ and $K_a$,  in
$ \mathfrak A$. Since $(\partial
D'(z_0)){\cal H}'(z_0)= \partial D'(z_0)$, we have   $D_a \Subset
K_a$ and $\partial D_a$ is smooth in $K_a$.
Similarly, since
$D'(z_0) \Sigma_{z_0} {\cal H}'(z_0)= D'(z_0)  $ and  $ K'(z_0)
\Sigma_{z_0} {\cal H}'(z_0)=  K'(z_0)$ by 2. (ii) in Lemma
\ref{liegroup-hz}, the quotients $D'(z_0)/ \Sigma_{z_0} {\cal
H}'(z_0)$ and  $K'(z_0)/ \Sigma_{z_0} {\cal H}'(z_0)$ define
domains, say
$D_0$ and $K_0$,
in $\mathfrak B$.  Since $g\in \partial D'(z_0)$ implies $g \Sigma_{z_0} {\cal
H}'(z_0) \subset \partial D'(z_0)$, it follows that  $D_0 \Subset K_0$ and
$\partial D_0$ is smooth in $K_0$. Thus,  $\pi_b^a(D_a)=D_0$  and
$\pi_b^a(K_a )=K_0 $.

\smallskip
On the other hand,  we have
$
D_{a} \approx D $ and $K_ {a} \approx K.
$ To see this, since ${\cal H}'(z_0) \subset H_{z_0}$, we have the
canonical projection
$$
\pi_1: g {\cal H}'(z_0) \in \mathfrak A \mapsto gH_{z_0}=g(z_0)\in M,
$$
 so that $( \mathfrak A, \pi_1)$ is an unramified  covering of $M$ (since $\widehat{
 H}_{w_0} \subset {\cal H}'(z_0) \subset H_{z_0}$ and $\widehat{
 H}_{w_0}$ is a normal subgroup of $H_{z_0}$).
Then the restriction $\pi_1:\,D_ a \ (K_a) \ \mapsto D \ (K) $ is
bijective. We only give the proof for $\pi_1:K_a \mapsto K$.
The surjectivity comes from 1. in Proposition \ref{eqn:elm} (applied
to $K$ instead of $D$). To prove
the injectivity, let $g_1, g_2 \in K'(z_0) $ with
$g_1(z_0)=g_2(z_0)$. We can find $h\in H_{z_0}$ with
$g_1=g_2h$. It follows from 3. in Proposition \ref{eqn:elm} that $h\in
H_{z_0} \cap K'(z_0)= {\cal H}'(z_0)$. Since ${\cal H}'(z_0)$ is a group, we have $g_1 {\cal H}'(z_0)= g_2 h {\cal H}'(z_0)= g_2
{\cal H}'(z_0)$, as required.

Setting ${ \pi}_0:= \pi_b^a|_{K_A}\circ(\pi_1) ^{-1}  |_K $, we have the holomorphic surjection
$$
{ \pi}_0: K \ \mapsto  K_a \ \mapsto  K_0,
$$
which satisfies $\pi_0(D)=D_0$. Moreover, let $z_1, \ z_2\in K $ and let
$g_1,g_2\in K'(z_0)$ with $z_1=g_1(z_0)$ and $z_2=g_2(z_0)$. Then
\begin{align}
 \label{eqn:visible3}
\pi_0(z_1)=\pi_0(z_2) \qquad \mbox{ iff }  \qquad  g_1 \sigma_{z_0} =
g_2 \sigma_{z_0}.
\end{align}

Assume first that $\pi_0(z_1)=\pi_0(z_2)$. Then
 by the definition of $\pi_0$ we have
$g_1 \Sigma_{z_0} {\cal H}'(z_0)=g_2 \Sigma_{z_0} {\cal H}'(z_0)$.
Since $g_i \Sigma_{z_0} {\cal H}'(z_0)(z_0)= g_i \sigma_{z_0}, \ i=1,2$,
we have $g_1 \sigma_{z_0} =g_2 \sigma_{z_0} $. Conversely,
assume that
$g_1 \sigma_{z_0} =g_2 \sigma_{z_0} $. Then we have $s_1,s_2\in
\Sigma_{z_0} $ and $h\in H_{z_0}$ with $g_1 s_1 =g_2 s_2 h$. By
 (3) in  Lemma \ref{liegroup-hz} (replacing $D'(z_0)$ by $K'(z_0)$) we have $g_is\in
K'(z_0),\  i=1,2$. It follows from 3. in Proposition \ref{eqn:elm} that
$h\in K'(z_0)\cap H_{z_0}= {\cal H}'(z_0)$. Since $\Sigma_{z_0} {\cal
H}'(z_0)$ is a group, we have $g _1 \Sigma_{z_0}
{\cal H}'(z_0)= g_2  \Sigma_{z_0} {\cal H}'(z_0)$; i.e.,
$\pi_0(\zeta_1)=\pi_0(\zeta_2)$. Thus (\ref{eqn:visible3}) is proved.

In particular, (\ref{eqn:visible3}) implies that $\pi_0  ^{-1}  (\zeta)\approx
\sigma_{z_0} $ for each $\zeta\in K$, and
$\pi_0  ^{-1}(D_0)=\cup _{g\in D'(z_0)} \, g \sigma_{z_0}= D$.

\smallskip
To complete the proof of $2-b.$ in
Theorem  \ref{thm:main-nonconnected}, it remains to verify that
$D_0$ is a Stein domain in $K_0$.

Fix
$\zeta\in D_0$. Take a point $z\in D$ with $ \pi_0(z) =\zeta$ and define
$
\lambda _0(\zeta ):= \lambda (z).
$
If we take another point $z'\in D$ with $\pi_0(z')=\zeta$,
there exists $g\in D'(z_0)$ with $z=g(z_0)$, and hence $z'\in g
\sigma_{z_0} $ by (\ref{eqn:visible3}). Since $\lambda(g \sigma_{z_0} )= const. =
\lambda (g(z_0))= \lambda (z)$,  it follows that
$\lambda(z')=\lambda(z)$; thus $\lambda_0$ is a well-defined real-valued function on
 $D_0$.  Since $\pi_0(D)=D_0$, it follows
that $-\lambda_0$ is plurisubharmonic exhaustion function on $D_0$.
 Noting  that assertion $\alpha.$ in the proof of Theorem
\ref{thm:mainth} holds under the conditions in Theorem
\ref{thm:main-nonconnected}, using the same method as in the proof of $1-b.$ of Theorem \ref{thm:mainth}, we see that $-\lambda_0$ is strictly plurisubharmonic on $D_0$; i.e., $D_0$ is Stein.
\hfill $\Box$

\end{document}